\documentclass[12pt,letterpaper]{article}

\usepackage{etoolbox}
\newtoggle{SUPPLEMENTAL}
\toggletrue{SUPPLEMENTAL}
\togglefalse{SUPPLEMENTAL}
\newcommand{\Supplemental}[2]{\iftoggle{SUPPLEMENTAL}{#1}{#2}}

\usepackage[T1]{fontenc}
\usepackage[margin=1in,letterpaper]{geometry}  
\usepackage{graphicx,grffile}
\usepackage{doi}
\usepackage{amsmath,amssymb,amsthm,cancel,staves}
\usepackage{mathtools}
\usepackage{epstopdf}
\usepackage{caption}
\usepackage[bottom]{footmisc}
\usepackage{alltt}
\usepackage[longnamesfirst]{natbib}
\usepackage{hypernat}
\usepackage[nodisplayskipstretch]{setspace}
\usepackage{enumerate}
\usepackage{hyperref}
\usepackage{cleveref}

\theoremstyle{plain}
\newtheorem{theorem}{Theorem}
\newtheorem{lemma}[theorem]{Lemma}
\newtheorem{proposition}[theorem]{Proposition}

\theoremstyle{definition}
\newtheorem{defn}{Definition}
\newtheorem{assumption}{Assumption}

\theoremstyle{remark}

\newenvironment{itemizecomp}%
  {\begin{itemize}%
    \setlength{\itemsep}{0pt}%
    \setlength{\parskip}{0pt}}%
  {\end{itemize}}

  {\begin{enumerate}%
    \setlength{\itemsep}{0pt}%
    \setlength{\parskip}{0pt}}%
  {\end{enumerate}}

\makeatletter
\newcommand{\vast}{\bBigg@{4}}
\newcommand{\Vast}{\bBigg@{5}}
\makeatother

\newcommand\independent{\protect\mathpalette{\protect\independenT}{\perp}} \def\independenT#1#2{\mathrel{\rlap{$#1#2$}\mkern2mu{#1#2}}}
\newcommand{\D}[2]{\frac{\partial #1}{\partial #2}}
\newcommand{\DD}[2]{\frac{\partial^2 #1}{\partial#2\partial#2'}} 
\newcommand{\Dz}[3][0]{\left.\frac{\partial}{\partial #3}#2\right|_{#3=#1}}
\newcommand{\DDz}[3][0]{\left.\frac{\partial^2}{\partial #3^2}#2\right|_{#3=#1}}

\DeclareMathOperator{\tr}{tr}
\DeclareMathOperator{\Var}{Var}

\def\citeposs#1{\citeauthor{#1}'s (\citeyear{#1})}

\newcommand{\trp}[1]{\tr\left(#1\right)}
\newcommand{\Varp}[1]{\Var\left(#1\right)}

\def\DHLhksqrt#1#2{\setbox0=\hbox{$#1\sqrt{#2\,}$}\dimen0=\ht0
  \advance\dimen0-0.2\ht0 
  \setbox2=\hbox{\vrule height\ht0 depth -\dimen0}%
  {\box0\lower0.4pt\box2}}

\let\originalleft\left
\let\originalright\right
\renewcommand{\left}{\mathopen{}\mathclose\bgroup\originalleft}
\renewcommand{\right}{\aftergroup\egroup\originalright}

  \renewenvironment{thebibliography}[1]{%
    \begin{oldthebibliography}{#1}%
      \setlength{\parskip}{0ex}%
      \setlength{\itemsep}{0ex}%
  }%
  {%
    \end{oldthebibliography}%
  }

\newcommand*\samethanks[1][\value{footnote}]{\footnotemark[#1]}

\newcommand{\mockalph}[1]{}  

\allowdisplaybreaks[3]

\title{\Supplemental{Supplemental appendix:\\}{} Fractional order statistic approximation for \\ nonparametric conditional quantile inference} 

\author{Matt Goldman\thanks{Goldman: \Supplemental{}{Microsoft, }mattgold@microsoft.com.  Kaplan (corresponding author): \Supplemental{}{Department of Economics, University of Missouri, }kaplandm@missouri.edu.  \Supplemental{}{We thank the co-editor (Oliver Linton), associate editor, referees, and Yixiao Sun for helpful comments and references, and Patrik Guggenberger and Andres Santos for feedback improving the clarity of presentation.  Thanks also to Brendan Beare, Karen Messer, and active audience members at seminars and conferences.  Thanks to Ruixuan Liu for providing and discussing code from \citet{FanLiu2016}. 
This paper was previously circulated as parts of ``IDEAL quantile inference via interpolated duals of exact analytic $L$-statistics'' and ``IDEAL inference on conditional quantiles.''
}} \and David M.\ Kaplan\samethanks} 
\date{\today} 

\begin{document}
\maketitle


\begin{abstract}
\Supplemental{This supplement includes longer proofs and additional details. 
}{%
Using and extending fractional order statistic theory, we characterize the $O(n^{-1})$ coverage probability error of the previously proposed confidence intervals for population quantiles using $L$-statistics as endpoints in \citet{Hutson1999}.  We derive an analytic expression for the $n^{-1}$ term, which may be used to calibrate the nominal coverage level to get $O\left(n^{-3/2}[\log(n)]^3\right)$ coverage error.  Asymptotic power is shown to be optimal.  Using kernel smoothing, we propose a related method for nonparametric inference on conditional quantiles.  This new method compares favorably with asymptotic normality and bootstrap methods in theory and in simulations.  Code is provided for both unconditional and conditional inference. 


}
\end{abstract}

\vfill

\begin{flushright}
\copyright\ 2016 by the authors. 
This manuscript version is made available under the CC-BY-NC-ND 4.0 license: \url{http://creativecommons.org/licenses/by-nc-nd/4.0/}
\end{flushright}

\doublespacing

\Supplemental{\vfill\pagebreak}{\newpage}

\section{Introduction\label{sec:intro}}

%
%
Quantiles contain information about a distribution's shape.  Complementing the mean, they capture heterogeneity, inequality, and other measures of economic interest.  Nonparametric conditional quantile models further allow arbitrary heterogeneity across regressor values. 
This paper concerns nonparametric inference on quantiles and conditional quantiles. 
In particular, we characterize the high-order accuracy of both \citeposs{Hutson1999} $L$-statistic-based confidence intervals (CIs) and our new conditional quantile CIs. 

Conditional quantiles appear across diverse topics because they are fundamental statistical objects.  
Such topics include wages \citep{Hogg1975,Chamberlain1994,Buchinsky1994}, 
infant birthweight \citep{Abrevaya2001}, demand for alcohol \citep{ManningEtAl1995}, 
and Engel curves \citep[][pp.\ 81--82]{AlanEtAl2005,Deaton1997}, which we examine in our empirical application.

%
%
We formally derive the coverage probability error (CPE) of the CIs from \citet{Hutson1999}, as well as asymptotic power of the corresponding hypothesis tests. 
\citet{Hutson1999} had proposed CIs for quantiles using $L$-statistics (interpolating between order statistics) as endpoints and found they performed well, but formal proofs were lacking. 
Using the analytic $n^{-1}$ term we derive in the CPE, we provide a new calibration to achieve $O\bigl(n^{-3/2}[\log(n)]^3\bigr)$ CPE, analogous to the \citet{HoLee2005a} analytic calibration of the CIs in \citet{BeranHall1993}. 

The theoretical results we develop contribute to the fractional order statistic literature and provide the basis for inference on other objects of interest explored in \citet{GoldmanKaplan2014b} and \citet{Kaplan2014}.  In particular, Theorem \ref{thm:cdferror} tightly links the distributions of $L$-statistics from the observed and `ideal' (unobserved) fractional order statistic processes. 
Additionally, Lemma \ref{lem:den-i} provides Dirichlet PDF and PDF derivative approximations. 

High-order accuracy is important for small samples (e.g.,\ for experiments) as well as nonparametric analysis with small \emph{local} sample sizes.  For example, if $n=1024$ and there are five binary regressors, then the smallest local sample size cannot exceed $1024/2^5=32$. 


For nonparametric conditional quantile inference, we apply the unconditional method to a local sample (similar to local constant kernel regression), smoothing over continuous covariates and also allowing discrete covariates.  
CPE is minimized by balancing the CPE of our unconditional method and the CPE from bias due to smoothing. 
We derive the optimal CPE and bandwidth rates, as well as a plug-in bandwidth when there is a single continuous covariate.  

Our $L$-statistic method has theoretical and computational advantages over methods based on normality or an unsmoothed bootstrap.  
The theoretical bottleneck for our approach is the need to use a uniform kernel.  
Nonetheless, even if normality or bootstrap methods assume an infinitely differentiable conditional quantile function (and hypothetically fit an infinite-degree local polynomial), our CPE is still of smaller order with one or two continuous covariates. 
Our method also computes more quickly than existing methods (of reasonable accuracy), handling even more challenging tasks in 10--15 seconds instead of minutes. 

Recent complementary work of \citet{FanLiu2016} also concerns a ``direct method'' of nonparametric inference on conditional quantiles.  They use a limiting Gaussian process to derive first-order accuracy in a general setting, whereas we use the finite-sample Dirichlet process to achieve high-order accuracy in an iid setting.  
\citet{FanLiu2016} also provide uniform (over $X$) confidence bands.  
We suggest a confidence band from interpolating a growing number of joint CIs (as in \citet{HorowitzLee2012}), although it will take additional work to rigorously justify.  A different, ad hoc confidence band described in Section \ref{sec:sim} generally outperformed others in our simulations. 

If applied to a local constant estimator with a uniform kernel and the same bandwidth, the \citet{FanLiu2016} approach is less accurate than ours due to the normal (instead of beta) reference distribution and integer (instead of interpolated) order statistics in their CI in equation (6). 
However, with other estimators like local polynomials or that in \citet{DonaldEtAl2012}, the \citet{FanLiu2016} method is not necessarily less accurate.  One limitation of our approach is that it cannot incorporate these other estimators, whereas Assumption GI(iii) in \citet{FanLiu2016} includes any estimator that weakly converges (over a range of quantiles) to a Gaussian process with a particular structure.  We compare further in our simulations.  One open question is whether using our beta reference and interpolation can improve accuracy for the general \citet{FanLiu2016} method beyond the local constant estimator with a uniform kernel; our Lemma \ref{lem:u1u2-approx} shows this at least retains first-order accuracy.  

%
%
The order statistic approach to quantile inference uses the idea of the probability integral transform, which dates back to R.\ A.\ \citet{Fisher1932}, Karl \citet{Pearson1933}, and \citet{Neyman1937}. 
For continuous $X_i\stackrel{iid}{\sim}F(\cdot)$, $F(X_i)\stackrel{iid}{\sim}\textrm{Unif}(0,1)$. 
Each order statistic from such an iid uniform sample has a known beta distribution for any sample size $n$. 
%
We show that the $L$-statistic linearly interpolating consecutive order statistics also follows an approximate beta distribution, with only $O(n^{-1})$ error in CDF.  
%
%
%
%
Although $O(n^{-1})$ is an asymptotic claim, the CPE of the CI using the $L$-statistic endpoint is bounded between the CPEs of the CIs using the two order statistics comprising the $L$-statistic, where one such CPE is too small and one is too big, for any sample size. 
This is an advantage over methods more sensitive to asymptotic approximation error.

%
%
Many other approaches to one-sample quantile inference have been explored. 
With Edgeworth expansions, \citet{HallSheather1988} and \citet{Kaplan2015} obtain two-sided $O(n^{-2/3})$ CPE. 
With bootstrap, smoothing is necessary for high-order accuracy.  This increases the computational burden and requires good bandwidth selection in practice.\footnote{For example, while achieving the impressive two-sided CPE of $O(n^{-3/2})$, \citet[][p.\ 833]{PolanskySchucany1997} admit, ``If this method is to be of any practical value, a better bandwidth estimation technique will certainly be required.''} See \citet[\S1]{HoLee2005b} for a review of bootstrap methods.  
Smoothed empirical likelihood \citep{ChenHall1993} also achieves nice theoretical properties, but with the same caveats. 

Other order statistic-based CIs dating back to \citet{Thompson1936} are surveyed in \citet[\S7.1]{DavidNagaraja2003}.  Most closely related to \citet{Hutson1999} is \citet{BeranHall1993}.  Like \citet{Hutson1999}, \citet{BeranHall1993} linearly interpolate order statistics for CI endpoints, but with an interpolation weight based on the binomial distribution. Although their proofs use expansions of the \citet{Renyi1953} representation instead of fractional order statistic theory, their $n^{-1}$ CPE term is identical to that for \citet{Hutson1999} other than the different weight. 
Prior work \citep[e.g.,][]{Bickel1967,Shorack1972} has established asymptotic normality of $L$-statistics and convergence of the sample quantile process to a Gaussian limit process, but without such high-order accuracy. 

The most apparent difference between the two-sided CIs of \citet{BeranHall1993} and \citet{Hutson1999} is that the former are symmetric in the order statistic index, whereas the latter are equal-tailed.  This allows \citet{Hutson1999} to be computed further into the tails. 
Additionally, our framework can be extended to CIs for interquantile ranges and two-sample quantile differences \citep{GoldmanKaplan2014b}, 
which has not been done in the R\'enyi representation framework. 
 
For nonparametric conditional quantile inference, in addition to the aforementioned \citet{FanLiu2016} approach, \citet{Chaudhuri1991} derives the pointwise asymptotic normal distribution of a local polynomial estimator.  
\citet{QuYoon2015} propose modified local linear estimators of the conditional quantile process that converge weakly to a Gaussian process, and they suggest using a type of bias correction that strictly enlarges a CI to deal with the first-order effect of asymptotic bias when using the MSE-optimal bandwidth rate.  

%
%

Section \ref{sec:cdf-err} contains our theoretical results on fractional order statistic approximation, which are applied to unconditional quantile inference in Section \ref{sec:inf-unconditional}.  
Section \ref{sec:inf-conditional} concerns our new conditional quantile inference method.  
An empirical application and simulation results are in Sections \ref{sec:empirical} and \ref{sec:sim}, respectively.  
Proof sketches are collected in Appendix \ref{sec:app-pfs}, while the supplemental appendix contains full proofs.  
The supplemental appendix also contains details of the plug-in bandwidth calculations, as well as additional empirical and simulation results. 

Notationally, 
$\phi(\cdot)$ and $\Phi(\cdot)$ are respectively the standard normal PDF and CDF, 
$\doteq$ should be read as ``is equal to, up to smaller-order terms'', 
$\asymp$ as ``has exact (asymptotic) rate/order of'', 
and $A_n=O(B_n)$ as usual. 
Acronyms used are those for cumulative distribution function (CDF), confidence interval (CI), coverage probability (CP), coverage probability error (CPE), and probability density function (PDF).

\section{Fractional order statistic theory}\label{sec:cdf-err}

In this section, we introduce notation and present our core theoretical results linking unobserved `ideal' fractional $L$-statistics with their observed counterparts.  

Given an iid sample $\{X_i\}_{i=1}^n$ of draws from a continuous CDF denoted\footnote{$F$ will often be used with a random variable subscript to denote the CDF of that particular random variable.  If no subscript is present, then $F(\cdot)$ refers to the CDF of $X$.  Similarly for the PDF $f(\cdot)$.} $F(\cdot)$, interest is in $Q(p) \equiv F^{-1}(p)$ for some $p\in(0,1)$, where $Q(\cdot)$ is the quantile function.  
For $u\in(0,1)$, the sample $L$-statistic commonly associated with $Q(u)$ is 
\begin{equation}
\label{eqn:QXLdef}
\hat Q^L_X(u)
  \equiv (1-\epsilon) X_{n:k}
         +\epsilon X_{n:k+1} , \quad
k \equiv \lfloor u(n+1)\rfloor, \quad
\epsilon \equiv u(n+1) - k ,
\end{equation}
where $\lfloor\cdot\rfloor$ is the floor function, $\epsilon$ is the interpolation weight, and $X_{n:k}$ denotes the $k$th order statistic (i.e.,\ $k$th smallest sample value).  While $Q(u)$ is latent and nonrandom, $\hat Q^L_X(u)$ is a random variable, and $\hat Q^L_X(\cdot)$ is a stochastic process, observed for arguments in $[1/(n+1),n/(n+1)]$.  

Let $\Xi_n\equiv\left\{k/(n+1)\right\}_{k=1}^n$ denote the set of quantiles corresponding to the observed order statistics.  If $u\in\Xi_n$, then no interpolation is necessary and $\hat Q^L_X(u)=X_{n:k}$.  As detailed in Section \ref{sec:inf-unconditional}, application of the probability integral transform yields exact coverage probability of a CI endpoint $X_{n:k}$ for $Q(p)$:  
$ P\left(X_{n:k}<F^{-1}(p)\right) = P\left( U_{n:k}<p\right) $, 
where $U_{n:k}\equiv F(X_{n:k})\sim\beta(k,n+1-k)$ is equal in distribution to the $k$th order statistic from $U_i\stackrel{iid}{\sim}\textrm{Unif}(0,1)$, $i=1,\ldots,n$ \citep[][\textbf{8.7.4}]{Wilks1962}. 
However, we also care about $u\notin\Xi_n$, in which case $k$ is fractional. 
To better handle such fractional order statistics, we will present a tight link between the marginal distributions of the stochastic process $\hat Q^L_X(\cdot)$ and those of the analogous `ideal' (I) process 
\begin{align}
\label{eqn:Qxidef}
\tilde Q^I_X(\cdot) &\equiv F^{-1}\left(\tilde Q^I_U(\cdot)\right) ,
\end{align}
where $\tilde Q^I_U(\cdot)$ is the ideal (I) uniform (U) fractional order ``statistic'' process.  We use a tilde in $\tilde Q^I_X(\cdot)$ and $\tilde Q^I_U(\cdot)$ instead of the hat like in $\hat Q^L_X(\cdot)$ to emphasize that the former are unobserved (hence not true statistics), whereas the latter is computable from the sample data.

This $\tilde Q^I_U(\cdot)$ in \eqref{eqn:Qxidef} is a Dirichlet process \citep{Ferguson1973,Stigler1977} on the unit interval with index measure $\nu\left([0,t]\right)=(n+1)t$.  Its univariate marginals are
\begin{align}\label{eqn:QUIdist}
\tilde Q^I_U(u) = U_{n:(n+1)u} \sim \beta\bigl((n+1)u,(n+1)(1-u)\bigr) .
\end{align}
The marginal distribution of $\left(\tilde Q^I_U(u_1),\tilde Q^I_U(u_2)-\tilde Q^I_U(u_1),\ldots,\tilde Q^I_U(u_k)-\tilde Q^I_U(u_{k-1})\right)$ for $u_1<\cdots<u_k$ is Dirichlet with parameters $\left(u_1(n+1),(u_2-u_1)(n+1),\ldots,(u_k-u_{k-1})(n+1)\right)$. 

For all $u\in\Xi_n$, $\tilde Q^I_X(u)$ coincides with $\hat Q^L_X(u)$; they differ only in their interpolation between these points.  Proposition \ref{prop:error-prob} shows $\tilde Q^I_X(\cdot)$ and $\hat Q^L_X(\cdot)$ to be closely linked in probability. 
\begin{proposition}\label{prop:error-prob}
For any fixed $\delta>0$ and $m>0$, define 
 $\mathcal{U}^\delta \equiv \{u\in(0,1) \mid \forall\,t\in(u-m,u+m), f\left(F^{-1}(t)\right) \ge \delta \}$ and 
 $\mathcal{U}^\delta_n \equiv \mathcal{U}^{\delta} \cap [\frac{1}{n+1},\frac{n}{n+1}]$; then, $\underset{u \in \mathcal{U}^\delta_n}{\sup}\left| \tilde Q^I_X(u) - \hat Q^L_X(u) \right| = O_p\left(n^{-1}\log(n)\right)$.
\end{proposition}
Although Proposition \ref{prop:error-prob} motivates approximating the distribution of $\hat Q^L_X(u)$ by that of $\tilde Q^I_X(u)$, it is not relevant to high-order accuracy.  In fact, its result is achieved by any interpolation between $X_{n:k}$ and $X_{n:k+1}$, not just $\hat Q^L_X(u)$; in contrast, the high-order accuracy we establish in Theorem \ref{thm:IDEAL-single} is only possible with precise interpolations like $\hat Q^L_X(u)$. 

Next, we consider marginal distributions of fixed dimension $J$. 
%
We also consider the Gaussian approximation to the sampling distribution of fractional order statistics. 
It is well known that the centered and scaled empirical process for standard uniform random variables converges to a Brownian bridge.  For standard Brownian bridge process $B(\cdot)$, we index by $u\in(0,1)$ the additional stochastic processes
\begin{align*}
\tilde Q^B_U(u) \equiv u + n^{-1/2} B(u) \qquad\textrm{and} \qquad \tilde Q^B_X(u) \equiv F^{-1}\bigl(\tilde Q^B_U(u)\bigr).
\end{align*}
The vector $\tilde Q^I_U(\mathbf u)$ has an ordered Dirichlet distribution (i.e.,\ the spacings between consecutive $\tilde Q^I_U(u_j)$ follow a joint Dirichlet distribution), while $\tilde Q^B_U(\mathbf u)$ is multivariate Gaussian. 
Lemma \ref{lem:den-i} in the appendix shows the close relationship between multivariate Dirichlet and Gaussian PDFs and PDF derivatives. 

Theorem \ref{thm:cdferror} shows the close distributional link among linear combinations of ideal, interpolated, and Gaussian-approximated fractional order statistics.  Specifically, for arbitrary weight vector $\boldsymbol{\psi} \in \mathbb{R}^J$, we (distributionally) approximate 
\begin{align}
\label{eqn:def-L}
L^L \equiv \sum_{j=1}^J \psi_j \hat Q^L_X(u_j)
\quad&\textrm{by}\quad
L^I \equiv \sum_{j=1}^J \psi_j \tilde Q^I_X(u_j)
,\\
\notag
\quad\textrm{or alternatively } &\textrm{by}\quad
L^B \equiv \sum_{j=1}^J \psi_j \tilde Q^B_X(u_j).
\end{align}

Our assumptions for this section are now presented, followed by the main theoretical result.  Assumption \ref{a:hut-pf} ensures that the first three derivatives of the quantile function are uniformly bounded in neighborhoods of the quantiles, $u_j$, which helps bound remainder terms in the proofs.  We use \textbf{bold} for vectors and \underline{underline} for matrices. 
\begin{assumption}\label{a:iid}
Sampling is iid: $X_i\stackrel{iid}{\sim} F$, $i=1,\ldots,n$. 
\end{assumption}
\begin{assumption}\label{a:hut-pf}
For each quantile $u_j$, the PDF $f(\cdot)$ (corresponding to CDF $F(\cdot)$ in \ref{a:iid}) satisfies 
(i)  $f(F^{-1}(u_j))>0$; 
(ii) $f''(\cdot)$ is continuous in some neighborhood of $F^{-1}(u_j)$, i.e.,\ $f\in C^2\left(U_\delta\left(F^{-1}(u_j)\right)\right)$ with $U_\delta(x)$ denoting some $\delta$-neighborhood of point $x\in\mathbb R$.
\end{assumption}
\begin{theorem}\label{thm:cdferror}
Define $\underline{\mathcal{V}}$ as the $J\times J$ matrix with row $i$, column $j$ entries 
$\underline{\mathcal{V}}_{i,j} = \min\{u_i,u_j\} - u_i u_j$.  
Let $\underline{\mathcal{A}}$ be the $J\times J$ matrix with main diagonal entries $\underline{\mathcal{A}}_{j,j}=f\left(F^{-1}(u_j)\right)$ and zeros elsewhere, and let
\begin{gather*}
\mathcal{V}_{\boldsymbol{\psi}} \equiv \boldsymbol{\psi}'\left(\underline{\mathcal{A}}^{-1}\underline{\mathcal{V}}\,\underline{\mathcal{A}}^{-1}\right)\boldsymbol{\psi}, \quad 
\mathbb{X}_0  \equiv \sum_{j=1}^J \psi_j F^{-1}(u_j).  
\end{gather*}
Let Assumption \ref{a:iid} hold, and let \ref{a:hut-pf} hold at $\mathbf{\bar u}$. 
Given the definitions in \eqref{eqn:QXLdef}, \eqref{eqn:Qxidef}, and \eqref{eqn:def-L}, the following results hold uniformly over $\mathbf u=\mathbf{\bar u}+o(1)$. 
\begin{enumerate}\def\theenumi{\roman{enumi}} 
\item \label{thm:cdferror-ptwise}For a given constant $K$,
\begin{align*}
P & \bigg(L^L<\mathbb{X}_{0} + n^{-1/2}K\bigg) 
   -P\bigg(L^I <\mathbb{X}_{0} + n^{-1/2}K\bigg)  \\
  &= \frac{K\exp\left\{-K^2/(2 \mathcal{V}_{\boldsymbol{\psi}})\right\}}{\sqrt{2\pi \mathcal{V}_{\boldsymbol{\psi}}^3}} 
     \left[\sum_{j=1}^J \left(\frac{\psi_j^2 \epsilon_j(1-\epsilon_j) }{f\left[F^{-1}(u_j)\right]^2}\right)\right] n^{-1}
    +O\left(n^{-3/2}[\log(n)]^3\right) , 
\end{align*}
where the remainder is uniform over all $K$. 
%
\item \label{thm:cdferror-unif} Uniformly over $K$, 
\begin{align*}
\sup_{K\in\mathbb R} 
  &\left[P\bigg(L^L<\mathbb{X}_{0} + n^{-1/2}K\bigg) 
        -P\bigg(L^I <\mathbb{X}_{0} + n^{-1/2}K\bigg)\right] \\
  &= \frac{e^{-1/2}}{\sqrt{2\pi \mathcal{V}_{\boldsymbol{\psi}}^2}} 
     \left[\sum_{j=1}^J \left(\frac{\psi_j^2 \epsilon_j(1-\epsilon_j) }{f\left[F^{-1}(u_j)\right]^2}\right)\right] n^{-1}
    +O\left(n^{-3/2}[\log(n)]^3\right) , \\
    \sup_{K\in\mathbb R} 
      &\left|P\bigg(L^L<\mathbb{X}_{0} + n^{-1/2}K\bigg) 
            -P\bigg(L^B <\mathbb{X}_{0} + n^{-1/2}K\bigg)\right| 
      = O\left(n^{-1/2}[\log(n)]^3\right) . 
\end{align*}
\end{enumerate}
\end{theorem}



\section{Quantile inference: unconditional}\label{sec:inf-unconditional}

For inference on $Q(p)$, we continue to maintain \ref{a:iid} and \ref{a:hut-pf}. 
For $p\in(0,1)$ and confidence level $1-\alpha$, define $u^h(\alpha)$ and $u^l(\alpha)$ to solve 
\begin{align}
\label{eqn:uhdef}
\alpha &= P\Bigl(\tilde Q^I_U\left(u^h(\alpha)\right)<p\Bigr)  ,  
\quad
\alpha = P\Bigl(\tilde Q^I_U\left(u^l(\alpha)\right)>p\Bigr)  ,  
\end{align}
with $\tilde Q^I_U(u)\sim\beta\bigl((n+1)u,(n+1)(1-u)\bigr)$ from \eqref{eqn:QUIdist}, parallel to (7) and (8) in \citet{Hutson1999}.  

One-sided CI endpoints for $Q(p)$ are $\hat Q^L_X(u^h)$ or $\hat Q^L_X(u^l)$. 
Two-sided CIs replace $\alpha$ with $\alpha/2$ in \eqref{eqn:uhdef} and use both endpoints.  This use of $\alpha/2$ yields the equal-tailed property; more generally, $t\alpha$ and $(1-t)\alpha$ can be used for $t\in(0,1)$.  

\begin{figure}[htbp]
  \centering
  \hspace*{\fill}%
  \includegraphics[clip=true,trim=5 45 25 55,width=0.45\textwidth]{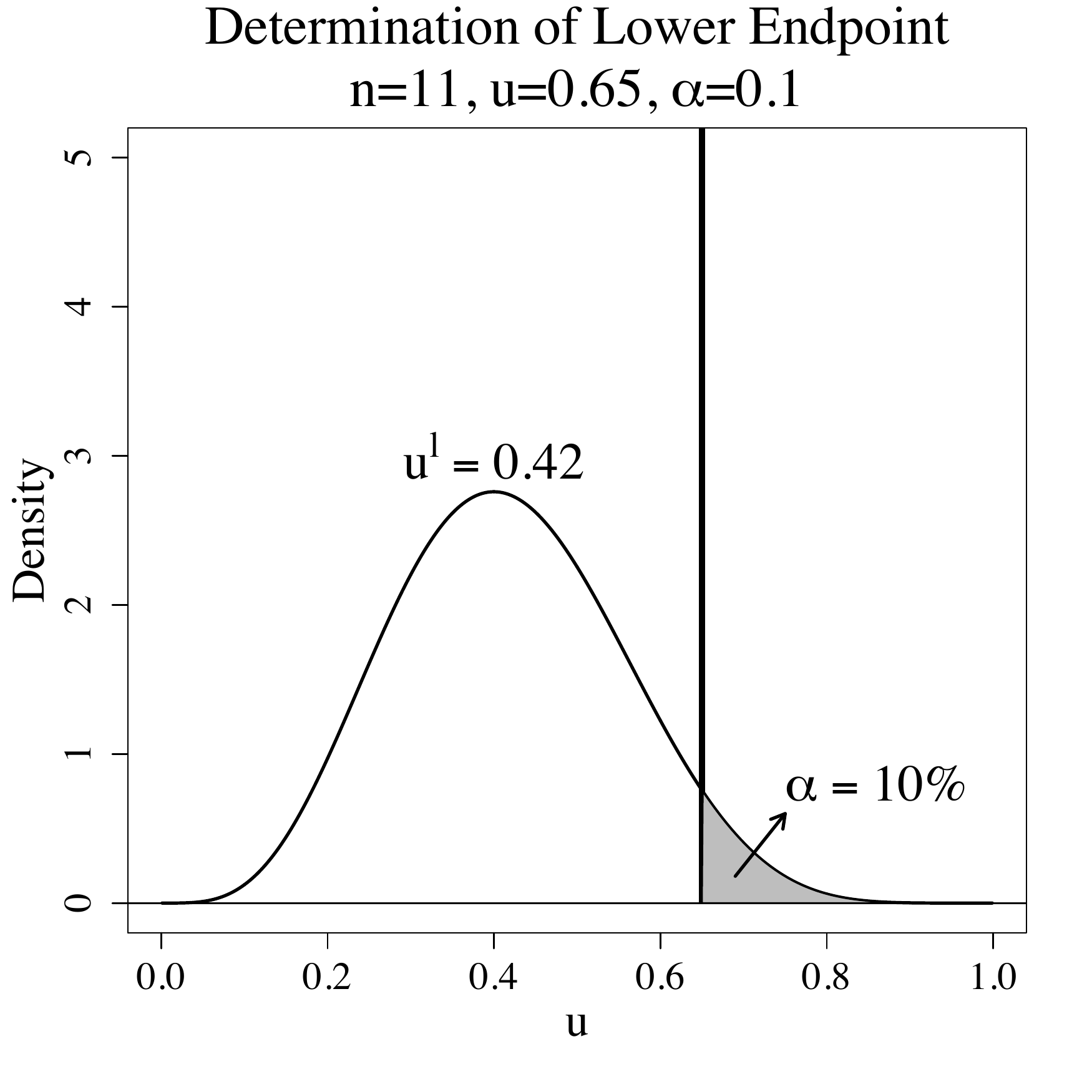}%
  \hfill%
  \includegraphics[clip=true,trim=5 45 25 55,width=0.45\textwidth]{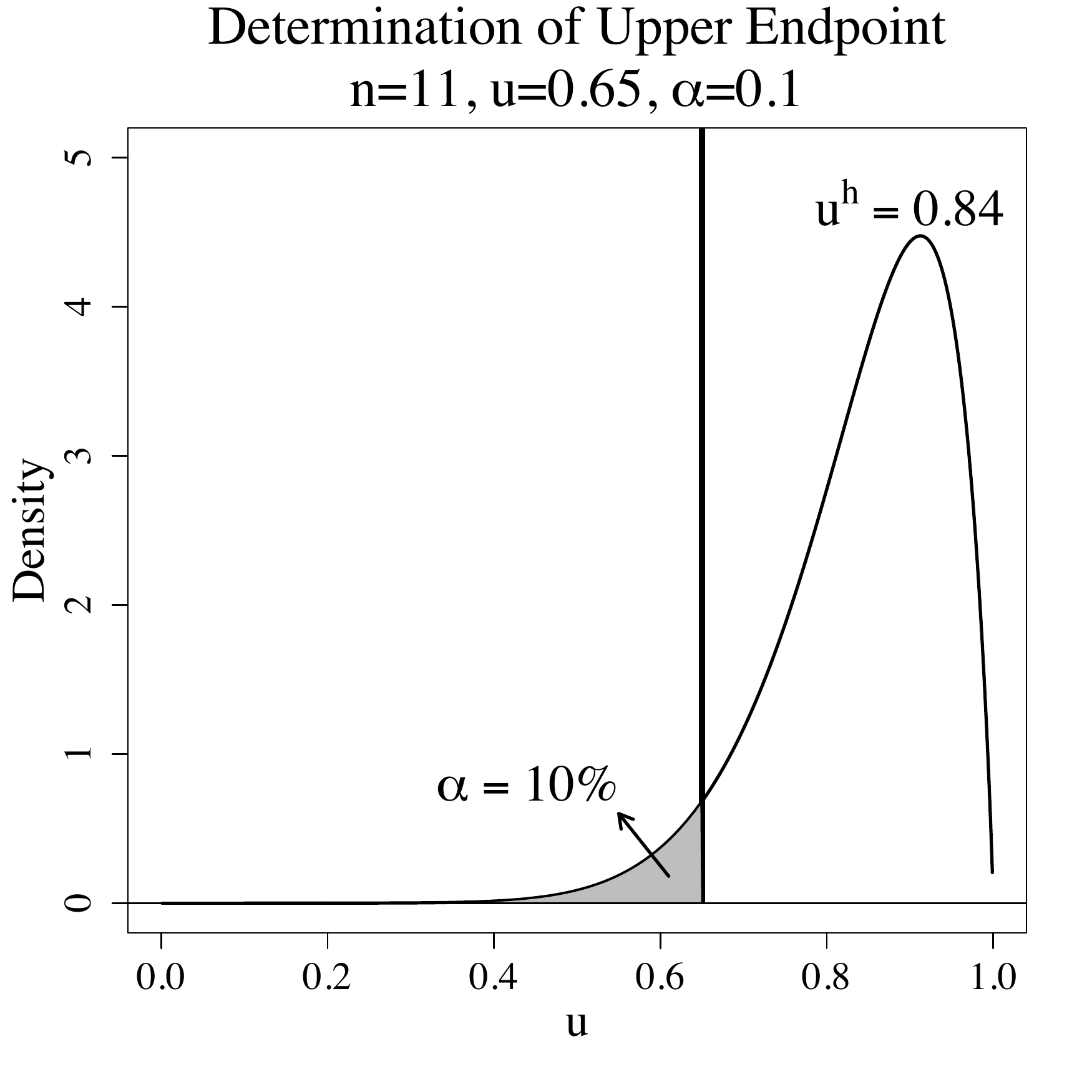}%
  \hspace*{\fill}%
  \caption{\label{fig:hutson-u-example}Example of one-sided CI endpoint determination, $n=11$, $p=0.65$, $\alpha=0.1$.  Left: $u^l$ makes the shaded region's area $P\bigl(\tilde Q^I_U(u^l)>p\bigr)=\alpha$.  Right: similarly, $u^h$ solves $P\bigl(\tilde Q^I_U(u^h)<p\bigr)=\alpha$.}
\end{figure}

Figure \ref{fig:hutson-u-example} visualizes an example.  The beta distribution's mean is $u^h$ (or $u^l$).  Decreasing $u^h$ increases the probability mass in the shaded region below $u$, while increasing $u^h$ decreases the shaded region, and vice-versa for $u^l$.  Solving \eqref{eqn:uhdef} is a simple numerical search problem. 

Lemma \ref{lem:u1u2-approx} shows the CI endpoint indices converge to $p$ at a $n^{-1/2}$ rate and may be approximated using quantiles of a normal distribution.  
\begin{lemma}\label{lem:u1u2-approx}
Let $z_{1-\alpha}$ denote the $(1-\alpha)$-quantile of a standard normal distribution, $z_{1-\alpha}\equiv \Phi^{-1}(1-\alpha)$.  From the definitions in \eqref{eqn:uhdef}, the values $u^l(\alpha)$ and $u^h(\alpha)$ can be approximated as
\begin{align*}
u^l(\alpha) 
  &= p - n^{-1/2}z_{1-\alpha}\sqrt{p(1-p)} - \frac{2p-1}{6n}(z_{1-\alpha}^2+2) +O(n^{-3/2}) , \\
u^h(\alpha) 
  &= p + n^{-1/2}z_{1-\alpha}\sqrt{p(1-p)} - \frac{2p-1}{6n}(z_{1-\alpha}^2+2) +O(n^{-3/2}) .
\end{align*}
\end{lemma}

For the lower one-sided CI, using \eqref{eqn:uhdef}, the $1-\alpha$ CI from \citet{Hutson1999} is
\begin{align}
\label{eqn:hutson-CI-lower}
\Bigl(-\infty,\hat Q^L_X\bigl(u^h(\alpha)\bigr)\Bigr) .
\end{align}
 Coverage probability is 
\begin{align*}
P & \left\{Q(p) \in \left(-\infty, \hat Q^L_X\bigl(u^h(\alpha)\bigr)\right)\right\}  
   = P\left(\hat Q^L_X\bigl(u^h(\alpha)\bigr)>Q(p)\right) \\
  &\!\!\!\!\stackrel{\textrm{Thm \ref{thm:cdferror}}}{=} 
     P\left(\tilde Q^I_X\bigl(u^h(\alpha)\bigr)>Q(p)\right) 
    +\frac{\epsilon_h(1-\epsilon_h)z_{1-\alpha}\exp\{-z_{1-\alpha}^2/2\}} {\sqrt{2\pi}u^h(\alpha)(1-u^h(\alpha))} n^{-1}
    +O\left(n^{-3/2}[\log(n)]^3\right)  \\
  &= 1-\alpha 
    +\frac{\epsilon_h(1-\epsilon_h)z_{1-\alpha}\phi(z_{1-\alpha})}{p(1-p)}n^{-1}
    +O\left(n^{-3/2}[\log(n)]^3\right),
\end{align*}
where $\phi(\cdot)$ is the standard normal PDF and the $n^{-1}$ term is non-negative. 
Similar to the \citet{HoLee2005a} calibration, we can remove the analytic $n^{-1}$ term with the calibrated CI
\begin{align}
\label{eqn:hutson-CI-lower-calibrated}
  \left(-\infty, \hat Q^L_X\left(u^h\left(\alpha+\frac{\epsilon_h(1-\epsilon_h)z_{1-\alpha}\phi(z_{1-\alpha})}{p(1-p)}n^{-1}\right)\right)\right),
\end{align}
which has CPE of order $O\left(n^{-3/2}[\log(n)]^3\right)$.  We follow convention and define $\textrm{CPE}\equiv\textrm{CP}-(1-\alpha)$, where CP is the actual coverage probability and $1-\alpha$ the desired confidence level. 

By parallel argument, \citeposs{Hutson1999} uncalibrated upper one-sided and two-sided CIs also have $O(n^{-1})$ CPE, or $O\left(n^{-3/2}[\log(n)]^3\right)$ with calibration. 
For the upper one-sided case, again using \eqref{eqn:uhdef}, the $1-\alpha$ Hutson CI and our calibrated CI are respectively given by 
\begin{equation}
\label{eqn:hutson-CI-upper}
\left(\hat Q^L_X\bigl(u^l(\alpha)\bigr) , \infty\right) , \quad
\left(\hat Q^L_X\left(u^l\left(\alpha + \frac{\epsilon_\ell(1-\epsilon_\ell)z_{1-\alpha}\phi(z_{1-\alpha})}{p(1-p)}n^{-1}\right)\right) , \infty\right) ,
\end{equation}
and for equal-tailed two-sided CIs,
\begin{align}
\label{eqn:hutson-CI-2s}
&\left(\hat Q^L_X\left[u^l\left(\alpha/2\right)\right],\hat Q^L_X\left(u^h\left(\alpha/2\right)\right)\right)   \quad\textrm{and}   \\
\label{eqn:hutson-CI-2s-calibrated}
\begin{split}
&\Bigg(\enspace \hat Q^L_X\left(u^l\left(\frac{\alpha}{2} + \frac{\epsilon_\ell(1-\epsilon_\ell)z_{1-\alpha/2}\phi(z_{1-\alpha/2})}{p(1-p)}n^{-1}\right)\right),  \\*
  &\qquad\quad  \hat Q^L_X\left(u^h\left(\frac{\alpha}{2} + \frac{\epsilon_h(1-\epsilon_h)z_{1-\alpha/2}\phi(z_{1-\alpha/2})}{p(1-p)}n^{-1}\right)\right)\enspace\Bigg).
\end{split}
\end{align}
Without calibration, in all cases the $n^{-1}$ CPE term is non-negative (indicating over-coverage). 

For relatively extreme quantiles $p$ (given $n$), the $L$-statistic method cannot be computed because the $(n+1)$th (or zeroth) order statistic is needed.  In such cases, our code uses the Edgeworth expansion-based CI in \citet{Kaplan2015}.  Alternatively, if bounds on $X$ are known a priori, they may be used in place of these ``missing'' order statistics to generate conservative CIs.  Regardless, as $n\to\infty$, the range of computable quantiles approaches $(0,1)$. 

The hypothesis tests corresponding to all the foregoing CIs achieve optimal asymptotic power against local alternatives.  The sample quantile is a semiparametric efficient estimator, so it suffices to show that power is asymptotically first-order equivalent to that of the test based on asymptotic normality.  
Theorem \ref{thm:IDEAL-single} collects all of our results on coverage and power. 

\begin{theorem}\label{thm:IDEAL-single}
Let $z_\alpha$ denote the $\alpha$-quantile of the standard normal distribution, and let $\epsilon_h=(n+1)u^h(\alpha)-\lfloor(n+1)u^h(\alpha)\rfloor$ and $\epsilon_\ell=(n+1)u^l(\alpha)-\lfloor(n+1)u^l(\alpha)\rfloor$.  
 Let Assumption \ref{a:iid} hold, and let \ref{a:hut-pf} hold at $p$.  Then, we have the following.
\begin{enumerate}\def\theenumi{\roman{enumi}}
 \item\label{cor:IDEAL-single-1s-CP} The one-sided lower and upper CIs in \eqref{eqn:hutson-CI-lower} and \eqref{eqn:hutson-CI-upper} have coverage probability  
 \[ 1-\alpha + \frac{\epsilon(1-\epsilon)z_{1-\alpha}\phi(z_{1-\alpha})}{p(1-p)}n^{-1} + O\left(n^{-3/2}[\log(n)]^3\right),\] 
 with $\epsilon=\epsilon_h$ for the former and $\epsilon=\epsilon_\ell$ for the latter. 
 \item\label{cor:IDEAL-single-2s-CP} The equal-tailed, two-sided CI in \eqref{eqn:hutson-CI-2s} has coverage probability 
 \[ 1-\alpha + \frac{[\epsilon_h(1-\epsilon_h)+\epsilon_\ell(1-\epsilon_\ell)]z_{1-\alpha/2}\phi(z_{1-\alpha/2})}{p(1-p)}n^{-1} + O\left(n^{-3/2}[\log(n)]^3\right). \] 
 \item\label{cor:IDEAL-single-calibrated} The calibrated one-sided lower, one-sided upper, and two-sided equal-tailed CIs given in \eqref{eqn:hutson-CI-lower-calibrated}, \eqref{eqn:hutson-CI-upper}, and \eqref{eqn:hutson-CI-2s-calibrated}, respectively, have $O\left(n^{-3/2}[\log(n)]^3\right)$ CPE. 
 \item\label{cor:IDEAL-single-power} 
 The asymptotic probabilities of excluding $D_n=Q(p)+\kappa n^{-1/2}$ from lower one-sided (l), upper one-sided (u), and equal-tailed two-sided (t) CIs (i.e.,\ asymptotic power of the corresponding hypothesis tests) are
 \begin{align*}
  \mathcal P_n^l(D_n) &\to \Phi\left(z_\alpha+S\right), \; 
  \mathcal P_n^u(D_n)  \to \Phi\left(z_\alpha-S\right), \;
  \mathcal P_n^t(D_n)  \to \Phi\left(z_{\alpha/2}+S\right)  
                          +\Phi\left(z_{\alpha/2}-S\right),
 \end{align*}
 where $S\equiv \kappa f(F^{-1}(p))/\sqrt{p(1-p)}$.
\end{enumerate}
\end{theorem}

The equal-tailed property of our two-sided CIs is a type of median-unbiasedness.  If $(\hat{L},\hat{H})$ is a CI for scalar $\theta$, then an equal-tailed CI is ``unbiased'' under loss function $L(\theta,\hat L,\hat H)=\max\{0,\theta-\hat H,\hat L-\theta\}$, as defined in (5) of \citet{Lehmann1951}.  
This median-unbiased property may be desirable \citep[e.g.,][footnote 11]{AndrewsGuggenberger2014}, although it is different than the usual ``unbiasedness'' where a CI is the inversion of an unbiased test.  
More generally, in \eqref{eqn:hutson-CI-2s}, we could replace $u^l(\alpha/2)$ and $u^h(\alpha/2)$ by $u^l(t\alpha)$ and $u^h((1-t)\alpha)$ for $t\in[0,1]$.  Different $t$ may achieve different optimal properties, which we leave to future work.

\section{Quantile inference: conditional}\label{sec:inf-conditional}

\subsection{Setup and bias}\label{sec:setup}

Let $Q_{Y|X}(u;x)$ be the conditional $u$-quantile function of scalar outcome $Y$ given conditioning vector $X\in\mathcal X\subset\mathbb R^d$, evaluated at $X=x$.  
The object of interest is $Q_{Y|X}(p;x_0)$, for $p\in(0,1)$ and interior point $x_0$. 
The sample $\{Y_i,X_i\}_{i=1}^n$ is drawn iid.  
Without loss of generality, let $x_0=0$. 

If $X$ is discrete so that $P(X=0)>0$, we can take the subsample with $X_i=0$ and compute a CI from the corresponding $Y_i$ values, using the method in Section \ref{sec:inf-unconditional}. 
Even with dependence like strong mixing among the $X_i$, CPE is the same $O(n^{-1})$ from Theorem \ref{thm:IDEAL-single} as long as the subsample's $Y_i$ are independent draws from the same $Q_{Y|X}(\cdot;0)$ and $N_n\stackrel{a.s.}{\asymp}n$. 

If $X$ is continuous, then $P(X_i=0)=0$, so observations with $X_i\ne0$ must be included.  
If $X$ contains mixed continuous and discrete components, then we can apply our method for continuous $X$ to each subsample corresponding to each unique value of the discrete subvector of $X$.  
The asymptotic rates are unaffected by the presence of discrete variables (although the finite-sample consequences may deserve more attention), so we focus on the case where all components of $X$ are continuous. 

We now present definitions and assumptions, continuing the normalization $x_0=0$. 

\begin{defn}[local smoothness]\label{def:smoothness}
Following \citet[pp.\ 762--3]{Chaudhuri1991}: if, in a neighborhood of the origin, function $g(\cdot)$ is continuously differentiable through order $k$, and its $k$th derivatives are uniformly H\"older continuous with exponent $\gamma\in(0,1]$, then $g(\cdot)$ has ``local smoothness'' of degree $s=k+\gamma$. 
\end{defn} 

\begin{assumption}\label{a:sampling}
Sampling of $(Y_i,X_i')'$ is iid, for continuous scalar $Y_i$ and continuous vector $X_i\in\mathcal X\subseteq\mathbb R^d$.  The point of interest $X=0$ is in the interior of $\mathcal X$, and the quantile of interest is $p\in(0,1)$. 
\end{assumption}
\begin{assumption}\label{a:f} 
The marginal density of $X$, denoted $f_X(\cdot)$, satisfies $0<f_X(0)<\infty$ and has local smoothness $s_X=k_X+\gamma_X>0$. 
\end{assumption}
\begin{assumption}\label{a:Q} 
For all $u$ in a neighborhood of $p$, $Q_{Y|X}(u;\cdot)$ (as a function of the second argument) has local smoothness\footnote{Our $s_Q$ corresponds to variable $p$ in \citet{Chaudhuri1991}; \citet{BhattacharyaGangopadhyay1990} use $s_Q=2$ and $d=1$.} $s_Q=k_Q+\gamma_Q>0$. 
\end{assumption}
\begin{assumption}\label{a:h} 
As $n\to\infty$, the bandwidth satisfies 
(i) $h\to0$, 
(i') $h^{b+d/2}\sqrt{n}\to0$ with $b\equiv\min\{s_Q,s_X+1,2\}$, 
(ii) $nh^d/[\log(n)]^2\to\infty$.
\end{assumption}
\begin{assumption}\label{a:GK1} 
For all $u$ in a neighborhood of $p$ and all $x$ in a neighborhood of the origin, $f_{Y|X}\left(Q_{Y|X}(u;x);x\right)$ is uniformly bounded away from zero.
\end{assumption}
\begin{assumption}\label{a:GK2} 
For all $y$ in a neighborhood of $Q_{Y|X}(p;0)$ and all $x$ in a neighborhood of the origin, $f_{Y|X}\left(y;x\right)$ has a second derivative in its first argument ($y$) that is uniformly bounded and continuous in $y$, having local smoothness $s_Y=k_Y+\gamma_Y>2$. 
\end{assumption}

Definition \ref{def:cube} refers to a window whose size depends on $h$: $C_h=[-h,h]$ if $d=1$, or more generally a hypercube as in \citet[pp.\ 763]{Chaudhuri1991}: 
letting $\|\cdot\|_\infty$ denote the $L_\infty$-norm, 
\begin{align}
\label{eqn:Ch}
C_h &\equiv \{x:x\in\mathbb R^d,\|x\|_\infty\le h\}, \quad
N_n \equiv \#\bigl(\{Y_i: X_i\in C_h, 1\le i\le n\}\bigr)  .
\end{align}
\begin{defn}[local sample]\label{def:cube}
Using $C_h$ and $N_n$ defined in \eqref{eqn:Ch}, the ``local sample'' consists of $Y_i$ values from observations with $X_i\in C_h\subset\mathbb R^d$, 
and the ``local sample size'' is $N_n$. 
Additionally, let the local quantile function $Q_{Y|X}(p;C_h)$ be the $p$-quantile of $Y$ given $X\in C_h$, satisfying 
$p = P\left(Y<Q_{Y|X}(p;C_h)\mid X\in C_h\right)$; similarly define the local CDF $F_{Y|X}(\cdot;C_h)$, local PDF $f_{Y|X}(\cdot;C_h)$, and derivatives thereof. 
\end{defn}

Given fixed values of $n$ and $h$, Assumption \ref{a:sampling} implies that the $Y_i$ in the local sample are independent and identically distributed,\footnote{This may be the case asymptotically even with substantial dependence, although we do not explore this point.  For example, \citet[p.\ 237]{PolonikYao2002} write, ``Only the observations with $X_t$ in a small neighbourhood of $x$ are effectively used\ldots [which] are not necessarily close with each other in the time space. Indeed, they could be regarded as asymptotically independent under appropriate conditions such as strong mixing\ldots.''} 
which is needed to apply Theorem \ref{thm:IDEAL-single}.  However, they do not have the quantile function of interest, $Q_{Y|X}(\cdot;0)$, but rather the biased $Q_{Y|X}(\cdot;C_h)$.  This is like drawing a global (any $X_i$) iid sample of wages, $Y_i$, and restricting it to observations in Japan ($X\in C_h$) when our interest is only in Tokyo ($X=0$): our restricted $Y_i$ constitute an iid sample from Japan, but the $p$-quantile wage in Japan may differ from that in Tokyo. 
Assumptions \ref{a:f}--\ref{a:h}(i) and \ref{a:GK2} are necessary for the calculation of this bias, $Q_{Y|X}(p;C_h)-Q_{Y|X}(p;0)$, in Lemma \ref{lem:bias}. 
Assumptions \ref{a:h}(ii) and \ref{a:GK1} (and \ref{a:sampling}) ensure $N_n\stackrel{a.s.}{\to}\infty$.  Assumptions \ref{a:GK1} and \ref{a:GK2} are conditional versions of Assumptions \ref{a:hut-pf}(i) and \ref{a:hut-pf}(ii), respectively. Their uniformity ensures uniformity of the remainder term in Theorem \ref{thm:IDEAL-single}, accounting for the fact that the local sample's distribution, $F_{Y|X}(\cdot;C_h)$, changes with $n$ (through $h$ and $C_h$). 

From \ref{a:h}(i), asymptotically $C_h$ is entirely contained within the neighborhoods implicit in \ref{a:f}, \ref{a:Q}, and \ref{a:GK2}.  This in turn allows us to examine only a local neighborhood around $p$ (e.g.,\ as in \ref{a:Q}) since the CI endpoints converge to the true value at a $N_n^{-1/2}$ rate.  

The $X_i$ being iid helps guarantee that $N_n$ is almost surely of order $nh^d$. 
The $h^d$ comes from the volume of $C_h$. 
Larger $h$ lowers CPE via $N_n$ but raises CPE via bias.  This tradeoff determines the optimal rate at which $h\to0$ as $n\to\infty$.  Using Theorem \ref{thm:IDEAL-single} and additional results on CPE from bias below, we determine the optimal value of $h$. 

\begin{defn}[steps to compute CI for $Q_{Y|X}(p;0)$]\label{def:method}
First, $C_h$ and $N_n$ are calculated as in Definition \ref{def:cube}.  
Second, using the $Y_i$ from observations with $X_i\in C_h$, a $p$-quantile CI is constructed as in \citet{Hutson1999}. 
If additional discrete conditioning variables exist, then repeat separately for each combination of discrete conditioning values.  
This procedure may be repeated for any number of $x_0$.  
For the bandwidth, we recommend the formulas in Section \ref{sec:opt-h}. 
\end{defn}

The bias characterized in Lemma \ref{lem:bias} is the difference between these two population conditional quantiles. 
\begin{lemma}\label{lem:bias}
Define $b$ as in \ref{a:h} and let $B_h\equiv Q_{Y|X}(p;C_h)-Q_{Y|X}(p;0)$. 
If Assumptions \ref{a:f}, \ref{a:Q}, \ref{a:h}(i), \ref{a:GK1}, and \ref{a:GK2} hold, then the bias is of order 
$|B_h| = O(h^b)$. 
%
Defining 
\begin{equation*}
\xi_p \equiv Q_{Y|X}(p;0) , \quad
F_{Y|X}^{(0,1)}(\xi_p;0) \equiv \Dz{F_{Y|X}(\xi_p;x)}{x} , \quad
F_{Y|X}^{(0,2)}(\xi_p;0) \equiv \DDz{F_{Y|X}(\xi_p;x)}{x} , 
\end{equation*}
with $d=1$, $k_X\ge1$, and $k_Q\ge2$, the bias is
\begin{equation}\label{eqn:Bh}
B_h
   = -h^2
      \frac{f_X(0) F_{Y|X}^{(0,2)}(\xi_p;0)
            +2 f_X'(0) F_{Y|X}^{(0,1)}(\xi_p;0)}
           {6 f_X(0) f_{Y|X}(\xi_p;0)}  
           +o(h^2) . 
\end{equation}
\end{lemma}
%
Equation \eqref{eqn:Bh} is the same as in \citet{BhattacharyaGangopadhyay1990}, who derive it using different arguments. 

\subsection{Optimal CPE order}\label{sec:opt-CPE}

The CPE-optimal bandwidth minimizes the sum of the two dominant high-order CPE terms.  It must be small enough to control the $O(h^b+N_nh^{2b})$ (two-sided) CPE from bias, but large enough to control the $O(N_n^{-1})$ CPE from applying the unconditional $L$-statistic method.  
The following theorem summarizes optimal bandwidth and CPE results. 

\begin{theorem}\label{thm:h-rate}
Let Assumptions \ref{a:sampling}--\ref{a:GK2} hold. 
The following results are for the method in Definition \ref{def:method}. 
%
For a one-sided CI, the bandwidth $h^*$ minimizing CPE has rate 
$ h^* \asymp n^{-3/(2b+3d)} $, 
corresponding to CPE of order $O(n^{-2b/(2b+3d)})$.  
For a two-sided CI, the optimal bandwidth rate is $h^*\asymp n^{-1/(b+d)}$, and the optimal CPE is $O(n^{-b/(b+d)})$.  
%
Using the calibration in Section \ref{sec:inf-unconditional}, if $p=1/2$, then the nearly (up to $\log(n)$) CPE-optimal two-sided bandwidth rate is $h^*\asymp n^{-5/(4b+5d)}$, yielding CPE of order $O\left(n^{-6b/(4b+5d)}[\log(n)]^3\right)$; if $p\ne1/2$, then $h^*\asymp n^{-3/(b+3d)}$ and CPE is $O\left( n^{-3b/(2b+6d)} [\log(n)]^3 \right)$. 
The nearly CPE-optimal calibrated one-sided bandwidth rate is $h^*\asymp n^{-2/(b+2d)}$, yielding CPE of order $O\left(n^{-3b/(2b+4d)}[\log(n)]^3\right)$.
\end{theorem}
%



As detailed in the supplemental appendix, Theorem \ref{thm:h-rate} implies that for the most common values of dimension $d$ and most plausible values of smoothness $s_Q$, even our uncalibrated method is more accurate than inference based on asymptotic normality with a local polynomial estimator.  
The same comparisons apply to basic bootstraps, which claim no refinement over asymptotic normality; in this (quantile) case, even Studentization does not improve theoretical CPE without the added complications of smoothed or $m$-out-of-$n$ bootstraps. 

The only opportunity for normality to yield smaller CPE is to greatly reduce bias by using a very large local polynomial if $s_Q$ is large; our approach implicitly uses a uniform kernel, so bias reduction beyond $O(h^2)$ is impossible.  Nonetheless, our method has smaller CPE when $d=1$ or $d=2$ even if $s_Q=\infty$, and in other cases the necessary local polynomial degree may be prohibitively large given common sample sizes. 



\begin{figure}[htbp]
\centering
\hfill
 \includegraphics[clip=true,trim=10 10 30 70,width=0.45\textwidth]{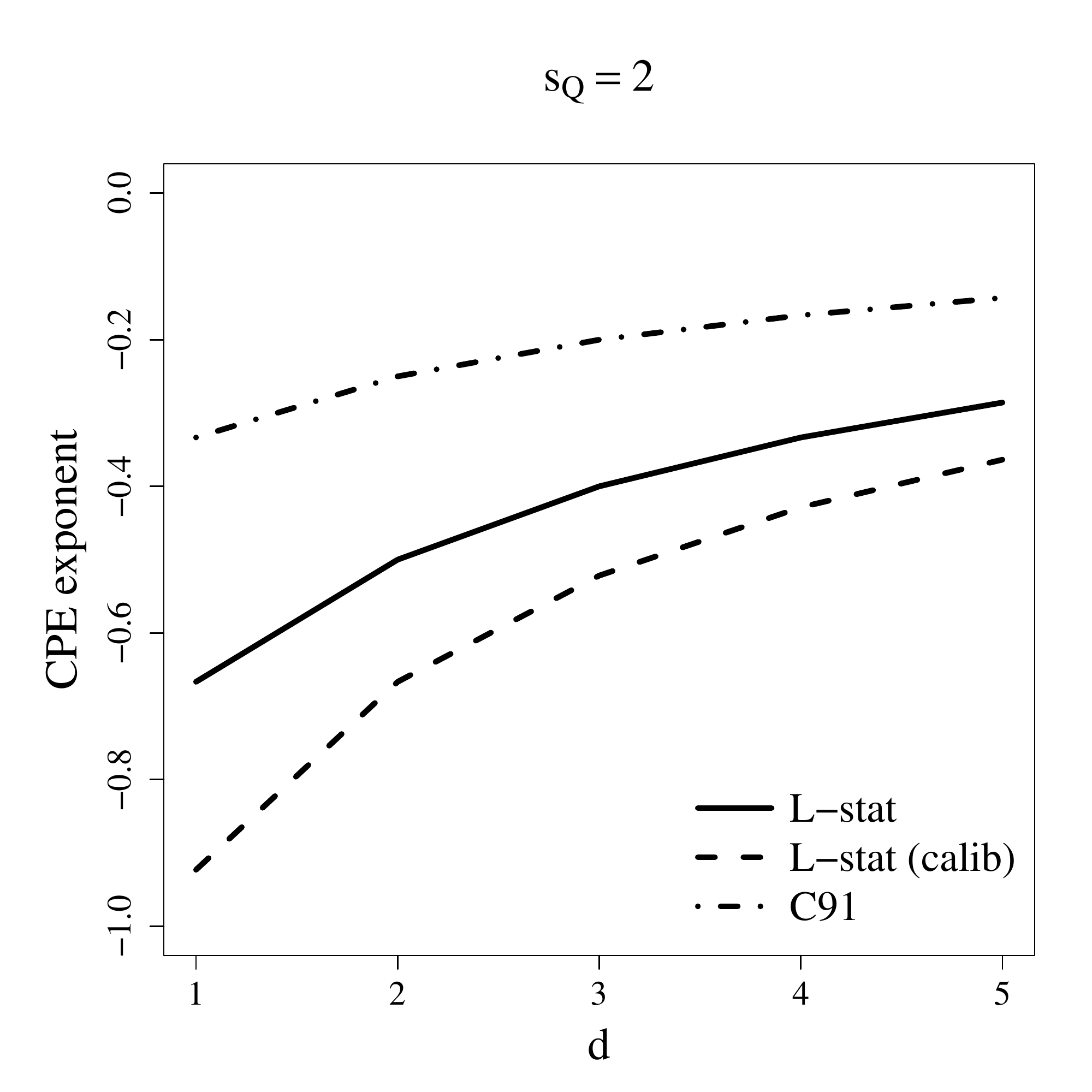}
 \hfill
 \includegraphics[clip=true,trim=10 10 30 70,width=0.45\textwidth]{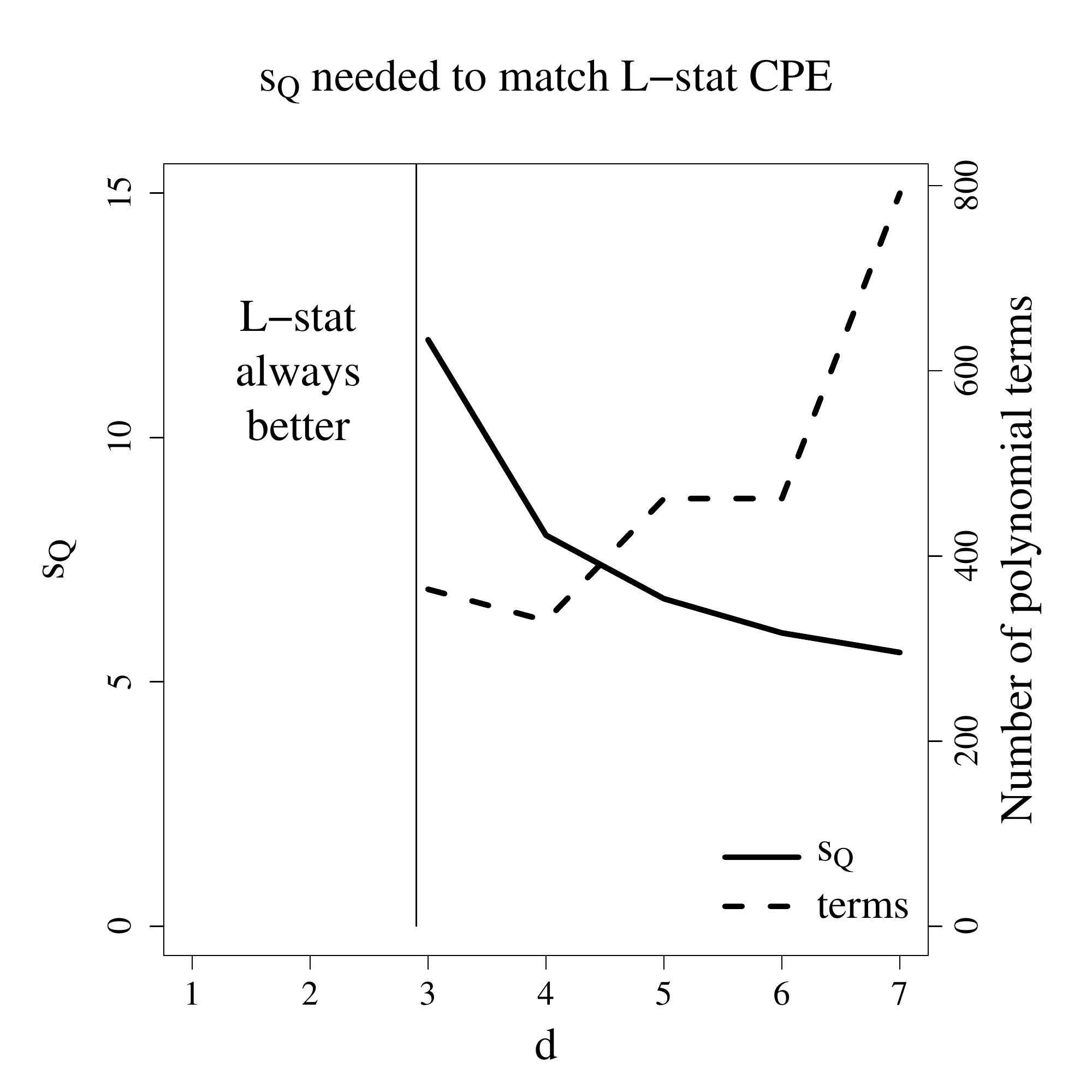}
 \hfill
 \null
 \caption{\label{fig:CPE-comp}Two-sided CPE comparison between new (``L-stat'') method and the local polynomial asymptotic normality method based on \citet{Chaudhuri1991}.  Left: with $s_Q=2$ and $s_X=1$, writing CPE as $n^\kappa$, comparison of $\kappa$ for different methods and different values of $d$.  Right: required smoothness $s_Q$ for the local polynomial normality-based CPE to match that of L-stat, as well as the corresponding number of terms in the local polynomial, for different $d$.}
\end{figure}
  
Figure \ref{fig:CPE-comp} (left panel) shows that if $s_Q=2$ and $s_X=1$, then the optimal CPE from asymptotic normality is always larger (worse) than our method's CPE.  
As shown in the supplement, CPE with normality is nearly $O\bigl(n^{-2/(4+2d)}\bigr)$. 
With $d=1$, this is $O\bigl(n^{-1/3}\bigr)$, much larger than our two-sided $O\bigl(n^{-2/3}\bigr)$.  
With $d=2$, $O\bigl(n^{-1/4}\bigr)$ is larger than our $O\bigl(n^{-1/2}\bigr)$.  
It remains larger for all $d$ since the bias is the same for both methods while the unconditional $L$-statistic inference is more accurate than normality. 

Figure \ref{fig:CPE-comp} (right panel) shows the required amount of smoothness and local polynomial degree for asymptotic normality to match our method's CPE.  
For the most common cases of $d=1$ and $d=2$, two-sided CPE with normality is larger even with infinite smoothness and a hypothetical infinite-degree polynomial. 
With $d=3$, to match our CPE, normality needs $s_Q\ge12$ and a local polynomial of degree $k_Q\ge11$.  
Since interaction terms are required, an $11$th-degree polynomial has $\sum_{T=d-1}^{k_Q+d-1} \binom{T}{d-1} = 364$ terms, which requires a large $N_n$ (and yet larger $n$). 
As $d\to\infty$, the required number of terms in the local polynomial only grows larger and may be prohibitive in realistic finite samples.  


\subsection{Plug-in bandwidth}\label{sec:opt-h}

We propose a feasible bandwidth value with the CPE-optimal rate. 
To avoid recursive dependence on $\epsilon$ (the interpolation weight), we fix its value.  This does not achieve the theoretical optimum, but it remains close even in small samples and seems to work well in practice.  The CPE-optimal bandwidth value derivation is shown for $d=1$ in the supplemental appendix; a plug-in version is implemented in our code.  
For reference, the plug-in bandwidth expressions are collected here.  
The $\alpha$-quantile of $N(0,1)$ is again denoted $z_\alpha$. 
We let $\hat B_h$ denote the estimator of bias term $B_h$; $\hat f_X$ the estimator of $f_X(x_0)$; $\hat f_X'$ the estimator of $f_X'(x_0)$; $\hat F_{Y|X}^{(0,1)}$ the estimator of $F_{Y|X}^{(0,1)}(\xi_p;x_0)$; and $\hat F_{Y|X}^{(0,2)}$ the estimator of $F_{Y|X}^{(0,2)}(\xi_p;x_0)$, with notation from Lemma \ref{lem:bias}. 

When $d=1$, the following are our CPE-optimal plug-in bandwidths.  
\begin{itemizecomp}
\item For one-sided inference, let
\begin{align}
\label{eqn:h-plugin-1s-pm}
\hat h_{+-}
  &= n^{-3/7} 
     \left(
      \frac{z_{1-\alpha}}
           {3 \left[p(1-p)\hat f_X\right]^{1/2} 
            \left[\hat f_X \hat F_{Y|X}^{(0,2)} 
                   +2 \hat f_X' \hat F_{Y|X}^{(0,1)} \right]}
     \right)^{2/7}  , \\
\label{eqn:h-plugin-1s-pp}
\hat h_{++} 
  &= -0.770\hat h_{+-}  .
\end{align}
%
  For lower one-sided inference, $\hat h_{+-}$ should be used if $\hat B_h<0$, and $\hat h_{++}$ otherwise.  
  For upper one-sided inference, $\hat h_{++}$ should be used if $\hat B_h<0$, and $\hat h_{+-}$ otherwise.  
%
\item For two-sided inference with general $p\in(0,1)$, 
\begin{align}
\label{eqn:h-plugin-2s-general}
\hat h
  &= n^{-1/3} \left( \frac
       { (\hat B_h/|\hat B_h|) (1-2p)
        +\sqrt{(1-2p)^2  +4} }
       {2 \left|   \hat f_X  \hat F_{Y|X}^{(0,2)}
                +2 \hat f_X' \hat F_{Y|X}^{(0,1)} \right| } 
     \right)^{1/3} , 
\end{align}
which simplifies to $\hat h = n^{-1/3} \bigl| \hat f_X \hat F_{Y|X}^{(0,2)} +2 \hat f_X' \hat F_{Y|X}^{(0,1)} \bigr|^{-1/3}$ with $p=0.5$. 
\end{itemizecomp}

While we suggest the CPE-optimal bandwidths for moderate $n$, we suggest shifting toward a larger bandwidth as $n\to\infty$.  Once CPE is small over a range of bandwidths, a larger bandwidth in that range is preferable since it yields shorter CIs.  
As an initial suggestion, we use a coefficient of $\max\{1,n/1000\}^{5/60}$ that keeps the CPE-optimal bandwidth for $n\le1000$ and then moves toward a $n^{-1/20}$ under-smoothing of the MSE-optimal bandwidth rate, as in \citet[p.\ 205]{FanLiu2016}.

\section{Empirical application}\label{sec:empirical}

We present an application of our $L$-statistic inference to \citet{Engel1857} curves.  Code is available from the latter author's website, and the data are publicly available. 

\citet{BanksEtAl1997} argue that a linear Engel curve is sufficient for certain categories of expenditure, while adding a quadratic term suffices for others.  Their Figure 1 shows nonparametrically estimated mean Engel curves (budget share $W$ against log total expenditure $\ln(X)$) with $95\%$ pointwise CIs at the deciles of the total expenditure distribution, using a subsample of 1980--1982 U.K.\ Family Expenditure Survey (FES) data.  

We present a similar examination, but for quantile Engel curves in the 2001--2012 U.K.\ Living Costs and Food Surveys \citep{UKLCFS2012}, which is a successor to the FES.  We examine the same four categories as in the original analysis: food; fuel, light, and power (``fuel''); clothing and footwear (``clothing''); and alcohol.  We use the subsample of households with one adult male and one adult female (and possibly children) living in London or the South East, leaving $8{,}528$ observations.  Expenditure amounts are adjusted to 2012 nominal values using annual CPI data.{\interfootnotelinepenalty10000\footnote{\url{http://www.ons.gov.uk/ons/datasets-and-tables/data-selector.html?cdid=D7BT&dataset=mm23&table-id=1.1}}}

\begin{table}[htbp]
\centering
\caption{\label{tab:emp}$L$-statistic $99\%$ CIs for various unconditional quantiles ($p$) of the budget share distribution, for different categories of expenditure described in the text.}
\begin{tabular}{lccc}
\multicolumn{1}{c}{Category} & $p=0.5$ & $p=0.75$ & $p=0.9$ \\
\hline
food     & (0.1532,0.1580) & (0.2095,0.2170) & (0.2724,0.2818) \\
fuel     & (0.0275,0.0289) & (0.0447,0.0470) & (0.0692,0.0741) \\
clothing & (0.0135,0.0152) & (0.0362,0.0397) & (0.0697,0.0761) \\
alcohol  & (0.0194,0.0226) & (0.0548,0.0603) & (0.1012,0.1111) \\
\hline
\hline
\end{tabular}
\end{table}
 
Table \ref{tab:emp} shows unconditional $L$-statistic CIs for various quantiles of the budget share distributions for the four expenditure categories.  (Due to the large sample size, calibrated CIs are identical at the precision shown.)  These capture some population features, but the conditional quantiles are of more interest. 

\begin{figure}[htbp]
 \centering
\hfill
 \includegraphics[clip=true,trim=10 40 30 70,width=0.49\textwidth]{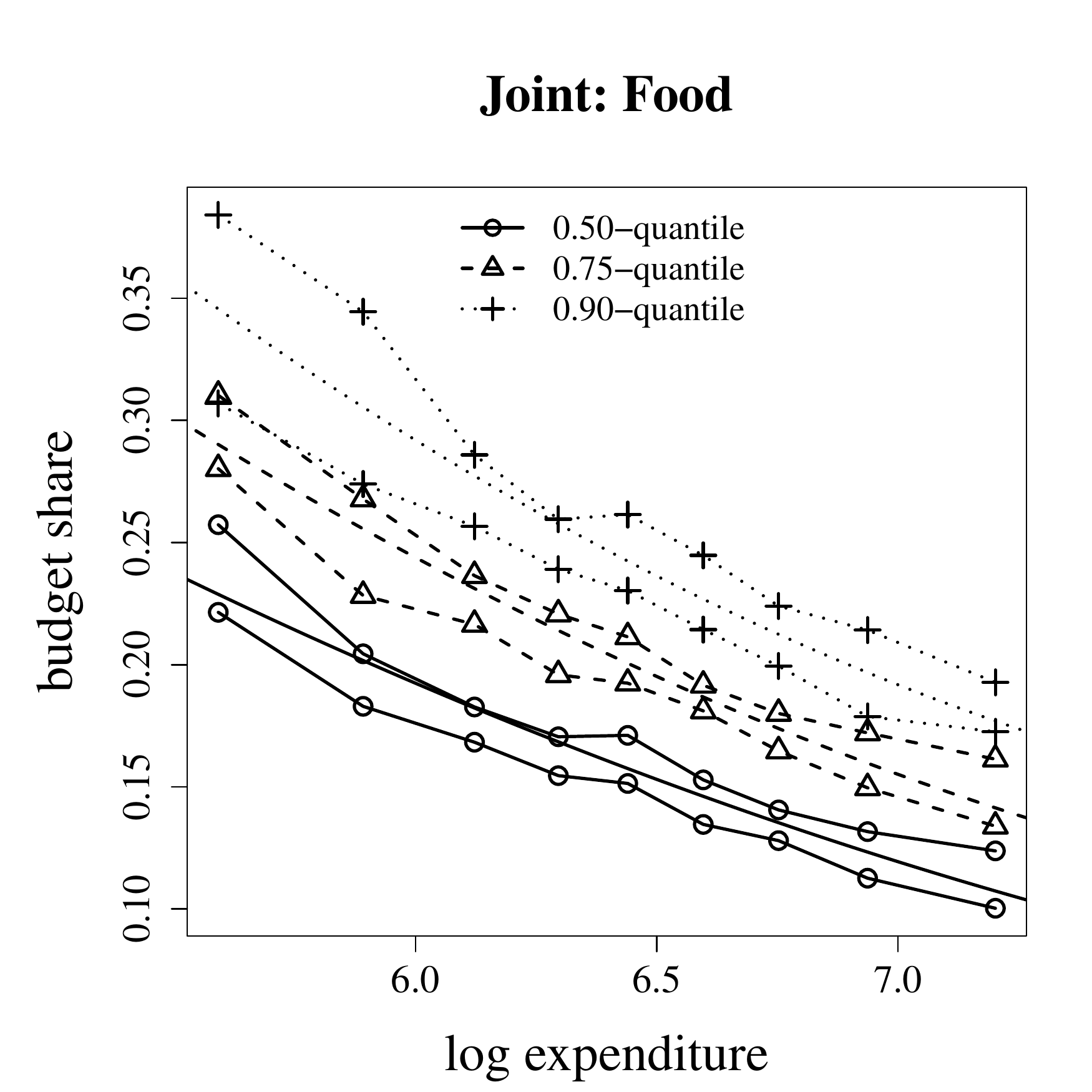}%
 \hfill%
 \includegraphics[clip=true,trim=40 40 0 70,width=0.49\textwidth]{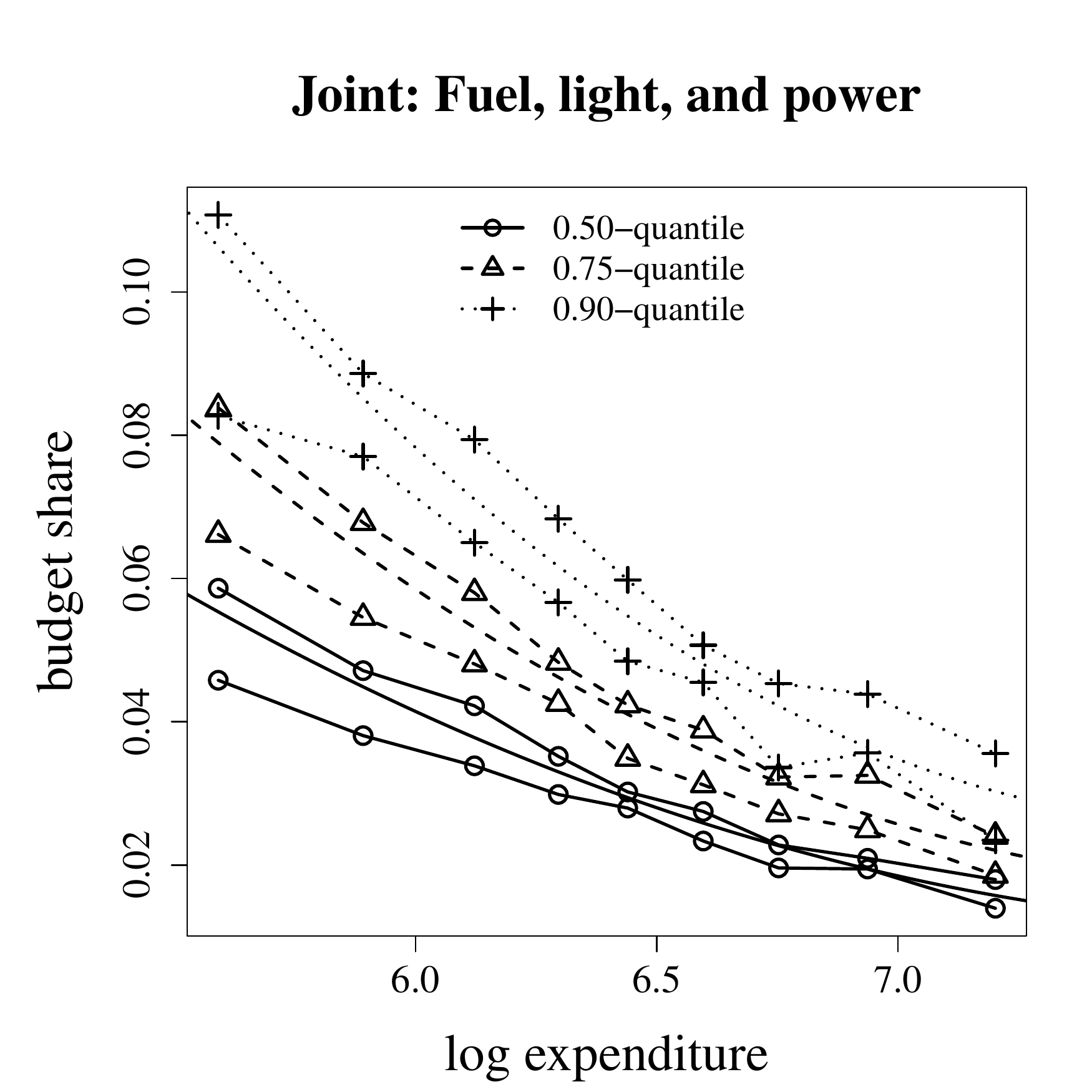}%
\hfill\null
 \\
\hfill
 \includegraphics[clip=true,trim=10 5 30 80,width=0.49\textwidth]{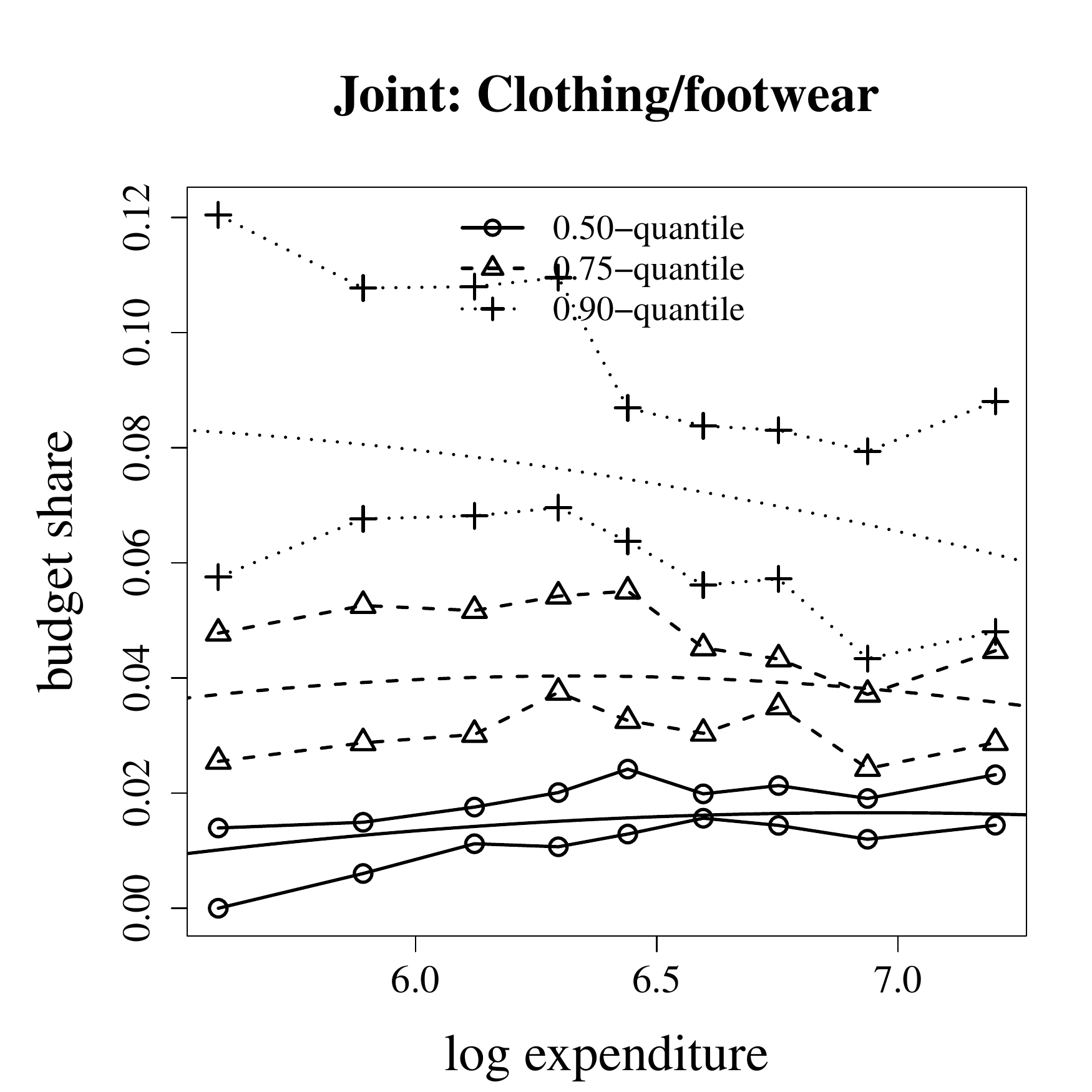}%
 \hfill%
 \includegraphics[clip=true,trim=40 5 0 80,width=0.49\textwidth]{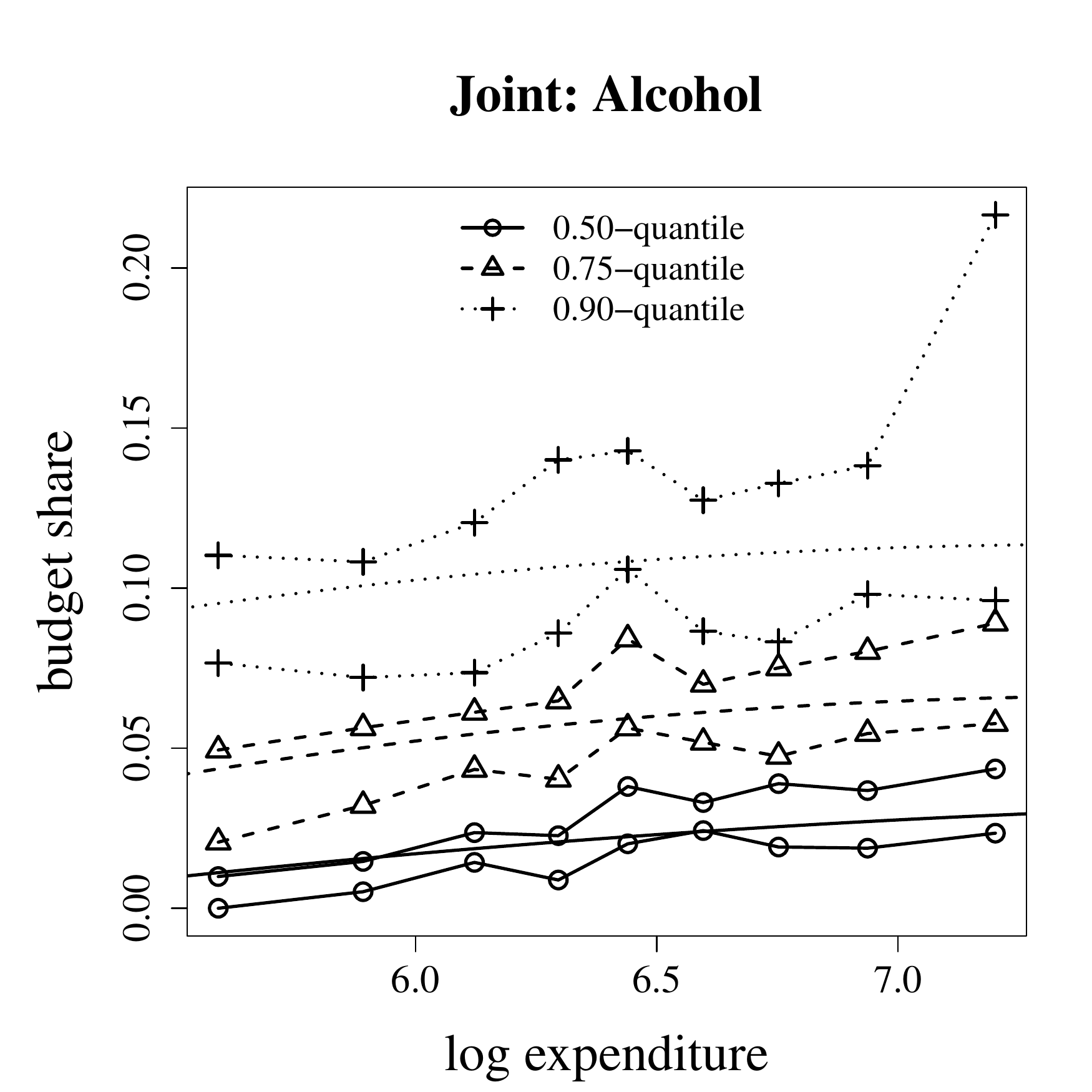}%
\hfill\null
 \caption{\label{fig:emp}Joint (over the nine expenditure levels) $90\%$ confidence intervals for quantile Engel curves: food (top left), fuel (top right), clothing (bottom left), and alcohol (bottom right).}
\end{figure}

Figure \ref{fig:emp} is comparable to Figure 1 of \citet{BanksEtAl1997} but with $90\%$ joint (over the nine expenditure levels) CIs instead of $95\%$ pointwise CIs, alongside quadratic quantile regression estimates.  (To get joint CIs, we simply use the Bonferroni adjustment and compute $1-\alpha/9$ pointwise CIs.)  Joint CIs are more intuitive for assessing the shape of a function since they jointly cover all corresponding points on the true curve with $90\%$ probability, rather than any given single point.  The CIs are interpolated only for visual convenience.  Although some of the joint CI shapes do not look quadratic at first glance, the only cases where the quadratic fit lies outside one of the intervals are for alcohol at the conditional median and clothing at the conditional upper quartile, and neither is a radical departure.  With a $90\%$ confidence level and 12 confidence sets, we would not be surprised if one or two did not cover the true quantile Engel curve completely.  Importantly, the CIs are relatively precise, too; the linear fit is rejected in 8 of 12 cases.  Altogether, this evidence suggests that the benefits of a quadratic (but not linear) approximation may outweigh the cost of approximation error. 

The supplemental appendix includes a similar figure but with a nonparametric (instead of quadratic) conditional quantile estimate along with joint CIs from \citet{FanLiu2016}.

\section{Simulation study}\label{sec:sim}

Code for our methods and simulations is available on the latter author's website.

\subsection{Unconditional simulations}\label{sec:sim-unconditional}

We compare two-sided unconditional CIs from the following methods: ``L-stat'' from Section \ref{sec:inf-unconditional}, originally in \citet{Hutson1999}; ``BH'' from \citet{BeranHall1993}; ``Norm'' using the sample quantile's asymptotic normality and kernel-estimated variance; ``K15'' from \citet{Kaplan2015}; and ``BStsym,'' a symmetric Studentized bootstrap (99 draws) with bootstrapped variance (100 draws).\footnote{Other bootstraps were consistently worse in terms of coverage: (asymmetric) Studentized bootstrap, and percentile bootstrap with and without symmetry.} 


Overall, L-stat and BH have the most accurate coverage probability (CP), avoiding under-coverage while maintaining shorter length than other methods achieving at least $95\%$ CP.  Near the median, L-stat and BH are nearly identical.  Away from the median, L-stat is closer to equal-tailed and often shorter than BH.  Farther into the tails, L-stat can be computed where BH cannot. 

\begin{table}[htbp]
\caption{\label{tab:sim-un1}CP and median CI length, $1-\alpha=0.95$; $n$, $p$, and distributions of $X_i$ ($F$) shown in table; $10{,}000$ replications.  ``Too high'' is the proportion of simulation draws in which the lower endpoint was above the true $F^{-1}(p)$, and ``too low'' is the proportion when the upper endpoint was below $F^{-1}(p)$.} 
\centering
\begin{tabular}[c]{cccccccc}
 $n$  &   $p$   &       $F$       & Method & CP & Too low & Too high & Length \\
\hline
$ 25$ & $0.5  $ &          Normal & L-stat & 0.953 & 0.022 & 0.025 & 0.99 \\
$ 25$ & $0.5  $ &          Normal & BH     & 0.955 & 0.021 & 0.024 & 1.00 \\
$ 25$ & $0.5  $ &          Normal & Norm   & 0.942 & 0.028 & 0.030 & 1.02 \\
$ 25$ & $0.5  $ &          Normal & K15    & 0.971 & 0.014 & 0.015 & 1.19 \\
$ 25$ & $0.5  $ &          Normal & BStsym & 0.942 & 0.028 & 0.030 & 1.13 \\[2pt]
$ 25$ & $0.5  $ &         Uniform & L-stat & 0.953 & 0.022 & 0.025 & 0.37 \\
$ 25$ & $0.5  $ &         Uniform & BH     & 0.954 & 0.021 & 0.025 & 0.37 \\
$ 25$ & $0.5  $ &         Uniform & Norm   & 0.908 & 0.046 & 0.046 & 0.35 \\
$ 25$ & $0.5  $ &         Uniform & K15    & 0.963 & 0.018 & 0.020 & 0.44 \\
$ 25$ & $0.5  $ &         Uniform & BStsym & 0.937 & 0.031 & 0.032 & 0.45 \\[2pt]
$ 25$ & $0.5  $ &     Exponential & L-stat & 0.953 & 0.024 & 0.023 & 0.79 \\
$ 25$ & $0.5  $ &     Exponential & BH     & 0.954 & 0.024 & 0.022 & 0.80 \\
$ 25$ & $0.5  $ &     Exponential & Norm   & 0.924 & 0.056 & 0.020 & 0.75 \\
$ 25$ & $0.5  $ &     Exponential & K15    & 0.968 & 0.022 & 0.010 & 0.96 \\
$ 25$ & $0.5  $ &     Exponential & BStsym & 0.941 & 0.039 & 0.020 & 0.93 \\
\hline
\hline
\end{tabular}
\end{table}

Table \ref{tab:sim-un1} shows nearly exact CP for both L-stat and BH when $n=25$ and $p=0.5$.  
``Norm'' can be slightly shorter, but it under-covers.  
The bootstrap has only slight under-coverage, and K15 none, but their CIs are longer than L-stat's. 
Additional results are in the supplemental appendix, but the qualitative points are the same. 

\begin{table}[htbp]
\caption{\label{tab:sim-un2}CP and median CI length, as in Table \ref{tab:sim-un1}.}
\centering
\begin{tabular}[c]{cccccccc}
 $n$  &   $p$   &       $F$       & Method & CP & Too low & Too high & Length \\
\hline
$ 99$ & $0.037$ &          Normal & L-stat & 0.951 & 0.023 & 0.026 & 1.02 \\
$ 99$ & $0.037$ &          Normal & BH     &   NA &   NA &   NA &   NA \\
$ 99$ & $0.037$ &          Normal & Norm   & 0.925 & 0.016 & 0.059 & 0.83 \\
$ 99$ & $0.037$ &          Normal & K15    & 0.970 & 0.009 & 0.021 & 1.55 \\
$ 99$ & $0.037$ &          Normal & BStsym & 0.950 & 0.020 & 0.030 & 1.20 \\[2pt]
$ 99$ & $0.037$ &          Cauchy & L-stat & 0.950 & 0.022 & 0.028 & 39.37 \\
$ 99$ & $0.037$ &          Cauchy & BH     &   NA &   NA &   NA &   NA \\
$ 99$ & $0.037$ &          Cauchy & Norm   & 0.784 & 0.082 & 0.134 & 18.90 \\
$ 99$ & $0.037$ &          Cauchy & K15    & 0.957 & 0.002 & 0.041 & 36.55 \\
$ 99$ & $0.037$ &          Cauchy & BStsym & 0.961 & 0.002 & 0.037 & 48.77 \\[2pt]
$ 99$ & $0.037$ &         Uniform & L-stat & 0.951 & 0.024 & 0.026 & 0.07 \\
$ 99$ & $0.037$ &         Uniform & BH     &   NA &   NA &   NA &   NA \\
$ 99$ & $0.037$ &         Uniform & Norm   & 0.990 & 0.000 & 0.010 & 0.12 \\
$ 99$ & $0.037$ &         Uniform & K15    & 0.963 & 0.028 & 0.009 & 0.11 \\
$ 99$ & $0.037$ &         Uniform & BStsym & 0.924 & 0.053 & 0.022 & 0.08 \\
\hline
\hline
\end{tabular}
\end{table}

Table \ref{tab:sim-un2} shows a case in the lower tail with $n=99$ where BH cannot be computed (because it needs the zeroth order statistic).  
Even then, L-stat's CP remains almost exact, and it is closest to equal-tailed.  
``Norm'' under-covers for two $F$ (severely for Cauchy) and is almost twice as long as L-stat for the third.  
BStsym has less under-coverage, and K15 none, but both are generally longer than L-stat. 
Again, additional results are in the supplemental appendix, with similar patterns. 

The supplemental appendix contains additional simulation results for $p\ne0.5$ but where BH is still computable.  
L-stat and BH both attain $95\%$ CP, but L-stat is much closer to equal-tailed and is shorter.  
The supplemental appendix also has results illustrating the effect of calibration.

Table \ref{tab:sim-beta-norm1} isolates the effects of using the beta distribution rather than the normal approximation, as well as the effects of interpolation.  Method ``Normal'' uses the normal approximation to determine $u^h$ and $u^l$ but still interpolates, while ``Norm/floor'' uses the normal approximation with no interpolation as in equations (5) and (6) of \citet[Ex.\ 2.1]{FanLiu2016}.  

\begin{table}[htbp]
\caption{\label{tab:sim-beta-norm1}CP and median CI length, $n=19$, $Y_i\stackrel{iid}{\sim}N(0,1)$, $1-\alpha=0.90$, $1{,}000$ replications, various $p$.  In parentheses below CP are probabilities of being too low or too high, as in Table \ref{tab:sim-un1}. 
Methods are described in the text.}
\centering
\begin{tabular}[c]{lcccccccc}
 && \multicolumn{3}{c}{Two-sided CP} && \multicolumn{3}{c}{} \\
 && \multicolumn{3}{c}{(Too low, Too high)} && \multicolumn{3}{c}{Median length} \\
\cline{3-5}\cline{7-9}
Method && $p=0.15$ & $p=0.25$ & $p=0.5$ && $p=0.15$ & $p=0.25$ & $p=0.5$ \\
\hline
L-stat      && 0.905 & 0.901 & 0.898 && 1.20 & 1.03 & 0.93 \\
            && (0.048,0.047) & (0.050,0.049) & (0.052,0.050) &&  &  &  \\
Normal      &&    NA & 0.926 & 0.912 &&   NA & 1.22 & 1.00 \\
            && (NA,NA) & (0.062,0.012) & (0.045,0.043) &&  &  &  \\
Norm/floor  &&    NA & 0.913 & 0.876 &&   NA & 1.47 & 0.91 \\
            && (NA,NA) & (0.083,0.004) & (0.087,0.037) &&  &  &  \\
\hline
\end{tabular}
\end{table}

Table \ref{tab:sim-beta-norm1} shows several advantages of L-stat. 
First, for $p=0.15$, Normal and Norm/floor cannot even be computed (hence ``NA'') because they require the zeroth order statistic, which does not exist, whereas L-stat is computable and has nearly exact CP ($0.905$). 
%
Second, with $p=0.25$ and $p=0.5$, the normal approximation (Normal) makes the CI needlessly longer than L-stat's CI.  
Third, additionally not interpolating (Norm/floor) makes the CI even longer for $p=0.25$ but leads to under-coverage for $p=0.5$. 
%
Fourth, whereas the L-stat CIs are almost exactly equal-tailed, the normal-based CIs are far from equal-tailed at $p=0.25$, where Norm/floor is essentially a one-sided CI. 

\subsection{Conditional simulations}\label{sec:sim-conditional}


For conditional quantile inference, we compare our $L$-statistic method 
(``L-stat'') with a variety of others. 
Implementation details may be seen in the supplemental appendix and available code. 
The first other method (``rqss'') is from the popular \texttt{quantreg} package in R \citep{R.quantreg}.  
%
The second (``boot'') is a local cubic method following \citet{Chaudhuri1991} but with bootstrapped standard errors; 
the bandwidth is L-stat's multiplied by $n^{1/12}$ to get the local cubic CPE-optimal rate. 
%
The third (``QYg'') uses the asymptotic normality of a local linear estimator with a Gaussian kernel, using results and ideas from \citet{QuYoon2015}, although they are more concerned with uniform (over quantiles) inference; they suggest using the MSE-optimal bandwidth (Corollary 1) and a particular type of bias correction (Remark 7).  
The fourth (``FLb'') is from Section 3.1 in \citet{FanLiu2016}, based on a symmetrized $k$-NN estimator using a bisquare kernel; we use the code from their simulations.\footnote{Graciously provided to us.  The code differs somewhat from the description in their text, most notably by an additional factor of $0.4$ in the bandwidth.}  Interestingly, although in principle they are just slightly undersmoothing the MSE-optimal bandwidth, their bandwidth is very close to the CPE-optimal bandwidth for the sample sizes considered. 
  
We now write $x_0$ as the point of interest, instead of $x_0=0$; we also take $d=1$, $b=2$, and focus on two-sided inference, both pointwise (single $x_0$) and joint (over multiple $x_0$). 
Joint CIs for all methods are computed using the Bonferroni approach.  
Uniform bands are also examined, with L-stat, QYg, and boot relying on the adjusted critical value from the \citet{Hotelling1939} tube computations in \texttt{plot.rqss}.  
Each simulation has $1{,}000$ replications unless otherwise noted.

\begin{figure}[thbp]
  \centering 
\hfill
  \includegraphics[clip=true,trim=15 15 10 58,width=0.445\textwidth]
    {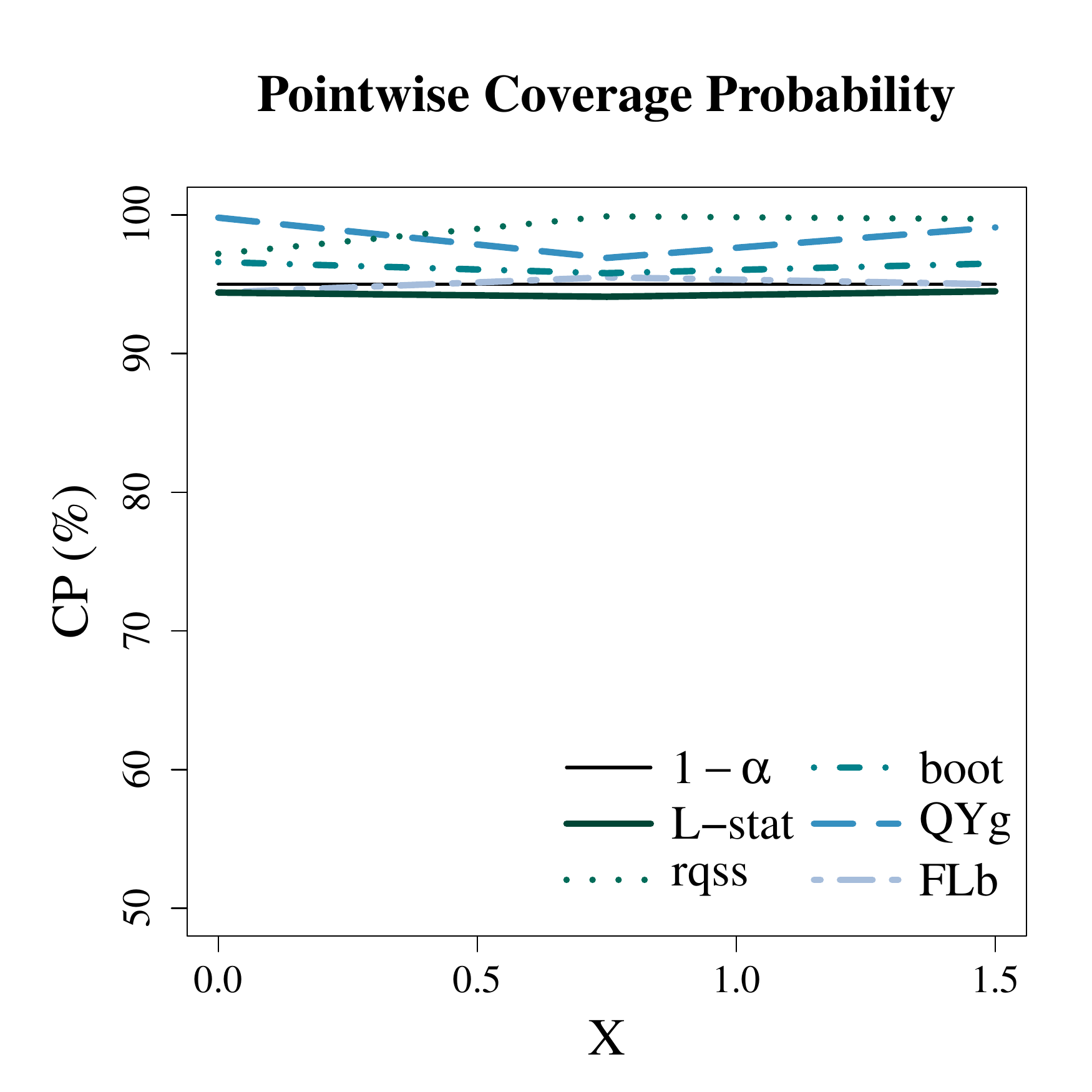}
    \hfill
  \includegraphics[clip=true,trim=15 15 10 58,width=0.445\textwidth]
    {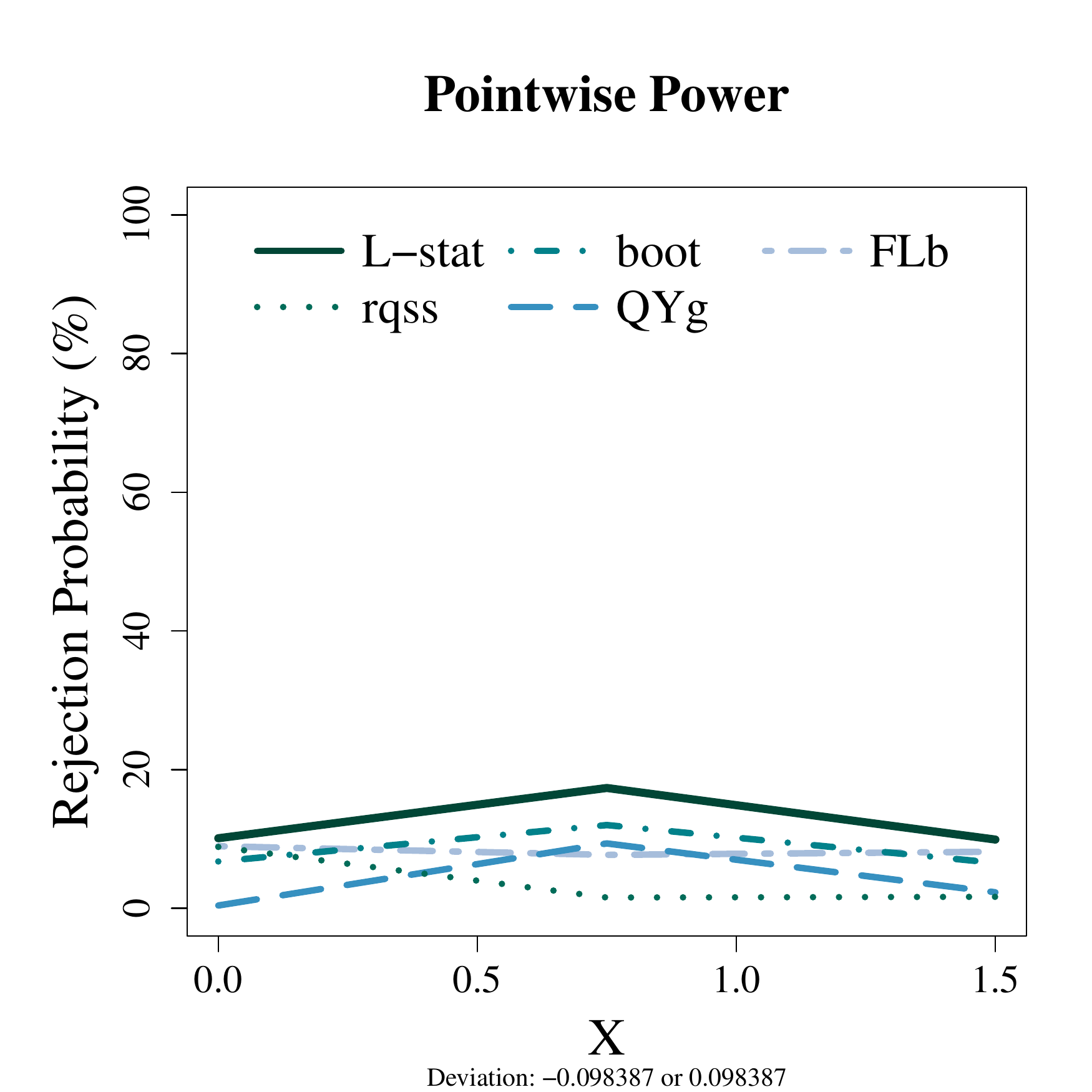}
\hfill
  \null
\\
\hfill
  \includegraphics[clip=true,trim=15 15 10 58,width=0.445\textwidth]
    {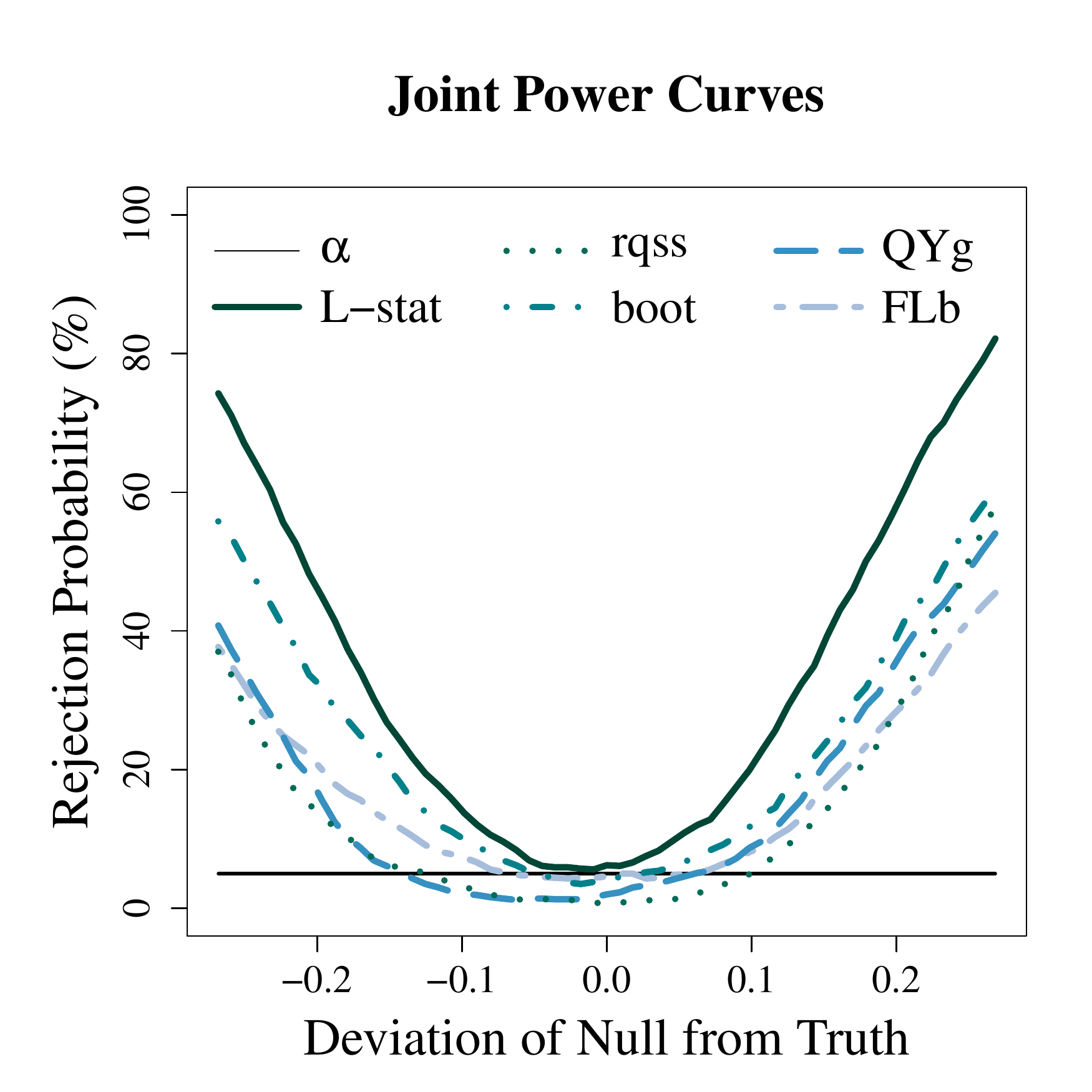}
    \hfill
  \includegraphics[clip=true,trim=15 15 10 58,width=0.445\textwidth]
    {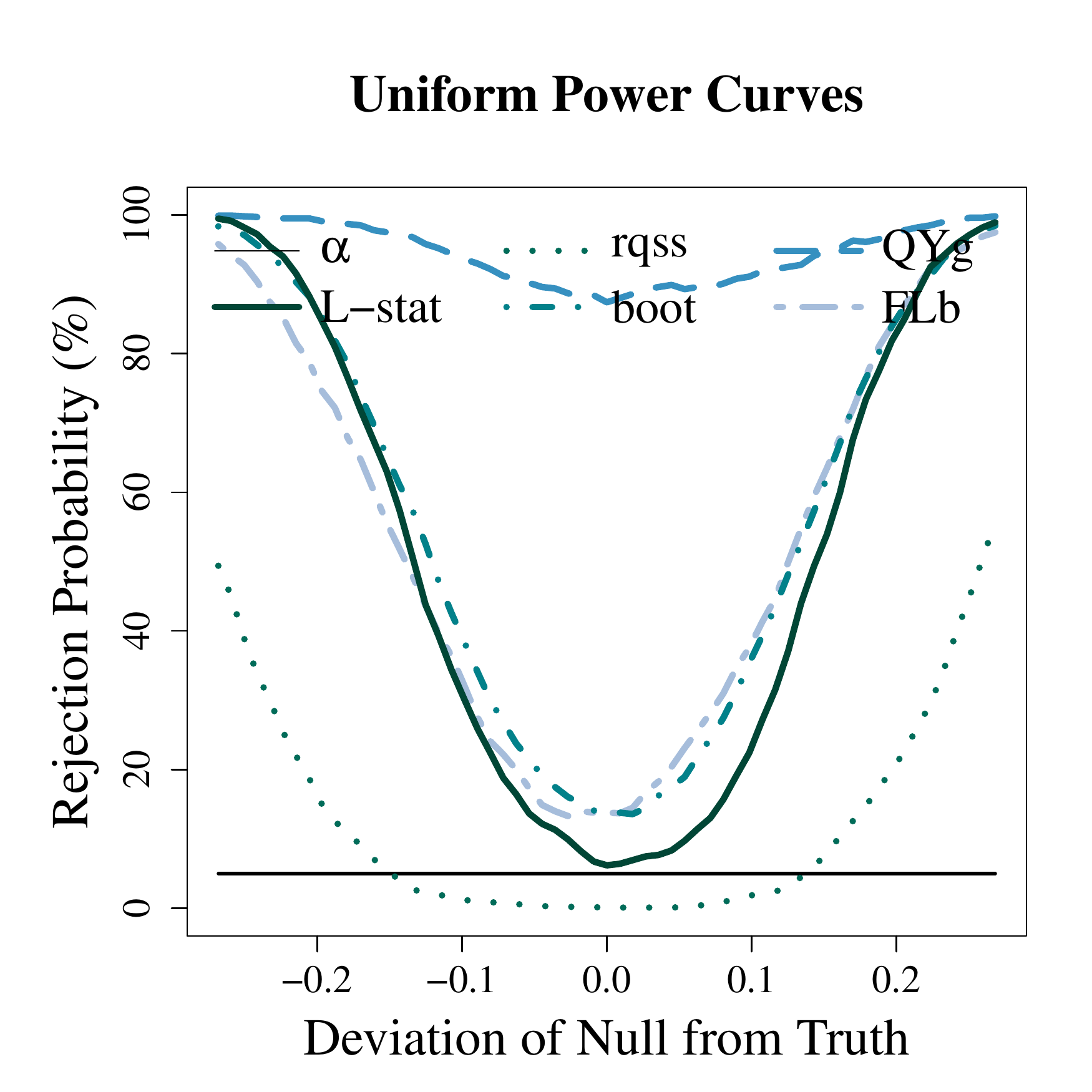}
\hfill
\null
  \caption{\label{fig:FLsim}Results from DGP in Model 1 of \citet{FanLiu2016}, $n=500$, $p=0.5$.  Top left: pointwise CP at $x_0\in\{0,0.75,1.5\}$, interpolated for visual ease.  Top right: pointwise power at the same $x_0$ against deviations of $\pm0.1$.  Bottom left: joint power curves.  Bottom right: uniform power curves.}
\end{figure}

Figure \ref{fig:FLsim} uses Model 1 from \citet[p.\ 205]{FanLiu2016}: 
$Y_i=2.5+\sin(2X_i)+2\exp\bigl(-16X_i^2\bigr)+0.5\epsilon_i$, 
$X_i\stackrel{iid}{\sim}N(0,1)$, 
$\epsilon_i\stackrel{iid}{\sim}N(0,1)$, 
$X_i\independent\epsilon_i$, 
$n=500$, $p=0.5$.  
The ``Direct'' method in their Table 1 is our FLb. 
All methods have good pointwise CP (top left). 
L-stat has the best pointwise power (top right). 

Figure \ref{fig:FLsim} (bottom left) shows power curves of the hypothesis tests corresponding to the joint (over $x_0\in\{0,0.75,1.5\}$) CIs, varying $H_0$ while maintaining the same DGP.  The deviations of $Q_{Y|X}(p;x_0)$ shown on the horizontal axis are the same at each $x_0$; zero deviation implies $H_0$ is true, in which case the rejection probability is the type I error rate.  
All methods have good type I error rates: L-stat's is 6.2\%, and other methods' are below the nominal 5\%.  L-stat has significantly better power, an advantage of 20--40\% at the larger deviations.  
The bottom right graph in Figure \ref{fig:FLsim} is similar, but based on uniform confidence bands evaluated at $231$ different $x_0$.  Only L-stat has nearly exact type I error rate and good power.

Next, we use the simulation setup of the \texttt{rqss} vignette in \citet{R.quantreg}, which in turn came in part from \citet[\S17.5.1]{RuppertEtAl2003}.  Here, $n=400$, $p=0.5$, $d=1$, $\alpha=0.05$, and
\begin{equation}\label{eqn:DGP-rqss}
X_i \stackrel{iid}{\sim}\textrm{Unif}(0,1), \quad
Y_i = \sqrt{X_i(1-X_i)}\sin\bigl(2\pi(1+2^{-7/5})/(X_i+2^{-7/5})\bigr) +\sigma(X_i)U_i,
\end{equation}
where the $U_i$ are iid $N(0,1)$, $t_3$, Cauchy, or centered $\chi^2_3$, and $\sigma(X)=0.2$ or $\sigma(X)=0.2(1+X)$.  
The conditional median function is graphed in the supplemental appendix.  
Although the function as a whole is not a common shape in economics (with multiple local maxima and minima), it provides insight into different types of functions at different points.  
For pointwise and joint CIs, we consider 47 equispaced points, $x_0=0.04,0.06,\ldots,0.96$; uniform confidence bands are evaluated at 231 equispaced values of $x_0$.

\begin{figure}[htbp]
  \centering 
  \hfill
  \includegraphics[clip=true,trim=20 45 30 25,width=0.32\textwidth]
    {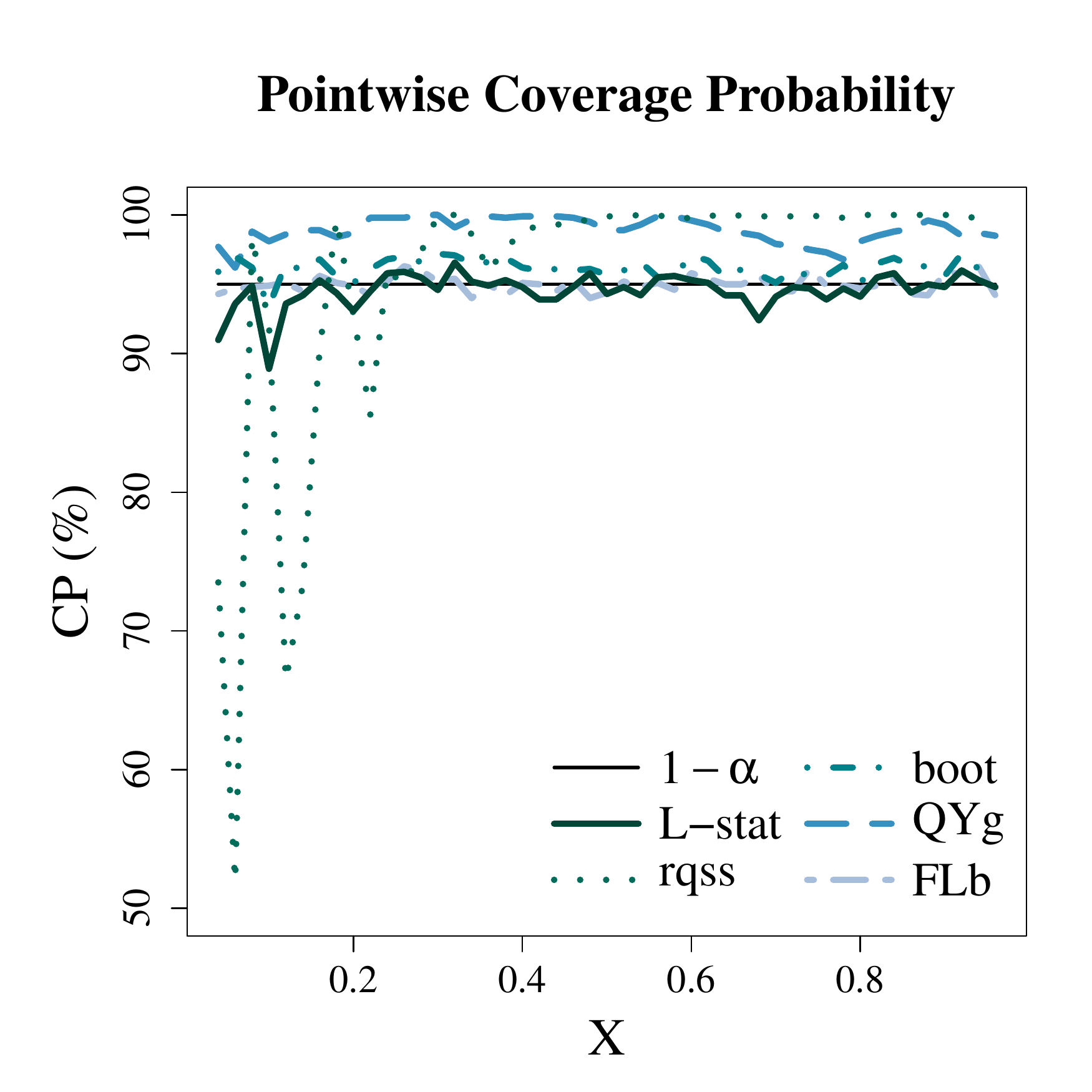}
  \hfill
  \includegraphics[clip=true,trim=51 45 30 25,width=0.30\textwidth]
    {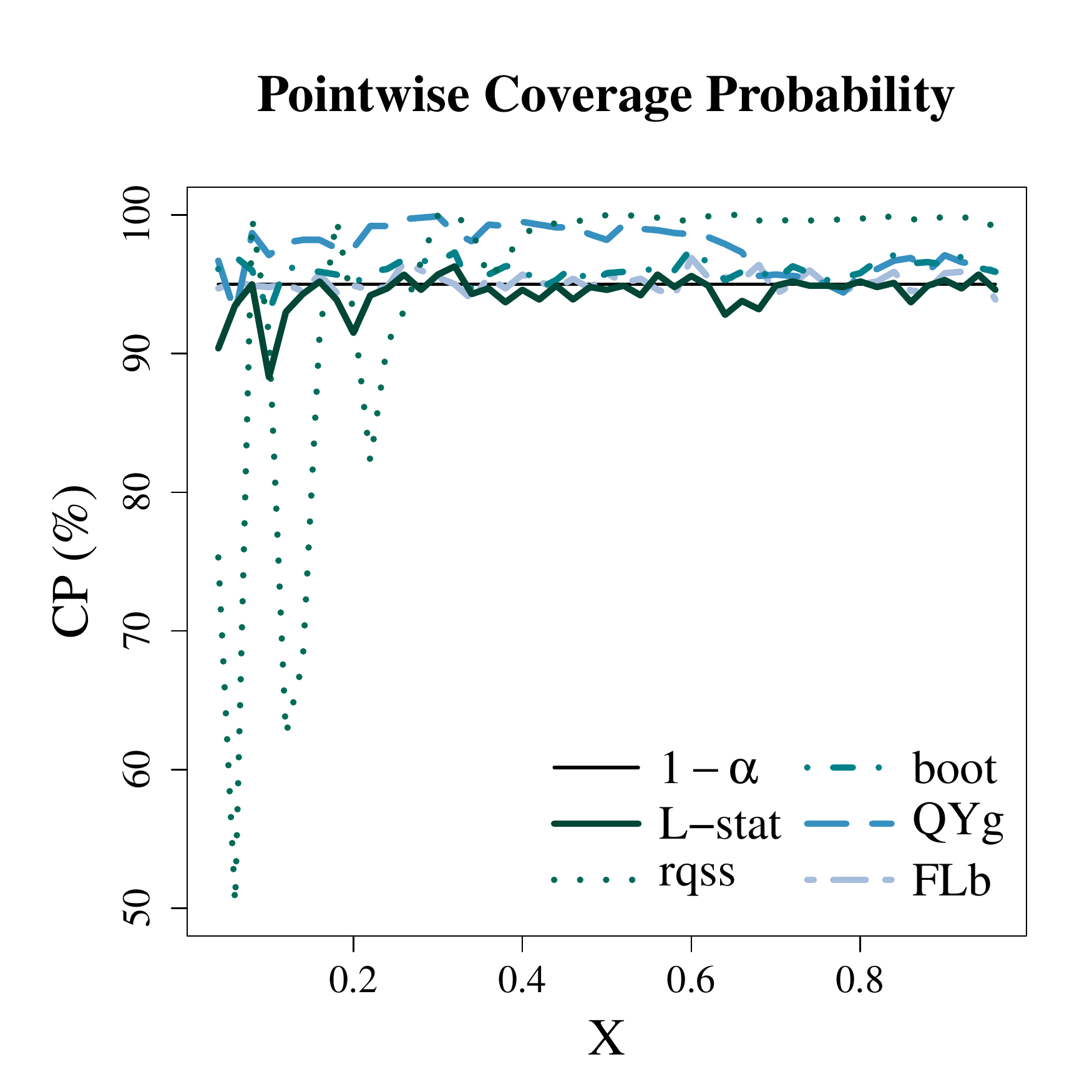}
  \hfill
  \includegraphics[clip=true,trim=10 45 30 25,width=0.329\textwidth]
    {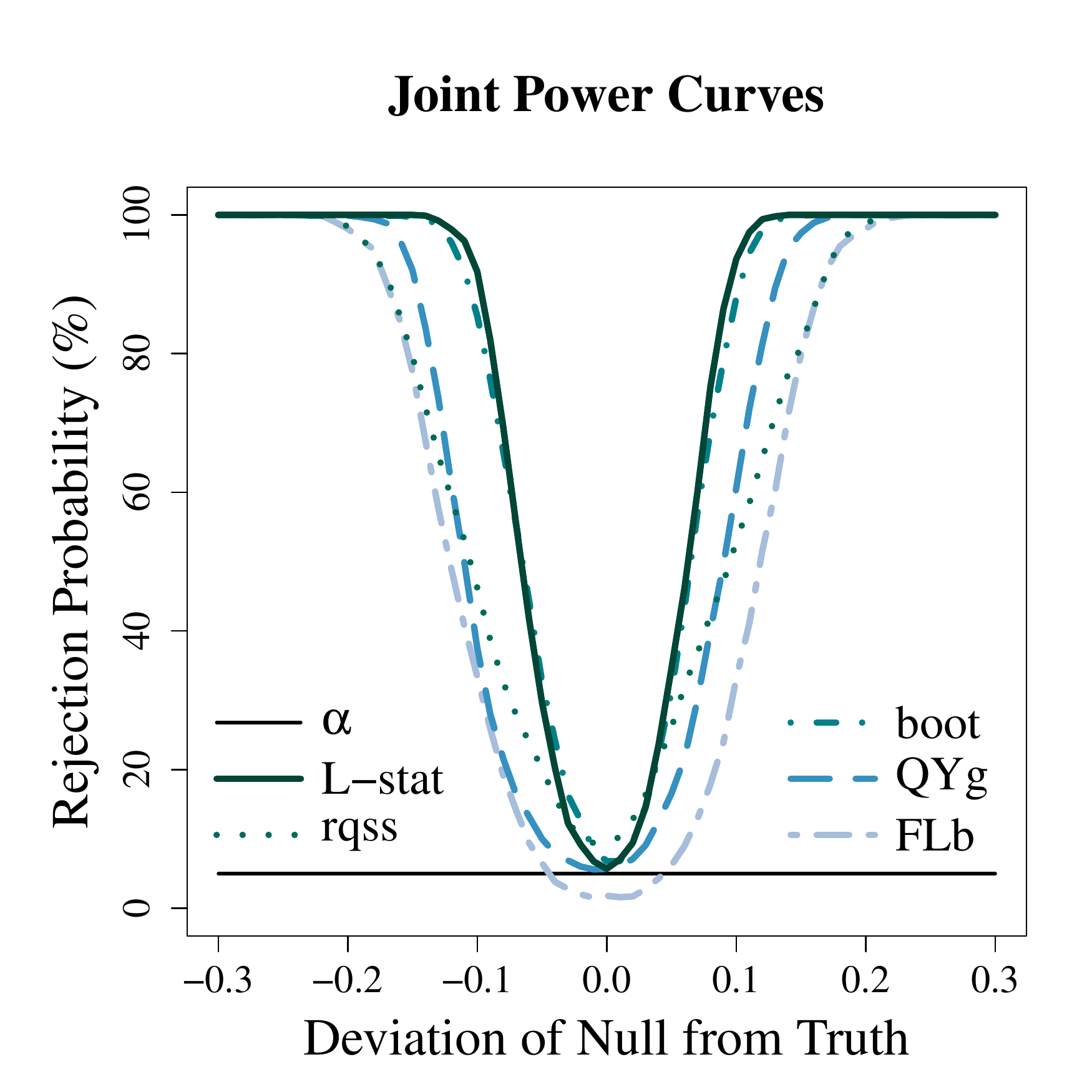}
  \hfill\null
\\
  \hfill
  \includegraphics[clip=true,trim=20 45 30 70,width=0.32\textwidth]
    {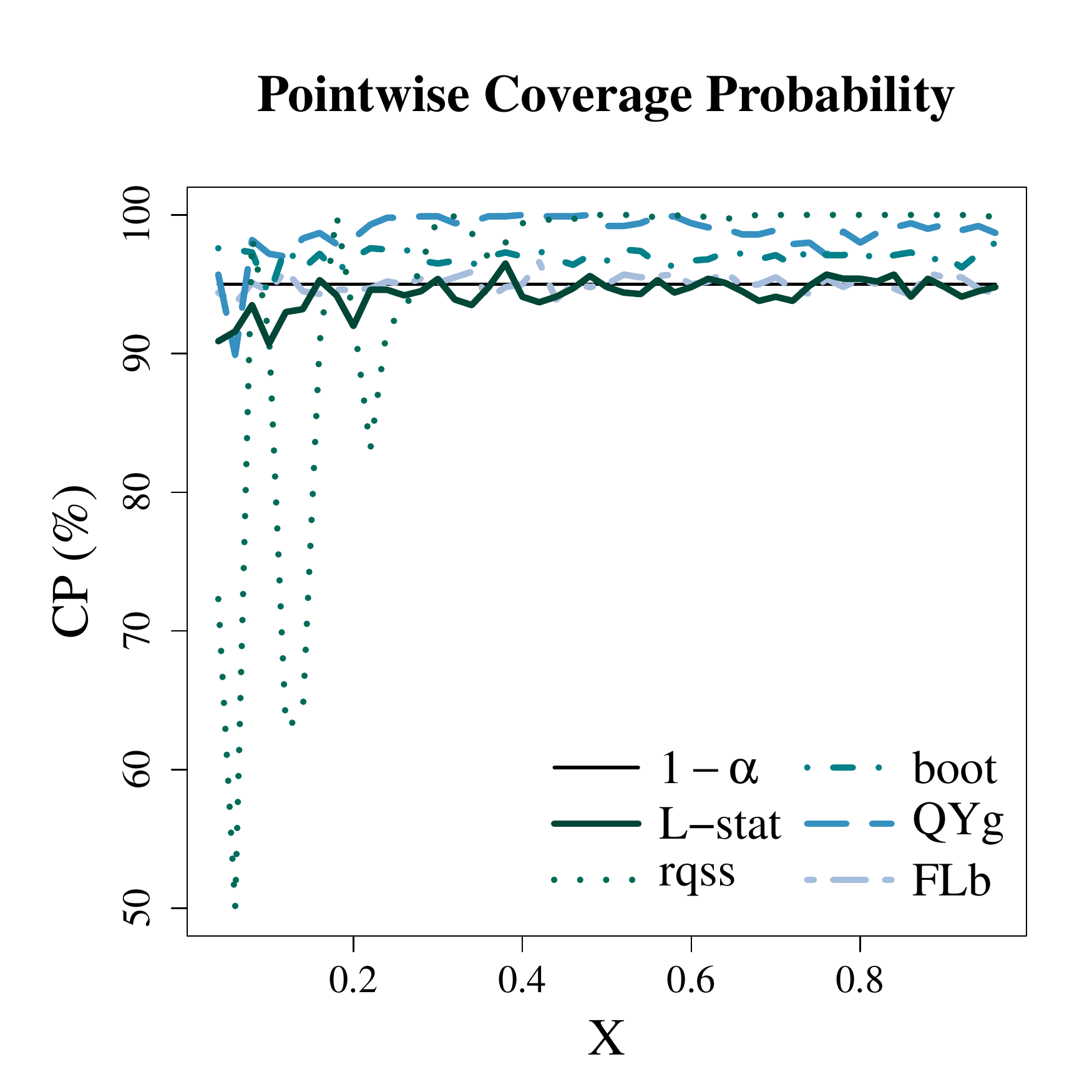}
  \hfill
  \includegraphics[clip=true,trim=51 45 30 70,width=0.30\textwidth]
    {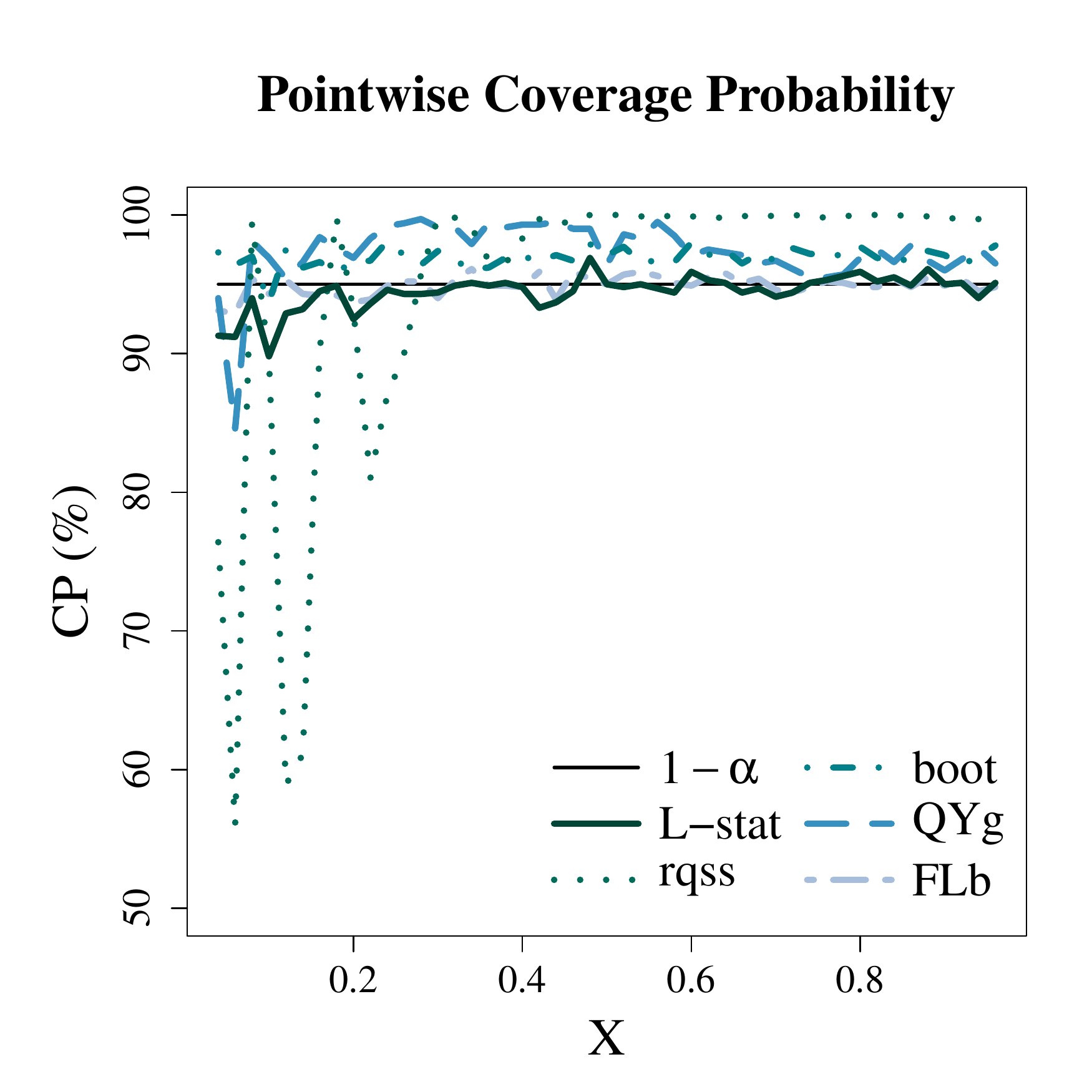}
  \hfill
  \includegraphics[clip=true,trim=10 45 30 70,width=0.329\textwidth]
    {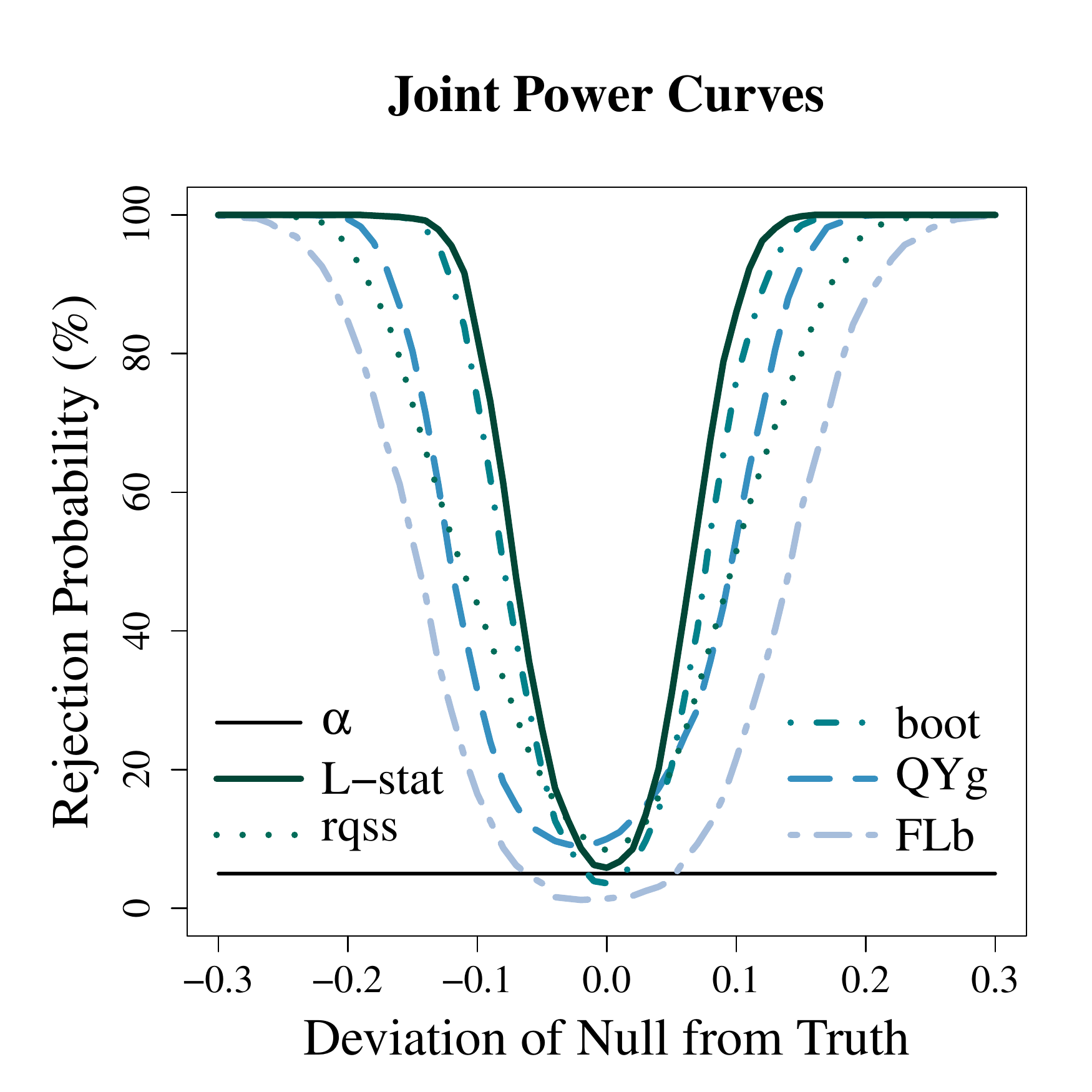}
  \hfill\null
\\
  \hfill
  \includegraphics[clip=true,trim=20 45 30 70,width=0.32\textwidth]
    {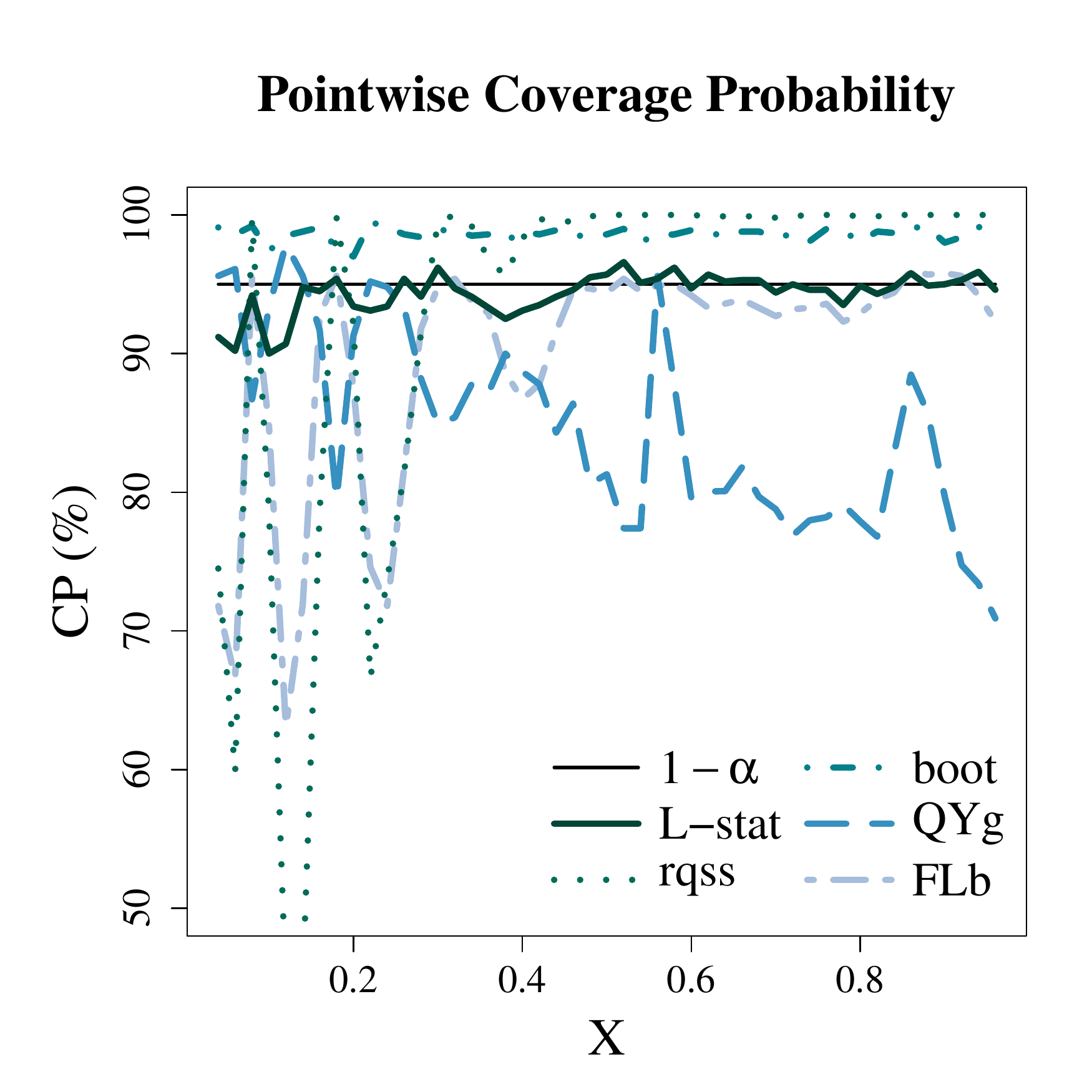}
  \hfill
  \includegraphics[clip=true,trim=51 45 30 70,width=0.30\textwidth]
    {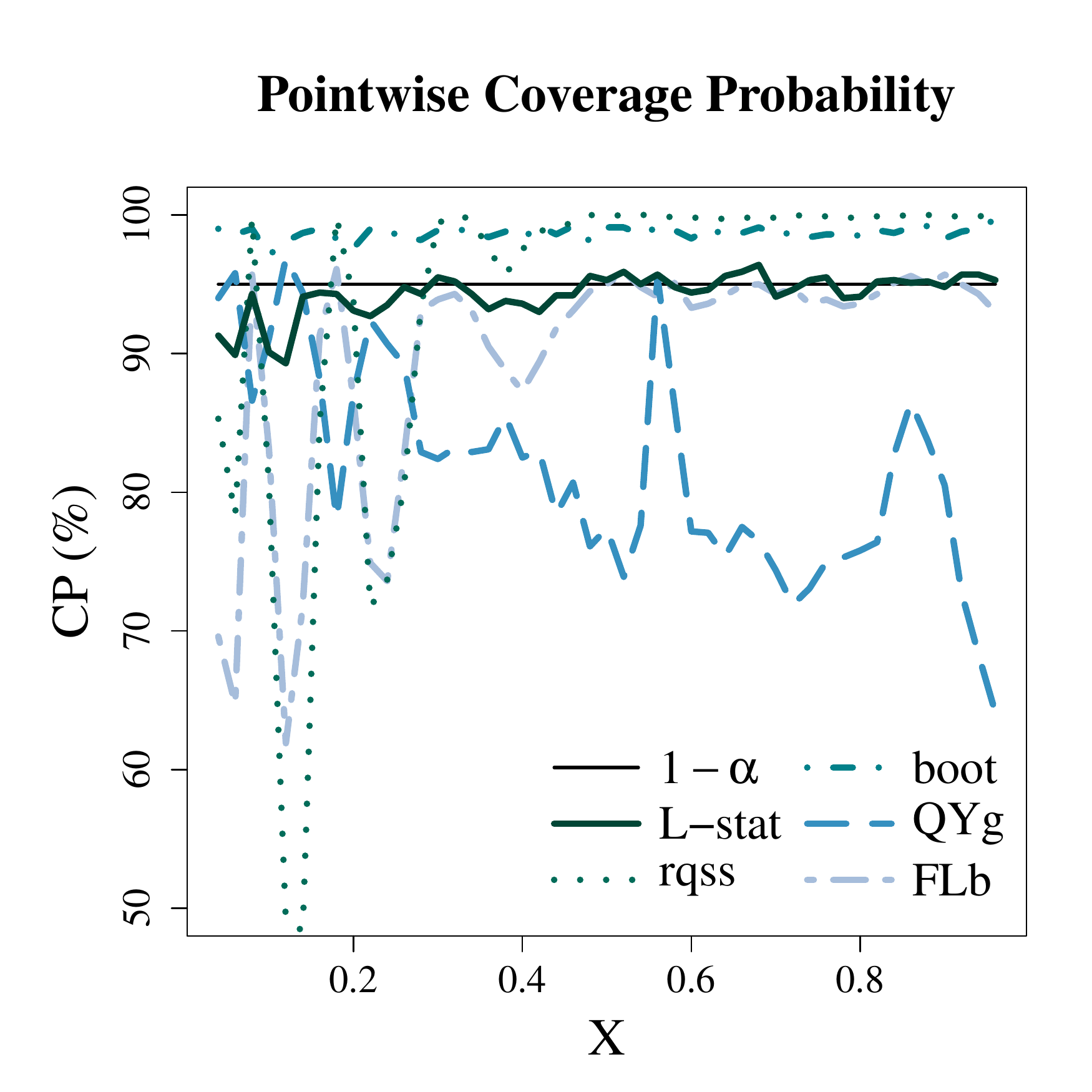}
  \hfill
  \includegraphics[clip=true,trim=10 45 30 70,width=0.329\textwidth]
    {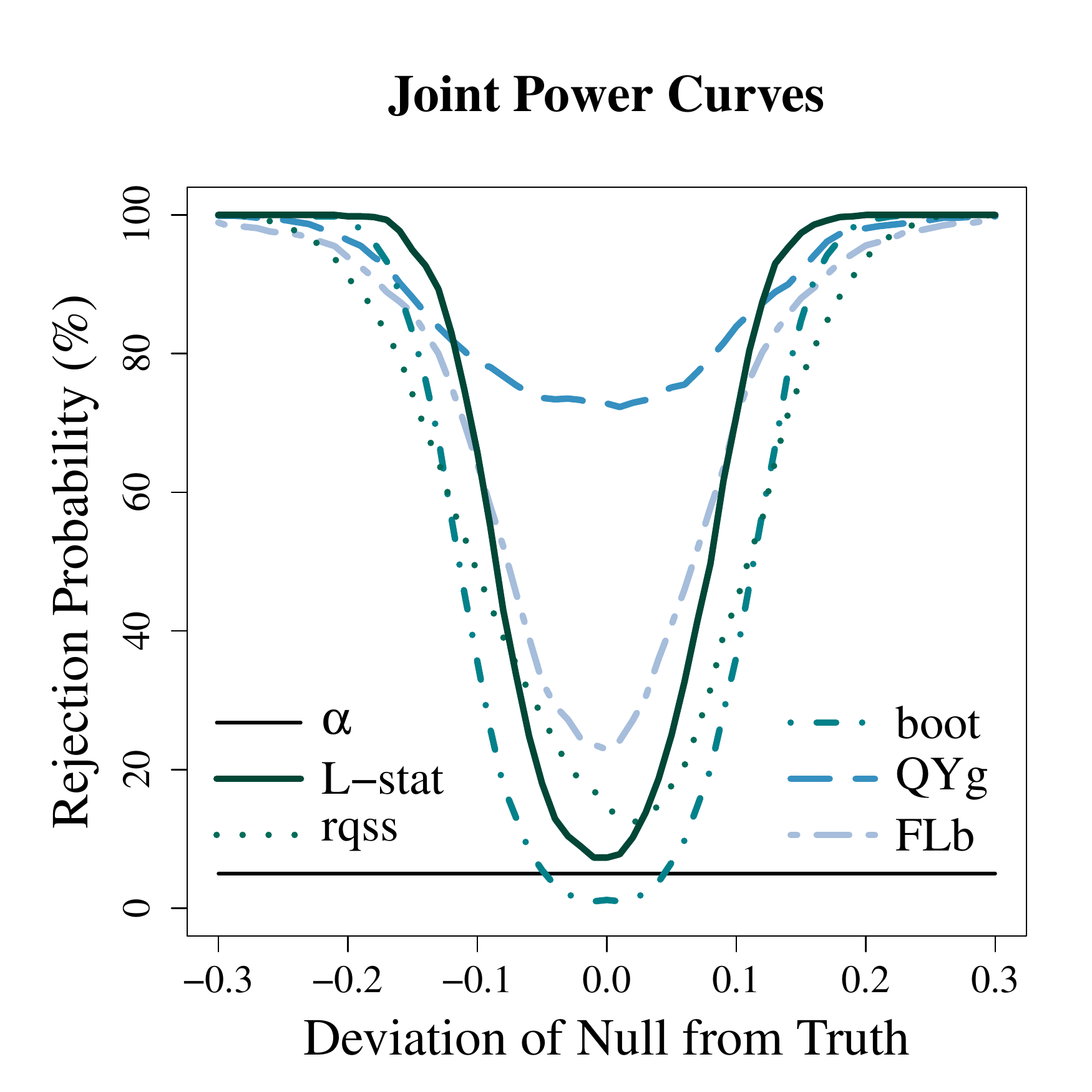}
  \hfill\null
\\
  \hfill
  \includegraphics[clip=true,trim=20 15 30 70,width=0.32\textwidth]
    {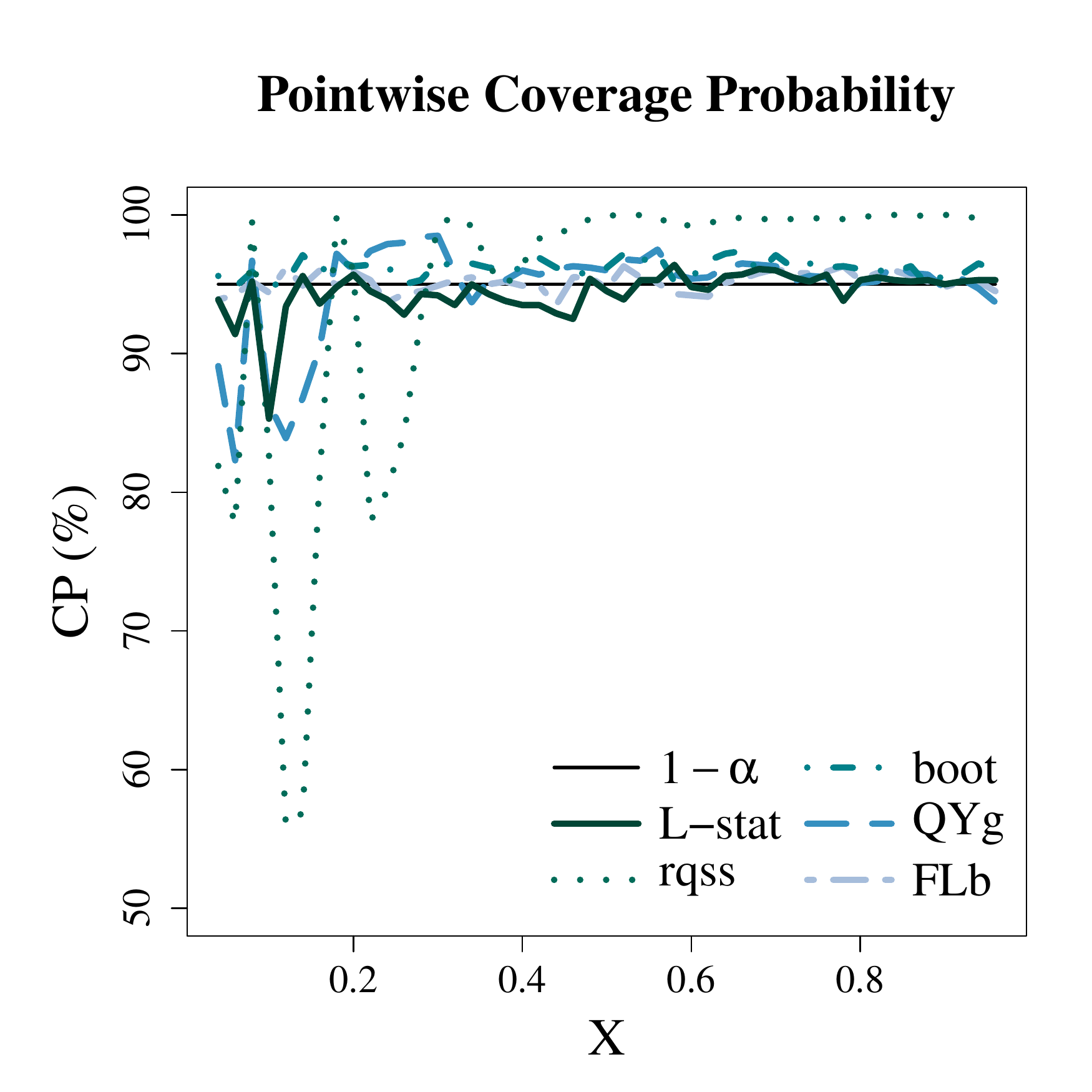}
  \hfill
  \includegraphics[clip=true,trim=51 15 30 70,width=0.30\textwidth]
    {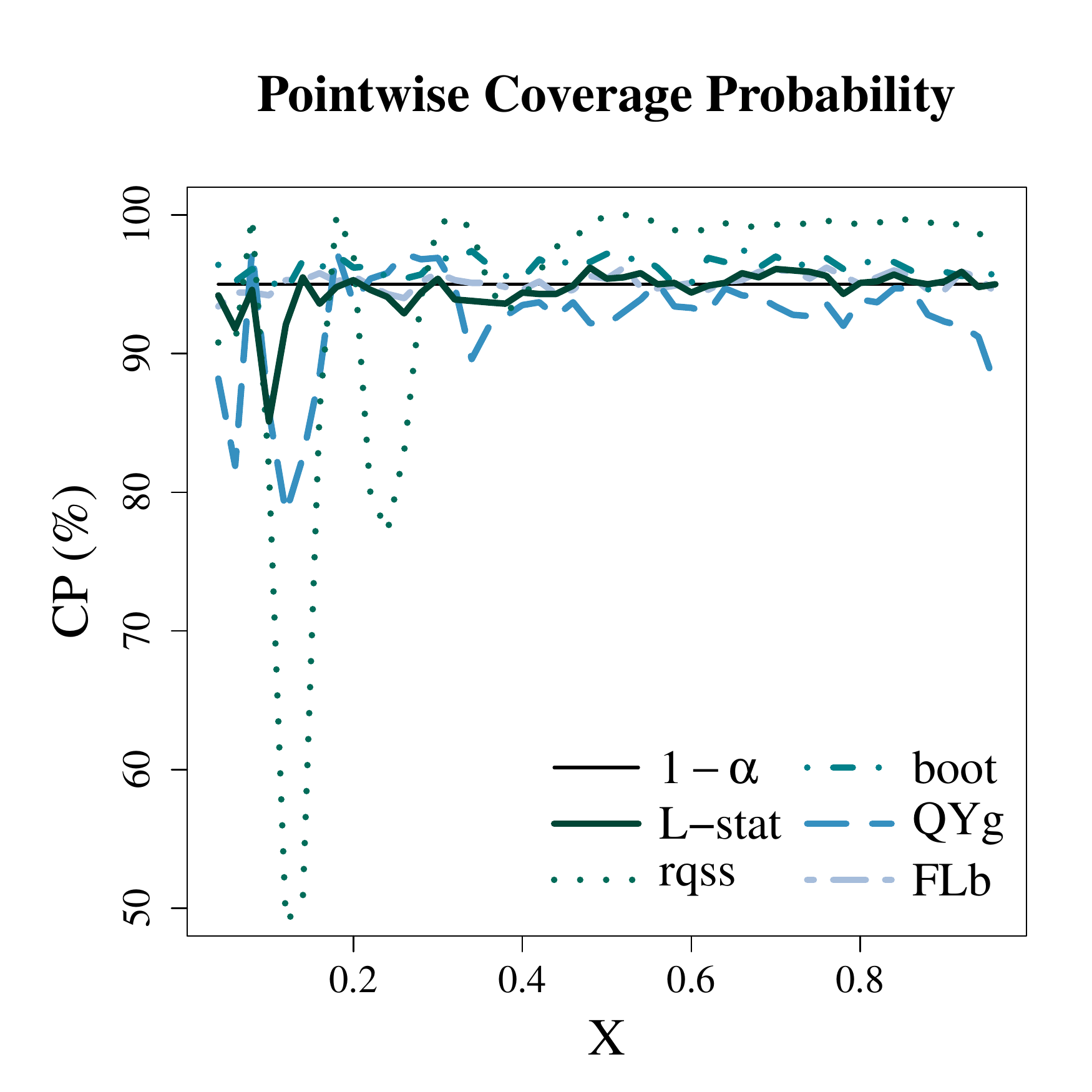}
  \hfill
  \includegraphics[clip=true,trim=10 15 30 70,width=0.329\textwidth]
    {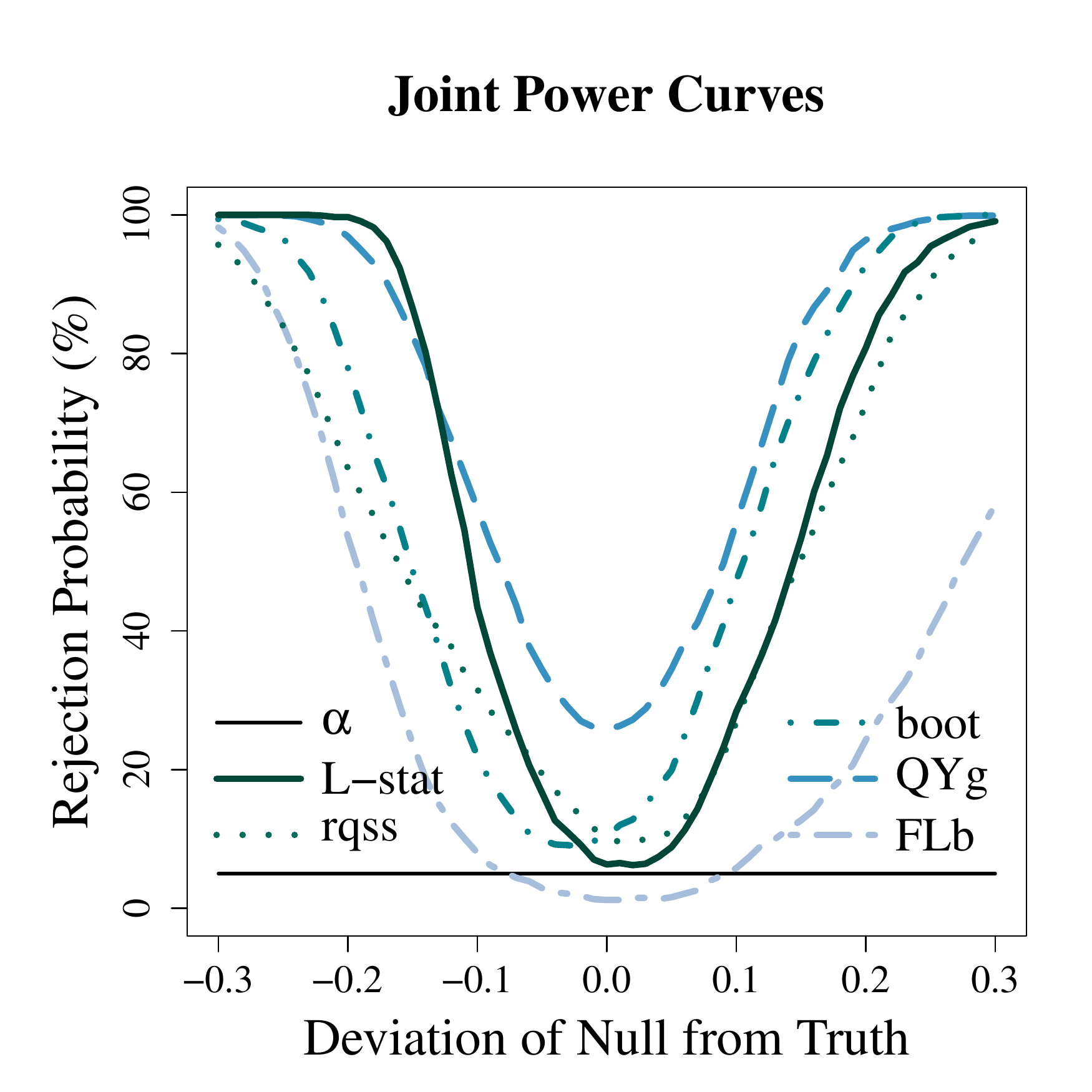}
  \hfill\null
  \caption{\label{fig:ptCP-1}Pointwise CP (first two columns) and joint power curves (third column), $1-\alpha=0.95$, $n=400$, $p=0.5$, DGP in \eqref{eqn:DGP-rqss}.  Distributions of $U_i$ are, top row to bottom row: $N(0,1)$, $t_3$, Cauchy, and centered $\chi^2_3$.  Columns 1 \& 3: $\sigma(x)=0.2$; Column 2: $\sigma(x)=(0.2)(1+x)$.}
\end{figure}

Figure \ref{fig:ptCP-1}'s first two columns show that across all eight DGPs (four error distributions, homoskedastic or heteroskedastic), L-stat has consistently accurate pointwise CP.  
At the most challenging points (smallest $x_0$), L-stat can under-cover by around five percentage points.  
Otherwise, CP is near $1-\alpha$ for all $x_0$ in all DGPs. 

In contrast, with the exception of boot, the other methods can have significant under-coverage. 
As seen in the first two columns of Figure \ref{fig:ptCP-1}, rqss has under-coverage (as low as 50--60\% CP) for $x_0$ closer to zero. 
QYg has under-coverage with the $\chi^2_3$ and (especially) Cauchy.  
FLb has good CP except with the Cauchy, where CP can dip below 70\%. 

Figure \ref{fig:ptCP-1}'s third column shows the joint power curves.  
The horizontal axis of the graphs indicates the deviation of $H_0$ from the true values.  For example, letting $\xi_{p,j}$ be the true conditional quantiles at the $j=1,\ldots,47$ values of $x_0$ (say, $x_j$), $-0.1$ deviation refers to $H_0:\{Q_{Y|X}(p;x_j)=\xi_{p,j}-0.1\textrm{ for }j=1,\ldots,47\}$ (which is false), and zero deviation means $H_0$ is true.  
Our method's type I error rate is close to $\alpha$ under all four $U_i$ distributions (5.7\%, 5.8\%, 7.3\%, 6.3\%).  
In contrast, other methods show size distortion under Cauchy and/or $\chi^2_3$ $U_i$; among them, boot is closest but still has $10.3\%$ type I error rate with the $\chi^2_3$. 
Next-best is rqss; size distortion for FLb and QYg is more serious. 
L-stat also has the steepest joint power curves among all methods. 
%
Beyond steepness, they are also the most robust to the underlying distribution.  
L-stat's type I error rate is near 5\% for all four distributions.  
In contrast, boot ranges from only 1.2\% for the Cauchy, leading to worse power, up to 10.3\% for the $\chi^2_3$. 

The supplemental appendix shows a comparison of hypothesis tests based on uniform confidence bands. 
The results are similar to the joint power curves, but with slightly higher rejection rates all around.  

\begin{figure}[thb]
  \centering 
  \hfill
  \includegraphics[clip=true,trim=15 15 10 58,width=0.445\textwidth]
    {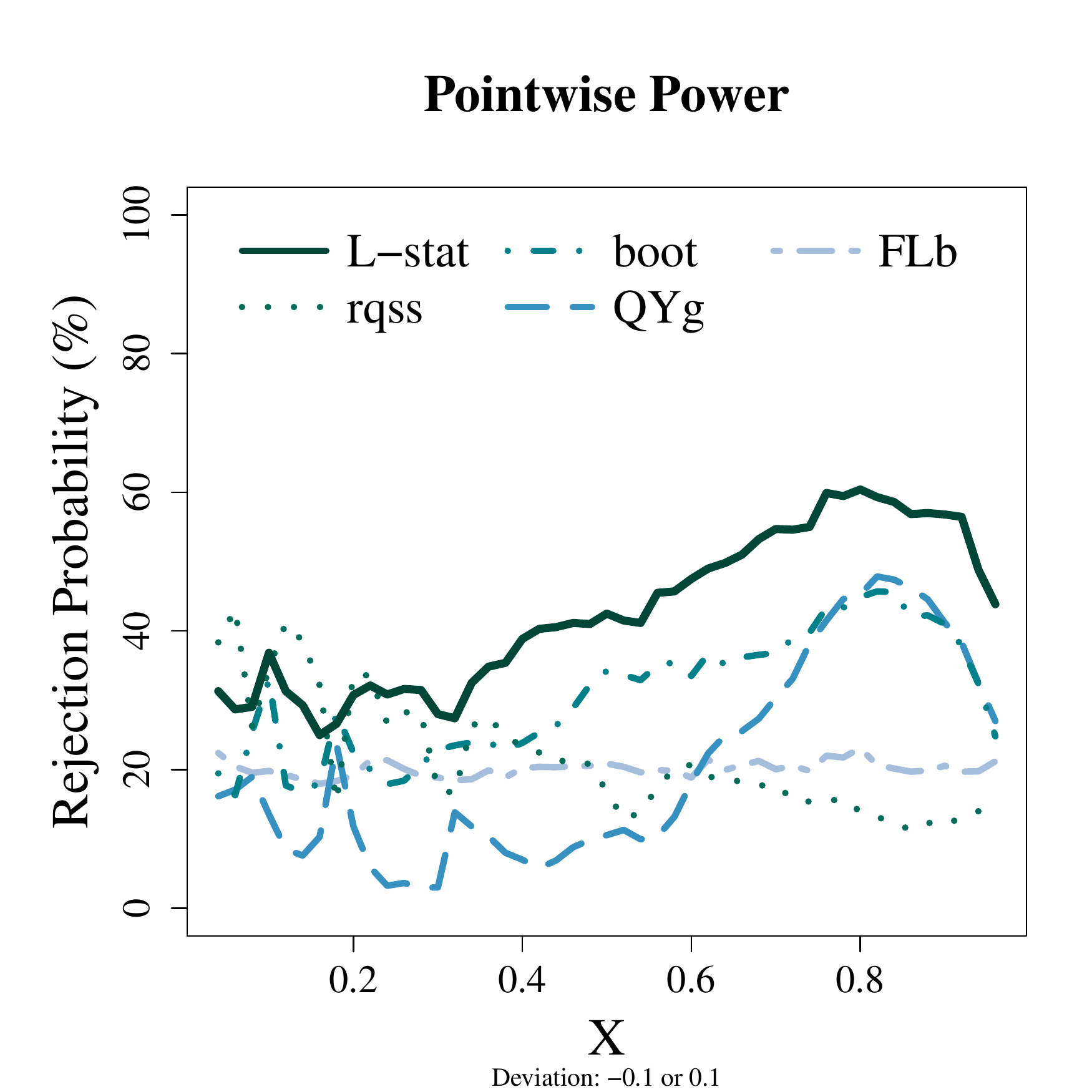}
  \hfill
  \includegraphics[clip=true,trim=51 15 10 58,width=0.406\textwidth]
    {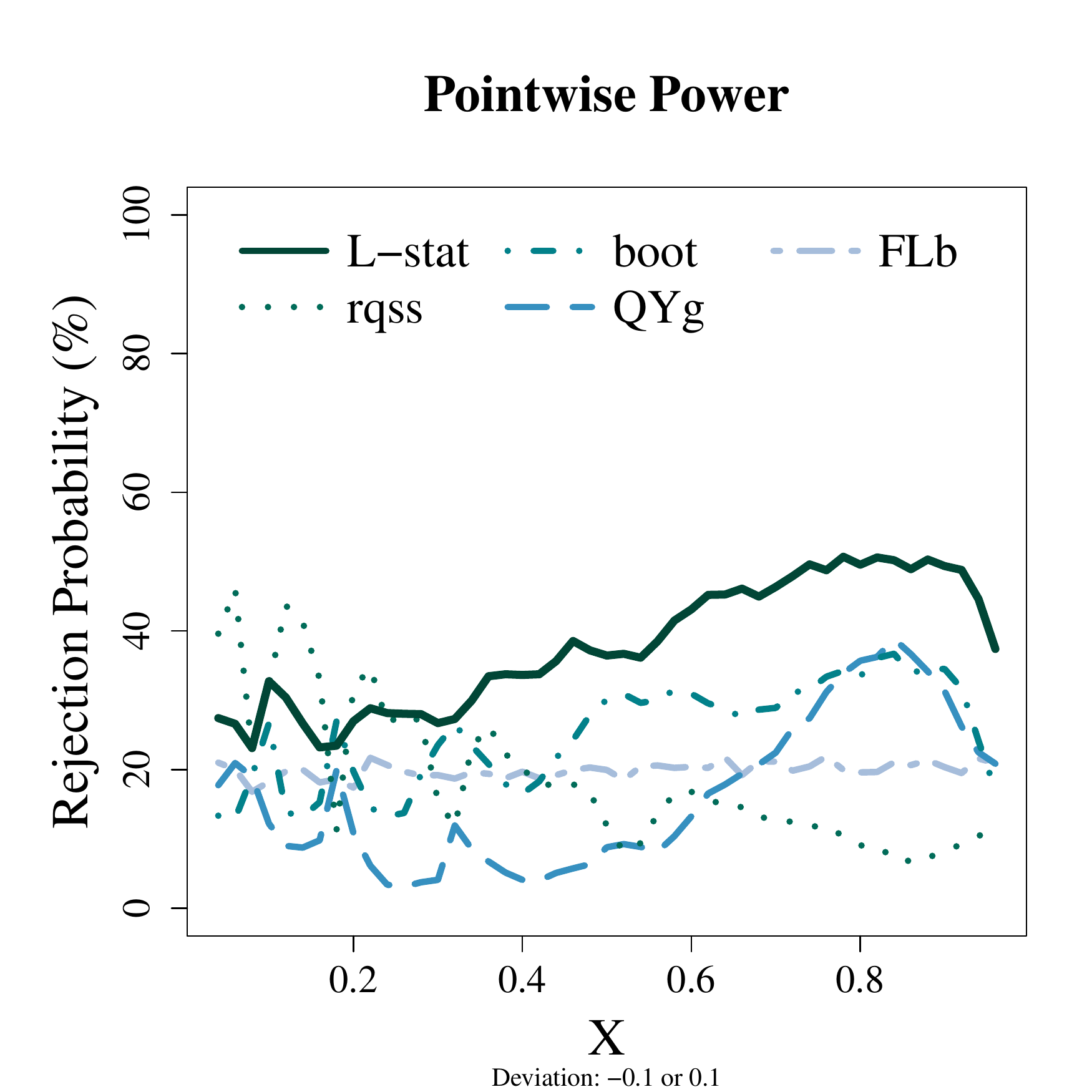}
  \hfill\null
\\
  \hfill
  \includegraphics[clip=true,trim=15 15 10 58,width=0.445\textwidth]
    {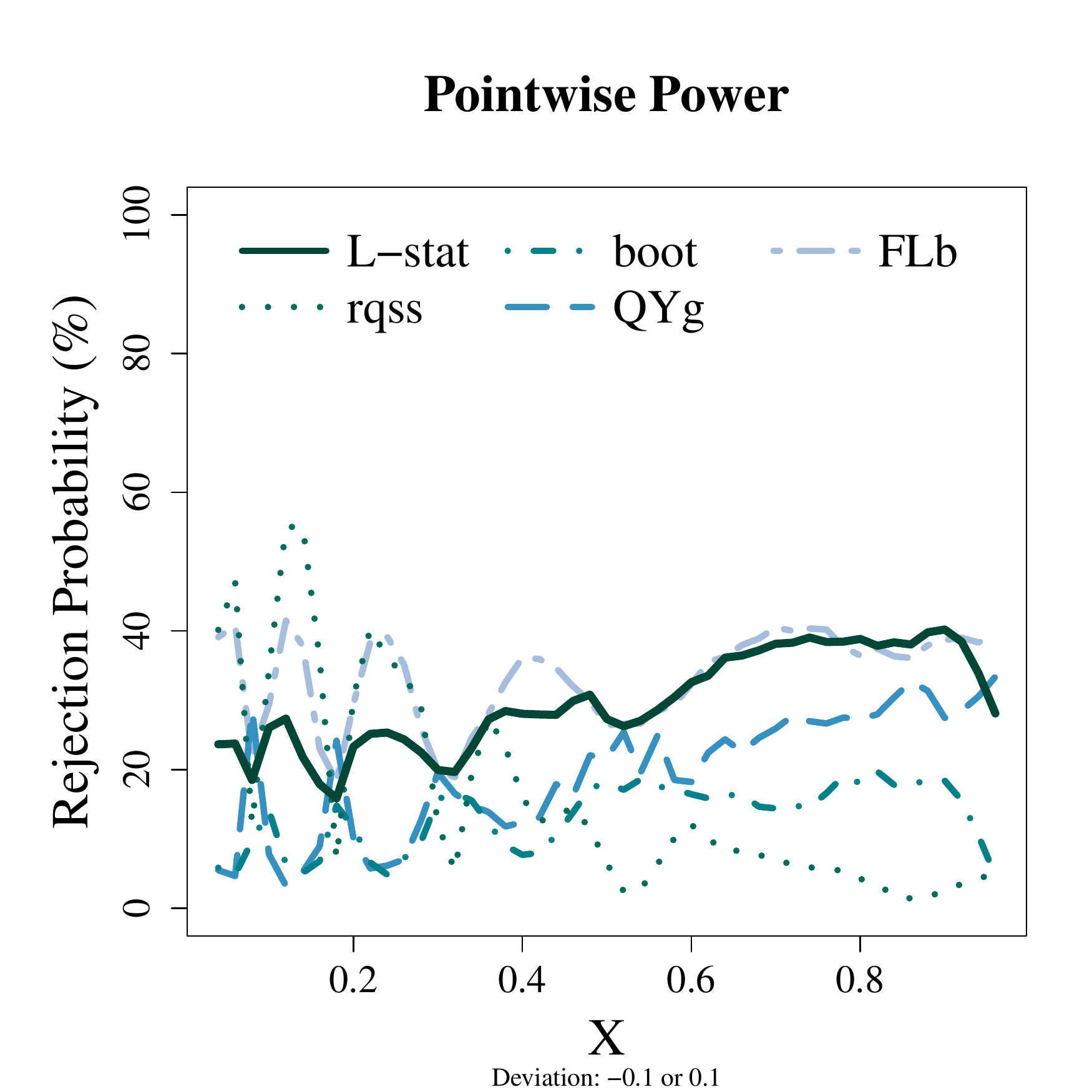}
  \hfill
  \includegraphics[clip=true,trim=51 15 10 58,width=0.406\textwidth]
    {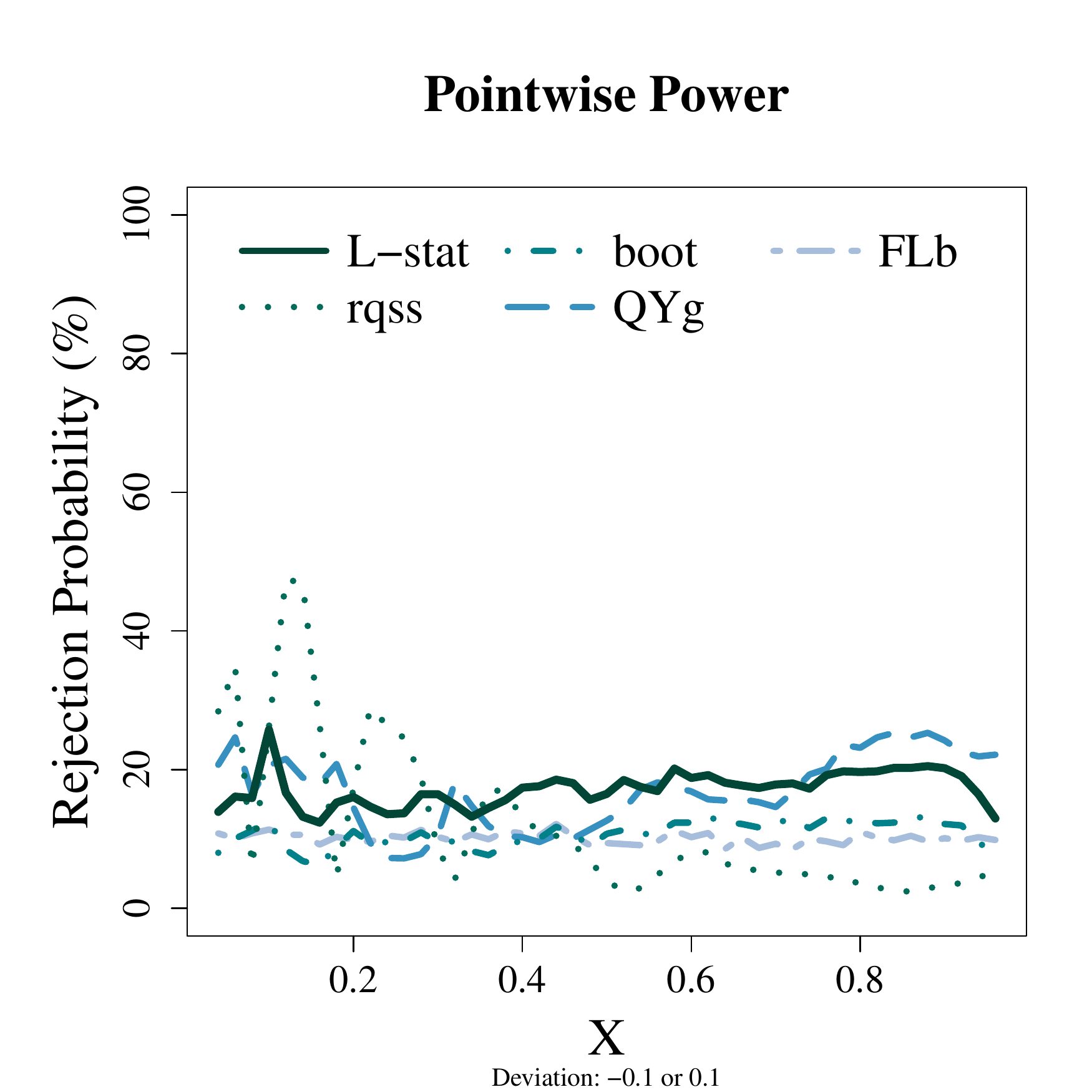}
  \hfill\null
  \caption{\label{fig:ptPWR2-1}Pointwise power (described in text), $1-\alpha=0.95$, $n=400$, $p=0.5$, DGP from \eqref{eqn:DGP-rqss}, $\sigma(x)=0.2$. The $U_i$ are $N(0,1)$ (top left), $t_3$ (top right), Cauchy (bottom left), and centered $\chi^2_3$ (bottom right).}
\end{figure}

Figure \ref{fig:ptPWR2-1} shows pointwise power.  Specifically, for a given $x_0$, this is the proportion of simulation draws in which $Q_{Y|X}(p;x_0)-0.1$ is excluded from the CI, averaged with the corresponding proportion for $Q_{Y|X}(p;x_0)+0.1$. 
L-stat generally has the best power among methods with correct CP (per first column of Figure \ref{fig:ptCP-1}).  

The supplemental appendix contains results for $p=0.25$, where L-stat continues to perform well.  
One additional advantage is that L-stat's joint test is nearly unbiased, whereas the other joint tests are all biased.  

The supplemental appendix also shows the computational advantage of our method.  For example, with $n=10^5$ and $100$ different $x_0$, L-stat takes only $10$ seconds, whereas the local cubic bootstrap takes $141$ seconds; rqss is even slower. 

Overall, the simulation results show the new L-stat method to be fast and accurate. 
Besides L-stat, the only method to avoid serious under-coverage is the local cubic with bootstrapped standard errors, perhaps due to its reliance on our newly proposed CPE-optimal bandwidth.  However, L-stat consistently has better power, greater robustness across different conditional distributions, and less bias of its joint hypothesis tests.

\section{Conclusion}

We derive a uniform $O(n^{-1})$ difference between the linearly interpolated and ideal fractional order statistic distributions.  
We generalize this to $L$-statistics to help justify quantile inference procedures. 
In particular, this translates to $O(n^{-1})$ CPE for the quantile CIs proposed by \citet{Hutson1999}, which we improve to $O\left(n^{-3/2}[\log(n)]^3\right)$ via calibration. 
We extend these results to a nonparametric conditional quantile model, with both theoretical and Monte Carlo success.  The derivation of an optimal bandwidth value (not just rate) and a fast approximation thereof are important practical advantages. 
 
Our results can be extended to other objects of interest, such as interquantile ranges and two-sample quantile differences \citep{GoldmanKaplan2014b}, quantile marginal effects \citep{Kaplan2014}, and entire distributions \citep{GoldmanKaplan2015c}.  

In ongoing work, we consider the connection with Bayesian bootstrap quantile inference, which may be a way to ``relax'' the iid assumption. 
%
%
Other future work may improve finite-sample performance, e.g.,\ by smoothing over discrete covariates \citep{LiRacine2007}. 

\singlespacing
\Supplemental{}{%
\bibliographystyle{chicago}

}


\doublespacing
\onehalfspacing
\singlespacing


\Supplemental{\pagebreak\setcounter{page}{2}}{}

\appendix

\allowdisplaybreaks[4]


\section{\Supplemental{Mathematical proofs}{Proof sketches} and additional lemmas\label{sec:app-pfs}}

\Supplemental{}{The following are only sketches of proofs.  The full proofs, with additional intermediate steps and explanations, may be found in the supplemental appendix.}

\Supplemental{%
Notationally, a tilde as in $\tilde{x}$ indicates a value determined by the mean value theorem (usually in the ``Lagrange form'' of the remainder from Taylor's theorem), and its value may change across equations.  For example, $f(x)=f(x_0)+(x-x_0)f'(\tilde{x})$, with $\tilde{x}$ between $x$ and $x_0$, may be followed by $g(x)=g(x_1)+(x-x_1)g'(\tilde{x})$, with $\tilde{x}$ between $x$ and $x_1$ and different than the first $\tilde{x}$. %
}{}

\subsection*{\Supplemental{Proof}{Sketch of proof} of Proposition \ref{prop:error-prob}}

\Supplemental{\begin{proof}}{}
For any $u$, let $k=\lfloor(n+1)u\rfloor$ and $\epsilon=(n+1)u-k\in[0,1)$.  If $\epsilon=0$, then the objects $\tilde Q^I_X(u)$, $\hat Q^L_X(u)$, and $F^{-1}\bigl(\tilde Q^L_U(u)\bigr)$ are identical and equal to $X_{n:k}$.  Otherwise, each lies in between $X_{n:k}$ and $X_{n:k+1}$ due to monotonicity of the quantile function and $k/(n+1)\le u<(k+1)/(n+1)$:
\begin{align*}
X_{n:k} &= \tilde Q^I_X\left(k/(n+1)\right) \le \tilde Q^I_X(u) \le \tilde Q^I_X\left((k+1)/(n+1)\right) = X_{n:k+1}, \\
X_{n:k} &\le \hat Q^L_X\left(u\right) = (1-\epsilon)X_{n:k} + \epsilon X_{n:k+1} \le X_{n:k+1}, \\
X_{n:k} &= F^{-1}\Bigl(\tilde Q^L_U\bigl(k/(n+1)\bigr)\Bigr) 
         \le F^{-1}\Bigl(\tilde Q^L_U\left(u\right)\Bigr)
         \le F^{-1}\Bigl(\tilde Q^L_U\bigl((k+1)/(n+1)\bigr)\Bigr) = X_{n:k+1} .
\end{align*}
\Supplemental{%
Thus, $\left|\tilde Q^I_X(u)-\hat Q^L_X(u)\right|\le X_{n:k+1}-X_{n:k}$ and $\left|\hat Q^L_X(u)-F^{-1}\left(\tilde Q^L_U(u)\right)\right|\le X_{n:k+1}-X_{n:k}$, so 
\begin{align}\label{eqn:thm1pf-Q-sup-max}
\sup_{u\in[1/(n+1),n/(n+1)]} \left|\tilde Q^I_X(u)-\hat Q^L_X(u)\right|
  &\le \max_{k\in\{1,\ldots,n-1\}} X_{n:k+1}-X_{n:k} ,
\end{align}
and similarly for $\left|\hat Q^L_X(u)-F^{-1}\left(\tilde Q^L_U(u)\right)\right|$. %
}{%
Thus, differences between the processes can be bounded by the maximum (over $k$) spacing $X_{n:k+1}-X_{n:k}$.  Using the assumption that the density is uniformly bounded away from zero over the interval of interest (and applying a maximal inequality from \citet[eqn.\ (3.7)]{Bickel1967}), this in turn can be bounded by a maximum of uniform order statistic spacings $U_{n:k+1}-U_{n:k}$. %
}%
\Supplemental{%

Taking a mean value expansion,
\begin{align}\label{eqn:thm1pf-MVT}
X_{n:k+1}-X_{n:k}
  &= F^{-1}(U_{n:k+1}) - F^{-1}(U_{n:k})
   = (U_{n:k+1}-U_{n:k}) / f\left(F^{-1}(\tilde u_k)\right) ,
\end{align}
where $\tilde u_k\in[U_{n:k},U_{n:k+1}]$.  

From (3.7) in \citet{Bickel1967}, for arbitrarily small $\eta>0$, 
\begin{align}\label{eqn:thm1pf-Bickel}
P\left[ \max_{k=1,\ldots,n} \left|U_{n:k}-k/(n+1)\right| \ge \eta\right] \to 0 ,
\end{align}
and since $(k+1)/(n+1)-k/(n+1)=1/(n+1)$, 
\begin{align*}
P\left[ \sup_{u\in[1/(n+1),n/(n+1))} \left|u-\tilde u_{\lfloor(n+1)u\rfloor}\right| \ge \eta\right] \to 0 .
\end{align*}
Combined with the definition of $\mathcal U^\delta$ in the statement of Proposition \ref{prop:error-prob},\footnote{As written, $\mathcal U^\delta$ cannot contain the whole unit interval $(0,1)$ if $X$ has an unbounded distribution.  However, the following equation could be modified to allow $f(F^{-1}(u))\to\infty$ at some rate as $u\to0$, and slowing $\mathcal U^\delta_n$ to a subset of $[n^{-a},1-n^{-a}]$ for $a<1$, at the expense of a larger right-hand side in the statement of the theorem.  We omit this since it is not helpful for our goal of high-order accurate confidence intervals.} this implies that 
\[ P\left(\max_{k=1,\ldots,n-1}f\left(F^{-1}(\tilde u_k)\right)\le\delta\right)\to1. \] 

Given $\max_{k=1,\ldots,n-1}f\left(F^{-1}(\tilde u_k)\right)\le\delta$, we can factor out $\delta$.  Using \eqref{eqn:thm1pf-Q-sup-max} and \eqref{eqn:thm1pf-MVT}, and since $\mathcal U^\delta_n\subseteq[1/(n+1),n/(n+1)]$, 
\begin{align}\label{eqn:thm1pf-central}
\sup_{u\in\mathcal U^\delta_n} \left|\tilde Q^I_X(u)-\hat Q^L_X(u)\right|
  &\le \max_{k\in\{1,\ldots,n-1\}} X_{n:k+1}-X_{n:k} 
  &\le \delta^{-1} \max_{k\in\{1,\ldots,n-1\}} (U_{n:k+1}-U_{n:k}) .
\end{align}
The marginal distribution of each uniform spacing is $(U_{n:k+1}-U_{n:k})\sim\beta(1,n)$, which is equivalent to a \citet{Kumaraswamy1980} distribution with the same parameters and CDF $F(x)=1-(1-x)^n$.  With $a_n$ such that $\lim_{n \to \infty} a_n/\left(n^{-1}[\log n]\right) = \infty$ and $a_n\to0$, %
\Supplemental{}{with $\delta$ and $\mathcal{U}_n^\delta$ as defined in the proposition and the first equality established formally in the supplemental appendix, }
\begin{align*}
P\biggl(\sup_{u \in \mathcal{U}^\delta_n}
\left|\tilde Q^I_X(u) - \hat Q^L_X(u)\right|>a_n\biggr) 
  &\le P\left( \max_{k\in\{1,\ldots,n-1\}}(U_{n:k+1}-U_{n:k}) > \delta a_n \right) \\
  &\le n P(\beta(1,n) > \delta a_n)  
   = n \max\{ (1-\delta a_n)^n,0\}  \\
  &= \max\left\{\exp\{\log(n) + n \log(1-\delta a_n)\}, 0\right\}  \\
  &\le \exp\{\log(n) -n\delta a_n\} \to 0.    \qedhere
\end{align*}
}{%
The marginal distribution $(U_{n:k+1}-U_{n:k})\sim\beta(1,n)$ can then be used to bound the probability as needed. 
}
\Supplemental{\end{proof}}{}

\subsection*{Lemma for {PDF} approximation}

\begin{lemma}\label{lem:den-i}
Let $\mathbf{\Delta k}$ be a positive $(J+1)$-vector of natural numbers such that $\sum_{j=1}^{J+1} \Delta k_j=n+1$, $\min_j\{\Delta k_j \}\to \infty$, and $\min_j\{n-\Delta k_j \}\to \infty$, and define $k_j \equiv \sum_{i=1}^j \Delta k_i$ and $\mathbf{k}\equiv(k_1,\ldots,k_J)'$.  Let $\mathbf X \equiv (X_1,\ldots,X_J)'$ be the random $J$-vector such that 
\[ \mathbf{\Delta X}\equiv (X_1,X_2-X_1,\ldots,1-X_J)' \sim \textrm{Dirichlet}(\mathbf{\Delta k}) . \]  Take any sequence $a_n$ that satisfies conditions 
  a) $a_n \to \infty$, 
  b) $a_n  n^{-1} \left[\max \{\Delta k_j\}\right]^{1/2}\to 0$, and
  c) $a_n^3 \left[\min \{\Delta k_j\}\right]^{-1/2} \to 0$.
Define Condition $\star(a_n)$ as satisfied by vector $\mathbf{x}$ if and only if
\begin{equation}\label{cond:star}
\max_{j} \left \{n \Delta k_j^{-1/2}\left| \Delta x_j - \Delta k_j/n\right|\right \} \leq a_n . \tag*{Condition $\star(a_n)$}
\end{equation}
Let $\|\mathbf{v}\|_\infty \equiv \max_{j\in\{1,\ldots,k\}}|v_j|$ denote the maximum norm of vector $\mathbf{v}=(v_1,\ldots,v_k)'$. 

\begin{enumerate}\def\theenumi{\roman{enumi}}
\item\label{lem:den-i-recenter} \ref{cond:star} implies
\begin{gather}\label{eqn:star-mean}
\max_{j} \left \{n \Delta k_j^{-1/2}\left| \Delta x_j - \Delta k_j/(n+1)\right|\right \} = O(a_n) , \\
\label{eqn:star-mode}
\max_{j} \left \{n \Delta k_j^{-1/2}\left| \Delta x_j - (\Delta k_j-1)/(n-J)\right|\right \} = O(a_n) , 
\end{gather}
where $\Delta k_j/(n+1)$ and $(\Delta k_j-1)/(n-J)$ are respectively the mean and mode of $\Delta X_j$. 
\item\label{lem:den-i-spacing-density} At any point of evaluation $\mathbf{\Delta x}$ satisfying \ref{cond:star}, the log Dirichlet PDF of $\mathbf{\Delta X}$ may be uniformly approximated as
\begin{align*}
\log f_{\mathbf{\Delta X}}(\mathbf{\Delta x}) = D 
     -\frac{(n-J)^2}{2} \sum_{j=1}^{J+1} \frac{\left(\Delta x_j - \frac{\Delta k_j -1}{n-J}\right)^2}{\Delta k_j -1} 
    + R_n,
\end{align*}
\Supplemental{where the constant }{}
\[ D\equiv \frac{J }{2}\log(n/2\pi) + \frac{1}{2}\sum_{j=1}^{J+1} \log\left(\frac{n}{\Delta k_j - 1}\right) , \]
and $R_n = O(a_n^3 \|\boldsymbol{\Delta k}^{-1/2}\|_{\infty})$ uniformly (over $\mathbf{\Delta x}$).  We also have the uniform (over $\mathbf{\Delta x}$) approximations
\begin{align*}
\D{\log [f_{\mathbf{\Delta X}}(\mathbf{\Delta x})]}{\Delta x_j} 
  &= (n-J) 
    -\frac{(n-J)^2}{\Delta k_j -1}
     \left(\Delta x_j - \frac{\Delta k_j -1}{n-J}\right)
    +O(a_n^2 n \|\mathbf{\Delta k}^{-1}\|_{\infty}), \\
\D{\log [f_{\mathbf{\Delta X}} (\mathbf{x})]}{\Delta k_j/n}   
  &= -\D{\log [f_{\mathbf{\Delta X}}(\mathbf{\Delta x})]}{\Delta x_j} + O\left(a_n^2 n \|\mathbf{\Delta k}^{-1}\|_{\infty} \right) .
\end{align*}
%
\item\label{lem:den-i-ordered-density} Uniformly over all $\mathbf x\in\mathbb{R}^{J}$ satisfying \ref{cond:star}, 
\begin{align*}
\log[f_{\mathbf{X}}(\mathbf{x})]
  &= D - \frac{1}{2}(\mathbf{x} - \mathbf{k}/(n+1))'\underline{H} (\mathbf{x} - \mathbf{k}/(n+1)) + O(a_n^3  \|\mathbf{\Delta k}^{-1/2}\|_{\infty})  , \\
\D{\log [f_{\mathbf{X}}(\mathbf{x})]}{\mathbf{x}} 
  &= -\underline{H} (\mathbf{x} - \mathbf{k}/(n+1)) + O\left(a_n^2 n \|\mathbf{\Delta k}^{-1}\|_{\infty} \right) , \\
\D{\log [f_{\mathbf{X}}(\mathbf{x})]}{\mathbf{k}/(n+1)} 
  &= \underline{H} (\mathbf{x} - \mathbf{k}/(n+1)) 
     +O(a_n^2 n\|\mathbf{\Delta k}^{-1}\|_{\infty}),
\end{align*}
where the constant $D$ is the same as in part (\ref{lem:den-i-spacing-density}), and the $J\times J$ matrix $\underline{H}$ has non-zero elements only on the diagonal $\underline{H}_{j,j} = n^2\left(\Delta k_j ^{-1} + \Delta k_{j+1}^{-1}\right)$ and one off the diagonal $\underline{H}_{j,j+1} = \underline{H}_{j+1,j} = -n^2 \Delta k_{j+1}^{-1}$. 
The covariance matrix for $\mathbf{x}$, $\underline{\mathcal V}/n \equiv \underline{H}^{-1}$, has row $i$, column $j$ elements
\begin{align}
\label{eqn:varMatrix}
\underline{\mathcal V}_{i,j}=\min(k_i,k_j)(n+1-\max(k_i,k_j))/[n(n+1)],
\end{align}
connecting the above with the conventional asymptotic normality results for sample quantiles.  That is,%
\Supplemental{\footnote{Compare \eqref{eqn:fX-phi} with equation (4.7.2) in \citet{Reiss1989}, which is very similar.  His equation is more precise for points of evaluation near the mean, but it only applies to $J=1$ and a particular version of our \ref{cond:star}.  His $b_{r,n}=k/(n+1)$ matches exactly, and his $a_{r,n}^2=(\mathcal{V}/n)[n^2/(n+1)^2]$ matches up to the $1+O(1/n)$ multiplicative difference.  His condition that $|x|\le\{[r(n-r)]/n\}^{1/6}$ is like our \ref{cond:star} with $a_n=\{[k(n-k)]/n\}^{1/6}$.  If $k\propto n$, then our condition (c) on $a_n$ becomes $a_n=o(k^{1/6})$, whereas his is (essentially) $a_n=O(k^{1/6})$.  With $k\propto n$, taking $|x|=a_n$, his PDF approximation then has a multiplicative error of $O\left(a_n^3 k^{-1/2}\right)$, matching \eqref{eqn:fX-phi}.}}{}
\begin{align}\label{eqn:fX-phi}
f_{\mathbf{X}}(\mathbf{x})
  &= \phi_{\underline{\mathcal{V}}/n} \left(\mathbf{x} - \mathbf{k}/(n+1) \right)
    \left[ 1 + O\left(a_n^3\|\mathbf{\Delta k}^{-1/2}\|_\infty\right) \right] , \\
\label{eqn:dfX-dphi}
\D{f_{\mathbf{X}}(\mathbf{x})}{\mathbf{x}}
  &= \D{}{\mathbf{x}} \phi_{\underline{\mathcal{V}}/n} \left(\mathbf{x} - \mathbf{k}/(n+1) \right)
    +O\left( a_n^4 n^{J/2} n \|\mathbf{\Delta k}^{-1}\|_\infty  \right) . 
\end{align}
\item\label{lem:den-i-star-prob} For the Dirichlet-distributed $\mathbf{\Delta X}$, \ref{cond:star} is violated with only exponentially decaying (in $n$) probability: 
\Supplemental{\begin{equation}\label{eqn:star-violation-prob}}{$}
1 - P\bigl(\star(a_n)\bigr) 
= O\bigl( a_n^{-1}e^{-a_n^2/2} \bigr) . 
\Supplemental{\end{equation}}{$}
%
\item\label{lem:den-i-star-prob-fixed} If instead there are asymptotically fixed components of the parameter vector, the largest of which is $\Delta k_j=M<\infty$, then 
with $M=1$, 
\Supplemental{\begin{equation}\label{eqn:star-violation-prob-fixed-1}}{$}
1-P\bigl( \star(a_n)\bigr) \le e^{-a_n-1} .
\Supplemental{\end{equation}}{$}
With $M\ge2$, for any $\eta>0$, 
\Supplemental{\begin{equation}\label{eqn:star-violation-prob-fixed-M}}{$}
1-P\bigl( \star(a_n)\bigr) = o\bigl( a_n^{-1} \exp\bigl\{ -a_n\sqrt{M}(1/2-\eta) \bigr\} \bigr) .
\Supplemental{\end{equation}}{$}
\end{enumerate}
\end{lemma}

\subsection*{\Supplemental{Proof}{Sketch of proof} of Lemma \ref{lem:den-i}\Supplemental{(\ref{lem:den-i-recenter})}{}}

\Supplemental{\begin{proof}}{}
\Supplemental{%
Since
\begin{gather}\label{eqn:mean-kn}
\Delta k_j/(n+1)=(\Delta k_j/n)(1-O(1/n))=(\Delta k_j/n)-O(1/n), \\
\label{eqn:mode-kn}
(\Delta k_j-1)/(n-J)=\Delta k_j/[n(1-J/n)]-O(1/n)=(\Delta k_j/n)+O(1/n), 
\end{gather}
both the mean and mode of $\Delta X_j$ are within $O(1/n)$ of $\Delta k_j/n$.  By the triangle inequality, 
\begin{align*}
|\Delta x_j-\Delta k_j/(n+1)|
  &\le |\Delta x_j-\Delta k_j/n|+\overbrace{|\Delta k_j/n - \Delta k_j/(n+1)|}^{\textrm{apply \eqref{eqn:mean-kn}}}
\\&= |\Delta x_j-\Delta k_j/n| + O(1/n) , \\
|\Delta x_j-(\Delta k_j-1)/(n-J)|
  &\le |\Delta x_j-\Delta k_j/n| + \overbrace{|\Delta k_j/n - (\Delta k_j-1)/(n-J)|}^{\textrm{apply \eqref{eqn:mode-kn}}}
\\&= |\Delta x_j-\Delta k_j/n| + O(1/n) , 
\end{align*}
so
\begin{align*}
n\Delta k_j^{-1/2} \left|\Delta x_j-\Delta k_j/(n+1)\right|
  &\le \overbrace{n\Delta k_j^{-1/2} |\Delta x_j-\Delta k_j/n|}^{=O(a_n)\textrm{ by \ref{cond:star}}}
     + \overbrace{O(\Delta k_j^{-1/2})}^{=o(a_n)} , \\
n\Delta k_j^{-1/2} \left|\Delta x_j-(\Delta k_j-1)/(n-J)\right|
  &\le \overbrace{n\Delta k_j^{-1/2} |\Delta x_j-\Delta k_j/n|}^{=O(a_n)\textrm{ by \ref{cond:star}}}
     + \overbrace{O(\Delta k_j^{-1/2})}^{=o(a_n)} , 
\end{align*}
where $\Delta k_j^{-1/2}=o(a_n)$ is by condition (c) on $a_n$. 
}{%
The proof of part (\ref{lem:den-i-recenter}) uses the triangle inequality and the fact that the mean and mode differ from $\Delta k_j/n$ by $O(1/n)$.}
\Supplemental{\end{proof}}{}

\Supplemental{\subsection*{Proof of Lemma \ref{lem:den-i}(\ref{lem:den-i-spacing-density})}}{}

\Supplemental{\begin{proof}}{}
\Supplemental{Since }{For part (\ref{lem:den-i-spacing-density}), since }%
$\mathbf{\Delta X} \sim \textrm{Dirichlet}(\mathbf{\Delta k})$, for any $\mathbf{\Delta x}$ that sums to one, 
\begin{align}\label{eqn:log-fDx-orig}
\log(f_{\mathbf{\Delta X}}(\mathbf{\Delta x})) 
  &= \log(\Gamma(n+1)) 
    +\sum_{j=1}^{J+1} \bigg[
      (\Delta k_j -1) \log(\Delta x_j) - \log(\Gamma(\Delta k_j)) \bigg] .
\end{align}
\Supplemental{We apply the Stirling-type bounds in \citet{Robbins1955} to the gamma functions: 
\begin{gather*}
\sqrt{2\pi}n^{n+1/2}e^{-n+1/(12n+1)} 
  < n! < \sqrt{2\pi}n^{n+1/2}e^{-n+1/(12n)} , \textrm{ so }  \\
\log(\Gamma(n))
  = (1/2)\log(2\pi) + (n-1/2)\log(n-1) + (1-n) + O\left(n^{-1}\right). 
\end{gather*}
Applying the above and using $\sum_{j=1}^{J+1}\Delta k_j=n+1$,}{%
Applying Stirling-type bounds in \citet{Robbins1955} to the gamma functions,} 
\begin{align}
\notag
\log(f_{\mathbf{\Delta X}}(\mathbf{\Delta x}))
\Supplemental{%
  &= (1/2)\log(2\pi) +(n+1/2)\log(n) - n + O(n^{-1}) \\ 
\notag
  &\quad + \sum_{j=1}^{J+1} \Bigl[(\Delta k_j -1) \log(\Delta x_j) 
               - (1/2)\log(2\pi) 
               - (\Delta k_j-1/2)\log(\Delta k_j-1) \\
\notag
  &\qquad\qquad
               - 1 + \Delta k_j
               + O(\Delta k_j^{-1}) \Bigr]  \\ 
\notag 
  &= -(J/2)\log(2\pi) - n +(n+1)-(J+1) +(n+1/2)\log(n) \\ 
\notag 	 	
  &\quad + \sum_{j=1}^{J+1} \left[(\Delta k_j -1) \log(\Delta x_j) - (\Delta k_j-1/2)\log(\Delta k_j-1) \right]  \\ 
\notag 	 	
  &\quad + O(n^{-1}+\max\{\Delta k_j^{-1}\}) \\  	 	
\notag
  &= \frac{J}{2}\log\left(\frac{1}{2\pi}\right) - J +\log(n)\left[(J/2)+(J+1)/2+(n-J)\right] \\
\notag
  &\quad + \sum_{j=1}^{J+1} \left[(\Delta k_j -1) \log(\Delta x_j) 
              + (\Delta k_j-1+1/2)\log\left(\frac{1}{\Delta k_j-1}\right) \right]  \\
\notag
  &\quad + O\left(\|\mathbf{\Delta k}^{-1}\|_\infty\right) \\  	 	
\notag 
  &= \frac{J}{2}\left[\log\left(\frac{1}{2\pi}\right)+\log(n)\right] 
    + \left[\sum_{j=1}^{J+1} (1/2)\log\left(\frac{1}{\Delta k_j-1}\right)\right]
    +\frac{1}{2}\overbrace{(J+1)\log(n)}^{=\sum_{j=1}^{J+1}\log(n)} \\
\notag
&\quad
  - J 
  +\overbrace{[(n+1)-(J+1)]}^{=\sum_{j=1}^{J+1}(\Delta k_j-1)} \log(n)
  + \sum_{j=1}^{J+1} (\Delta k_j -1) \log\left(\frac{\Delta x_j}{\Delta k_j-1}\right) 
  \\
\notag
  &\quad + O\left(\|\mathbf{\Delta k}^{-1}\|_\infty\right) \\  	 	
\notag 
}{}
  &= \overbrace{\frac{J}{2}\log(n/(2\pi)) + \frac{1}{2}\sum_{j =1}^{J+1} \log\left(\frac{n}{\Delta k_j - 1}\right)}^D + \overbrace{\sum_{j=1}^{J+1} (\Delta k_j-1) \log\left(\frac{n \Delta x_j}{\Delta k_j -1}\right) - J}^{h(\mathbf{\Delta x})}  \\
\label{eqn:betaRepStep}
  &\quad +O(\|\mathbf{\Delta k}^{-1}\|_\infty) ,
\end{align}
where $D$ is the same constant as in the statement of the lemma.  
\Supplemental{Since the remainder terms only come from the Stirling-type approximations of $\log(\Gamma(n+1))$ and $\log(\Gamma(\Delta k_j))$, they do not depend on $\mathbf{\Delta x}$. %
}{}
 
\Supplemental{%
To expand $h(\cdot)$, we calculate derivatives
\begin{align}
\label{eqn:hj}
h_j(\mathbf{\Delta x}) &= \frac{\Delta k_j - 1}{\Delta x_j},\\
\label{eqn:hjj}
h_{j,j}(\mathbf{\Delta x}) &= -\frac{\Delta k_j - 1}{\Delta x_j ^2}, \qquad j\neq k \implies  h_{j,k}(\mathbf{\Delta x}) =0,\\
\label{eqn:hjjj}
h_{j,j,j}(\mathbf{\Delta x})&=2\frac{\Delta k_j - 1}{\Delta x_j^3} , \qquad (j\neq k) \cup (i\neq j) \implies h_{i,j,k}(\mathbf{\Delta x})=0,\\
\label{eqn:hjjjj}
h_{j,j,j,j}(\mathbf{\Delta x})&=-6\frac{\Delta k_j - 1}{\Delta x_j^4} , \qquad (j\neq k) \cup (i\neq j) \cup (k\neq l) \implies h_{i,j,k,l}(\mathbf{\Delta x})=0 .
\end{align}
The mode of the Dirichlet distribution is $\mathbf{\Delta x_0}\equiv(\mathbf{\Delta k}-1)/(n-J)$.  Using \eqref{eqn:hj}--\eqref{eqn:hjjjj},
\begin{align}
\notag
h(\mathbf{\Delta x_0})
  &= \log\left (\frac{n}{n-J}\right ) \sum_{j=1}^{J+1} (\Delta k_j -1)- J 
   = -\left[0-J/n+O(n^{-2})\right](n-J) - J  \\
\label{eqn:h0}
  &= J+O(n^{-1})-J 
   = O(n^{-1}),\\
\label{eqn:hj0}
h_j(\mathbf{\Delta x_0})&=n-J,\\
\label{eqn:hjj0}
h_{j,j}(\mathbf{\Delta x_0}) &= -(n-J)^2(\Delta k_j-1)^{-1},\\
\label{eqn:hjjj0}
h_{j,j,j}(\mathbf{\Delta x_0}) &= 2(n-J)^3(\Delta k_j-1)^{-2},\\
\notag
h_{j,j,j,j}(\mathbf{\Delta x}) &= -6\frac{\Delta k_j - 1}{\Delta x_j^4}\\
\notag
  &\ge -6\frac{\Delta k_j - 1}{\left(\Delta x_{0,j} - O(a_n n^{-1} \Delta k_j^{1/2})\right)^4} \\
\notag
  &= -6\frac{\Delta k_j}{\left(\Delta k_j/n+O(1/n) - O(a_n n^{-1} \Delta k_j^{1/2})\right)^4} +O(n^4/\Delta k_j^4) \\
\label{eqn:hjjjj0}
  &= -6n^4\Delta k_j^{-3}[1 +O(a_n\Delta k_j^{-1/2})] +O(n^4/\Delta k_j^4)
    = -O(n^4 \Delta k_j^{-3}),
\end{align}
where the rate on the fourth derivative follows by plugging in the smallest value of $\Delta x_j$ that satisfies \ref{cond:star}, and $O(a_n\Delta k_j^{-1/2})=o(1)$ by the stated condition (c) on $a_n$.  Since the remainder in \eqref{eqn:hjjjj0} relies only on \ref{cond:star}, it is uniform over all $\mathbf{\Delta x}$ satisfying \ref{cond:star}; the remainder in \eqref{eqn:h0} does not depend on $\mathbf{\Delta x}$ at all.  The fourth derivative is always negative since $\Delta k_j-1>0$. 

Since $h_j(\mathbf{\Delta x_0})$ is a constant, 
\begin{align*}
\sum_{j=1}^{J+1}h_j(\mathbf{\Delta x_0})(\Delta x_j-\Delta x_{0j})
  &= (n-J)\left[\left(\sum_{j=1}^{J+1} \Delta x_j\right)-\left(\sum_{j=1}^{J+1}\Delta x_{0j}\right)\right] \Supplemental{\\
  &= (n-J)[(n+1)-(n+1)]           }{}
   = 0  .
\end{align*}
Applying this and other results from \eqref{eqn:h0}--\eqref{eqn:hjjjj0}, %
}{%
We then expand $h(\cdot)$ around the Dirichlet mode, $\mathbf{\Delta x_0}$.  The cross partials are zero, the first derivative terms sum to zero, and the fourth derivative is smaller-order uniformly over $\mathbf{\Delta x}$ satisfying \ref{cond:star}:}
\begin{align}
\label{eqn:def-Rn} 
\begin{split}
\Supplemental{%
h(\mathbf{\Delta x})
  &= \overbrace{h(\mathbf{\Delta x_0})}^{\equiv R_{1n}}
  +\sum_{j=1}^{J+1} h_j(\mathbf{\Delta x_0})(\Delta x_j-\Delta x_{0j}) 	 	
  +\frac{1}{2}\sum_{j=1}^{J+1} h_{j,j}(\mathbf{\Delta x_0})(\Delta x_j-\Delta x_{0j})^2  \\ 	 	
  &\quad  	 	
  +\overbrace{\frac{1}{6}\sum_{j=1}^{J+1} h_{j,j,j}(\mathbf{\Delta x_0})(\Delta x_j-\Delta x_{0j})^3}^{\equiv R_{2n}}
  +\overbrace{\frac{1}{24}\sum_{j=1}^{J+1} h_{j,j,j,j}(\mathbf{\Delta \tilde x})(\Delta x_j-\Delta x_{0j})^4}^{\equiv R_{3n}}
}{%
h(\mathbf{\Delta x})
  &= \overbrace{h(\mathbf{\Delta x_0})}^{\equiv R_{1n}=O(n^{-1})}
  +\sum_{j=1}^{J+1} h_j(\mathbf{\Delta x_0})(\Delta x_j-\Delta x_{0j}) 	 	
  +\frac{1}{2}\sum_{j=1}^{J+1} h_{j,j}(\mathbf{\Delta x_0})(\Delta x_j-\Delta x_{0j})^2  \\ 	 	
  &\quad  	 	
  +\overbrace{\frac{1}{6}\sum_{j=1}^{J+1} h_{j,j,j}(\mathbf{\Delta x_0})(\Delta x_j-\Delta x_{0j})^3}^{\equiv R_{2n}=O(a_n^3 \|\mathbf{\Delta k}^{-1/2}\|_{\infty})}
  +\overbrace{\frac{1}{24}\sum_{j=1}^{J+1} h_{j,j,j,j}(\mathbf{\Delta \tilde x})(\Delta x_j-\Delta x_{0j})^4}^{\equiv R_{3n}=O(a_n^4 \|\mathbf{\Delta k}^{-1}\|_{\infty})} 
  ,
}
\end{split}
\Supplemental{%
\\
\notag
  &= \overbrace{O(n^{-1})}^{=R_{1n}} + 0  
    -\frac{(n-J)^2}{2} \sum_{j=1}^{J+1} \frac{\left(\Delta x_j - \frac{\Delta k_j -1}{n-J}\right)^2}{\Delta k_j -1} 
    +\overbrace{\frac{(n-J)^3}{3}\sum_{j=1}^{J+1} \frac{\left(\Delta x_j - \frac{\Delta k_j -1}{n-J}\right)^3}{(\Delta k_j -1)^2}}^{=R_{2n}} \\
  &\quad
\label{eqn:hTaylor2}
    + \overbrace{O(a_n^4 \|\mathbf{\Delta k}^{-1}\|_{\infty})}^{=R_{3n}} .
}{}
\end{align}
\Supplemental{%
For the quartic term $R_{3n}$, the rate is uniform over all $\mathbf{\Delta x}$ satisfying \ref{cond:star}.  
For the cubic term, uniformly over all $\mathbf{\Delta x}$ satisfying \ref{cond:star}, plugging in \eqref{eqn:star-mode} yields
\begin{equation*}
R_{2n} = O(n^3)O\left((a_nn^{-1}\Delta k_j^{1/2})^3/\Delta k_j^2\right)=O(a_n^3 \|\mathbf{\Delta k}^{-1/2}\|_{\infty}) , 
\end{equation*}
which is bigger than both $R_{1n}$ and $R_{3n}$. 
Thus, the remainder is
\begin{equation*}
R_n = \overbrace{R_{1n}+R_{2n}+R_{3n}}^{\textrm{defined in \eqref{eqn:def-Rn}}}
= O\left( n^{-1} + a_n^4\|\mathbf{\Delta k}^{-1}\|_{\infty} + a_n^3 \|\mathbf{\Delta k}^{-1/2}\|_{\infty} \right)
= O\left(a_n^3 \|\mathbf{\Delta k}^{-1/2}\|_{\infty}\right) . 
\end{equation*}
}{%
where the quadratic term expands to the form in the statement of the lemma.}

\Supplemental{%
We now consider the derivative of \eqref{eqn:betaRepStep} with respect to $\mathbf{\Delta x}$. As noted, the $O(\|\mathbf{\Delta k}^{-1}\|_\infty)$ remainder does not depend on $\mathbf{\Delta x}$ and thus has a derivative of zero.  Similarly, $D$ is a constant with respect to $\mathbf{\Delta x}$ and also has a zero derivative. 
The first derivative result is then an expansion of $h_j(\mathbf{\Delta x})=(\Delta k_j-1)\Delta x_j^{-1}\equiv g_j(\Delta x_j)$ in \eqref{eqn:hj} around the mode, $\Delta x_{0j}=(\Delta k_j-1)/(n-J)$:
\[ g_j(\Delta x_j) 
 = g_j(\Delta x_{0j}) 
  +g_j'(\Delta x_{0j})(\Delta x_j-\Delta x_{0j})
  +(1/2)g_j''(\Delta\tilde{x}_j)(\Delta x_j-\Delta x_{0j})^2 .
\]
Using an expansion of the form $(x+\eta)^{-1}=x^{-1}-x^{-2}\eta+\tilde{x}^{-3}\eta^2$ for small $\eta$ and $\tilde{x}$ between $x$ and $x+\eta$, with $x=\Delta x_{0j}$ and $\eta=\Delta x_j-\Delta x_{0j}$, 
\begin{align}\notag
\begin{split}
(\Delta k_j-1)\Delta x_j^{-1}
  &= (\Delta k_j-1)
     \biggl[ \frac{n-J}{\Delta k_j-1} 
           -\frac{(n-J)^2}{(\Delta k_j-1)^2}\left(\Delta x_j-\frac{\Delta k_j-1}{n-J}\right) \\
  &\qquad\qquad\qquad
           +\Delta\tilde{x}_j^{-3}\left(\Delta x_j-\frac{\Delta k_j-1}{n-J}\right)^2 \biggr] \\
  &= (n-J) -\frac{(n-J)^2}{\Delta k_j-1}\left(\Delta x_j-\frac{\Delta k_j-1}{n-J}\right)
    +O\left(\Delta k_j   n^3 \Delta k_j^{-3}   a_n^2 \Delta k_j n^{-2}\right) 
\end{split}
\\&= (n-J) -\frac{(n-J)^2}{\Delta k_j-1}\left(\Delta x_j-\frac{\Delta k_j-1}{n-J}\right)
    +O\left(n \Delta k_j^{-1} a_n^2 \right) .
\label{eqn:gj-exp}
\end{align}
The first term is order $n$, the second is $O(a_n\Delta k_j^{-1/2}n)$ by \ref{cond:star}, and since $a_n=o(\Delta k_j^{1/2})$ by condition (c) on $a_n$, the remainder is $o(a_n\Delta k_j^{-1/2}n)$.  
The remainder uses two applications of \eqref{eqn:star-mode}.  First, directly, $\left(\Delta x_j-\Delta x_{0j}\right)^2=O(a_n^2\Delta k_j n^{-2})$.  Second, the smallest possible $\Delta\tilde{x}_j$ satisfying \ref{cond:star} is $\Delta x_{0j}-O(a_n\Delta k_j^{1/2}n^{-1})$.  Using this, the fact that condition (c) on $a_n$ implies 
\begin{equation*}
a_n\Delta k_j^{1/2}n^{-1}
= a_n\Delta k_j^{-1/2} \Delta k_j/n 
= o(\Delta k_j/n) 
\end{equation*}
for all $j$, and $\Delta x_{0j}\asymp\Delta k_j n^{-1}$, then 
\begin{equation*}
\Delta\tilde{x}_j^{-3}
  = \frac{1}{\left(\Delta x_{0j}-o(\Delta k_j/n)\right)^3}
  =\frac{1}{\Delta x_{0j}^3\left[1-o(1)\right]}
  =O(n^3 \Delta k_j^{-3}) .
\end{equation*}
The overall remainder in \eqref{eqn:gj-exp} is thus uniform over all $\Delta x_j$ satisfying \ref{cond:star}.  Taking the maximum remainder over $j$ yields the uniform remainder $O\left(n \|\mathbf{\Delta k}^{-1}\|_\infty a_n^2 \right)$. 

For the derivative with respect to $\Delta k_j/n$,
\begin{align*}
\D{\log[f_{\mathbf{\Delta X}}(\mathbf{\Delta x})]}{\Delta k_j/n} 
  &= n\D{h(\mathbf{\Delta x})}{\Delta k_j}
    +n\D{D}{\Delta k_j}  ,  \\
\D{h(\mathbf{\Delta x})}{\Delta k_j/n}
\Supplemental{%
  &= \D{h(\mathbf{\Delta x_0})}{\Delta k_j/n}
    +\D{[h_j(\mathbf{\Delta x_0})(\Delta x_j-\Delta x_{0j})]}{\Delta k_j/n} \\
  &\quad
    +(1/2)\D{[h_{j,j}(\mathbf{\Delta x_0})(\Delta x_j-\Delta x_{0j})^2]}{\Delta k_j/n} \\
  &\quad
    +(1/6)\D{[h_{j,j,j}(\mathbf{\Delta x_0})(\Delta x_j-\Delta x_{0j})^3]}{\Delta k_j/n} \\
  &\quad
    +(1/24)\D{[h_{j,j,j,j}(\mathbf{\Delta \tilde x})(\Delta x_j-\Delta x_{0j})^4]}{\Delta k_j/n} \\
  &= n\log[n/(n-J)]
    +(n-J)(-1) \\
  &\quad
    +(1/2)\Bigl[ (n-J)^2n\Delta k_j^{-2} (\Delta x_j-\Delta x_{0j})^2 \\
  &\qquad\qquad\quad
                -(n-J)^2\Delta k_j^{-1} (-2)(\Delta x_j-\Delta x_{0j}) \Bigr] \\
  &\quad
    +(1/6)\Bigl[-4(n-J)^3n\Delta k_j^{-3} (\Delta x_j-\Delta x_{0j})^3 \\
  &\qquad\qquad\quad
                +2(n-J)^3\Delta k_j^{-2} (-3)(\Delta x_j-\Delta x_{0j})^2\Bigr] \\
  &\quad
    +(1/24)\Bigl[ O(n\Delta\tilde x_j^{-4}) (\Delta x_j-\Delta x_{0j})^4 \\
  &\qquad\qquad\quad
                 +O(n^4\Delta k_j^{-3}) (-4)(\Delta x_j-\Delta x_{0j})^3\Bigr] \\
  &= [J+O(n^{-2})] -(n-J) \\
  &\quad
    +O(n^3\Delta k_j^{-2} a_n^2 n^{-2}\Delta k_j)
      +(n-J)^2\Delta k_j^{-1} (\Delta x_j-\Delta x_{0j})  \\
  &\quad
    +O(n^4\Delta k_j^{-3} a_n^3 n^{-3} \Delta k_j^{3/2})
      +O(n^3\Delta k_j^{-2} a_n^2 n^{-2} \Delta k_j) \\
  &\quad
    +O(n^5\Delta k_j^{-4} a_n^4 n^{-4}\Delta k_j^2)
      +O(n^4\Delta k_j^{-3} a_n^3n^{-3}\Delta k_j^{3/2}) \\ %
}{}
  &= -(n-2J) +(n-J)^2 \Delta k_j^{-1}(\Delta x_j-\Delta x_{0j}) \\
  &\quad
    +O\left( a_n^2 n \Delta k_j^{-1} 
            +a_n^3 n \Delta k_j^{-3/2}
            +a_n^2 n \Delta k_j^{-1}
            +a_n^4 n \Delta k_j^{-2}
            +a_n^3 n \Delta k_j^{-3/2}\right)  , \\
  &= -(n-J) +J +(n-J)^2 (\Delta k_j-1)^{-1}(1-O(\Delta k_j^{-1}))(\Delta x_j-\Delta x_{0j}) \\
  &\quad
    +O\left( a_n^2 n \Delta k_j^{-1} \right)  , \\
\Supplemental{%
\D{D}{\Delta k_j}
  &= -(1/2)(\Delta k_j-1)^{-1}
   = -\frac{1}{2\Delta k_j(1-1/(2\Delta k_j))}
   = -(2\Delta k_j)^{-1}(1+o(1)) , \\ %
}{}
\D{D}{\Delta k_j/n}
  &= O(n\Delta k_j^{-1})
   = o(a_n^2 n \Delta k_j^{-1}). 
\end{align*}
Since $J$ is fixed asymptotically while $a_n^2 n \Delta k_j^{-1}\to\infty$, $J=O(a_n^2 n \Delta k_j^{-1})$ can go in the remainder, as can the $O(n^2\Delta k_j^{-2} a_n\Delta k_j^{1/2}n^{-1})$ from replacing $\Delta k_j^{-1}$ with $(\Delta k_j-1)^{-1}$ since $a_n\Delta k_j^{-1/2}=o(1)$ by condition (c) on $a_n$.  Again, the remainders are uniform over all $\mathbf{\Delta x}$ satisfying \ref{cond:star}.%
}{%
The derivative with respect to $\mathbf{\Delta x}$ is computed by expanding $h_j(\mathbf{\Delta x})=(\Delta k_j-1)\Delta x_j^{-1}$ around the mode, and then simplifying with \ref{cond:star} and the fact that $\sum_{j=1}^{J+1}(\Delta x_j-\Delta x_{0,j})=1-1=0$. 
The derivative with respect to $\Delta k_j/n$ is computed from an expansion of $h(\mathbf{\Delta x})$ around the mode, reusing many results from the original computation of the PDF. 
}
\Supplemental{\end{proof}}{}

\Supplemental{\subsection*{Proof of Lemma \ref{lem:den-i}(\ref{lem:den-i-ordered-density})}}{}

\Supplemental{\begin{proof}}{}
\Supplemental{%
The density for $\mathbf X$ may be obtained from the Lemma \ref{lem:den-i}(\ref{lem:den-i-spacing-density}) density for $\mathbf{\Delta X}$ simply by plugging in $\Delta x_j=x_j-x_{j-1}$ since the transformation is unimodular\Supplemental{:
\begin{equation*}
\left(\begin{array}{c}\Delta x_1 \\ \Delta x_2 \\ \Delta x_3 \\ \vdots \\ \Delta x_J\end{array}\right)
  = \left(\begin{array}{ccccc}
          1 & 0 & 0 & \cdots & 0 \\
          -1& 1 & 0 & \cdots & 0 \\
          0 & -1& 1 & \ddots & \vdots \\
          \vdots & \ddots & \ddots & \ddots & \\
          0 & \cdots & 0 & -1 & 1
          \end{array}\right)
    \left(\begin{array}{c}x_1 \\ x_2 \\ x_3 \\ \vdots \\ x_J\end{array}\right) ,
\end{equation*}
where the transformation matrix has determinant equal to one}{}. 
Centering at the mean instead of the mode introduces only smaller-order error\Supplemental{: from \eqref{eqn:mean-kn} and \eqref{eqn:mode-kn}, the mean and mode differ by $O(1/n)$, so 
\begin{align*}
\left[\Delta x_j-\Delta k_j/(n+1)\right]^2
  &= \left[\Delta x_j-(\Delta k_j-1)/(n-J)\right]^2
  +O\left(n^{-1}[\Delta x_j-\Delta k_j/(n+1)]\right) \\
  &= \left[\Delta x_j-(\Delta k_j-1)/(n-J)\right]^2
  +O(n^{-1} a_n \Delta k_j^{1/2} n^{-1}) ,
\end{align*}
using \eqref{eqn:star-mean} for the last line.}{.} 
We also replace $n-J=n[1-O(1/n)]$ and $\Delta k_j-1=\Delta k_j[1-O(\Delta k_j^{-1})]$.  The constant $D$ does not change.  
From the Lemma \ref{lem:den-i}(\ref{lem:den-i-spacing-density}) density for $\mathbf{\Delta x}$, plugging in $\Delta x_j=x_j-x_{j-1}$ and simplifying,
\begin{align*}
\log[f_{\mathbf X}(\mathbf x)] 
  &= \left[ \sum_{j=1}^{J+1} \Delta k_j^{-1}[1+O(\Delta k_j^{-1})] 
  \left\{ \left[ x_j-x_{j-1} - \left(\frac{k_j-k_{j-1}}{n+1}\right) \right]^2 
         +O(a_n \Delta k_j^{1/2} n^{-2}) \right\} \right] \\
  &\quad \times \left(- \frac{n^2[1-O(1/n)]}{2} \right) \\
  &\quad + D + O(a_n^3 \|\mathbf{\Delta k}^{-1/2}\|_\infty) \\
  &= D - \frac{n^2}{2} \sum_{j=1}^{J+1}\Delta k_j^{-1}\left\{\left[x_j-k_j/(n+1)\right]-\left[x_{j-1}-k_{j-1}/(n+1)\right]\right\}^2 \\
  &\quad +O(a_n^3 \|\mathbf{\Delta k}^{-1/2}\|_\infty) .
\end{align*}
For the $j=0$ entries, $\Delta x_1-\Delta k_1/(n+1)=x_1-k_1/(n+1)$, so implicitly $x_0\equiv0$ and $k_0\equiv0$.  For $j=J+1$,
\begin{equation*}
\Delta x_{J+1}-\Delta k_{J+1}/(n+1)
  = (1-x_J) - (n+1-k_J)/(n+1) , 
\end{equation*}
so $x_{J+1}\equiv1$ and $k_{J+1}\equiv n+1$. 
The diagonal terms $[x_j-k_j/(n+1)]^2$ have coefficients $-(n^2/2)[\Delta k_{j+1}^{-1}+\Delta k_j^{-1}]$. 
The one-off-diagonal terms $[x_j-k_j/(n+1)][x_{j+1}-k_{j+1}/(n+1)]$ have coefficients $n^2\Delta k_{j+1}^{-1}$.  Since we factor out $-1/2$, but have both $(j,j+1)$ and $(j+1,j)$ entries in the symmetric matrix $\underline{H}$, the diagonal elements of $\underline{H}$ become $\underline{H}_{j,j}=n^2[\Delta k_{j+1}^{-1}+\Delta k_j^{-1}]$ while the off-diagonal elements become $\underline{H}_{j,j+1}=\underline{H}_{j+1,j}=-n^2\Delta k_{j+1}^{-1}$.%
}{%
For part (\ref{lem:den-i-ordered-density}), the results are intuitive given part (\ref{lem:den-i-spacing-density}), so we defer to the supplemental appendix.  
It is helpful that the transformation from the values $X_j$ to spacings $\Delta X_j$ is unimodular. %
} 
\Supplemental{


For the derivative with respect to $x_j$, we differentiate \eqref{eqn:log-fDx-orig} with respect to $x_j$, 
\begin{equation}\label{eqn:df-dx-approx}
\D{\log[f_{\mathbf X}(\mathbf x)]}{x_j}
= \D{\log[f_{\mathbf{\Delta X}}(\mathbf{\Delta x})]}{x_j}
= \frac{\Delta k_j - 1}{\Delta x_j} - \frac{\Delta k_{j+1} - 1}{\Delta x_{j+1}} . 
\end{equation}
\Supplemental{%
We again expand around the modes: 
\begin{align} \notag
\frac{\Delta k_j - 1}{\Delta x_j} 
 &= (\Delta k_j - 1)\left[ \Delta x_{0,j}^{-1} - \Delta x_{0,j}^{-2}(\Delta x_j - \Delta x_{0,j}) + (1/2)\Delta \tilde{x}_j^{-3}\overbrace{\left(\Delta x_j - \Delta x_{0,j}\right)^2}^{=O(a_n^2\Delta k_j n^{-2})} \right] \\ \notag
 &= n-J - \frac{(n-J)^2}{\Delta k_j-1}\left(\Delta x_j - \Delta x_{0,j}\right) \\ \notag
 &\quad
  +(\Delta k_j - 1)\left[ \Delta x_{0,j}-O(a_n n^{-1} \Delta k_j^{1/2}) \right]^{-3} O(a_n^2\Delta k_j n^{-2}) \\ \notag
 &= n-J - \frac{n^2}{\Delta k_j-1}\left(\Delta x_j - \Delta x_{0,j}\right) +O(a_n n^{-1}\Delta k_j^{-1}) \\ \notag
 &\quad
  +\frac{(n-J)^3}{(\Delta k_j-1)^2}\left[1+O\left(\overbrace{a_n \Delta k_j^{-1/2}}^{\to0}\right)\right] O(a_n^2 \Delta k_j n^{-2}) \\
 &= n-J - \frac{n^2}{\Delta k_j}\left[1+O(\Delta k_j^{-1})\right]\left(\Delta x_j - \Delta x_{0,j}\right) +O(a_n^2 n\Delta k_j^{-1}) , \label{eqn:Dkj-Dxj}
\end{align}
using \ref{cond:star} as well as conditions (a) and (c) on $a_n$ and $\min_j \Delta k_j\to\infty$.  Continuing from \eqref{eqn:df-dx-approx}, plugging in \eqref{eqn:Dkj-Dxj} with both $j$ and $j+1$,
\begin{align*}
\D{\log[f_{\mathbf X}(\mathbf x)]}{x_j}
  &= \frac{\Delta k_j - 1}{\Delta x_j} - \frac{\Delta k_{j+1} - 1}{\Delta x_{j+1}} \\
  &= n-J - \frac{n^2}{\Delta k_j} \left(\Delta x_j - \Delta x_{0,j}\right) 
    -\left\{n-J - \frac{n^2}{\Delta k_{j+1}} \left(\Delta x_{j+1} - \Delta x_{0,j+1}\right) \right\} \\
  &\quad
    +O(a_n^2 n\|\mathbf{\Delta k}^{-1}\|_\infty) \\
  &= -n^2 \Delta k_j^{-1} \left[ x_j - x_{0,j} - (x_{j-1}-x_{0,j-1})\right] \\
  &\quad
     +n^2 \Delta k_{j+1}^{-1}\left[ x_{j+1} - x_{0,j+1} -(x_j-x_{0,j})\right]
     +O(a_n^2 n \|\mathbf{\Delta k}^{-1}\|_\infty) \\
  &= -n^2(\Delta k_j^{-1}+\Delta k_{j+1}^{-1})(x_j-x_{0,j})
     +n^2 \Delta k_j^{-1}(x_{j-1}-x_{0,j-1}) \\
  &\quad
     +n^2 \Delta k_{j+1}^{-1} (x_{j+1}-x_{0,j+1})
     +O(a_n^2 n \|\mathbf{\Delta k}^{-1}\|_\infty) . 
\end{align*}

For $j=1$, $\Delta x_1=x_1$ and $\Delta k_1=k_1$, so $\Delta x_1-\Delta x_{0,1}=x_1-(k_1-1)/(n-J)$; that is, there is no $x_{j-1}=x_0$ term.  
For $j=J$, $\Delta x_{J+1}=1-x_J$ and $\Delta k_{J+1}=n+1-k_J$, so 
\begin{gather}\notag
\Delta x_{0,J+1} = \frac{n+1-k_J-1}{n-J}
  = \frac{(n-J)-(k_J-1)+(J-1)}{n-J}
  = 1 - x_{0,J} +O(1/n) , \\
\label{eqn:DxJ1}
\Delta x_{J+1} - \Delta x_{0,J+1}
  = (1-x_J) - \left[1 - x_{0,J} +O(1/n)\right]
  = -(x_J-x_{0,J}) -O(1/n) , 
\end{gather}
so when $j=J$ there is no $x_{j+1}=x_{J+1}$ term, up to the maintained order of magnitude. 

Finally, from \eqref{eqn:mean-kn} and \eqref{eqn:mode-kn}, $(k_j-1)/(n-J)=k_j/(n+1)+O(1/n)$, and $\underline{H}O(1/n)=O(n^2\|\mathbf{\Delta k}^{-1}\|_\infty n^{-1}) = O(n\|\mathbf{\Delta k}^{-1}\|_\infty) = o(a_n^2 n \|\mathbf{\Delta k}^{-1}\|_\infty)$ by condition (a) of $a_n$, so we can recenter at the mean instead of mode to get the final result,
\begin{equation*}
\D{\log[f_{\mathbf{X}}(\mathbf x)]}{\mathbf{x}} 
  = -\underline{H}(\mathbf{x}-\mathbf{k}/(n+1)) + O(a_n^2 n \|\mathbf{\Delta k}^{-1}\|_\infty) .
\end{equation*}

For the derivative with respect to $k_j/(n+1)$, since (as shown) the transformation from $\mathbf{\Delta x}$ to $\mathbf{x}$ is unimodular, we may differentiate \eqref{eqn:log-fDx-orig} with respect to $k_j$ and later multiply by $n+1$:
\begin{align*}
\D{\log[f_{\mathbf X}(\mathbf x)]}{k_j}
  &= \D{\log[f_{\mathbf{\Delta X}}(\mathbf{\Delta x})]}{k_j}
   = \D{\log[f_{\mathbf{\Delta X}}(\mathbf{\Delta x})]}{\Delta k_j}
    -\D{\log[f_{\mathbf{\Delta X}}(\mathbf{\Delta x})]}{\Delta k_{j+1}}, \\ 
\D{\log[f_{\mathbf{\Delta X}}(\mathbf{\Delta x})]}{\Delta k_j}
  &= 0 + \log(\Delta x_j) - \Gamma'(\Delta k_j)/\Gamma(\Delta k_j) .
\end{align*}
Since $\Delta k_j$ is a natural number, the gamma function derivative is 
\begin{equation*}
\Gamma'(\Delta k_j)
  = (\Delta k_j-1)! \left[ -\gamma + \sum_{k=1}^{\Delta k_j-1}(1/k)\right] ,
\end{equation*}
where $\gamma$ is the Euler--Mascheroni constant.  By the asymptotic approximation for the $m$th harmonic number, $H_m=\sum_{k=1}^{m}(1/k)=\log(m)+\gamma+O(1/m)$, this simplifies to
\begin{equation*}
\frac{\Gamma'(\Delta k_j)}{\Gamma(\Delta k_j)}
  = \frac{(\Delta k_j-1)! \left[ \log(\Delta k_j-1) + O(\Delta k_j^{-1})\right]}
         {(\Delta k_j-1)!}
  = \log(\Delta k_j-1) + O(\Delta k_j^{-1}) .
\end{equation*}
Altogether so far,
\begin{align}\label{eqn:dfDx-dDk}
\D{\log[f_{\mathbf{\Delta X}}(\mathbf{\Delta x})]}{\Delta k_j}
  &= \log(\Delta x_j/(\Delta k_j-1)) + O(\Delta k_j^{-1}) , \\
\label{eqn:dfDx-dk}
\D{\log[f_{\mathbf{\Delta X}}(\mathbf{\Delta x})]}{k_j}
  &= \log(\Delta x_j/(\Delta k_j-1)) - \log(\Delta x_{j+1}/(\Delta k_{j+1}-1))
    +O(\|\mathbf{\Delta k}^{-1}\|_\infty) .
\end{align}

Expanding \eqref{eqn:dfDx-dDk} around the mode,
\begin{align}
\begin{split}
& \log(\Delta x_j/(\Delta k_j-1)) \\
&\quad= \log\left( \frac{\Delta x_{0,j} +(\Delta x_j-\Delta x_{0,j})}{\Delta k_j-1} \right) \\
&\quad= \log\left( \frac{1}{n-J} + \frac{\Delta x_j - \Delta x_{0,j}}{\Delta k_j-1} \right) \\
&\quad= \log(1/(n-J))
    +(n-J)\frac{\Delta x_j-\Delta x_{0,j}}{\Delta k_j-1} \\
&\quad\quad
    -(1/2)\left(\frac{1}{n-J}+\frac{\Delta \tilde{x}_j}{\Delta k_j-1}\right)^{-2}\left(\frac{\Delta x_j-\Delta x_{0,j}}{\Delta k_j-1}\right)^2 \\
&\quad= -\log(n-J) + \frac{n-J}{\Delta k_j-1}(\Delta x_j-\Delta x_{0,j}) 
     -O(n^2) O(\Delta k_j^{-2} a_n^2 n^{-2} \Delta k_j) \\
&\quad= -\log(n-J) + n \Delta k_j^{-1}[1+O(\Delta k_j^{-1})](\Delta x_j-\Delta x_{0,j}) \\
&\qquad
     +\overbrace{O(a_n n^{-1}\Delta k_j^{1/2} \Delta k_j^{-1})}^{=o(\Delta k_j^{-1})}
     +O(a_n^2 \Delta k_j^{-1}) 
\end{split}\\
\label{eqn:log-x-over-k-exp}
&\quad= -\log(n-J) + n \Delta k_j^{-1} (\Delta x_j-\Delta x_{0,j}) 
     +O(a_n n^{-1} \Delta k_j^{-3/2}) + O(a_n^2 \Delta k_j^{-1}) ,
\end{align}
where as before $\Delta \tilde{x}$ is from the mean value theorem, lying between zero and $\Delta x_j-\Delta x_{0,j}$, and therefore $\Delta \tilde{x}=O(a_n n^{-1} \Delta k_j^{1/2})=o(1)$ by \ref{cond:star} and condition (c) on $a_n$.  
Plugging \eqref{eqn:log-x-over-k-exp} into \eqref{eqn:dfDx-dk},
\begin{align*}
\D{\log[f_{\mathbf{\Delta X}}(\mathbf{\Delta x})]}{k_j}
  &= \left\{ -\log(n-J) +n \Delta k_{j}^{-1} (\Delta x_{j}-\Delta x_{0,j}) \right\} 
  \\
  &\quad
    -\left\{ -\log(n-J) +n \Delta k_{j+1}^{-1} (\Delta x_{j+1}-\Delta x_{0,j+1}) \right\} 
    +O(a_n^2\|\mathbf{\Delta k}^{-1}\|_\infty) \\
  &= n \Delta k_{j}^{-1} [x_{j}-x_{0,j} - (x_{j-1}-x_{0,j-1})] 
    -n \Delta k_{j+1}^{-1} [x_{j+1}-x_{0,j+1} - (x_{j}-x_{0,j})] \\
  &\quad
    +O(a_n^2\|\mathbf{\Delta k}^{-1}\|_\infty) \\
  &= n[\Delta k_{j}^{-1}+\Delta k_{j+1}^{-1}](x_j-x_{0,j})
    -n\Delta k_j^{-1} (x_{j-1}-x_{0,j-1}) \\
  &\quad
    -n\Delta k_{j+1}^{-1}(x_{j+1}-x_{0,j+1})
    +O(a_n^2\|\mathbf{\Delta k}^{-1}\|_\infty) , \\
\D{\log[f_{\mathbf{\Delta X}}(\mathbf{\Delta x})]}{k_j/(n+1)}
  &= n^2[\Delta k_{j}^{-1}+\Delta k_{j+1}^{-1}](x_j-x_{0,j})
    -n^2\Delta k_j^{-1} (x_{j-1}-x_{0,j-1}) \\
  &\quad
    -n^2\Delta k_{j+1}^{-1}(x_{j+1}-x_{0,j+1}) \\
  &\quad
    +O(n\|\mathbf{\Delta k}^{-1}\|_\infty a_n n^{-1} \|\mathbf{\Delta k}^{1/2}\|_\infty)
    +O(a_n^2 n \|\mathbf{\Delta k}^{-1}\|_\infty) . 
\end{align*}

As before, for $j=1$, $\Delta x_1=x_1$ and $\Delta k_1=k_1$, so $\Delta x_1-\Delta x_{0,1}=x_1-(k_1-1)/(n-J)$; that is, there is no $x_{j-1}=x_0$ term.  
As seen in \eqref{eqn:DxJ1}, $\Delta x_{J+1} - \Delta x_{0,J+1} = -(x_J-x_{0,J}) -O(1/n)$, so when $j=J$ there is no $x_{j+1}=x_{J+1}$ term, up to the maintained order of magnitude. 
Thus, the result can be written in the matrix form using $\underline{H}$ seen in the statement of the lemma. 
}{Expanding around the mode and simplifying leads to the result.  Similar calculations are done for the derivative with respect to $k_j/(n+1)$.}

\Supplemental{%
We now show that $\underline{H}\,\underline{\mathcal V}/n=\underline{I}_J$, the $J\times J$ identity matrix.  Immediately the $n^2$ in $\underline{H}$ cancels with the $1/n$ and the additional $1/n$ in each $\underline{\mathcal V}_{i,j}$.  We enumerate the different cases of multiplying row $i$ of $\underline{H}$ by row $j$ of $\underline{\mathcal V}/n$ below.  We omit the $n+1$ denominator in $\underline{\mathcal V}$, so the below results should be $n+1$ if $i=j$ and zero if $i\ne j$. 

If $i=j=1$, 
\begin{align*}
&\left( \Delta k_1^{-1} + \Delta k_2^{-1}\right) k_1(n+1-k_1)
- \Delta k_2^{-1} k_1(n+1-k_2) \\
&\quad= n+1-k_1 + \Delta k_2^{-1} k_1\left[n+1-k_1-(n+1-k_2)\right] \\
&\quad= n+1-k_1 + k_1 \Delta k_2^{-1} \Delta k_2
= n+1 .
\end{align*}

If $1<i<J$, and $i=j$,
\begin{align*}
&(-\Delta k_i^{-1}) k_{i-1}(n+1-k_i)
 + (\Delta k_i^{-1} + \Delta k_{i+1}^{-1}) k_i(n+1-k_i)
 - \Delta k_{i+1}^{-1} k-i(n+1-k_{i+1}) \\
&\quad= \Delta k_i^{-1}\left[ k_i(n+1-k_i) - k_{i-1}(n+1-k_i)\right]
 +\Delta k_{i+1}^{-1}\left[ k_i(n+1-k_i) - k_i(n+1-k_{i+1})\right] \\
&\quad= \Delta k_i^{-1} \Delta k_i(n+1-k_i) + \Delta k_{i+1}^{-1} k_i(n+1-k_i-n-1+k_{i+1}) \\
&\quad= n+1 - k_i +\Delta k_{i+1}^{-1} k_i \Delta k_{i+1}^{-1}
= n+1 .
\end{align*}

If $i=j=J$,
\begin{align*}
&-\Delta k_J^{-1} k_{J-1}(n+1-k_J) + (\Delta k_J^{-1} + \Delta k_{J+1}^{-1})k_J(n+1-k_J) \\
&\quad= \Delta k_J^{-1}\left[ k_J(n+1-k_J) - k_{J-1}(n+1-k_J)\right] + \Delta k_{J+1}^{-1} k_J\overbrace{(n+1-k_J)}^{=\Delta k_{J+1}} \\
&\quad= \Delta k_J^{-1} \Delta k_J (n+1-k_J) + k_J
= n+1 .
\end{align*}

If $i=1$ and $j>1$,
\begin{align*}
&(\Delta k_1^{-1} + \Delta k_2^{-1}) k_1 (n+1-k_j)
- \Delta k_2^{-1} k_2(n+1-k_j) \\
&\quad= \Delta k_2^{-1} (n+1-k_j)(k_1-k_2) + \Delta k_1^{-1} k_1(n+1-k_j) \\
&\quad= -(n+1-k_j) + (n+1-k_j)
= 0. 
\end{align*}

If $1<i<J$, and $j<i$,
\begin{align*}
&(-\Delta k_i^{-1}) k_j(n+1-k_{i-1}) + (\Delta k_i^{-1} + \Delta k_{i+1}^{-1}) k_j(n+1-k_i)
- \Delta k_{i+1}^{-1} k_j(n+1-k_{i+1}) \\
&\quad= \Delta k_i^{-1}\left[ k_j(n+1-k_i) - k_j(n+1-k_{i-1})\right]
+ \Delta k_{i+1}^{-1}\left[ k_j(n+1-k_i) - k_j(n+1-k_{i+1})\right] \\
&\quad= k_j \Delta k_i^{-1}(-\Delta k_i) + k_j \Delta k_{i+1}^{-1} \Delta k_{i+1}
= 0 .
\end{align*}

If $1<i<J$, and $j>i$,
\begin{align*}
&(-\Delta k_j^{-1}) k_{i-1}(n+1-k_j) + (\Delta k_i^{-1} + \Delta k_{i+1}^{-1}) k_i(n+1-k_j)
- \Delta k_{i+1}^{-1} k_{i+1}(n+1-k_j) \\
&\quad= \Delta k_i^{-1}(n+1-k_j)(k_i-k_{i-1}) + \Delta k_{i+1}^{-1}(n+1-k_j)(k_i-k_{i+1}) \\
&\quad= n+1-k_j - (n+1-k_j)
= 0 .
\end{align*}

If $i=J$ and $j<J$,
\begin{align*}
&-\Delta k_J^{-1} k_j(n+1-k_{J-1}) + (\Delta k_J^{-1} + \Delta k_{J+1}^{-1}) k_j(n+1-k_J) \\
&\quad= \Delta k_J^{-1}\left[ k_j \Delta k_{J+1} - k_j(n+1-k_{J-1})\right] + k_j\Delta k_{J+1}^{-1}\Delta k_{J+1} \\
&\quad= \Delta k_J^{-1} k_j(-\Delta k_J) + k_j
= 0 .  %
\end{align*}
}{%
For the result on $\underline{\mathcal{V}}/n$, it can simply be checked that $\underline{H}\,\underline{\mathcal V}/n=\underline{I}_J$, the $J\times J$ identity matrix, as shown explicitly in the supplement. %
}

Last, to establish the multivariate normal PDF representation, we show that
\begin{equation*}
e^D = (2\pi)^{-J/2} \left| \underline{\mathcal{V}}/n \right|^{-1/2} \left[ 1 + O(1/n)\right] ,
\end{equation*}
\Supplemental{%
which immediately reduces to showing
\begin{equation*}
\begin{split}
\left| \underline{\mathcal{V}}/n \right|^{-1/2} \left[ 1+O(1/n)\right]
 &= n^{J/2}
    \exp \left\{ (1/2)\sum_{j=1}^{J+1}\log\left(\frac{n}{\Delta k_j -1}\right) \right\} 
  = n^{J/2}
    \left[ \prod_{j=1}^{J+1} \frac{n}{\Delta k_j -1} \right]^{1/2} \\
 &= n^{J/2} \left[ \prod_{j=1}^{J+1} \Delta u_j^{-1}[1+O(1/n)] \right]^{1/2} ,
\end{split}
\end{equation*}
where $|\underline{\mathcal{V}}/n|$ is the determinant of matrix $\underline{\mathcal{V}}/n$. 
Since $\underline{\mathcal{V}}/n=\underline{H}^{-1}$, $|\underline{\mathcal{V}}/n|=1/|\underline{H}|$. Thus, we wish to show 
\begin{equation}\label{eqn:det-H-WTS}
|\underline{H}| = n^{J} \prod_{j=1}^{J+1}\Delta u_j^{-1} \left[1+O(1/n)\right] . 
\end{equation}

When $J=2$, 
\begin{align}\notag
|\underline{H}|
 &= n^2\left(\Delta k_1^{-1}+\Delta k_2^{-1}\right)
     n^2\left(\Delta k_2^{-1}+\Delta k_3^{-1}\right)
     -\left( -n^2\Delta k_2^{-1}\right)^2 \\ \notag
 &= n^2 \left[ \left(\frac{n}{\Delta k_1}+\frac{n}{\Delta k_2}\right)
               \left(\frac{n}{\Delta k_2}+\frac{n}{\Delta k_3}\right)
              -\left(\frac{n}{\Delta k_2}\right)^2 \right] \\ \notag
 &= n^2 \Bigl[ \left(\Delta u_1^{-1}[1+O(1/n)]+\Delta u_2^{-1}[1+O(1/n)]\right)
               \left(\Delta u_2^{-1}[1+O(1/n)]+\Delta u_3^{-1}[1+O(1/n)]\right) \\ \notag
 &\qquad\quad
              -\Delta u_2^{-2}[1+O(1/n)] \Bigr] \\ 
\Supplemental{%
\notag
 &= n^2 \left[ \Delta u_1^{-1}\Delta u_2^{-1}+\Delta u_2^{-1}\Delta u_2^{-1}
              +\Delta u_1^{-1}\Delta u_3^{-1}+\Delta u_2^{-1}\Delta u_3^{-1}
              -\Delta u_2^{-2} \right]
    \left[1+O(1/n)\right] \\ \notag
 &= n^2 \left[ \Delta u_1^{-1}\Delta u_2^{-1}
              +\Delta u_1^{-1}\Delta u_3^{-1}
              +\Delta u_2^{-1}\Delta u_3^{-1} \right]
    \left[1+O(1/n)\right] \\ \notag
 &= n^2 \frac{\Delta u_1+\Delta u_2+\Delta u_3}{\Delta u_1 \Delta u_2 \Delta u_3}
    \left[1+O(1/n)\right] \\ 
}{}
\label{eqn:det-H-J2}
 &= n^2 \Delta u_1^{-1} \Delta u_2^{-1} \Delta u_3^{-1}
    \left[1+O(1/n)\right] ,
\end{align}
which matches \eqref{eqn:det-H-WTS} as desired. 

More generally, since $\underline{H}$ is a tridiagonal matrix, its determinant may be calculated as the $J$th term of the continuant.  The continuant refers to the recursive sequence \citep[e.g.,][Ch.\ XIII, eqn.\ (4), p.\ 518]{Muir1960}
\begin{equation}\label{eqn:continuant}
f_{-1} = 0, \quad
f_{0} = 1, \quad
f_{1} = \underline{H}_{1,1}, \quad
f_{j} = \underline{H}_{j,j} f_{j-1} - \underline{H}_{j-1,j} \underline{H}_{j,j-1} f_{j-2} .
\end{equation}
Since $\underline{H}_{j,j+1}=\underline{H}_{j+1,j}=-n^2\Delta k_{j+1}^{-1}$ and $\underline{H}_{j,j}=n^2\left(\Delta k_j^{-1}+\Delta k_{j+1}^{-1}\right)$, 
\begin{align}\notag
f_{j} 
 &= \left( n \Delta u_{j}^{-1} [1+O(1/n)] + n \Delta u_{j+1}^{-1} [1+O(1/n)]\right) f_{j-1}
       -\left( n^2 \Delta u_{j}^{-2} [1+O(1/n)] \right) f_{j-2} \\
 &= \left( \Delta u_{j}^{-1} + \Delta u_{j+1}^{-1} \right) n f_{j-1} [1+O(1/n)]
       -\Delta u_{j}^{-2} n^2 f_{j-2} [1+O(1/n)] .
\label{eqn:continuant-recursion}
\end{align}
Up to the $[1+O(1/n)]$ coefficient,
\begin{gather*}
\begin{split}
f_{2} [1+O(1/n)]
     &= \left( \Delta u_{2}^{-1} + \Delta u_{3}^{-1} \right) n f_{1}
       -\Delta u_{2}^{-2} n^2 f_{0} \\
     &= n^2 \left( \Delta u_{1}^{-1} + \Delta u_{2}^{-1} \right) \left( \Delta u_{2}^{-1} + \Delta u_{3}^{-1} \right)
       -n^2 \Delta u_{2}^{-2} \\
     &= n^2 \left(
              \frac{\Delta u_1+\Delta u_2}{\Delta u_1\Delta u_2}
              \frac{\Delta u_2+\Delta u_3}{\Delta u_2\Delta u_3}
             -\frac{1}{\Delta u_2^2}
            \right) \\
     &= n^2 \frac{(\Delta u_1+\Delta u_2)(\Delta u_2+\Delta u_3) - \Delta u_1\Delta u_3}
                 {\Delta u_1\Delta u_2^2 \Delta u_3} \\
     &= n^2 \frac{\Delta u_2(\Delta u_1+\Delta u_2+\Delta u_3) +\Delta u_1\Delta u_3- \Delta u_1\Delta u_3}
                 {\Delta u_1\Delta u_2^2 \Delta u_3} \\
     &= n^2 \frac{\Delta u_1+\Delta u_2+\Delta u_3}
                 {\Delta u_1\Delta u_2\Delta u_3} ,
\end{split}
\end{gather*}
so for $i=2$, 
\begin{equation}\label{eqn:continuant-i}
f_i[1+O(1/n)] = n^{i} \left(\sum_{t=1}^{i+1} \Delta u_t\right) \prod_{t=1}^{i+1} \Delta u_t^{-1} .
\end{equation}
By induction, assuming \eqref{eqn:continuant-i} holds for all $i<j$, 
using \eqref{eqn:continuant-recursion}, 
\begin{align*}
f_{j} [1+O(1/n)]
 &= \left( \Delta u_{j}^{-1} + \Delta u_{j+1}^{-1} \right) n f_{j-1}
   -\Delta u_{j}^{-2} n^2 f_{j-2} \\
 &= n\frac{\Delta u_{j}+\Delta u_{j+1}}{\Delta u_{j}\Delta u_{j+1}}
    n^{j-1} \left(\sum_{t=1}^{j} \Delta u_{t}\right) \prod_{t=1}^{j} \Delta u_t^{-1} \\
 &\quad
   -n^2\frac{1}{\Delta u_{j}^2}
    n^{j-2} \left(\sum_{t=1}^{j-1} \Delta u_{t}\right) \prod_{t=1}^{j-1} \Delta u_t^{-1} \\
 &= n^{j}
    \frac{\Delta u_{j}+\Delta u_{j+1}}{\Delta u_{j}\Delta u_{j+1}}
    \left(\sum_{t=1}^{j} \Delta u_{t}\right) \prod_{t=1}^{j} \Delta u_t^{-1}
   -n^j
    \frac{1}{\Delta u_{j}^2}
    \left(\sum_{t=1}^{j-1} \Delta u_{t}\right) \prod_{t=1}^{j-1} \Delta u_t^{-1} \\
 &= n^j \Biggl\{
      \left( \sum_{t=1}^{j-1} \Delta u_t \right)
      \left( \prod_{t=1}^{j} \Delta u_t^{-1} \right)
      \overbrace{%
       \left[ \frac{\Delta u_j+\Delta u_{j+1}}{\Delta u_j\Delta u_{j+1}} - \frac{1}{\Delta u_j} \right]%
      }^{=(\Delta u_j+\Delta u_{j+1}-\Delta u_{j+1})/(\Delta u_j\Delta u_{j+1})=\Delta u_{j+1}^{-1}} \\
 &\qquad\quad
     +\frac{\Delta u_j+\Delta u_{j+1}}{\Delta u_j\Delta u_{j+1}} 
      \Delta u_j
      \prod_{t=1}^{j} \Delta u_t^{-1}
    \Biggr\} \\
 &= n^j \left\{
      \left( \sum_{t=1}^{j-1} \Delta u_t \right)
      \left( \prod_{t=1}^{j+1} \Delta u_t^{-1} \right)
     +(\Delta u_j+\Delta u_{j+1})
      \left( \prod_{t=1}^{j+1} \Delta u_t^{-1} \right)
    \right\} \\
 &= n^j
      \left( \sum_{t=1}^{j+1} \Delta u_t \right)
      \left( \prod_{t=1}^{j+1} \Delta u_t^{-1} \right) .
\end{align*}
Thus, $|\underline{H}|$ is
\begin{align*}
f_J
 &= n^J \left( \sum_{t=1}^{J+1} \Delta u_t \right)
      \left( \prod_{t=1}^{J+1} \Delta u_t^{-1} \right) [1+O(1/n)] 
  = n^J \prod_{t=1}^{J+1} \Delta u_t^{-1} [1+O(1/n)] ,
\end{align*}
which is the desired expression from \eqref{eqn:det-H-WTS}. %
}{%
which is simplified with $|\underline{\mathcal V}/n|=1/|\underline{H}|$.  The matrix $\underline{H}$ is tridiagonal, so its determinant may be calculated by a recursive argument using results on the continuant of a tridiagonal matrix from \citet{Muir1960}. 
}

\Supplemental{%
For the derivative of $f_{\mathbf{X}}(\mathbf{x})$ in \eqref{eqn:dfX-dphi}, 
\begin{align*}
\underline{H} 
 &= O(n^2 \|\mathbf{\Delta k}^{-1}\|_\infty), \\
\mathbf{x}-\mathbf{k}/(n+1) 
 &= O(a_n n^{-1}\|\mathbf{\Delta k}^{1/2}\|_\infty), \\
\D{\log[f_{\mathbf{X}}(\mathbf{x})]}{\mathbf{x}} 
 &= -\underline{H}\left(\mathbf{x}-\mathbf{k}/(n+1)\right) +O(a_n^2 n \|\mathbf{\Delta k}^{-1}\|_\infty ) \\
\D{\log[f_{\mathbf{X}}(\mathbf{x})]}{\mathbf{x}} 
 &= \frac{1}
          {f_{\mathbf{X}}(\mathbf{x})}
     \D{f_{\mathbf{X}}(\mathbf{x})}{\mathbf{x}} , \\
\D{f_{\mathbf{X}}(\mathbf{x})}{\mathbf{x}}
 &= f_{\mathbf{X}}(\mathbf{x})
     \D{\log[f_{\mathbf{X}}(\mathbf{x})]}{\mathbf{x}} \\
 &= \overbrace{\phi_{\underline{\mathcal{V}}/n}\left(\mathbf{x}-\mathbf{k}/(n+1)\right)}^{O(n^{J/2})}
     \left[ 1 + O\left( a_n^3 \|\mathbf{\Delta k}^{-1/2}\|_\infty \right) \right] \\
 &\quad\times
    \bigl[ \overbrace{-\underline{H}\left(\mathbf{x}-\mathbf{k}/(n+1)\right)}^{O\left(a_n n \|\mathbf{\Delta k}^{-1/2}\|_\infty \right)} 
          +O(a_n^2 n \|\mathbf{\Delta k}^{-1}\|_\infty ) \bigr] \\
 &= \overbrace{%
    \D{}{\mathbf{x}} 
    \phi_{\underline{\mathcal{V}}/n}\left(\mathbf{x}-\mathbf{k}/(n+1)\right) %
    }^{O\left(a_n n^{J/2} n\|\mathbf{\Delta k}^{-1/2}\|_\infty \right)}
   +O\left( a_n^4 n^{J/2} n \|\mathbf{\Delta k}^{-1}\|_\infty  \right) 
    ,
\end{align*}
where the remainder is smaller than the main term by condition (c) on $a_n$. 
}{}
\Supplemental{\end{proof}}{}
}{}

\Supplemental{\subsection*{Proof of Lemma \ref{lem:den-i}(\ref{lem:den-i-star-prob})}}{}

\Supplemental{\begin{proof}}{}
\Supplemental{
\Supplemental{Let }{For part (\ref{lem:den-i-star-prob}), let }%
\begin{equation*}
l_j \equiv (\Delta k_j/n) - a_n \Delta k_j^{1/2}n^{-1} , 
\quad
u_j \equiv (\Delta k_j/n) + a_n \Delta k_j^{1/2}n^{-1} . 
\end{equation*}
The probability that \ref{cond:star} is violated is bounded (by Boole's inequality) by
\begin{equation}\label{eqn:P-star-violated}
P\left( \cup_{j=1}^J \left\{\Delta X_j\not\in[l_j,u_j]\right\} \right) 
  \le J \max_j P\left( \Delta X_j\not\in[l_j,u_j] \right) .
\end{equation}
\Supplemental{%

For a single $j$, $(\Delta X_j,1-\Delta X_j)\sim\textrm{Dirichlet}(\Delta k_j,n+1-\Delta k_j)$, so $\Delta X_j\sim\beta(\Delta k_j,n+1-\Delta k_j)$.  
For $X\sim\beta(k,n+1-k)$, equation (2.17) from \citet{DasGupta2000} gives
\begin{equation}\label{eqn:beta-tail-bound}
P(X>u) \le \frac{u(1-u)}{u[k+(n+1-k)-1]-k}f_X(u) 
       = \frac{1-u}{n-k/u}f_X(u) 
\end{equation}
if $n+1-k>1$ and $u>k/n$, where $f_X(\cdot)$ is the PDF of $X$.  Plugging in $k=\Delta k_j$, $X=\Delta X_j$, and $u=u_j$, 
\begin{align}\notag
\begin{split}
n-\Delta k_j/u_j
  &= n-\frac{\Delta k_j}%
            {\Delta k_j/n + a_n \Delta k_j^{1/2}n^{-1}} 
   = \frac{\Delta k_j+a_n\Delta k_j^{1/2}-\Delta k_j}%
          {\Delta k_j/n + a_n \Delta k_j^{1/2}n^{-1}} \\
  &= \frac{a_n\Delta k_j^{1/2}}%
          {n^{-1}\Delta k_j^{1/2}[\Delta k_j^{1/2} + a_n]} 
   = \frac{a_n n}{a_n + \Delta k_j^{1/2}}, \\
P(\Delta X_j > u_j) 
  &\le \frac{1-u_j}%
            {n-\Delta k_j/u_j}
       f_{\Delta X_j}(u_j) \\
  &= \frac{1-\Delta k_j/n - a_n\Delta k_j^{1/2}/n}%
          {a_n n/(a_n + \Delta k_j^{1/2})}
       f_{\Delta X_j}(u_j) \\
  &\le \frac{(1-\Delta k_j/n)\left(a_n + \Delta k_j^{1/2}\right)}{a_n n}
     f_{\Delta X_j}(u_j) 
\end{split}
\\&= a_n^{-1} n^{-1} (1-\Delta k_j/n)\Delta k_j^{1/2}\left[1+\overbrace{O\left(a_n/\Delta k_j^{1/2}\right)}^{=o(1)}\right]
     f_{\Delta X_j}(u_j) .
\label{eqn:Xj-upper-int}
\end{align}
Using the PDF approximation (and notation) from Lemma \ref{lem:den-i}(\ref{lem:den-i-spacing-density}) for the marginal (beta) distribution of $\Delta X_j$, 
\begin{align*}
f_{\Delta X_j}(u_j) 
  &= e^D 
     \exp\left\{ -(1/2)\frac{(u_j-\Delta k_j/n)^2 n^2}{\Delta k_j(n+1-\Delta k_j)/(n+1)} + R_n \right\} \\
  &= \sqrt{\frac{n^3}{2\pi (\Delta k_j-1)(n-\Delta k_j)}}
     \exp\left\{ -(1/2)\frac{(a_n\Delta k_j^{1/2}/n)^2 n^2}{\Delta k_j[1-\Delta k_j/(n+1)]} 
                 +\overbrace{O(a_n^3 \Delta k_j^{-1/2})}^{=o(1)\textrm{ by (c) on $a_n$}} \right\} \\
  &= \frac{n}{\sqrt{2\pi (\Delta k_j-1)(1-\Delta k_j/n)}}
     \exp\left\{ -(1/2)\frac{a_n^2}{1-\Delta k_j/(n+1)} 
                 +o(1) \right\} \\
  &\le \frac{n}{\sqrt{\Delta k_j-1}\sqrt{1-\Delta k_j/n}}
     \exp\left\{ -a_n^2/2 \right\} \left[1+o(1)\right] . 
\end{align*}
Plugging this into \eqref{eqn:Xj-upper-int},
\begin{align*}
P(\Delta X_j > u_j) 
  &\le a_n^{-1} n^{-1} (1-\Delta k_j/n)\Delta k_j^{1/2}\left[1+\overbrace{O\left(a_n/\Delta k_j^{1/2}\right)}^{=o(1)}\right] \\
  &\quad\times
       \frac{n}{\Delta k_j^{1/2}[1-O(1/\Delta k_j)]\sqrt{1-\Delta k_j/n}}
       \exp\left\{ -a_n^2/2 \right\} \left[1+o(1)\right] \\
  &= a_n^{-1} (1-\Delta k_j/n)^{1/2}
     \left[1+O\left(\Delta k_j^{-1}\right)\right]
     \left[1+o(1)\right] 
     \exp\left\{ -a_n^2/2 \right\} \\ 
  &\le a_n^{-1} \left[1+o(1)\right]
     \exp\left\{ -a_n^2/2 \right\} ,
\end{align*}
using $O(a_n\Delta k_j^{-1/2})=o(1)$ by condition (c) on $a_n$. 

For $P(\Delta X_j<l_j)$, $1-\Delta X_j \sim \beta(n+1-\Delta k_j,\Delta k_j)$, so
\begin{align*}
P(\Delta X_j<l_j)
  &= P(-\Delta X_j > -l_j)
   = P(1-\Delta X_j > 1-l_j)  \\
  &= P\left(\beta(n+1-\Delta k_j) > (n+1-\Delta k_j)/n -1/n +a_n \Delta k_j^{1/2}/n\right) .
\end{align*}
This is similar to $P(\Delta X_k>u_k)$ for $\Delta X_k\sim\beta(\Delta k_k,n+1-\Delta k_k)$ with $\Delta k_k=n+1-\Delta k_j$, but now 
\begin{equation}\label{eqn:uk-def}
u_k = \Delta k_k/n -1/n +a_n(n+1-\Delta k_k)^{1/2}/n . 
\end{equation}
If $\Delta k_k$ is enough below $n/2$ so that $(n+1-\Delta k_k)^{1/2} \ge \Delta k_k^{1/2} + a_n^{-1}$, then $u_k\ge \Delta k_k/n +a_n\Delta k_k^{1/2}/n$ and the prior result for $P(\Delta X_j>u_j)$ applies directly; but if $\Delta k_k>n/2$, the result is not immediate. 
Applying \eqref{eqn:beta-tail-bound} as before, and using Lemma \ref{lem:den-i}(\ref{lem:den-i-recenter}) for the density, 
\begin{align*}
P(\Delta X_j<l_j)
  &= P\left(\Delta X_k > \Delta k_k/n -1/n +a_n(n+1-\Delta k_k)^{1/2}/n\right) \\
  &= P\left(\Delta X_k > u_k\right) \\
  &\le \frac{1 - u_k}{n-\Delta k_k/u_k} f_{\Delta X_k}(u_k) , \\
n-\Delta k_k/u_k
  &= n - \frac{\Delta k_k}{\Delta k_k/n -1/n +a_n (n+1-\Delta k_k)^{1/2}/n} \\
  &= \frac{\Delta k_k - 1 +a_n (n+1-\Delta k_k)^{1/2} - \Delta k_k}{\Delta k_k/n -1/n +a_n (n+1-\Delta k_k)^{1/2}/n} \\
  &= \frac{a_n (n+1-\Delta k_k)^{1/2}n - n}{\Delta k_k -1 +a_n (n+1-\Delta k_k)^{1/2}} , \\
P(\Delta X_k>u_k)
  &\le \left(1 - \left[\Delta k_k/n -1/n +a_n (n+1-\Delta k_k)^{1/2}/n\right]\right) \\
  &\quad\times 
       \frac{\Delta k_k -1 +a_n (n+1-\Delta k_k)^{1/2} }%
            {a_n (n+1-\Delta k_k)^{1/2}n - n} 
       f_{\Delta X_k}(u_k) \\
  &= n^{-1} \left(n+1-\Delta k_k\right)
     \left[1 - a_n (n+1-\Delta k_k)^{-1/2} \right] \\
  &\quad\times 
      n^{-1}\frac{a_n^{-1}(n+1-\Delta k_k)^{-1/2}\left(\Delta k_k -1\right) +1 }%
                 {1 - a_n^{-1}(n+1-\Delta k_k)^{-1/2}} 
     \sqrt{\frac{n^3}{2\pi (\Delta k_k-1)(n-\Delta k_k)}} \\
  &\quad\times
     \exp\left\{ -(1/2)\frac{(u_k-\Delta k_k/n)^2 n^2}{\Delta k_k(n+1-\Delta k_k)/(n+1)} + O(a_n^3\Delta k_k^{-1/2}) \right\} \\
  &\le n^{-1} \left(n+1-\Delta k_k\right)
       \left[1 - a_n (n+1-\Delta k_k)^{-1/2} \right] \\
  &\quad\times 
      \frac{a_n^{-1}(n+1-\Delta k_k)^{-1/2}\left(\Delta k_k -1\right) +1 }%
                 {1 - a_n^{-1}(n+1-\Delta k_k)^{-1/2}} \\
  &\quad\times
     \frac{1}{\sqrt{\Delta k_k-1}\sqrt{1-\Delta k_k/n}} \\
  &\quad\times
     \exp\left\{ -(1/2)\frac{(-1/n +a_n(n+1-\Delta k_k)^{1/2}/n)^2 n^2}{\Delta k_k(n+1-\Delta k_k)/(n+1)} + O(a_n^3\Delta k_k^{-1/2}) \right\} \\
  &= (1-\Delta k_k/n)\left[1+O(1/n)\right]\left[1-o(1)\right]
     \left[ 1 + a_n^{-1}\left(n+1-\Delta k_k\right)^{-1/2}(\Delta k_k-1) \right] \\
  &\quad\times 
      \left[1+O\left(a_n^{-1}(n+1-\Delta k_k)^{-1/2}\right)\right] 
      (\Delta k_k-1)^{-1/2} (1-\Delta k_k/n)^{-1/2} \\
  &\quad\times
     \exp\left\{ -(1/2)\frac{a_n^2(n+1-\Delta k_k)-O(a_n n^{1/2})}{(n+1-\Delta k_k)\Delta k_k/(n+1)} \right\} \\
  &\quad\times
     \overbrace{\left[1+ O(a_n^3\Delta k_k^{-1/2})\right]}^{=1+o(1)\textrm{ by condition (c) on $a_n$}} \\
  &= (1-\Delta k_k/n)^{1/2}(\Delta k_k-1)^{-1/2}[1+o(1)] \\
  &\quad\times
     \left[ 1 + (\Delta k_k-1) a_n^{-1} (n+1-\Delta k_k)^{-1/2}\right] \\
  &\quad\times
     \exp\left\{ -(1/2)a_n^2 (n+1)/\Delta k_k -O(a_n n^{-1/2}) \right\}  \\
  &\le \left[ \left(\frac{1-\Delta k_k/n}{\Delta k_k-1}\right)^{1/2}
             +a_n^{-1} (\Delta k_k-1)^{1/2} \left(\frac{1-\Delta k_k/n}{1-\Delta k_k/n +1/n}\right)^{1/2} n^{-1/2} \right] \\
  &\quad\times
     \exp\left\{ -a_n^2/2 \right\}
     \left[ 1 + O(a_n n^{-1/2}) \right]
     \left[ 1 + o(1) \right] \\
  &= \overbrace{\left[ O\left(\Delta k_k^{-1/2}\right) + O\left(a_n^{-1} \Delta k_k^{1/2} n^{-1/2}\right) \right]}^{=O(a_n^{-1})\textrm{ by condition (c) on $a_n$ and $\Delta k_k/n<1$}}
       \left[ 1 + o(1) \right] 
       e^{-a_n^2/2} \\
  &= O\left( a_n^{-1} e^{-a_n^2/2} \right) . 
\end{align*}
This is a slight improvement over using Inequality 11.1.1 in \citet[p.\ 440]{ShorackWellner1986}, as seen leading up to \eqref{eqn:star-violation-prob-SW} below. %

For the special case $a_n=2\log(n)$ and all $\Delta k_j\asymp n$, 
\begin{equation*}
\max_j \left\{ n\Delta k_j^{-1/2} |\Delta x_j - \Delta k_j(n+1)^{-1}| \right\}
= o_{a.s.}(a_n) , 
\end{equation*}
following from equation (1.6) in \citet{CsorgoRevesz1978} and the law of iterated logarithm for the corresponding Kiefer process mentioned in their Remark 2. 
(Their result is even stronger: it applies to the $\sup$ over the uniform empirical quantile process and has $\sqrt{\log[\log(n)]}$ instead of $\log(n)$.) 
That is, with probability one, \ref{cond:star} holds for all large enough $n$. 
}{%
Since $\Delta X_j$ has a beta distribution, we can use the beta tail probability bounds from equation (2.17) in \citet{DasGupta2000} and plug in our beta PDF approximation from Lemma \ref{lem:den-i}(\ref{lem:den-i-spacing-density}).  This approach yields a slight improvement over using Inequality 11.1.1 in \citet[p.\ 440]{ShorackWellner1986}, as detailed in the supplemental appendix. 
}
\Supplemental{\end{proof}}{}
}{%
For part (\ref{lem:den-i-star-prob}), we use Boole's inequality along with the beta tail probability bounds from \citet{DasGupta2000} and our beta PDF approximation from Lemma \ref{lem:den-i}(\ref{lem:den-i-spacing-density}). 
}

\Supplemental{\subsection*{Proof of Lemma \ref{lem:den-i}(\ref{lem:den-i-star-prob-fixed})}}{}

\Supplemental{\begin{proof}}{}
\Supplemental{Since }{For part (\ref{lem:den-i-star-prob-fixed}), since }%
$\Delta k_j=M<\infty$ is a fixed natural number, we can write $\Delta X_j = \sum_{i=1}^M \delta_i$, where each $\delta_i$ is a spacing between consecutive uniform order statistics.  The marginal distribution of each $\delta_i$ is $\beta(1,n)$%
\Supplemental{%
, with closed form CDF \citep{Kumaraswamy1980}
\[ F_{\delta_i}(d) = 1 - \left(1-d\right)^{n} \textrm{ for }0\le d\le1 . \]
Using this, and using Boole's inequality in the second line,%
}{%
. Using the corresponding CDF formula and Boole's inequality leads to the result. %
}%
\Supplemental{
\begin{align*}
P\left(\Delta X_j > M/n + a_n n^{-1} M^{1/2}\right) 
\Supplemental{%
  &\leq P\left(\exists i : \delta_i > \frac{1+a_n M^{-1/2}}{n} \right)\\
}{}
  &\leq M P\left (\delta_i > \frac{1+a_n M^{-1/2}}{n} \right )
\Supplemental{\\
  &= M \left(1-\frac{1+a_n M^{-1/2}}{n}\right)^n \\
  &= M \exp\left\{ n \log\left(1-\frac{1+a_n M^{-1/2}}{n}\right) \right\} \\
  &\le M \exp\left\{ n \left(-\frac{1+a_n M^{-1/2}}{n}\right) \right\} \\
  &= }{\le}
  M e^{-1} \exp\{-a_n M^{-1/2}\} .
\end{align*}
}{}%
\Supplemental{%
The final inequality follows because $\forall x \in (0,1)$, $\log(1-x) <-x$. 
If there are multiple fixed spacings in the parameter vector, this argument can be repeated and only the largest (fixed) spacing will enter the asymptotic rate\Supplemental{ since the probability of the union of all violations is bounded above by the sum of the individual violation probabilities (i.e.,\ Boole's inequality)}{}. 
}{}
\Supplemental{

The above bound is quite good for $M=1$, in which the first two inequalities become equalities, but it can be somewhat improved for $M>2$. 
Consider the Hoeffding bound for binomial random variables in equation (6) in Inequality 11.1.1 in \citet[p.\ 440]{ShorackWellner1986}.  The distribution of $U_{n:k}\sim\beta(k,n+1-k)$ is the same as that of any spacing $\Delta X_j$ with $\Delta k_j=k$.  In their notation for $\mathbb{G}_n$ (the empirical distribution) and $\mathbb{U}_n$ (the empirical process) for iid Uniform$(0,1)$ random variables $U_i$, for a fixed $k$, 
\begin{align*}
P\left( U_{n:k} > k/n + a_n\sqrt{k}/n\right)
  &= P\left( \mathbb{G}_n\left(\frac{k+a_n\sqrt{k}}{n}\right) \le (k-1)/n\right) \\
  &= P\left( \sqrt{n}\left(\mathbb{G}_n\left(\frac{k+a_n\sqrt{k}}{n}\right)-\frac{k+a_n\sqrt{k}}{n}\right) \le -\frac{1+a_n\sqrt{k}}{\sqrt{n}} \right) \\
  &= P\left( \mathbb{U}_n\left(\frac{k+a_n\sqrt{k}}{n}\right) \le -\frac{1+a_n\sqrt{k}}{\sqrt{n}} \right) \\
  &= P\left( -\mathbb{U}_n\left(\frac{k+a_n\sqrt{k}}{n}\right) \ge \frac{1+a_n\sqrt{k}}{\sqrt{n}} \right) \\
  &\le \exp\left\{ -\frac{(k+a_n\sqrt{k})^2/n}{2\left(\frac{1+a_n\sqrt{k}}{n}\right)\left(\frac{n-k-a_n\sqrt{k}}{n}\right)} \right\} \\
  &\le \exp\left\{ -\frac{a_n^2/2}{(1+a_n/\sqrt{k})(1-a_n\sqrt{k}/n -k/n)}\right\} \\
  &\le \exp\left\{ -\frac{a_n^2/2}{1+a_n/\sqrt{k}}\right\} ,
\end{align*}
where the first inequality is the application of the aforementioned inequality from \citet{ShorackWellner1986}, and the last inequality is from $1-a_n\sqrt{k}/n -k/n\le 1$. 

If $k\to\infty$ and $a_n/\sqrt{k}\to0$ like in Lemma \ref{lem:den-i}(\ref{lem:den-i-star-prob}), then we get the further bound
\begin{equation}\label{eqn:star-violation-prob-SW}
\exp\left\{ -\frac{a_n^2/2}{1+a_n/\sqrt{k}}\right\}
   \le \exp\left\{ -a_n^2(1/2)[1-o(1)] \right\}
\end{equation}
since $1+a_n/\sqrt{k}=1+o(1)$.  This is similar to Lemma \ref{lem:den-i}(\ref{lem:den-i-star-prob}) but without the leading $O(a_n^{-1})$ term and with the $o(1)$ inside the $\exp\{\cdot\}$, so it is slightly looser. 

If instead $k$ is fixed, then $a_n/\sqrt{k}\to\infty$ and $a_n^{-1}\sqrt{k}\to0$, unlike in Lemma \ref{lem:den-i}(\ref{lem:den-i-star-prob}) where $a_n/\sqrt{k}\to0$.  Consequently, for any $\eta>0$, there exists $N_\eta<\infty$ such that $a_n^{-1}\sqrt{k}<\eta$ for all $n>N_\eta$. 
Then, we may write the bound as
\begin{align*}
\exp\left\{ -\frac{a_n^2/2}{a_n k^{-1/2}(1+a_n^{-1}\sqrt{k})} \right\} 
  &= \exp\left\{ -(1/2)a_n\sqrt{k}[1-o(1)] \right\} \\
  &= o\left( \exp\left\{ -a_n\sqrt{k}(1/2-\eta) \right\} \right) 
\end{align*}
for any $\eta>0$. 
If $k=1$, then we have $-a_n(1/2-\eta)$ inside the exponential function, not as good as $-a_n$.  If $k=2$, then they are almost identical: $-(1/2-\eta)a_n\sqrt{2}=-a_n(1/\sqrt{2}-\eta\sqrt{2})$ compared to $-a_n/\sqrt{2}$.  If $k>2$, then the Hoeffding bound is tighter. 
Since we only use $k=1$ for our quantile inference application, the first result is more helpful. 
}{%
The other bound may be derived using equation (6) in Inequality 11.1.1 in \citet[p.\ 440]{ShorackWellner1986}, as seen in the supplemental appendix, but for our quantile inference application we only use the first result since it is better for $M=1$. %
}%

For violations of \ref{cond:star} in the other direction, the probability is zero for large enough $n$%
\Supplemental{:
\begin{align*}
P\left( \Delta X_j < M/n - a_n n^{-1} M^{1/2}\right) 
  &= P\left( \Delta X_j < n^{-1}M^{1/2}[M^{1/2}-a_n]\right) . 
\end{align*}
Since $M$ is fixed and $a_n\to\infty$, $M^{1/2}-a_n<0$ for large enough $n$, and $P\left(\Delta X_j < 0\right)=0$.  %
\Supplemental{\end{proof}}{}
}{ 
since $P(\Delta X_j<0)=0$ and $M^{1/2}-a_n<0$ for large enough $n$. 
}

\subsection*{Lemma for proving Theorem \ref{thm:cdferror}}

First, we introduce notation.  
\Supplemental{Vectors are always column vectors (unless otherwise noted) and in bold, e.g.,\ $\mathbf{u}=(u_1,\ldots,u_J)'$, and matrices (later) are underlined, e.g.,\ $\underline{\Omega}$. }{}%
From earlier, $u_0\equiv 0$ and $u_{J+1} \equiv 1$.  
\Supplemental{For all $j \in \{1,2,\ldots,J\}$, 
}{
For all $j$, $k_j \equiv \lfloor(n+1)u_j \rfloor$, $\epsilon_j \equiv (n+1)u_j - k_j$. 
Let $\mathbf{\Delta k}$ denote the $(J+1)$-vector such that $\Delta k_j=k_j-k_{j-1}$, let $\boldsymbol\psi=(\psi_1,\ldots,\psi_J)'$ be the fixed weight vector from \eqref{eqn:def-L}, and 
}
\begin{align}
\Supplemental{
\notag
k_j &\equiv \lfloor(n+1)u_j \rfloor ,
&  \epsilon_j &\equiv (n+1)u_j - k_j  ,
\intertext{where the $\epsilon_j \in [0,1)$ are interpolation weights as in \eqref{eqn:QXLdef}.  Let $\mathbf{\Delta k}$ denote the $(J+1)$-vector such that $\Delta k_j=k_j-k_{j-1}$, let $\boldsymbol\psi=(\psi_1,\ldots,\psi_J)'$ be the fixed weight vector from \eqref{eqn:def-L}, and}
}{}
\notag
Y_j &\equiv U_{n:k_j} \sim \beta(k_j,n+1-k_j) , 
& \mathbf{\Delta Y}
 &\equiv(Y_1,Y_2-Y_1,\ldots,1-Y_J) \sim \textrm{Dirichlet}(\mathbf{\Delta k}) , \\
\notag
  \Lambda_j &\equiv U_{n:k_j+1} - U_{n:k_j} \sim \beta(1,n) ,
 & \\
\label{eqn:def-lots}
Z_j &\equiv \sqrt{n}\left(Y_j-u_j\right) , 
&  V_j &\equiv \sqrt{n}\left[F^{-1}(Y_j)-F^{-1}(u_j)\right] , \\
\notag
\mathbb{X} &\equiv \sum_{j=1}^J \psi_j F^{-1}(Y_j) , 
&  \mathbb{X}_0 &\equiv \sum_{j=1}^J \psi_j F^{-1}(u_j) ,  \\
\notag
\mathbb{W} &\equiv \sqrt{n}\left(\mathbb{X} - \mathbb{X}_0\right) = \boldsymbol{\psi}'\mathbf{V} , 
&\mathbb{W}_{\boldsymbol\epsilon,\boldsymbol\Lambda}
  \equiv \mathrlap{\mathbb W + n^{1/2} \sum_{j=1}^J \epsilon_j \psi_j \Lambda_j \left[Q'(u_j) + Q''(u_j)(Y_j - u_j)\right] , } 
\end{align}
where 
\Supplemental{%
$Q(\cdot)=F^{-1}(\cdot)$ is the quantile function of interest, with first and second derivatives $Q'(\cdot)$ and $Q''(\cdot)$.  
The values and distributions of %
}{%
}the preceding variables are all understood to vary with $n$.  
 
Let $\phi_{\underline\Sigma}(\cdot)$ be the PDF of a mean-zero multivariate normal distribution with covariance $\underline\Sigma$. 
\begin{lemma}\label{lem:den}
\begin{enumerate}    \def\theenumi{\roman{enumi}} 
Let Assumption \ref{a:hut-pf} hold at $\mathbf{\bar u}$, and let each element of $\mathbf Y$ and $\boldsymbol\Lambda$ satisfy \ref{cond:star} (as defined in Lemma \ref{lem:den-i}) with $a_n=2\log(n)$.  The following results hold uniformly over any $\mathbf u=\mathbf{\bar u}+o(1)$. 
\item \label{lem:den-bd2}%
Let $\mathbf{C}$ be a $J$-vector of random interpolation coefficients as defined in \citet{Jones2002}: each $C_j \sim \beta(\epsilon_j,1-\epsilon_j)$, and they are mutually independent and independent of all other random variables.  Then,
\begin{equation}\label{eqn:LL-LI-approx}
\begin{split}
\left| n^{1/2}(L^L - \mathbb X_0) - \mathbb{W}_{\boldsymbol\epsilon,\boldsymbol\Lambda} \right|
  &= O\left(n^{-3/2}[\log(n)]^3\right) , \\
\left| n^{1/2}(L^I - \mathbb X_0) - \mathbb{W}_{\mathbf C,\boldsymbol\Lambda} \right| 
  &= O\left(n^{-3/2}[\log(n)]^3\right) .
\end{split}
\end{equation}

\item \label{lem:den-iii} Define $\underline{\mathcal{V}}$ as the $J\times J$ matrix with row $i$, column $j$ elements $\underline{\mathcal{V}}_{i,j} = \min\{u_i,u_j\}(1-\max\{u_i,u_j\})$, and define $\underline{\mathcal{A}}=\textrm{diag}\{f(F^{-1}(\mathbf{u}))\}$, i.e.,\ $\underline{\mathcal{A}}_{i,j}=f(F^{-1}(u_i))$ if $i=j$ and zero if $i\ne j$.  
Define $\mathcal{V}_{\boldsymbol{\psi}} \equiv \boldsymbol{\psi}'\left(\underline{\mathcal{A}}^{-1}\underline{\mathcal{V}}\,\underline{\mathcal{A}}^{-1}\right)\boldsymbol{\psi} \in \mathbb{R}$. 
For any realization $\boldsymbol{\lambda}$ of $\boldsymbol{\Lambda}=(\Lambda_1,\ldots,\Lambda_J)$ satisfying Condition $\star(2\log(n))$,
\begin{gather*}
\sup_{\{w:\star(2\log(n))\textrm{ holds}\}} \left|
\frac{f_{\mathbb{W}_{\mathbb \epsilon,\boldsymbol{\Lambda}}\mid \boldsymbol{\Lambda}}(w\mid \boldsymbol{\lambda})}
     {\phi_{\mathcal{V}_{\boldsymbol{\psi}}}(w)}
     - 1
\right|
= O\left(n^{-1/2}[\log(n)]^3\right), \\
\sup_{\{w:\star(2\log(n))\textrm{ holds}\}} \left|
    \D{f_{\mathbb{W}_{\mathbb \epsilon,\boldsymbol{\Lambda}}\mid \boldsymbol{\Lambda}}(w\mid \boldsymbol{\lambda})}{w}
    - \D{\phi_{\mathcal{V}_{\boldsymbol{\psi}}}(w)}{w}\right|
= O\left(n^{-1/2}[\log(n)]^{3+J}\right) , 
\end{gather*}
where the notation $\phi_{\mathcal{V}_{\boldsymbol{\psi}}}(\cdot)$ denotes the PDF of a normal random variable with mean zero and variance $\mathcal{V}_{\boldsymbol{\psi}}$. 
For any value $\boldsymbol{\tilde \epsilon}\in [0,1)^J$, 
uniformly over $K$ satisfying \ref{cond:star}, 
\begin{equation*}
\left. \frac{\partial^2 F_{\mathbb W_{\boldsymbol \epsilon,\boldsymbol\Lambda}|\boldsymbol{\Lambda}} ( K  \mid \boldsymbol\lambda)}{\partial\epsilon_j^2} \right|_{\boldsymbol\epsilon=\boldsymbol{\tilde \epsilon}} 
  = n \psi_j^2 Q'(u_j)^2  \lambda_j^2 \left[\left.\D{\phi_{\mathcal{V}_{\boldsymbol{\psi}}}(w)}{w} \right|_{w=K}\right] + O\left(n^{-3/2}[\log(n)]^{5+J}\right) . 
\end{equation*} 
\end{enumerate}
\end{lemma}

\subsection*{\Supplemental{Proof}{Sketch of proof} of Lemma \ref{lem:den}\Supplemental{(\ref{lem:den-bd2})}{}}

\Supplemental{\begin{proof}}{}
\Supplemental{With }{For part (\ref{lem:den-bd2}), with }%
a Taylor expansion, the object $L^L$ may be rewritten as
\begin{align}\notag
 L^L 
\Supplemental{%
    &=\sum_{j=1}^J \psi_j\bigg(Q(Y_j)+ \epsilon_j \left[Q(Y_j + \Lambda_j) - Q(Y_j)\right]\bigg)\\ \notag
    &=\mathbb{X}_0 + \sum_{j=1}^J \psi_j\bigg(Q(Y_j) - Q(u_j) +  \epsilon_j \left[Q(Y_j + \Lambda_j) - Q(Y_j)\right]\bigg)\\ \notag
    &=\mathbb{X}_0 + n^{-1/2}\mathbb{W} + \sum_{j=1}^J \psi_j \epsilon_j \left[Q'(Y_j)\Lambda_j+\frac{Q''(\tilde y_j)}{2}\Lambda_j^2\right]\\ 
\begin{split}\label{eqn:LL-MVT}
    &=\mathbb{X}_0 + n^{-1/2}\mathbb{W} + \sum_{j=1}^J \psi_j \epsilon_j \bigg[Q'(u_j)\Lambda_j + Q''(u_j)(Y_j - u_j)\Lambda_j \\*
    &\qquad\qquad\qquad\qquad\qquad\qquad+ \underbrace{\frac{Q'''(\tilde u_j)}{2}(Y_j - u_j)^2\Lambda_j}_{\nu^L_{j,1}} 
       +\underbrace{\frac{Q''(\tilde y_j)}{2}\Lambda_j^2}_{\nu_{j,2}^L}\bigg]
\end{split}\\ \notag
}{}
    &= \mathbb{X}_0 + n^{-1/2}\mathbb{W}_{\boldsymbol{\epsilon,\Lambda}} + \sum_{j=1}^J \psi_j \epsilon_j \left(\nu^L_{j,1} +\nu_{j,2}^L\right) ,
\Supplemental{}{%
\\ \nu^L_{j,1} &\equiv \frac{Q'''(\tilde u_j)}{2}[Y_j - u_j]^2\Lambda_j, \quad
   \nu_{j,2}^L \equiv \frac{Q''(\tilde y_j)}{2}\Lambda_j^2 ,
\label{eqn:LL-MVT}
}
\end{align}
where $\forall j$, $\tilde y_j \in (Y_j,Y_j+\Lambda_j)$ and $\tilde u_j$ is between $u_j$ and $Y_j$. 
\Supplemental{%
We want to show that the the final term (the summation) is $O\bigl(n^{-3/2}[\log(n)]^3\bigr)$. 
Since $J$ is fixed and finite, there exists $n_0$ such that for all $n>n_0$, all $u_j$ lie inside the neighborhoods in Assumption \ref{a:hut-pf}.  Since the third derivative of the quantile function is 
\[ Q'''(u)=[3f'(Q(u))^2-f(Q(u))f''(Q(u))]/f(Q(u))^5 , \]
$Q'''(\cdot)$ is uniformly bounded by \ref{a:hut-pf}, which states that $f(Q(u_j))$ is uniformly bounded away from zero and $f''(Q(u_j))$ is uniformly bounded.  This helps uniformly bound the remainder terms. 
Applying \ref{a:hut-pf} and Condition $\star(2\log(n))$ gives the desired rate on the $\nu^L$ terms. 
}{%
The remainder is $O\bigl(n^{-3/2}[\log(n)]^3\bigr)$ by applying \ref{cond:star}, noting that \ref{a:hut-pf} uniformly bounds the quantile function derivatives for large enough $n$ under \ref{cond:star}.  The argument for $L^I$ is essentially the same. 
}

\Supplemental{%
To be more explicit about the $\nu^L$ terms, we show that Condition $\star(2\log(n))$ implies $\tilde y_j\to u_j$ and $\tilde u_j\to u_j$ so that \ref{a:hut-pf} is applicable, which uniformly bounds the respective quantile function derivatives.  
Alternatively, without invoking \ref{cond:star}, $\tilde{y}_j\stackrel{a.s.}{\to}u_j$ and $\tilde{u}_j\stackrel{a.s.}{\to}u_j$, following from the result in (e.g.)\ \citet[pf.\ of Lem.\ 3.1, p.\ 583]{Bickel1967} that $\max_{1\le k\le n}|U_{n:k}-k(n+1)^{-1}|\stackrel{a.s.}{\to}0$, so \ref{a:hut-pf} holds almost surely (for large enough $n$). 

In \eqref{eqn:LL-MVT}, $\tilde y_j\in(Y_j,Y_j+\Lambda_j)$ and $\tilde u_j$ is between $u_j$ and $Y_j$, so it suffices (for both) to show that $Y_j\to u_j$ and $Y_j+\Lambda_j\to u_j$ under Condition $\star(2\log(n))$.  Since $J$ is fixed, convergence of each $\tilde y_j$ and $\tilde u_j$ implies uniform (over $j$) convergence. 
Using the triangle inequality and $u_j=(k_j+\epsilon_j)/(n+1)$,
\begin{align*}
Y_j  &= \sum_{t=1}^j \Delta Y_t, \\
| Y_j - u_j |  &=  \left| \sum_{t=1}^j \left( 
        \Delta Y_t - \Delta k_t/n \right) + 
        k_j/n - (k_j+\epsilon_j)/(n+1)
                    \right| \\
  &\le \sum_{t=1}^j \left|\Delta Y_t - \Delta k_t/n \right|
       + \left| \frac{k_j/n}{n+1} - \epsilon_j/(n+1) \right| \\
\label{eqn:Yj-cond-star}\refstepcounter{equation}\tag{\theequation}
  &\le \sum_{t=1}^j a_n \Delta k_t^{1/2} n^{-1}
       + O(1/n) \\
  &= \sum_{t=1}^j o(1) + O(1/n) 
   \to 0 , 
\end{align*}
so \ref{cond:star} implies $Y_j\to u_j$ as long as condition (b) on $a_n$ holds, which it does for $a_n=2\log(n)$ and $k_j\asymp n$. 
(The $o(1)$ is sufficient for this purpose but quite conservative; with $a_n=2\log(n)$ and all $\Delta k_t\asymp n$ as for quantile inference, $a_n\Delta k_t^{1/2}n^{-1}\asymp n^{-1/2}\log(n)$.) 
For $\Lambda_j$, \ref{cond:star} implies $|\Lambda_j-1/n|\le a_n / n$, and $a_n n^{-1}\to0$ by condition (b) on $a_n$, so $\Lambda_j\to0$. 
Thus, Condition $\star(2\log(n))$ implies 
\begin{equation}\label{eqn:tilde-u-limit}
\tilde y_j\to u_j, \quad \tilde u_j\to u_j, 
\end{equation}
so for large enough $n$, Assumption \ref{a:hut-pf} is satisfied at all $\tilde y_j$ and $\tilde u_j$, i.e.,\ $Q''(\tilde y_j)$ and $Q'''(\tilde u_j)$ are uniformly bounded. 

From Condition $\star(2\log(n))$ and \eqref{eqn:star-mean}, $|\Lambda_j-1/(n+1)|\le 2\log(n) / n$, so $\Lambda_j=O(\log(n)/n)$ and $\Lambda_j^2=O\left([\log(n)]^2 n^{-2}\right)$.  Also by Condition $\star(2\log(n))$, as just seen, 
\[ | Y_j - u_j |=O\left(a_n n^{-1} \|\mathbf{\Delta k}^{1/2}\|_\infty\right)=O\left(\log(n)n^{-1}n^{1/2}\right)=O\left(\log(n)n^{-1/2}\right) , \]
so $(Y_j - u_j)^2=O\left( [\log(n)]^2 n^{-1} \right)$, and $(Y_j-u_j)^2 \Lambda_j = O\left( [\log(n)]^3 n^{-2} \right)$. 
}{}

\Supplemental{%
The argument is similar for $L^I$ but also depends on the random, mutually independent interpolation weights $C_{j} \sim \beta(\epsilon_j, 1-\epsilon_j)$ that replace $\epsilon_j$%
\Supplemental{.  From \citet{Jones2002}, 
\begin{gather*}
\tilde Q^I_U(u_j) \stackrel{d}{=} C_{j} U_{n:k_j+1} + (1-C_{j}) U_{n:k_j} = Y_j + C_{j} \Lambda_j,\\
C_{j} \sim \beta(\epsilon_j, 1-\epsilon_j), \quad  C_{j} \independent (\mathbf{Y},\boldsymbol{\Lambda},\mathbf C_{{-j}}), 
\end{gather*}
so}{:}
\begin{align*}
L^I 
\Supplemental{%
  &= \sum_{j=1}^J \psi_j\bigg(Q(Y_j)+ \left[Q(Y_j + C_j \Lambda_j) - Q(Y_j)\right]\bigg)\\   
  &= \mathbb{X}_0 + \sum_{j=1}^J \psi_j\bigg(Q(Y_j) - Q(u_j) +  \left[Q(Y_j + C_j \Lambda_j) - Q(Y_j)\right]\bigg)\\ 	 	
  &= \mathbb{X}_0 + n^{-1/2}\mathbb{W}  	 	
    +\sum_{j=1}^J \psi_j \left[Q'(Y_j)C_j\Lambda_j+\frac{Q''(\tilde y_j)}{2} (C_j \Lambda_j)^2\right]\\%
  &= \mathbb{X}_0 + n^{-1/2}\mathbb{W} 
    +\sum_{j=1}^J \psi_j \bigg[Q'(u_j) C_j \Lambda_j + Q''(u_j)[Y_j - u_j]C_j \Lambda_j  \\
  &\qquad\qquad\qquad\qquad\qquad\quad + \underbrace{\frac{Q'''(\tilde u_j)}{2}[Y_j - u_j]^2C_j\Lambda_j}_{\nu^I_{j,1}} 
       +\underbrace{\frac{Q''(\tilde y_j)}{2}(C_j\Lambda_j)^2}_{\nu_{j,2}^I}\bigg]\\
}{}
  &= \mathbb{X}_0 + n^{-1/2}\mathbb{W}_{\boldsymbol{C,\Lambda}} 
    +\sum_{j=1}^J \psi_j \left[\nu^I_{j,1} +\nu_{j,2}^I\right] 
\Supplemental{.  }{, \\
   \nu^I_{j,1} &\equiv \frac{Q'''(\tilde u_j)}{2}[Y_j - u_j]^2 C_j\Lambda_j, \quad
   \nu^I_{j,2}  \equiv \frac{Q''(\tilde y_j)}{2}(C_j\Lambda_j)^2 .
}
\end{align*}
\Supplemental{%
We want to show that
\begin{equation*}
\sum_{j=1}^J \psi_j \left[\nu^I_{j,1} +\nu_{j,2}^I\right]
= O\left(n^{-3/2}[\log(n)]^3\right) . 
\end{equation*}
Since $C_j\in(0,1)$, $\tilde y_j\in(Y_j,Y_j+C_j\Lambda_j)\subset(Y_j,Y_j+\Lambda_j)$ and the $\tilde{u}_j$ are still between $Y_j$ and $u_j$, so \eqref{eqn:tilde-u-limit} holds.  Also, $0\le C_j\Lambda_j\le\Lambda_j$, so the remainder terms are bounded by the same rates as before. %
}{%
Since $C_j\in(0,1)$, the same arguments from before lead to the same rates. 
}
}{}
\Supplemental{\end{proof}}{}

\Supplemental{\subsection*{Proof of Lemma \ref{lem:den}(\ref{lem:den-iii})}}{}

\Supplemental{\begin{proof}}{}
\Supplemental{Since }{For part (\ref{lem:den-iii}), since }%
$\boldsymbol\Lambda$ contains finite spacings, we cannot apply Lemma \ref{lem:den-i}(\ref{lem:den-i-ordered-density}).  Instead, we use the result from \textbf{8.7.5} in \citet[p.\ 238]{Wilks1962},
\begin{equation*}
(\Lambda_1,\ldots,\Lambda_J,1-\Lambda_1-\cdots-\Lambda_J)
  \sim \textrm{Dirichlet}(1,\ldots,1,n+1-J) . 
\end{equation*}
\Supplemental{%
Using the Dirichlet PDF formula,
\begin{align*}
\begin{split}
f_{\boldsymbol{\Lambda}}(\boldsymbol{\lambda})
  &= \frac{\Gamma(n+1)}{\Gamma(1)\cdots\Gamma(1)\Gamma(n+1-J)}
     \lambda_1^{1-1}\cdots\lambda_J^{1-1}(1-\lambda_1-\cdots-\lambda_J)^{n+1-J-1} \\
  &= n(n-1)\cdots(n-J+1)\left(1-\lambda_1-\cdots-\lambda_J\right)^{n-J}, 
\end{split}
\\
\log f_{\boldsymbol{\Lambda}}(\boldsymbol{\lambda})
  &= J \log(n) 
   +\sum_{j=1}^{J-1} \log\left(\frac{n-j}{n}\right) 
   +(n-J)\log\left(1-\sum_{j=1}^{J}\lambda_j\right) .
\end{align*}
Now
\[ \log\left(\frac{n-j}{n}\right) = \log(1-j/n) = -j/n + O(n^{-2}) , \]
so $\sum_{j=1}^{J}\log((n-j)/n)=O(1/n)$, and 
\[   \log\left(1-\sum_{j=1}^{J}\lambda_j\right)
  = -\sum_{j=1}^{J}\lambda_j + O\left(\left(\sum_{j=1}^{J}\lambda_j\right)^2\right) . \]
Under \ref{cond:star} with $a_n=2\log(n)$, $\lambda_j=O(n^{-1}\log(n))$, so 
\[ \left(\sum_{j=1}^{J}\lambda_j\right)^2 = O\left(n^{-2}[\log(n)]^2\right), \quad
   \sum_{j=1}^{J}\lambda_j = O(n^{-1}\log(n)) , \]
since $J$ is fixed.  Altogether,
}{%
Directly approximating the corresponding PDF yields
}
\begin{align}
\log f_{\boldsymbol{\Lambda}}(\boldsymbol{\lambda})
\Supplemental{\notag
  &= J \log(n) + O(1/n) - n\sum_{j=1}^{J}\lambda_j +O(n^{-1}\log(n)) + O(n^{-2}[\log(n)]^2) \\
}{}
\label{eqn:LambdaExp}
  &= J \log(n) - n \sum_{j=1}^J \lambda_j + O(n^{-1}\log(n)) . 
\end{align}

For the joint density of $\{\mathbf{Y},\boldsymbol\Lambda\}$, define 
\Supplemental{%
\begin{gather*}
\mathbf{T} 
  \equiv \left(\begin{array}{c}\Delta Y_1 \\ \Lambda_1 \\ \Delta Y_2-\Lambda_1 \\ \vdots \\ \Lambda_J \\ \Delta Y_{J+1}-\Lambda_J \end{array} \right)
  = \underline{T} \left(\begin{array}{c} Y_1 \\ \vdots \\ Y_{J+1} \\ \Lambda_1 \\ \vdots \\ \Lambda_J \end{array} \right)
  = \underline{T}\left(\mathbf{Y}',\boldsymbol\Lambda'\right)' ,  
\Supplemental{\\
\underline{T} 
  \equiv \left(\begin{array}{rrrrrrrrrrrr}
       1  & 0  & 0 & 0 & \cdots & 0 & 0 & 0  & 0  & 0 & \cdots & 0 \\
       0  & 0  & 0 & 0 & \cdots & 0 & 0 & 1  & 0  & 0 & \cdots & 0 \\
       -1 & 1  & 0 & 0 & \cdots & 0 & 0 & -1 & 0  & 0 & \cdots & 0 \\
       0  & 0  & 0 & 0 & \cdots & 0 & 0 & 0  & 1  & 0 & \cdots & 0 \\
       0  & -1 & 1 & 0 & \cdots & 0 & 0 & 0  & -1 & 0 & \cdots & 0 \\
       \vdots  & \vdots  & \vdots & \vdots & \vdots & \vdots & \vdots  & \vdots  & \vdots & \vdots & \vdots & \vdots \\
       0  & 0  & 0 & 0 & \cdots & 0  & 0 & 0  & 0  & 0 & \cdots & 1 \\
       0  & 0  & 0 & 0 & \cdots & -1 & 1 & 0  & 0  & 0 & \cdots & -1 \\
  \end{array} \right) , 
}{}
\end{gather*}
}{%
$\mathbf{T}\equiv\bigl(\Delta Y_1 , \Lambda_1 , \Delta Y_2-\Lambda_1 , \ldots , \Lambda_J , \Delta Y_{J+1}-\Lambda_J\bigr)=\underline{T}\bigl(\mathbf{Y}',\boldsymbol\Lambda'\bigr)'$, %
}
\Supplemental{%
so even-numbered row $2i$ has a single non-zero entry of $1$ in column $J+1+i$ while odd-numbered row $2i+1$ has entries $-1$ in columns $i$ and $J+1+i$ and entry $1$ in column $i+1$, and row $i=1$ is simply $(1,0,\ldots,0)$. 
The determinant of $\underline{T}$ is invariant to adding scalar multiples of one row to another row.  A series of such row operations on $\underline{T}$ yields a matrix with only a single entry of $1$ in each row, whose determinant is thus one: first, add row 1 to row 3; second, add row 2 to row 3; continue adding rows $2i-1$ and $2i$ to row $2i+1$ for $i=2,\ldots,J$.  Thus, $\det(\underline{T})=1$, i.e.,\ $\underline{T}$ is unimodular. 
}{where $\det(\underline{T})=1$ can be shown. }
Now
\begin{equation*}
\mathbf{T} \sim\textrm{Dirichlet}(\Delta k_1,1,\Delta k_2 -1,\ldots,1,\Delta k_{J+1}-1) . 
\end{equation*}
\Supplemental{%
Unfortunately, we still may not apply Lemma \ref{lem:den-i}(\ref{lem:den-i-ordered-density}).  However, using }{%
Using }the formula for the PDF of a transformed vector and plugging in the Dirichlet PDF formula for $\mathbf{T}$ %
\Supplemental{%
\begin{align}
\notag
\log(f_{\mathbf Y,\boldsymbol\Lambda} (\mathbf y,\boldsymbol\lambda )) 
  &= \log\left(f_{\mathbf{T}}\left(\underline{T}\left(\mathbf{y}',\boldsymbol\lambda'\right)'\right)\right) 
     \overbrace{\det(\underline{T})}^{=1} \\
\begin{split}\label{eqn:Ylambda}
  &= \log(\Gamma(n+1)) + (\Delta k_1 - 1) \log(\Delta y_1)  - \log\left(\Gamma(\Delta k_1)\right) \\
  &\quad + \sum_{j=2}^{J+1} \bigg[(\Delta k_j - 2) \log(\Delta y_j - \lambda_{j-1})  - \log\left(\Gamma(\Delta k_j -1)\right) \bigg] \\
  &\quad + \sum_{j=1}^{J} (1-1)\log(\lambda_j) . 
\end{split}
\end{align}

Combining \eqref{eqn:LambdaExp} and \eqref{eqn:Ylambda}, 
\begin{align}
\notag
\log f_{\mathbf Y \mid \boldsymbol{\Lambda}} \left( \mathbf{y} \mid \boldsymbol{\lambda} \right)
  &= \log(f_{\mathbf Y,\boldsymbol\Lambda} (\mathbf y,\boldsymbol\lambda )) 
    -\log f_{\boldsymbol{\Lambda}}(\boldsymbol \lambda)\\
\Supplemental{%
\notag 	 	
  &= \log(\Gamma(n+1)) + (\Delta k_1 - 1) \log(\Delta y_1)  - \log\left(\Gamma(\Delta k_1)\right) \\ 	 	
\notag 	 	
  &\quad + \sum_{j=2}^{J+1} \bigg [(\Delta k_j - 2) \log(\Delta y_j - \lambda_{j-1})  - \log\left(\Gamma(\Delta k_j - 1)\right) - \log n + n \lambda_{j-1} \bigg ]\\ 	 	
\notag 	 	
  &= \log f_{\mathbf Y}(\mathbf y) +  \sum_{j=2}^{J+1} \bigg [(\Delta k_j - 2) \log(\Delta y_j - \lambda_{j-1}) - (\Delta k_j - 1) \log(\Delta y_j)\\ 	 	
\notag 	 	
  &\quad+ \left[ \log\left(\Gamma(\Delta k_j)\right) - \log\left(\Gamma(\Delta k_j - 1)\right) \right] 
        - \log n + n \lambda_{j-1} \bigg ]\\ %
}{}
\begin{split}\label{eqn:YconLamb}
  &= \log f_{\mathbf Y}(\mathbf y) \\
  &\quad
    +\sum_{j=2}^{J+1} \bigg [\bigg\{(\Delta k_j - 2) [\log(\Delta y_j - \lambda_{j-1}) -  \log(\Delta y_j)] + n  \lambda_{j-1}\bigg\}\\
  &\qquad\qquad\;+ \bigg\{\left[\log\left(\Gamma(\Delta k_j)\right) - \log\left(\Gamma(\Delta k_j - 1)\right)\right] - \log n - \log\Delta y_j\bigg\}  \bigg ] 
\end{split}
\\
\Supplemental{%
\notag
  &= \log f_{\mathbf Y}(\mathbf y) +  \sum_{j=2}^{J+1} \bigg [\bigg\{(\Delta k_j - 2) [-\lambda_{j-1}/\Delta y_j + \overbrace{O(\lambda_{j-1}^2/\Delta y_j^2)}^{=O(n^{-2}[\log(n)]^2)}] + n \lambda_{j-1}\bigg\} \\
\notag 	 	
  &\qquad\qquad\qquad\qquad\;+ \bigg\{[\log(\Delta k_j) + O(\Delta k_j^{-1})] - \log n - \log\Delta y_j\bigg\}  \bigg ] \\ \notag
  &= \log f_{\mathbf Y}(\mathbf y) +  \sum_{j=2}^{J+1} \bigg[
        \lambda_{j-1}\left(n-\Delta k_j/\Delta y_j +2/\Delta y_j\right) + O\left(n^{-1}[\log(n)]^2\right) \\
\notag
  &\qquad\qquad\qquad\qquad\;+ O(n^{-1}) +\log\left(\frac{\Delta k_j/n}{\Delta y_j}\right) \bigg ] \\  \notag
  &= \log f_{\mathbf Y}(\mathbf y) +  \sum_{j=2}^{J+1} \bigg[
        \lambda_{j-1}\frac{n\Delta y_j-\Delta k_j}{\Delta y_j} +O\left(n^{-1}\log(n)\right) + O\left(n^{-1}[\log(n)]^2\right) \\
\notag
  &\qquad\qquad\qquad\qquad\;+ \log\left(\frac{\Delta k_j/n}{\Delta k_j/n\left[1+O\left(n^{-1/2}\log(n)\right)\right]}\right) \bigg ] \\  \notag
  &= \log f_{\mathbf Y}(\mathbf y) +  \sum_{j=2}^{J+1} \bigg[
        O\left(n^{-1}\log(n)\right)\frac{O\left(n^{-1/2}\log(n)\right)}{\Delta y_j}  + O\left(n^{-1}[\log(n)]^2\right) \\
\notag
  &\qquad\qquad\qquad\qquad\;+ \log\left( 1 + O\left(n^{-1/2}\log(n)\right)\right) \bigg ] \\  
\label{eqn:Y-PDF-cond-log-approx}
  &= \log f_{\mathbf Y}(\mathbf y) +  \sum_{j=2}^{J+1} \bigg[
        O\left(n^{-3/2}[\log(n)]^2 + n^{-1}[\log(n)]^2\right) \\
\notag
  &\qquad\qquad\qquad\qquad\;+ O\left(n^{-1/2}\log(n)\right) \bigg ] \\ 
}{}
\label{eqn:Y-cond-L-final}
  &= \log \phi_{\underline{\mathcal V}/n}\left(\mathbf y - \mathbf k/(n+1)\right)  
    +O\left(n^{-1/2}[\log(n)]^3\right) , 
\end{align}
repeatedly applying \ref{cond:star} with $a_n=2\log(n)$ and
applying Lemma \ref{lem:den-i}(\ref{lem:den-i-ordered-density}) to the log density of $\mathbf Y$, and since $\Delta k_j\asymp n$. 
The entries $\underline{\mathcal{V}}_{i,j}$ in matrix $\underline{\mathcal{V}}$ come from Lemma \ref{lem:den-i}(\ref{lem:den-i-ordered-density}), 
\Supplemental{%
\begin{align}
\notag
\begin{split}
\underline{\mathcal{V}}_{i,j}
  &= \overbrace{\min\{k_i,k_j\}(n+1-\max\{k_i,k_j\})/[n(n+1)]}^{\textrm{Lemma \ref{lem:den-i}(\ref{lem:den-i-ordered-density})}} \\
  &= \min\{k_i/(n+1),k_j/(n+1)\}(1-\max\{k_i/(n+1),k_j/(n+1)\}) \frac{n+1}{n} \\
  &= \min\{u_i+\epsilon_i/(n+1),u_j+\epsilon_j/(n+1)\} \\
  &\quad \times (1-\max\{u_i+\epsilon_i/(n+1),u_j+\epsilon_j/(n+1)\}) 
                \left[ 1 + 1/n \right] \\
  &= \left[ \min\{u_i,u_j\} + O(1/n) \right]
     \left[ 1 - \max\{u_i,u_j\} + O(1/n) \right]
     \left[ 1 + O(1/n) \right] \\
\end{split}\\
\label{eqn:Vn-alt}
  &= \overbrace{\min\{u_i,u_j\}(1-\max\{u_i,u_j\})}^{\textrm{Lemma \ref{lem:den}(\ref{lem:den-iii})}} + O(1/n)
\end{align}
since all $u_j\in(0,1)$ and $u_j=k_j/(n+1)+\epsilon_j/(n+1)$ with $\epsilon_j\in[0,1)$.  The error in the overall log density in \eqref{eqn:Y-cond-L-final} resulting from the altered covariance matrix $\underline{\mathcal{V}}/n$ in \eqref{eqn:Vn-alt} is 
\begin{equation*}
O\left( (1/n) n\|(\mathbf{y}-\mathbf{k}/(n+1))^2\|_\infty \right)
  = O\left( a_n^2 n^{-2} \|\mathbf{\Delta k}\|_\infty \right)
  = O\left( n^{-1} [\log n]^2 \right),
\end{equation*}
using $\Delta k_j\asymp n$ for all $j$ and \ref{cond:star} with $a_n=2\log(n)$.  This $O\left( n^{-1} [\log n]^2 \right)$ is smaller than the current remainder of $O\left(n^{-1/2}[\log(n)]^3\right)$.  Thus, using $\underline{\mathcal{V}}_{i,j}=\min\{u_i,u_j\}(1-\max\{u_i,u_j\})$ as in Lemma \ref{lem:den}(\ref{lem:den-iii}) provides a valid approximation. %
}{%
and they equal the entries stated in this lemma up to a smaller-order adjustment. 
}
}{%
yields the joint log PDF of $\mathbf{Y}$ and $\boldsymbol\Lambda$.  Combining this with the marginal log PDF of $\boldsymbol\Lambda$ in \eqref{eqn:LambdaExp} yields the log conditional PDF. %
}

\Supplemental{%
Note that \eqref{eqn:Y-cond-L-final} implies
\begin{align}\notag
f_{\mathbf Y \mid \boldsymbol{\Lambda}} \left( \mathbf{y} \mid \boldsymbol{\lambda} \right)
 &= \phi_{\underline{\mathcal V}/n}\left(\mathbf y - \mathbf k/(n+1)\right) 
    \exp\left\{ O\left(n^{-1/2}[\log(n)]^3\right) \right\}  \\
\label{eqn:Y-PDF-cond-nolog}
 &= \phi_{\underline{\mathcal V}/n}\left(\mathbf y - \mathbf k/(n+1)\right) 
    \left[ 1 + O\left(n^{-1/2}[\log(n)]^3\right) \right] . 
\end{align}
Equation \eqref{eqn:Y-cond-L-final} shows that, given \ref{cond:star} with $a_n=2\log(n)$, conditioning on $\boldsymbol{\Lambda}$ adds only 
$O\left(n^{-1/2}[\log(n)]^3\right)$  
to the unconditional log PDF of $\mathbf{Y}$.  

}{}
For the PDF of $\mathbb{W}$, we can use the formula for the PDF of a transformed random vector and then expand around the $u_j$.  
\Supplemental{%
Let $F^{-1}(\cdot)=Q(\cdot)$ denote the quantile function of interest, and let $F^{-1}(\mathbf{t})\equiv\left(F^{-1}(t_1),\ldots,F^{-1}(t_J)\right)'$ for any vector $\mathbf{t}=(t_1,\ldots,t_J)'$.  First, let $\mathbf{Z}=F^{-1}(\mathbf{Y})$, so $\mathbf{Y}=F(\mathbf{Z})$.  The Jacobian of the (inverse) transformation has elements 
\begin{equation}\label{eqn:FY-J}
\frac{\partial Y_i}{\partial Z_k} = \left\{ \begin{array}{cl} f(Z_k) & \textrm{if }i=k \\ 
                                                0      & \textrm{if }i\ne k
                      \end{array} \right. ,
\end{equation}
making its determinant 
\begin{equation*}
\left| d\mathbf{Y}/d\mathbf{Z}' \right|
  = \prod_{j=1}^{J} f(Z_j) .
\end{equation*}

We now expand both $F(Z_j)$ and $f(Z_j)$ around $Q(u_j)$.  For $F(Z_j)=F(Q(Y_j))$, 
\begin{align*}
F(Z_j)
  &= F(Q(u_j)) + f(Q(u_j))\left[Q(Y_j)-Q(u_j)\right] +(1/2)f'(\tilde{z}_j)\left[Q(Y_j)-Q(u_j)\right]^2 \\
  &= u_j + f(Q(u_j))\left[Q(Y_j)-Q(u_j)\right] 
    +(1/2)f'(\tilde{z}_j)\left[ Q'(\tilde{u}_j)(Y_j-u_j) \right]^2  .
\end{align*}
The argument in the proof of Lemma \ref{lem:den}(\ref{lem:den-bd2}) shows $\tilde{u}_j\to u_j$ in \eqref{eqn:tilde-u-limit}, and thus \ref{a:hut-pf} guarantees uniformly bounded $Q'(\tilde{u}_j)$.  Similarly, \eqref{eqn:tilde-u-limit} shows $\tilde{y}_j\to u_j$, so $\tilde{z}_j\to Q(u_j)$, and \ref{a:hut-pf} guarantees $f'(\tilde{z}_j)$ is uniformly bounded.  Finally, $Y_j-u_j=O\left(n^{-1/2}\log(n)\right)$ by \ref{cond:star}.  
Altogether, uniformly, 
\begin{equation}\label{eqn:F-Zj-exp}
F(Z_j) = u_j + f(Q(u_j)) \left[Z_j-Q(u_j)\right] 
        +O\left(n^{-1}[\log(n)]^2\right) .
\end{equation}

For the $f(Z_j)$ expansion,
\begin{align}\notag
f(Z_j)
  &= f\left(Q(u_j)\right) + f'(\tilde{z}_j) \left[ Z_j-Q(u_j)\right]
   = f\left(Q(u_j)\right) + f'(\tilde{z}_j) Q'(\tilde{u_j})(Y_j-u_j) \\
\label{eqn:fZj-exp}
  &= f\left(Q(u_j)\right) + O\left(n^{-1/2}\log(n)\right)
   = f\left(Q(u_j)\right) \left[ 1 + O\left(n^{-1/2}\log(n)\right) \right] ,
\end{align}
where \ref{a:hut-pf} uniformly bounds $f'(\tilde{z}_j)$ and bounds $f(Q(u_j))$ away from zero, and the arguments for $\tilde{z}_j\to Q(u_j)$, $Q'(\tilde{u}_j)=O(1)$ uniformly, and $Y_j-u_j=O(n^{-1/2}\log(n))$ are identical to those preceding \eqref{eqn:F-Zj-exp}. 

Returning to the transformed PDF formula, using the PDF for $\mathbf{Y}$ from Lemma \ref{lem:den-i}(\ref{lem:den-i-ordered-density}), 
\begin{align*}
\Supplemental{&}{} 
f_{\mathbf{Z}|\boldsymbol{\Lambda}}\left(z_1,\ldots,z_J \mid \boldsymbol{\lambda} \right) 
\Supplemental{\\}{} 
 &= 
 f_{\mathbf{Y}|\boldsymbol{\Lambda}}\left(F(z_1),\ldots,F(z_J) \mid \boldsymbol{\lambda}\right)
     \prod_{j=1}^{J} f(z_j) 
\label{eqn:Z-PDF-exp}\refstepcounter{equation}\tag{\theequation}
     \\
\Supplemental{%
 &= f_{\mathbf{Y}|\boldsymbol{\Lambda}}
     \left(\ldots,
           u_j + f(Q(u_j)) \left[z_j-Q(u_j)\right] +O\left(n^{-1}[\log(n)]^2\right),
           \ldots \mid \boldsymbol{\lambda}
     \right) \\
 &\quad\times 
     \prod_{j=1}^{J} \left[ f\left(Q(u_j)\right) + O\left(n^{-1/2}\log(n)\right) \right] \\
 &= e^D 
    \exp\Biggl\{ -\frac{1}{2} 
        \left( \ldots, \overbrace{u_j - k_j/(n+1)}^{=\epsilon_j/(n+1)=O(n^{-1})} + f(Q(u_j)) \left[z_j-Q(u_j)\right] +O\left(n^{-1}[\log(n)]^2\right), \ldots \right)
        \underline{H} \\
 &\qquad\qquad\qquad\times
        \left( \ldots, u_j - k_j/(n+1) 
                      +f(Q(u_j)) \left[z_j-Q(u_j)\right] 
                      +O\left(n^{-1}[\log(n)]^2\right) , \ldots \right)'
    \Biggr\} \\
 &\quad\times
    \left[ 1 + O\left( n^{-1/2}[\log(n)]^3\right) \right] \\
  &\quad\times
     \left[ \left( \prod_{j=1}^{J}  f\left(Q(u_j)\right) \right) + O\left(n^{-1/2}\log(n)\right) \right] \\
 &= e^D 
    \exp\Biggl\{ -\frac{1}{2} 
        \left( \ldots, f(Q(u_j)) \left[z_j-Q(u_j)\right] , \ldots \right)
        \underline{H} 
        \left( \ldots, f(Q(u_j)) \left[z_j-Q(u_j)\right] , \ldots \right)' \\
 &\qquad\qquad\quad
        +\overbrace{O\left(n^{-1}[\log(n)]^2\right)\overbrace{O(\underline{H})}^{=O(n)}\overbrace{O(z_j-Q(u_j))}^{=O(n^{-1/2}\log(n))}}^{O\left(n^{-1/2}[\log(n)]^3\right)}
    \Biggr\} \\
  &\quad\times
     \left[ \left( \prod_{j=1}^{J}  f\left(Q(u_j)\right) \right) + O\left(n^{-1/2}\log(n)\right) \right]
     \left[ 1 + O\left( n^{-1/2}[\log(n)]^3\right) \right] \\
}{}
\notag
 &= (2\pi)^{-J/2} |\underline{H}|^{1/2} 
    \left[ \prod_{j=1}^{J}  f\left(Q(u_j)\right) \right]
    \exp\left\{ -\frac{ 
        [\mathbf{z}-Q(\mathbf{u})]'
        \underline{\mathcal{A}}\,\underline{H}\,\underline{\mathcal{A}} 
        [\mathbf{z}-Q(\mathbf{u})]
        }{2} \right\} \\
\label{eqn:Z-PDF}\refstepcounter{equation}\tag{\theequation}
 &\quad\times
     \left[ 1 + O\left( n^{-1/2}[\log(n)]^3\right) \right] .
\end{align*}
\Supplemental{Since the determinant of a matrix product is the product of determinants, and $|\underline{\mathcal{A}}|=\left( \prod_{j=1}^{J}  f\left(Q(u_j)\right) \right)$, then 
\begin{equation*}
| \underline{\mathcal{A}}\,\underline{H}\,\underline{\mathcal{A}} | ^{1/2}
  = \left( \prod_{j=1}^{J}  f\left(Q(u_j)\right) \right) |\underline{H}|^{1/2} ,
\end{equation*}
so }{As shown in the supplemental appendix, }\eqref{eqn:Z-PDF} is (up to the remainder shown) a multivariate normal PDF with mean $Q(\mathbf{u})$ and covariance matrix 
\begin{equation*}
\left(\underline{\mathcal{A}}\,\underline{H}\,\underline{\mathcal{A}}\right)^{-1} 
  = \underline{\mathcal{A}}^{-1}
    (\underline{\mathcal{V}}/n)
    \underline{\mathcal{A}}^{-1} .
\end{equation*}

}{%
}%
The transformation from $\mathbf{Z}\equiv Q(\mathbf{Y})$ to $\mathbf{V}\equiv\sqrt{n}\left[Q(\mathbf{Y})-Q(\mathbf{u})\right]$ (from \eqref{eqn:def-lots}) is straightforward centering and $\sqrt{n}$-scaling. %
\Supplemental{%
Using the PDF for $Q(\mathbf{Y})$ in \eqref{eqn:Z-PDF}, the (conditional on $\boldsymbol{\Lambda}$) PDF of $\mathbf{V}$ evaluated at $\mathbf{v}=(v_1,\ldots,v_J)'$ is
\begin{align}
\Supplemental{\label{eqn:V-PDF-transform}
\begin{split}
& f_{\mathbf{V}|\boldsymbol{\Lambda}}(\mathbf{v} \mid \boldsymbol{\lambda}) \\
 &= f_{\mathbf{Z}|\boldsymbol{\Lambda}}\left(\mathbf{v}n^{-1/2}+Q(\mathbf{u}) \mid \boldsymbol{\lambda} \right) 
    \overbrace{\left| \frac{d\left[\mathbf{v}n^{-1/2}+Q(\mathbf{u})\right]}{d\mathbf{v}'} \right|}^{=|n^{-1/2}\underline{I}_J|=n^{-J/2}} 
\end{split}
\\
\notag
\begin{split}
 &= (2\pi)^{-J/2} 
    |\underline{\mathcal{A}}^{-1} (\underline{\mathcal{V}}/n) \underline{\mathcal{A}}^{-1}|^{-1/2} 
    n^{-J/2} \\
 &\quad\times
    \exp\left\{ -\frac{1}{2} 
        \left[ \mathbf{v}n^{-1/2}+Q(\mathbf{u}) - Q(\mathbf{u}) \right]'
        \left[ \underline{\mathcal{A}}^{-1} (\underline{\mathcal{V}}/n) \, \underline{\mathcal{A}}^{-1} \right]^{-1}
        \left[ \mathbf{v}n^{-1/2}+Q(\mathbf{u}) - Q(\mathbf{u}) \right] \right\} \\
 &\quad\times
     \left[ 1 + O\left( n^{-1/2}[\log(n)]^3\right) \right] ,
\end{split}
\\
}{}
\begin{split}\label{eqn:V-PDF}
\Supplemental{}{f_{\mathbf{V}|\boldsymbol{\Lambda}}(\mathbf{v} \mid \boldsymbol{\lambda})}
 &= (2\pi)^{-J/2} 
    |\underline{\mathcal{A}}^{-1} \underline{\mathcal{V}} \underline{\mathcal{A}}^{-1}|^{-1/2} 
    \exp\left\{ -\frac{1}{2} 
        \mathbf{v}'
        \left( \underline{\mathcal{A}}^{-1} \underline{\mathcal{V}} \, \underline{\mathcal{A}}^{-1} \right)^{-1}
        \mathbf{v} \right\} 
\Supplemental{\\&\quad\times}{}
     \left[ 1 + O\left( n^{-1/2}[\log(n)]^3\right) \right] ,
\end{split}
\end{align}
i.e.,\ multivariate normal with mean zero and covariance $\underline{\mathcal{A}}^{-1} \underline{\mathcal{V}} \, \underline{\mathcal{A}}^{-1}$. 

}{}%
The last transformation is from $\mathbf{V}$ to $\mathbb{W}=\boldsymbol{\psi}'\mathbf{V}$, as defined in \eqref{eqn:def-lots}.  
For the special case $J=1$, as in our quantile inference application, this step is trivial since $\mathbb{W}=\mathbf{V}$%
. %
\Supplemental{%
For this transformation, we integrate over $\mathbf{V}$ to get $\mathbb{W}$.  Since the lemma results and proofs are all conditional on \ref{cond:star}, we integrate the conditional PDF over values satisfying \ref{cond:star}. 
With some abuse of notation, let $1\{\mathbf{v}\in\star(a_n)\}$ be the indicator function of whether vector $\mathbf{v}$ satisfies \ref{cond:star}. 
In terms of the unconditional PDF $f_{\mathbf{V}}(\cdot)$, the PDF of $\mathbf{V}$ conditional on \ref{cond:star} is
\begin{align}\notag
f_{\mathbf{V}|\star(a_n)}(\mathbf{v})
  &= \frac{f_{\mathbf{V}}(\mathbf{v})}
          {P\left(\star(a_n)\right)} 
     1\{\mathbf{v}\in\star(a_n)\}
   = \frac{f_{\mathbf{V}}(\mathbf{v}) 1\{\mathbf{v}\in\star(a_n)\}}
          {1-\left[1-P\left(\star(a_n)\right)\right]} 
   = \frac{f_{\mathbf{V}}(\mathbf{v}) 1\{\mathbf{v}\in\star(a_n)\}}
          {1-O(n^{-2})} 
\\&
\label{eqn:V-PDF-star}
   = f_{\mathbf{V}}(\mathbf{v})
     1\{\mathbf{v}\in\star(a_n)\}
     \left[ 1 + O(n^{-2}) \right] , 
\end{align}
with $a_n=2\log(n)$, where the $O(n^{-2})$ probability is from Lemma \ref{lem:den-i}(\ref{lem:den-i-star-prob},\ref{lem:den-i-star-prob-fixed}) and is a constant (i.e.,\ independent of the point of evaluation, $\mathbf{v}$).  

Additional conditioning on $\boldsymbol{\Lambda}$ may be added to the left-hand side and right-hand side of \eqref{eqn:V-PDF-star}.  Since \eqref{eqn:V-PDF-star} has a smaller-order multiplicative error than the $\left[1+O\left(n^{-1/2}[\log(n)]^3\right)\right]$ approximation error from \eqref{eqn:V-PDF},  
\begin{equation*}
\left[ 1 + O(n^{-2}) \right] \left[1+O\left(n^{-1/2}[\log(n)]^3\right)\right]
  = \left[1+O\left(n^{-1/2}[\log(n)]^3\right)\right] . 
\end{equation*}
That is, we can additionally condition on \ref{cond:star} in the left-hand side of \eqref{eqn:V-PDF} without affecting the right-hand side beyond the addition of $1\{\mathbf{v}\in\star(a_n)\}$. 

The PDF of $\mathbb{W}$ conditional on $\boldsymbol{\Lambda}$ and $\star(2\log(n))$ is
\begin{align*}
\label{eqn:W-PDF-first}\refstepcounter{equation}\tag{\theequation}
\begin{split}
& f_{\mathbb{W}|\boldsymbol{\Lambda},\star(2\log(n))}(w \mid \boldsymbol{\lambda}) \\
  &= \idotsint\limits_{\star(2\log(n))} 
     f_{\mathbf{V}|\boldsymbol{\Lambda},\star(2\log(n))}\left(v_1,\ldots,v_{J-1},\frac{w-\psi_1 v_1 - \cdots - \psi_{J-1} v_{J-1}}{\psi_J} \mid \boldsymbol{\lambda} \right) dv_1 \cdots dv_{J-1} \\
\end{split}
\\
  &= \idotsint\limits_{\star(2\log(n))} 
     \phi_{\underline{\mathcal{A}}^{-1} \underline{\mathcal{V}} \, \underline{\mathcal{A}}^{-1}}%
     \left(v_1,\ldots,v_{J-1},\frac{w-\psi_1 v_1 - \cdots - \psi_{J-1} v_{J-1}}{\psi_J} \right) \\
 &\qquad\quad\times
     \left[ 1 + O\left( n^{-1/2}[\log(n)]^3\right) \right]
     dv_1 \cdots dv_{J-1} \\
  &= \idotsint\limits_{\star(2\log(n))} 
     \phi_{\underline{\mathcal{A}}^{-1} \underline{\mathcal{V}} \, \underline{\mathcal{A}}^{-1}}%
     \left(v_1,\ldots,v_{J-1},\frac{w-\psi_1 v_1 - \cdots - \psi_{J-1} v_{J-1}}{\psi_J} \right) 
     dv_1 \cdots dv_{J-1} \\
 &\quad\times
     \left[ 1 + O\left( n^{-1/2}[\log(n)]^3\right) \right] \\
\label{eqn:W-PDF}\refstepcounter{equation}\tag{\theequation}
  &= \phi_{\mathcal{V}_{\psi}}(w)
     \left[ 1 + O\left( n^{-1/2}[\log(n)]^3\right) \right] 
\end{align*}
since the $O(\cdot)$ term is uniform over $\mathbf{v}$ satisfying $\star(2\log(n))$, where $\mathcal{V}_{\psi}=\boldsymbol{\psi}'\underline{\mathcal{A}}^{-1} \underline{\mathcal{V}} \, \underline{\mathcal{A}}^{-1}\boldsymbol{\psi}$ as in the statement of the lemma. 

Going back to \eqref{eqn:Y-PDF-cond-nolog}, so far we have 
}{%
Altogether, up to this point, 
}
\begin{equation}\label{eqn:W-PDF-cond}
f_{\mathbb{W}|\boldsymbol{\Lambda}}(w\mid \boldsymbol{\lambda})
  = \phi_{\mathcal{V}_{\psi}}(w)
     \left[ 1 + O\left( n^{-1/2}[\log(n)]^3\right) \right] ,
\end{equation}
and it remains to account for the difference between $\mathbb{W}$ and $\mathbb{W}_{\boldsymbol{\epsilon},\boldsymbol{\Lambda}}$. 
\Supplemental{For now, we turn to the derivative. }{}

\Supplemental{%
Differentiating \eqref{eqn:YconLamb} with respect to $y_j$, %
\Supplemental{noting $\D{}{y_j}\Delta y_j=1=-\D{}{y_j}\Delta y_{j+1}$,}{it can be shown that}
\begin{align*}
\D{}{y_j}
\log f_{\mathbf{Y}|\boldsymbol{\Lambda}}(\mathbf{y}\mid\boldsymbol{\lambda})
\Supplemental{%
  &= \overbrace{\D{}{y_j} \log f_{\mathbf{Y}}(\mathbf{y})}^{O\left(a_n n \|\mathbf{\Delta k}^{-1/2}\|_\infty \right)}
    -\frac{1}{\Delta y_j} +\frac{1}{\Delta y_{j+1}} \\
 &\quad
    +\frac{\Delta k_j-2}{\Delta y_j-\lambda_{j-1}}
    -\frac{\Delta k_j-2}{\Delta y_j}
    -\frac{\Delta k_{j+1}-2}{\Delta y_{j+1}-\lambda_{j}}
    +\frac{\Delta k_{j+1}-2}{\Delta y_{j+1}} \\
 &= \D{}{y_j} \log f_{\mathbf{Y}}(\mathbf{y})
    +O(1) \\
 &\quad
    +\Delta k_j \frac{\lambda_{j-1}}{\Delta y_j(\Delta y_j-\lambda_{j-1})}
    +\Delta k_{j+1} \frac{-\lambda_{j}}{\Delta y_{j+1}(\Delta y_{j+1}-\lambda_{j})} \\
}{}
 &= \D{}{y_j} \log f_{\mathbf{Y}}(\mathbf{y})
    +\overbrace{O\left(a_n n^{-1} \|\mathbf{\Delta k}\|_\infty\right)}^{O(\log(n))\textrm{ here}} , \\
\Supplemental{%
\D{}{\mathbf{y}}
\log f_{\mathbf{Y}|\boldsymbol{\Lambda}}(\mathbf{y}\mid\boldsymbol{\lambda})
 &= \D{}{\mathbf{y}} \log f_{\mathbf{Y}}(\mathbf{y})
    +O\left(\log(n)\right) \\
 &= \frac{1}{f_{\mathbf{Y}}(\mathbf{y})}
    \D{}{\mathbf{y}} f_{\mathbf{Y}}(\mathbf{y})
    +O\left(\log(n)\right) , \\
}{}
\D{}{\mathbf{y}}
f_{\mathbf{Y}|\boldsymbol{\Lambda}}(\mathbf{y}\mid\boldsymbol{\lambda})
\Supplemental{%
 &= f_{\mathbf{Y}|\boldsymbol{\Lambda}}(\mathbf{y}\mid\boldsymbol{\lambda})
    \D{}{\mathbf{y}}
\log f_{\mathbf{Y}|\boldsymbol{\Lambda}}(\mathbf{y}\mid\boldsymbol{\lambda})
    \\
 &= \overbrace{f_{\mathbf{Y}}(\mathbf{y}) \left[ 1 + O\left(n^{-1/2}\log(n)\right) \right]}^{\textrm{from \eqref{eqn:Y-PDF-cond-log-approx}}} \\
 &\quad\times
    \left[ 
      \frac{1}{f_{\mathbf{Y}}(\mathbf{y})}
      \D{}{\mathbf{y}} f_{\mathbf{Y}}(\mathbf{y})
      +O\left(\log(n)\right)
    \right] \\
 &= \left[ 1 + O\left(n^{-1/2}\log(n)\right) \right]
    \biggl[ 
      \overbrace{\D{}{\mathbf{y}} f_{\mathbf{Y}}(\mathbf{y})}^{=O(n^{(J+1)/2}\log(n))}
      +O\left( n^{J/2} \log(n) \right)
    \biggr] \\
 &= \D{}{\mathbf{y}} f_{\mathbf{Y}}(\mathbf{y})
   +O\left( n^{J/2} \log(n) \right) \\
}{}
\label{eqn:Y-PDF-cond-deriv}\refstepcounter{equation}\tag{\theequation}
 &= \D{}{\mathbf{y}} \phi_{\underline{\mathcal{V}}/n} \left(\mathbf{y} - \mathbf{k}/(n+1) \right)
    +O\left( n^{J/2} [\log(n)]^4 \right) 
\Supplemental{,}{.}
\end{align*}
\Supplemental{where the final normal approximation is from \eqref{eqn:dfX-dphi}. 

}{}
Differentiating the conditional (on $\boldsymbol{\Lambda}$) version of \eqref{eqn:Z-PDF-exp} with respect to $\mathbf{z}$, continuing the notation $F(\mathbf{z})\equiv(F(z_1),\ldots,F(z_J))'$ (etc.),
\begin{align*}
\Supplemental{&}{} 
\D{}{\mathbf{z}} f_{\mathbf{Z}|\boldsymbol{\Lambda}}(z_1,\ldots,z_J \mid \boldsymbol{\lambda}) 
\Supplemental{%
\\&= \D{}{\mathbf{z}} \left[
     f_{\mathbf{Y}|\boldsymbol{\Lambda}}\left(F(z_1),\ldots,F(z_J) \mid \boldsymbol{\lambda}\right)
     \prod_{j=1}^{J} f(z_j) 
    \right] \\
 &= \left[ \D{}{\mathbf{z}} 
     f_{\mathbf{Y}|\boldsymbol{\Lambda}}\left(F(z_1),\ldots,F(z_J) \mid \boldsymbol{\lambda}\right)
    \right] 
     \prod_{j=1}^{J} f(z_j) 
  + f_{\mathbf{Y}|\boldsymbol{\Lambda}}\left(F(z_1),\ldots,F(z_J) \mid \boldsymbol{\lambda}\right)
     \D{}{\mathbf{z}} \left[ \prod_{j=1}^{J} f(z_j) \right] \\
 &=  \left( \ldots, \D{}{\mathbf{y}} f_{\mathbf{Y}|\boldsymbol{\Lambda}}\left(F(z_1),\ldots,F(z_J) \mid \boldsymbol{\lambda}\right) f(z_j), \ldots \right)'
     \prod_{j=1}^{J} f(z_j) \\
 &\quad
  + f_{\mathbf{Y}|\boldsymbol{\Lambda}}\left(F(z_1),\ldots,F(z_J) \mid \boldsymbol{\lambda}\right)
     \left[ \prod_{j=1}^{J} f(z_j) \right]
     \overbrace{\left( f'(z_1)/f(z_1), \ldots, f'(z_J)/f(z_J) \right)'}^{=O(1)\textrm{ by \ref{a:hut-pf}, $z_j\to Q(u_j)$}} \\
 &= \overbrace{\underline{\mathcal{A}}\left[ 1+O\left(n^{-1/2}\log(n)\right) \right]}^{\textrm{from \eqref{eqn:fZj-exp}}}
    \left[ \overbrace{\D{}{\mathbf{y}} \phi_{\underline{\mathcal{V}}/n} \left( F(\mathbf{z}) - \mathbf{k}/(n+1) \right)
     + O\left(n^{J/2}[\log(n)]^4\right) }^{\textrm{from \eqref{eqn:Y-PDF-cond-deriv}}} 
     \right]
      \\
 &\quad\times
     \overbrace{\prod_{j=1}^{J} f(Q(u_j))\left[ 1 + O\left(n^{-1/2}\log(n)\right) \right]}^{\textrm{from \eqref{eqn:Z-PDF}}} \\
 &\quad
  + \overbrace{\phi_{\underline{\mathcal{A}}^{-1}(\underline{\mathcal{V}}/n)\underline{\mathcal{A}}^{-1}}%
         \left( \mathbf{z} - Q(\mathbf{u}) \right) 
     \left[ 1 + O\left( n^{-1/2}[\log(n)]^3\right) \right]}^{\textrm{from \eqref{eqn:Z-PDF-exp}, \eqref{eqn:Z-PDF}}} 
     O(1) \\
 &= -\underline{\mathcal{A}}\,\underline{H}\left( F(\mathbf{z}) - \mathbf{u} +O(1/n) \right)
    \phi_{\mathcal{V}/n}%
         \left( F(\mathbf{z}) - \mathbf{u} +O(1/n) \right) 
    \left| \underline{\mathcal{A}} \right|
    \left[ 1 + O\left(n^{-1/2}[\log(n)]^3\right) \right] 
     \\
 &\quad
   + O\left( n^{J/2}[\log(n)]^4 \right) 
   + O(n^{J/2})
     \\
 &= -\underline{\mathcal{A}}\,\underline{H}
     \overbrace{\underline{\mathcal{A}} \left( \mathbf{z}-Q(\mathbf{u}) \right) \left[ 1+O\left(n^{-1/2}\log(n)\right) \right]}^{\textrm{from \eqref{eqn:F-Zj-exp2}}}
    \phi_{\underline{\mathcal{V}}/n}%
         \left( F(\mathbf{z}) - \mathbf{u} \right) \\
 &\qquad\times
    \left[ 1 + \overbrace{O\left([F(\mathbf{z})-\mathbf{u}]\underline{H}n^{-1}\right)}^{=O\left(n^{-1/2}\log(n)\right)\textrm{ by }\star(a_n)} \right] 
    \left| \underline{\mathcal{A}} \right|
    \left[ 1 + O\left(n^{-1/2}[\log(n)]^3\right) \right] 
     \\
 &\quad
  + O\left( n^{J/2} [\log(n)]^4 \right) \\
 &= -\underline{\mathcal{A}}\,\underline{H}\,\underline{\mathcal{A}}
     \left( \mathbf{z}-Q(\mathbf{u}) \right) 
     \left| \underline{\mathcal{A}} \right|
     \phi_{\underline{\mathcal{V}}/n}%
         \left( \overbrace{\underline{\mathcal{A}} \left( \mathbf{z}-Q(\mathbf{u}) \right) +O\left(n^{-1}[\log(n)]^2\right)}^{\textrm{from \eqref{eqn:F-Zj-exp}}} \right) \\
 &\qquad\times
    \left[ 1 + O\left(n^{-1/2}[\log(n)]^3\right) \right] 
     \\
 &\quad
  + O\left( n^{J/2} [\log(n)]^4 \right) \\
 &= -\underline{\mathcal{A}}\,\underline{H}\,\underline{\mathcal{A}}
     \left( \mathbf{z}-Q(\mathbf{u}) \right) 
     \left| \underline{\mathcal{A}} \right|
     \phi_{\underline{\mathcal{V}}/n}%
         \left( \underline{\mathcal{A}} \left[ \mathbf{z}-Q(\mathbf{u}) \right] \right) 
  + O\left( n^{J/2}[\log(n)]^2 \right) + O\left( n^{J/2}[\log(n)]^4 \right) \\
\label{eqn:dfZ-dz}\refstepcounter{equation}\tag{\theequation}
 &= -\underline{\mathcal{A}}\,\underline{H}\,\underline{\mathcal{A}}
     \left( \mathbf{z}-Q(\mathbf{u}) \right) 
     \phi_{\underline{\mathcal{A}}^{-1}(\underline{\mathcal{V}}/n)\underline{\mathcal{A}}^{-1}}%
         \left( \mathbf{z}-Q(\mathbf{u}) \right) 
   +O\left( n^{J/2}[\log(n)]^4 \right) \\
}{}
 &= \D{}{\mathbf{z}} 
     \phi_{\underline{\mathcal{A}}^{-1}(\underline{\mathcal{V}}/n)\underline{\mathcal{A}}^{-1}}%
         \left( \mathbf{z}-Q(\mathbf{u}) \right) 
   +O\left( n^{J/2}[\log(n)]^4 \right) .
\end{align*}
\Supplemental{%
The result $f'(z_j)/f(z_j)=O(1)$ is from the combination of $z_j\to Q(u_j)$ under \ref{cond:star} with Assumption \ref{a:hut-pf}. 
From \eqref{eqn:F-Zj-exp}, noting that the original remainder term was $(1/2)f'(\tilde z_j)[Q(Y_j)-Q(u_j)]^2$, 
\begin{equation}\label{eqn:F-Zj-exp2}
F(Z_j) - u_j = f\left(Q(u_j)\right) \left[ Z_j - Q(u_j)\right] \left[ 1 + O\left(n^{-1/2}\log(n)\right)\right] .
\end{equation}
}{}

As before, the transformation from $\mathbf{Z}\equiv Q(\mathbf{Y})$ to $\mathbf{V}\equiv\sqrt{n}\left[Q(\mathbf{Y})-Q(\mathbf{u})\right]$ is simply centering and scaling%
\Supplemental{, so 
\Supplemental{using \eqref{eqn:V-PDF-transform} and \eqref{eqn:dfZ-dz}, }{}%
the (conditional on $\boldsymbol{\Lambda}$) PDF derivative becomes
\begin{align*}
& \D{}{\mathbf{v}} f_{\mathbf{V}|\boldsymbol{\Lambda}}(\mathbf{v}\mid\boldsymbol{\lambda}) \\
\Supplemental{%
 &= \D{}{\mathbf{v}} f_{\mathbf{Z}|\boldsymbol{\Lambda}}\left(\mathbf{v}n^{-1/2}+Q(\mathbf{u}) \mid \boldsymbol{\lambda}\right) n^{-J/2} \\
 &= n^{-1/2} 
    \left. \D{}{\mathbf{z}} f_{\mathbf{Z}|\boldsymbol{\Lambda}}\left(\mathbf{z} \mid \boldsymbol{\lambda}\right) \right|_{\mathbf{z}=\mathbf{v}n^{-1/2}+Q(\mathbf{u})}
    n^{-J/2} \\
 &= -n^{-1/2} n^{-J/2} 
     \underline{\mathcal{A}}\,\underline{H}\,\underline{\mathcal{A}}
     \left[ \mathbf{v}n^{-1/2}+Q(\mathbf{u})-Q(\mathbf{u}) \right] \\
 &\quad\times
     \phi_{\underline{\mathcal{A}}^{-1}(\underline{\mathcal{V}}/n)\underline{\mathcal{A}}^{-1}}%
         \left( \mathbf{v}n^{-1/2}+Q(\mathbf{u})-Q(\mathbf{u}) \right) \\
 &\quad
   +O\left( n^{-(J+1)/2} n^{J/2}[\log(n)]^4 \right) \\
 &= -n^{-1} 
    \left[ \underline{\mathcal{A}}^{-1}(\underline{\mathcal{V}}/n)\underline{\mathcal{A}}^{-1} \right]^{-1}
    \mathbf{v}
    n^{-J/2} \\
 &\quad\times
   (2\pi)^{-J/2}
   \left| \left[ \underline{\mathcal{A}}^{-1}(\underline{\mathcal{V}}/n)\underline{\mathcal{A}}^{-1} \right]^{-1} \right| 
   \exp\left\{ -(1/2)
      \mathbf{v}'n^{-1/2} 
      \left[ \underline{\mathcal{A}}^{-1}(\underline{\mathcal{V}}/n)\underline{\mathcal{A}}^{-1} \right]^{-1}
      n^{-1/2}\mathbf{v} \right\} \\
 &\quad
   +O\left( n^{-1/2} [\log(n)]^4 \right) \\
 &= -\left[ \underline{\mathcal{A}}^{-1} \underline{\mathcal{V}}\,\underline{\mathcal{A}}^{-1} \right]^{-1}
    \mathbf{v} 
    (2\pi)^{-J/2}
   \left| \left( \underline{\mathcal{A}}^{-1} \underline{\mathcal{V}} \, \underline{\mathcal{A}}^{-1} \right)^{-1} \right| 
   \exp\left\{ -(1/2)
      \mathbf{v}'
      \left( \underline{\mathcal{A}}^{-1} \underline{\mathcal{V}} \, \underline{\mathcal{A}}^{-1} \right)^{-1}
      \mathbf{v} \right\} \\
 &\quad
   +O\left( n^{-1/2} [\log(n)]^4 \right) \\
}{}
 &= \D{}{\mathbf{v}}
    \phi_{\underline{\mathcal{A}}^{-1} \underline{\mathcal{V}} \, \underline{\mathcal{A}}^{-1}}%
         \left( \mathbf{v} \right) 
   +O\left( n^{-1/2}[\log(n)]^4 \right) .
\end{align*}
}{.}

\Supplemental{%
Finally, taking the linear combination $\mathbb{W}\equiv \boldsymbol{\psi}'\mathbf{V}$, we need to invoke \ref{cond:star} more directly to account for the error in $\D{}{\mathbf{v}}f_{\mathbf{V}}(\mathbf{v})$ being additive rather than multiplicative. 
First, from \eqref{eqn:V-PDF-star}, the effect on the PDF of $\mathbf{V}$ of conditioning on $\star(2\log(n))$ is a multiplicative $[1+O(n^{-2})]$ error. 
Second, taking a derivative of \eqref{eqn:V-PDF-star}, since the $O(n^{-2})$ does not depend on $\mathbf{v}$,
\begin{align*}
\D{}{\mathbf{v}}f_{\mathbf{V}|\star(a_n)}(\mathbf{v})
  &= \D{}{\mathbf{v}}f_{\mathbf{V}}(\mathbf{v})
     1\{\mathbf{v}\in\star(a_n)\}
     \left[ 1 + O(n^{-2}) \right] . 
\end{align*}
Third, using \eqref{eqn:Yj-cond-star}, \ref{cond:star} with $a_n=2\log(n)$ and all $\Delta k_j\asymp n^{1/2}$ implies 
\begin{equation*}
Y_j = u_j +O(n^{-1/2}\log(n)) . 
\end{equation*}
Since $Y_j\to u_j$, for large enough $n$, Assumption \ref{a:hut-pf} ensures bounded quantile function derivatives, so
\begin{equation*}
Q(Y_j) = Q(u_j) + O(n^{-1/2}\log(n)) , \quad
\sqrt{n}[Q(Y_j)-Q(u_j)] = O\left(\log(n)\right) , 
\end{equation*}
uniformly over $j$ since $J$ is fixed.  Thus, integrating over all values of $\mathbf{v}$ satisfying \ref{cond:star} with $a_n=2\log(n)$ is integrating over a hypercube with volume $O\left([\log(n)]^J\right)$.  Below, integrating over $J-1$ of the vector components corresponds to a hypercube with volume $O\left([\log(n)]^{J-1}\right)$. 
}{%
With $J=1$, the transformation to $\mathbb{W}\equiv \boldsymbol{\psi}'\mathbf{V}$ is again trivial; details for $J>1$ are in the supplemental appendix.
}
}{}

\Supplemental{Combining the preceding points, following \eqref{eqn:W-PDF-first}, }{Altogether, it can be shown that }%
the PDF of $\mathbb{W}$ conditional on \ref{cond:star} and $\boldsymbol{\Lambda}$ is
\begin{align*}
& f_{\mathbb{W}|\star(a_n),\boldsymbol{\Lambda}}(w \mid \boldsymbol{\lambda}) 
\Supplemental{\\&}{}
= \idotsint\limits_{\star(a_n)} f_{\mathbf{V}|\star(a_n),\boldsymbol{\Lambda}}\left(v_1,\ldots,v_{J-1},\frac{w-\psi_1 v_1 - \cdots - \psi_{J-1} v_{J-1}}{\psi_J} \mid \boldsymbol{\lambda} \right) 
  \,dv_1 \cdots dv_{J-1} 
\Supplemental{, }{.}
\Supplemental{%
\\&
\D{}{w}
f_{\mathbb{W}|\star(a_n),\boldsymbol{\Lambda}}(w\mid\boldsymbol{\lambda}) 
\\&= \D{}{w}
     \idotsint\limits_{\star(a_n)} f_{\mathbf{V}|\star(a_n),\boldsymbol{\Lambda}}\left(v_1,\ldots,v_{J-1},\frac{w-\psi_1 v_1 - \cdots - \psi_{J-1} v_{J-1}}{\psi_J} \mid\boldsymbol{\lambda}\right) 
     \,dv_1 \cdots dv_{J-1} \\
  &= \idotsint\limits_{\star(a_n)} 
     \D{}{v_J}
      f_{\mathbf{V}|\star(a_n),\boldsymbol{\Lambda}}
        \left(v_1,\ldots,v_{J-1},\frac{w-\psi_1 v_1 - \cdots - \psi_{J-1} v_{J-1}}{\psi_J} \mid\boldsymbol{\lambda} \right) 
     \psi_J^{-1} 
     \,dv_1 \cdots dv_{J-1} \\
  &= \idotsint\limits_{\star(a_n)} \left[ 
  \D{}{v_J}
    \phi_{\underline{\mathcal{A}}^{-1} \underline{\mathcal{V}} \, \underline{\mathcal{A}}^{-1}}%
         \left( \mathbf{v} \right) 
   +O\left( n^{-1/2}[\log(n)]^4 \right)
  \right] \psi_J^{-1} 
    \,dv_1 \cdots dv_{J-1} \\
  &= \idotsint\limits_{\star(a_n)} \left[ 
  \D{}{v_J}
    \phi_{\underline{\mathcal{A}}^{-1} \underline{\mathcal{V}} \, \underline{\mathcal{A}}^{-1}}%
         \left( \mathbf{v} \right) 
  \right] \psi_J^{-1} 
    \,dv_1 \cdots dv_{J-1}
    + O\left( n^{-1/2}[\log(n)]^{4+J-1} \right) \\
\label{eqn:W-PDF-derivative}\refstepcounter{equation}\tag{\theequation}
  \Supplemental{&}{}
  = \phi'_{\mathcal{V}_\psi}(w)
    + O\left( n^{-1/2}[\log(n)]^{4+J-1} \right) , 
}{}
\end{align*}
\Supplemental{%
using $a_n=2\log(n)$ and the normal distribution of a linear combination of jointly normal random variables. %
}{}


To transition to $\mathbb{W}_{\epsilon,\Lambda}$, 
\Supplemental{we use an extension of the law of iterated expectations to PDFs and PDF derivatives as shown in, e.g.,\ \citet[p.\ 45]{KaplanSun2016}:
\begin{equation}\label{eqn:LIE-PDF}
\begin{split}
 E\left[ f_{U|\mathbf{Z},\mathbf{X}}(u\mid \mathbf{Z},\mathbf{X}) \mid \mathbf{Z}=\mathbf{z} \right] = f_{U|\mathbf{Z}}(u\mid \mathbf{z}) , \\
 E\left[ f'_{U|\mathbf{Z},\mathbf{X}}(u\mid \mathbf{Z},\mathbf{X}) \mid \mathbf{Z}=\mathbf{z} \right] = f'_{U|\mathbf{Z}}(u\mid \mathbf{z}) . 
\end{split}
\end{equation}
Define}{define}
\Supplemental{\begin{equation*}}{$}
\eta = \sqrt{n} \sum_{j=1}^{J} \epsilon_j \psi_j \Lambda_j \left[ Q'(u_j) + Q''(u_j)(Y_j-u_j)\right] 
\Supplemental{,\end{equation*}}{$,}
so $\mathbb{W}_{\epsilon,\Lambda}=\mathbb{W}+\eta$. 
Conditional on $\mathbb{W}=w$, $Y_1=y_1$, \ldots, $Y_{J-1}=y_{J-1}$, the value of $Y_J$ is fully determined\Supplemental{ as
\begin{equation*}
Y_J = F\left( \frac{w-\sqrt{n}\sum_{j=1}^{J-1}\psi_j\left[Q(y_j)-Q(u_j)\right]}{\sqrt{n}\psi_J} + Q(u_J) \right) .
\end{equation*}
}{.}
Additionally conditioning on $\boldsymbol{\Lambda}=\boldsymbol{\lambda}$, the value of $\eta$ is fully determined.  %
\Supplemental{%
That is, the conditional distribution of $\eta$ given $(\mathbb{W},\boldsymbol{\Lambda},Y_1,\ldots,Y_{J-1})$ is degenerate, with 
$P(\eta=\tilde\eta)=1$ for 
\begin{equation}\label{eqn:eta-tilde}
\begin{split}
\tilde{\eta}
  &= 
   \sqrt{n}\epsilon_J\psi_J\lambda_J \\
  &\quad\times
   \left\{
     Q'(u_J)+Q''(u_J)\left[ 
       \overbrace{%
       F\left( \frac{w-\sqrt{n}\sum_{j=1}^{J-1}\psi_j\left[Q(y_j)-Q(u_j)\right]}{\sqrt{n}\psi_J} + Q(u_J) \right)%
       }^{=y_J}
       - u_J
     \right]
    \right\} \\
  &\quad+
  \sqrt{n}\sum_{j=1}^{J-1}\epsilon_j\psi_j\lambda_j\left[Q'(u_j)+Q''(u_j)(y_j-u_j)\right] .
\end{split}
\end{equation}
Applying \eqref{eqn:LIE-PDF},
\begin{align*}
f_{\mathbb{W}_{\epsilon,\Lambda}|\boldsymbol{\Lambda}}(\tilde{w}\mid\boldsymbol{\Lambda}=\boldsymbol{\lambda})
  &= E\left[ f_{\mathbb{W}_{\epsilon,\Lambda}|\boldsymbol{\Lambda},Y_1,\ldots,Y_{J-1}}(\tilde{w}\mid \boldsymbol{\lambda}, y_1,\ldots, y_{J-1}) \mid \boldsymbol{\lambda} \right] \\
  &= E\Biggl[\; \int\limits_{w+\eta=\tilde{w}}\! f_{\mathbb{W}|\boldsymbol{\Lambda},Y_1,\ldots,Y_{J-1}}(w \mid \boldsymbol{\lambda},y_1,\ldots,y_{J-1}) \\ 
  &\qquad\qquad\qquad
    \overbrace{dF_{\eta|\mathbb{W},\Lambda,Y_1,\ldots,Y_{J-1}}(\eta\mid w,\boldsymbol{\lambda},y_1,\ldots,y_{J-1})}^{\textrm{degenerate: \eqref{eqn:eta-tilde} with }w=\tilde{w}-\tilde{\eta}} \mid \boldsymbol{\lambda} \Biggr] \\
  &= E\left[ f_{\mathbb{W}|\boldsymbol{\Lambda},Y_1,\ldots,Y_{J-1}}(\tilde{w}-\tilde{\eta}\mid\boldsymbol{\lambda},y_1,\ldots,y_{J-1}) \mid \boldsymbol{\lambda} \right] , 
\end{align*}
where $\tilde{\eta}$ solves \eqref{eqn:eta-tilde} after plugging in $w=\tilde{w}-\tilde\eta$, and using \eqref{eqn:W-PDF-cond} for the normal approximation.  Since $\tilde{\eta}$ depends on the $y_j$ values, we treat it before integrating over $(Y_1,\ldots,Y_{J-1})$, using the fact that $\mathbb{W}$ has an approximately normal PDF given $\boldsymbol{\Lambda}$ and $Y_1,\ldots,Y_{J-1}$ by application of \eqref{eqn:W-PDF-cond} with $J=1$ to the term $\sqrt{n}\psi_J\left[Q(Y_J)-Q(u_J)\right]$; the other terms in $\mathbb{W}$ only serve to recenter it since they are (conditionally) fixed.  That is, given $Y_j=y_j$ for $j=1,\ldots,J-1$,
\begin{equation}\label{eqn:W-singleJ}
\mathbb{W}
  = \psi_J \overbrace{\sqrt{n}\left[Q(Y_J)-Q(u_J)\right]}^{\textrm{random; PDF in \eqref{eqn:W-PDF-cond}}}
   +\overbrace{\sqrt{n} \sum_{j=1}^{J-1}\psi_j \left[Q(y_j)-Q(u_j)\right]}^{\equiv\mu\textrm{; fixed, }O(\log(n))} . 
\end{equation}
By \ref{cond:star}, 
\begin{equation}\label{eqn:eta-order}
\tilde{\eta}=O\left(n^{-1/2}\log(n)\right), \quad
\tilde{w}=O\left(\log(n)\right) , \quad
\mu=O(\log(n)) . 
\end{equation}
Letting $\mu$ (as defined in \eqref{eqn:W-singleJ}) and $\sigma^2$ denote the conditional normal mean and variance, 
\begin{align*}
& f_{\mathbb{W}|\boldsymbol{\Lambda},Y_1,\ldots,Y_{J-1}}(\tilde{w}-\tilde{\eta}\mid\boldsymbol{\lambda},y_1,\ldots,y_{J-1}) \\
&= (2\pi\sigma^2)^{-1/2} 
     \exp\left\{ -(1/2)\sigma^{-2}(\tilde{w}-\tilde{\eta}-\mu)^2 \right\}
     \left[1+O\left(n^{-1/2}[\log(n)]^3\right)\right] \\
&= (2\pi\sigma^2)^{-1/2} 
     \exp\left\{ -(1/2)\sigma^{-2}(\tilde{w}-\mu)^2 +O\left(\tilde{\eta}(\tilde{w}-\mu)\right) \right\}
     \left[1+O\left(n^{-1/2}[\log(n)]^3\right)\right] \\
&= (2\pi\sigma^2)^{-1/2} 
     \exp\left\{ -(1/2)\sigma^{-2}(\tilde{w}-\mu)^2 \right\}
     [1+\overbrace{O\left(\tilde{\eta}(\tilde{w}-\mu)\right)}^{O\left(n^{-1/2}[\log(n)]^2\right)}]
     \left[1+O\left(n^{-1/2}[\log(n)]^3\right)\right] \\
&= f_{\mathbb{W}|\boldsymbol{\Lambda},Y_1,\ldots,Y_{J-1}}(\tilde{w}\mid\boldsymbol{\lambda},y_1,\ldots,y_{J-1})
     \left[1+O\left(n^{-1/2}[\log(n)]^3\right)\right] , \\
& f_{\mathbb{W}_{\epsilon,\Lambda}|\boldsymbol{\Lambda}}(\tilde{w} \mid \boldsymbol{\lambda}) \\
&= E\left[ f_{\mathbb{W}|\boldsymbol{\Lambda},Y_1,\ldots,Y_{J-1}}(\tilde{w}-\tilde{\eta}\mid\boldsymbol{\lambda},y_1,\ldots,y_{J-1}) \mid \boldsymbol{\lambda} \right] \\
&= E\left[ f_{\mathbb{W}|\boldsymbol{\Lambda},Y_1,\ldots,Y_{J-1}}(\tilde{w}\mid\boldsymbol{\lambda},y_1,\ldots,y_{J-1})
     \left[1+O\left(n^{-1/2}[\log(n)]^3\right)\right] \mid \boldsymbol{\lambda} \right] \\
&= f_{\mathbb{W}|\boldsymbol{\Lambda}}(\tilde{w}\mid \boldsymbol{\lambda})
     \left[1+O\left(n^{-1/2}[\log(n)]^3\right)\right] \\
\label{eqn:W-EPS-PDF-cond}\refstepcounter{equation}\tag{\theequation}
&= \phi_{\mathcal{V}_\psi}(\tilde{w})
     \left[ 1 + O\left(n^{-1/2}[\log(n)]^3\right) \right] , 
\end{align*}
by applying \eqref{eqn:LIE-PDF} and \eqref{eqn:W-PDF-cond}. 
}{%
Along with the implicit function theorem, this can be used to derive the final normal approximation. %
}

\Supplemental{%
For the derivative, using the same notation and applying \eqref{eqn:LIE-PDF},
\begin{align*}
f'_{\mathbb{W}_{\epsilon,\Lambda}|\boldsymbol{\Lambda}}(\tilde{w}\mid\boldsymbol{\Lambda}=\boldsymbol{\lambda})
  &= E\left[ f'_{\mathbb{W}_{\epsilon,\Lambda}|\boldsymbol{\Lambda},Y_1,\ldots,Y_{J-1}}(\tilde{w}\mid \boldsymbol{\lambda}, y_1,\ldots, y_{J-1} \mid \boldsymbol{\lambda} \right] \\
  &= E\Biggl[\D{}{\tilde{w}} \; \int\limits_{w+\eta=\tilde{w}}\! f_{\mathbb{W}|\boldsymbol{\Lambda},Y_1,\ldots,Y_{J-1}}(w \mid \boldsymbol{\lambda},y_1,\ldots,y_{J-1}) \\ 
  &\qquad\qquad\qquad\qquad
    \overbrace{dF_{\eta|\mathbb{W},\Lambda,Y_1,\ldots,Y_{J-1}}(\eta\mid w,\boldsymbol{\lambda},y_1,\ldots,y_{J-1})}^{\textrm{degenerate: \eqref{eqn:eta-tilde} with }w=\tilde{w}-\tilde{\eta}} \mid \boldsymbol{\lambda} \Biggr] \\
  &= E\left\{ f'_{\mathbb{W}|\boldsymbol{\Lambda},Y_1,\ldots,Y_{J-1}}(\tilde{w}-\tilde{\eta}\mid\boldsymbol{\lambda},y_1,\ldots,y_{J-1})
              \D{\tilde{w}-\tilde{\eta}}{\tilde{w}} 
      \mid \boldsymbol{\lambda} \right\} \\
\label{eqn:W-EPS-PDF-cond-deriv-int}\refstepcounter{equation}\tag{\theequation}
  &= E\left\{ f'_{\mathbb{W}|\boldsymbol{\Lambda},Y_1,\ldots,Y_{J-1}}(\tilde{w}-\tilde{\eta}\mid\boldsymbol{\lambda},y_1,\ldots,y_{J-1})
            \left[ 1 + O\left(n^{-1}\log(n)\right) \right]
      \mid \boldsymbol{\lambda} \right\} . 
\end{align*}
The last line is from applying the implicit function theorem to \eqref{eqn:eta-tilde} with $w=\tilde{w}-\tilde{\eta}$:
\begin{align*}
0 &= h(\tilde\eta,\tilde{w}) \\
  &\equiv    \sqrt{n}\epsilon_J\psi_J\lambda_J \\
  &\quad\times
   \Biggl\{
     Q''(u_J)\left[ 
       F\left( \frac{\tilde{w}-\tilde{\eta}-\sqrt{n}\sum_{j=1}^{J-1}\psi_j\left[Q(y_j)-Q(u_j)\right]}{\sqrt{n}\psi_J} + Q(u_J) \right)
       - u_J
     \right] 
\\&\qquad\quad+Q'(u_J)
    \Biggr\} \\
  &\quad+
  \sqrt{n}\sum_{j=1}^{J-1}\epsilon_j\psi_j\lambda_j\left[Q'(u_j)+Q''(u_j)(y_j-u_j)\right] 
  -\tilde{\eta} , \\
\D{h}{\tilde{w}}
  &= \sqrt{n}\epsilon_J\psi_J\lambda_J Q''(u_J) f\left(\frac{\tilde{w}-\tilde{\eta}-\sqrt{n}\sum_{j=1}^{J-1}\psi_j[Q(y_j)-Q(u_j)]}{\sqrt{n}\psi_J} + Q(u_J)\right)\frac{1}{\sqrt{n}\psi_J} \\
  &= O\left(n^{-1}\log(n)\right) , \\
\D{h}{\tilde{\eta}}
  &= -\D{h}{\tilde{w}} - 1 , \\
\D{\tilde{\eta}}{\tilde{w}}
  &= -\frac{\partial h/\partial\tilde{w}}{\partial h/\partial\tilde{\eta}}
   = O\left(n^{-1}\log(n)\right) , \\
\D{(\tilde{w}-\tilde{\eta})}{\tilde{w}}
  &= 1 + O\left(n^{-1}\log(n)\right) , 
\end{align*}
having applied \ref{a:hut-pf} to bound $Q''(u_J)$ and \ref{cond:star} to get $\lambda_J=O(n^{-1}\log(n))$. 

Using the decomposition in \eqref{eqn:W-singleJ} and the same $\mu$ and $\sigma^2$ notation, applying to it the normal PDF derivative approximation in \eqref{eqn:W-PDF-derivative} with $J=1$, and continuing from \eqref{eqn:W-EPS-PDF-cond-deriv-int},
\begin{align*}
& f'_{\mathbb{W}_{\epsilon,\Lambda}|\boldsymbol{\Lambda}}(\tilde{w}\mid\boldsymbol{\lambda}) 
\Supplemental{%
\\&= E\left\{ f'_{\mathbb{W}|\boldsymbol{\Lambda},Y_1,\ldots,Y_{J-1}}(\tilde{w}-\tilde{\eta}\mid\boldsymbol{\lambda},y_1,\ldots,y_{J-1})
      \mid \boldsymbol{\lambda} \right\}
    \left[ 1 + O\left(n^{-1}\log(n)\right) \right] 
    \\
&= E\left\{ \left[ \phi'_{\sigma^2}(\tilde{w}-\tilde{\eta}-\mu) +O\left(n^{-1/2}[\log(n)]^{4}\right) \right]
      \mid \boldsymbol{\lambda} \right\}
    \left[ 1 + O\left(n^{-1}\log(n)\right) \right] 
    \\
&= E\left\{ 
      -\sigma^{-2}(\tilde{w}-\tilde{\eta}-\mu)
       \exp\left\{ -(1/2)(\tilde{w}-\tilde{\eta}-\mu)^2 / \sigma^2 \right\}
      \mid \boldsymbol{\lambda} \right\}
    +O\left(n^{-1/2}[\log(n)]^{4}\right)
    \\
&= E\left\{ 
      -\sigma^{-2}(\tilde{w}-\mu)
       \exp\bigl\{ -(1/2)(\tilde{w}-\mu)^2 / \sigma^2 
                   +\overbrace{O\left(\log(n)n^{-1/2}\log(n)\right)}^{\textrm{from \eqref{eqn:eta-order}}} \bigr\}
      \mid \boldsymbol{\lambda} \right\} \\
&\quad
    +\overbrace{O\left(n^{-1/2}\log(n)\right)}^{\textrm{from \eqref{eqn:eta-order}}} 
    +O\left(n^{-1/2}[\log(n)]^{4}\right)
    \\
&= E\left\{ 
      -\sigma^{-2}(\tilde{w}-\mu)
       \exp\left\{ -(1/2)(\tilde{w}-\mu)^2 / \sigma^2 \right\}
       \left[ 1 + O\left(n^{-1/2}[\log(n)]^2\right) \right]
      \mid \boldsymbol{\lambda} \right\} \\
&\quad
    +O\left(n^{-1/2}[\log(n)]^{4}\right)
    \\
&= E\left\{ 
      f'_{\mathbb{W}|\boldsymbol{\Lambda},Y_1,\ldots,Y_{J-1}}(\tilde{w}\mid\boldsymbol{\lambda},y_1,\ldots,y_{J-1})
      \mid \boldsymbol{\lambda} \right\}
    +O\left(n^{-1/2}[\log(n)]^{4}\right)
    \\
&= f'_{\mathbb{W}|\boldsymbol{\Lambda}}(\tilde{w}\mid\boldsymbol{\lambda})
    +O\left(n^{-1/2}[\log(n)]^{4}\right) \\
}{}
\label{eqn:W-EPS-PDF-cond-derivative}\refstepcounter{equation}\tag{\theequation}
\Supplemental{&}{}
= \phi'_{\mathcal{V}_\psi}(\tilde{w})
    +O\left(n^{-1/2}[\log(n)]^{4+J-1}\right) 
\Supplemental{,}{.} 
\end{align*}
\Supplemental{using \eqref{eqn:W-PDF-derivative} and \eqref{eqn:LIE-PDF}. }{}
}{%
Results for the PDF derivative follow the same sequence of transformations; details are left to the supplemental appendix. 
}

For the last result in this part of the lemma, %
\Supplemental{}{in addition to Condition $\star(2\log(n))$ and \ref{a:hut-pf}, }%
we use the law of iterated expectations for CDFs%
\Supplemental{,
\begin{align}\notag
E\left[ F_{U|Z,X}(u\mid Z,X) \mid Z=z\right]
&= E\left\{ E\left[ 1\{U<u\} \mid Z,X\right] \mid Z=z\right\} 
= E\left[ 1\{U<u\} \mid Z=z\right] \\
\label{eqn:LIE-CDF}
&= F_{U|Z}(u\mid z) . 
\end{align}
}{.}
\Supplemental{%
Defining 
\begin{equation*}
\mathbb{W}_{\boldsymbol{\epsilon},\boldsymbol{\Lambda},-j}
 \equiv 
 \mathbb{W} 
+\sqrt{n} \sum_{t\ne j} \epsilon_t \psi_t \lambda_t \left[ Q'(u_t)+Q''(u_t)(Y_t-u_t)\right] , 
\end{equation*}
we have 
\begin{align*}
& \frac{\partial^2}{\partial \epsilon_j^2}
  F_{\mathbb{W}_{\boldsymbol{\epsilon},\boldsymbol{\Lambda}}|\boldsymbol{\Lambda}} ( K \mid \boldsymbol{\lambda} ) \\
  &= \D{}{\epsilon_j} \left[
       \D{}{\epsilon_j} P\left( \mathbb{W}_{\boldsymbol{\epsilon},\boldsymbol{\Lambda}} < K \mid \boldsymbol{\lambda} \right)
     \right] \\
  &= \D{}{\epsilon_j} 
       \D{}{\epsilon_j} E\left[ 1\left\{\mathbb{W}_{\boldsymbol{\epsilon},\boldsymbol{\Lambda}} < K\right\} \mid \boldsymbol{\lambda} \right]
      \\
  &= \D{}{\epsilon_j} 
       \D{}{\epsilon_j} E\left\{ E\left[ 1\left\{\mathbb{W}_{\boldsymbol{\epsilon},\boldsymbol{\Lambda},-j} < K - \sqrt{n}\epsilon_j\psi_j\lambda_j\left[Q'(u_j)+Q''(u_j)(y_j-u_j)\right] \right\} \mid \boldsymbol{\lambda}, y_j \right] \mid \boldsymbol{\lambda} \right\}
      \\
  &= E\left\{ \D{}{\epsilon_j} \D{}{\epsilon_j} E\left[ 1\left\{\mathbb{W}_{\boldsymbol{\epsilon},\boldsymbol{\Lambda},-j} < K - \sqrt{n}\epsilon_j\psi_j\lambda_j\left[Q'(u_j)+Q''(u_j)(y_j-u_j)\right] \right\} \mid \boldsymbol{\lambda}, y_j \right] \mid \boldsymbol{\lambda} \right\}
      \\
  &= E\biggl\{ 
      \D{}{\epsilon_j} 
    \D{}{K} E\left[ 1\left\{\mathbb{W}_{\boldsymbol{\epsilon},\boldsymbol{\Lambda},-j} < K - \sqrt{n}\epsilon_j\psi_j\lambda_j\left[Q'(u_j)+Q''(u_j)(y_j-u_j)\right] \right\} \mid \boldsymbol{\lambda}, y_j \right] 
     \\
  &\qquad\quad\times
     \Bigl(-\sqrt{n}\psi_j\lambda_j\bigl[Q'(u_j)+\overbrace{Q''(u_j)(y_j-u_j)}^{\star(2\log(n))\textrm{, \ref{a:hut-pf}}}\bigr]\Bigr) 
     \mid \boldsymbol{\lambda} 
     \biggr\}
      \\
  &= E\biggl\{ 
       n \psi_j^2 \lambda_j^2 \left[Q'(u_j)^2+O\left(n^{-1/2}\log(n)\right)\right] \\
  &\qquad\quad\times
       \frac{\partial^2}{\partial K^2} E\left[ 1\left\{\mathbb{W}_{\boldsymbol{\epsilon},\boldsymbol{\Lambda},-j} < K - \sqrt{n}\epsilon_j\psi_j\lambda_j\left[Q'(u_j)+Q''(u_j)(y_j-u_j)\right] \right\} \mid \boldsymbol{\lambda}, y_j \right] 
     \mid \boldsymbol{\lambda} \biggr\}
      \\
  &= n \psi_j^2 \lambda_j^2 \left[Q'(u_j)\right]^2
     \frac{\partial^2}{\partial K^2} 
     E\left\{ E\left[ 1\left\{\mathbb{W}_{\boldsymbol{\epsilon},\boldsymbol{\Lambda}} < K \right\} \mid \boldsymbol{\lambda}, y_j \right] \mid \boldsymbol{\lambda} \right\}
     +O\left(n^{-3/2}[\log(n)]^3\right) 
      \\
  &= n \psi_j^2 \lambda_j^2 \left[Q'(u_j)\right]^2
     \overbrace{%
     \frac{\partial^2}{\partial K^2} 
     F_{\mathbb{W}_{\boldsymbol{\epsilon},\boldsymbol{\Lambda}}|\boldsymbol{\Lambda}}( K \mid \boldsymbol{\lambda})%
     }^{\textrm{apply \eqref{eqn:W-EPS-PDF-cond-derivative}}}
     +O\left(n^{-3/2}[\log(n)]^3\right) 
      \\
  &= n \psi_j^2 \lambda_j^2 \left[Q'(u_j)\right]^2
     \left[ \phi'_{\mathcal{V}_\psi}(K) + O\left(n^{-1/2}[\log(n)]^{3+J}\right) \right]
     +O\left(n^{-3/2}[\log(n)]^3\right) 
      \\
  &= n \psi_j^2 \lambda_j^2 \left[Q'(u_j)\right]^2
     \phi'_{\mathcal{V}_\psi}(K)
     +O\left(n^{-3/2}[\log(n)]^{5+J}\right) . 
\end{align*}
In addition to the chain rule, the above steps use $Q'(u_j)$ and $Q''(u_j)$ being uniformly bounded from Assumption \ref{a:hut-pf}, and $\lambda_j=O(n^{-1}\log(n))$ and $Y_j-u_j=O(n^{-1/2}\log(n))$ from \ref{cond:star} with $a_n=2\log(n)$. %
}{%
%
}
\Supplemental{\end{proof}}{}

\subsection*{\Supplemental{Proof}{Sketch of proof} of Theorem \ref{thm:cdferror}\Supplemental{(\ref{thm:cdferror-ptwise})}{}}

\Supplemental{The following is brief intuition behind the proof. }{}
\Supplemental{We }{For part (\ref{thm:cdferror-ptwise}), we }%
start by restricting attention to cases where the largest of the $J$ spacings between relevant uniform order statistics, $U_{n:\lfloor(n+1)u_j\rfloor+1}-U_{n:\lfloor(n+1)u_j\rfloor}$, and the largest difference between the $U_{n:\lfloor(n+1)u_j\rfloor}$ and $u_j$ satisfy Condition $\star(2\log(n))$ as in Lemma \ref{lem:den-i}. 
By Lemma \ref{lem:den-i}(\ref{lem:den-i-star-prob},\ref{lem:den-i-star-prob-fixed}), the error from this restriction is smaller-order.  We then use the representation of ideal uniform fractional order statistics from \citet{Jones2002}, which is equal in distribution to the linearly interpolated form but with random interpolation weights $C_j\sim\beta(\epsilon_j,1-\epsilon_j)$ instead of fixed $\epsilon_j$, where each $C_j$ is independent of every other random variable we have.  The leading term in the error is due to $\Varp{C_j}$, and by plugging in other calculations from Lemma \ref{lem:den}, we see that it is uniformly $O(n^{-1})$ and can be calculated analytically.  

\Supplemental{\begin{proof}}{}
\Supplemental{
We assume that the realized values of $\mathbf{Y}$ and $\boldsymbol{\Lambda}$ all satisfy Condition $\star (2\log(n))$.  By application of Lemma \ref{lem:den-i}(\ref{lem:den-i-star-prob},\ref{lem:den-i-star-prob-fixed}), this induces at most $O(n^{-2})$ error.  For any event $A$, given $a_n=2\log(n)$, 
\begin{align*}
P(A)
&= P(A\mid\star(a_n))P(\star(a_n)) + P(A\mid\textrm{not }\star(a_n))[1-P(\star(a_n))] \\
&= P(A\mid\star(a_n))[1-O(n^{-2})] + P(A\mid\textrm{not }\star(a_n))O(n^{-2}) \\
&= P(A\mid\star(a_n)) + O(n^{-2}) . 
\label{eqn:star-error}\refstepcounter{equation}\tag{\theequation}
\end{align*}
This $O(n^{-2})$ remainder is asymptotically negligible (with respect to the nearly $O(n^{-3/2})$ remainder in the theorem), so we generally omit it going forward. 

Applying Lemma \ref{lem:den}(\ref{lem:den-bd2}) and Lemma \ref{lem:den}(\ref{lem:den-iii}) to the CDF of $L^L$,
\begin{align*}
P \left(L^L < \mathbb{X}_{0} + n^{-1/2} K \right)
 &= P\left( \sqrt{n}(L^L-\mathbb{X}_{0}) < K \right) \\
 &= P\left( \mathbb{W}_{\boldsymbol{\epsilon},\boldsymbol{\Lambda}} 
           +\overbrace{O\left(n^{-3/2}[\log(n)]^3\right)}^{\textrm{\ref{lem:den}(\ref{lem:den-bd2})}}
           < K \right) \\
\label{eqn:W-EL-CDF-K}\refstepcounter{equation}\tag{\theequation}
 &= P\left(\mathbb{W}_{\boldsymbol{\epsilon,\Lambda}}<K\right) 
   +\overbrace{\overbrace{f_{\mathbb{W}_{\boldsymbol{\epsilon},\boldsymbol{\Lambda}}}(\tilde{K})}^{=O(1)\textrm{ by \ref{lem:den}(\ref{lem:den-iii})}}
    O\left(n^{-3/2}[\log(n)]^3\right)}^{\textrm{mean value theorem}} \\ 
 &= \int\limits_{\underbrace{[0,2n^{-1}\log(n)]^J}_{\star(2\log(n))}}
   P\left(\mathbb W_{\boldsymbol{\epsilon},\boldsymbol{\Lambda}} < K  \mid \boldsymbol\lambda\right) 
   \,dF_{\boldsymbol{\Lambda}}(\boldsymbol \lambda)
\\&\quad +\overbrace{O(n^{-2})}^{\textrm{\eqref{eqn:star-error}}} + O(n^{-3/2}\log(n)^3).
\end{align*}
By a similar series of manipulations, and then using $\mathbf{C}\independent\boldsymbol{\Lambda}$ from \citet{Jones2002}, 
\begin{align*}
P \left( L^I < \mathbb{X}_{0} + n^{-1/2} K \right)
\Supplemental{%
 &= P\left( \sqrt{n}(L^I-\mathbb{X}_{0}) < K \right) \\
 &= P\left( \mathbb{W}_{\boldsymbol{C,\Lambda}} 
           +\overbrace{O\left(n^{-3/2}[\log(n)]^3\right)}^{\textrm{\ref{lem:den}(\ref{lem:den-bd2})}}
           < K \right) \\
 &= P(\mathbb{W}_{\boldsymbol{C,\Lambda}} < K) 
   +\overbrace{\overbrace{f_{\mathbb{W}_{\boldsymbol{C,\Lambda}}}(\tilde{K})}^{O(1)\textrm{ by \ref{lem:den}(\ref{lem:den-iii})}}
    O\left(n^{-3/2}[\log(n)]^3\right)}^{\textrm{mean value theorem}} \\
}{}
 &= \int\limits_{[0,2n^{-1}\log(n)]^J} 
    \; \int\limits_{[0,1]^J}
    P\left(\mathbb W_{\mathbf{c},\boldsymbol\Lambda} < K \mid \boldsymbol\lambda\right) 
    \,dF_{\mathbf{C}}(\mathbf{c})
    \,dF_{\boldsymbol\Lambda}(\boldsymbol\lambda) 
\\&\quad +\overbrace{O(n^{-2})}^{\textrm{\eqref{eqn:star-error}}} + O(n^{-3/2}\log(n)^3) . 
\end{align*}
The CDF difference between the two distributions is
\begin{align*}
\notag
P & \left(L^I<\mathbb{X}_{0} + n^{-1/2}K\right) 
    - P\left(L^L <\mathbb{X}_{0} + n^{-1/2}K\right) \\
\notag
&= \int\limits_{[0,2n^{-1}\log(n)]^J}
    \;\int\limits_{[0,1]^J}\left[F_{\mathbb{W}_{\mathbf{c},\boldsymbol{\Lambda}}|\boldsymbol{\Lambda}} ( K  \mid \boldsymbol\lambda) - F_{\mathbb{W}_{\boldsymbol\epsilon,\boldsymbol{\Lambda}}|\boldsymbol{\Lambda}} ( K  \mid \boldsymbol\lambda)\right] 
    \,dF_{\mathbf C}(\mathbf{c}) 
    \,dF_{\boldsymbol{\Lambda}}(\boldsymbol{\lambda})  
\\\notag&\qquad+ O(n^{-3/2}\log(n)^3)\\
\notag
&= \int\limits_{[0,2n^{-1}\log(n)]^J}
   \;\int\limits_{[0,1]^J} \biggl[
   \left. (\mathbf c - \boldsymbol\epsilon)' \frac{\partial F_{\mathbb{W}_{\mathbf c,\boldsymbol{\Lambda}}|\boldsymbol{\Lambda}} ( K  \mid \boldsymbol\lambda)}{\partial \mathbf{c}} \right|_{\mathbf c = \mathbf \epsilon}\\
\notag
  &\qquad\qquad\qquad\qquad\quad+\frac{1}{2} (\mathbf c - \boldsymbol\epsilon)'
  \left. \frac{\partial^2 F_{\mathbb{W}_{\mathbf c,\boldsymbol{\Lambda}}|\boldsymbol{\Lambda}} ( K  \mid \boldsymbol\lambda)}{\partial \mathbf{c} \partial \mathbf{c}'} \right|_{\mathbf c = \tilde{\mathbf c}\in [0,1]^J} (\mathbf c - \boldsymbol\epsilon)\biggr] 
  \,dF_{\mathbf C}(\mathbf{c}) 
  \,dF_{\boldsymbol{\Lambda}}(\boldsymbol{\lambda}) \\
\notag
  &\quad + O\left(n^{-3/2}[\log(n)]^3\right). \\
\intertext{Since $E(\mathbf C) = \boldsymbol\epsilon$ and $\mathbf{C}\independent\boldsymbol{\Lambda}$, the first term zeroes out.  
Since (additionally) the elements of $\mathbf C$ are mutually independent, the off-diagonal elements of the Hessian in the quadratic term also zero out%
\Supplemental{%
: for $j\ne k$,
\begin{equation*}
E[(C_j-\epsilon_j)(C_k-\epsilon_k)]
  = E(C_j-\epsilon_j)E(C_k-\epsilon_k)
  = 0 .
\end{equation*}
}{. }
Then we apply Lemma \ref{lem:den}(\ref{lem:den-iii}) to the Hessian 
and use $\Varp{C_j}=\epsilon_j(1-\epsilon_j)/2$, leaving}
  &= \int\limits_{[0,2n^{-1}\log(n)]^J}
    \frac{1}{2} 
    \sum_{j=1}^{J} 
     \overbrace{\frac{\epsilon_j(1-\epsilon_j)}{2}}^{\Varp{C_j}}
     \overbrace{\left[ n \psi_j^2 [Q'(u_j)]^2 \lambda_j^2 \phi'_{\mathcal{V}_\psi}(K) + O\left(n^{-3/2}[\log(n)]^{5+J}\right)\right]}^{\textrm{Lemma \ref{lem:den}(\ref{lem:den-iii})}}
     \,dF_{\boldsymbol{\Lambda}}(\boldsymbol{\lambda}) \\
  &\quad + O\left(n^{-3/2}[\log(n)]^3\right) \\
  &= \frac{n}{2} \phi_{\mathcal{V}_{\boldsymbol{\psi}}}'(K)
    \sum_{j=1}^J \left[[\psi_j Q'(u_j)]^2\frac{\epsilon_j(1-\epsilon_j)}{2} 
      \int\limits_{[0,2n^{-1}\log(n)]^J}
      \lambda_j^2 dF_{\boldsymbol{\Lambda}}(\boldsymbol{\lambda})\right] \\
  &\quad
   +O\left(n^{-5/2}[\log(n)]^{6+J}\right)
    + O\left(n^{-3/2}[\log(n)]^3\right) . 
\end{align*} 
The proof is completed by using $\Lambda_j\sim\beta(1,n)$, 
\Supplemental{%
whose PDF is $f_{\Lambda_j}(\lambda_j)=n(1-\lambda_j)^{n-1}$, so 
\begin{align*}
\int_{0}^{2n^{-1}\log(n)} \lambda_j^2 n(1-\lambda_j)^{n-1} \,d\lambda_j
  &= \left[ -\frac{(1-\lambda_j)^{n} \left(n^2\lambda_j^2 +n\lambda_j(\lambda_j+2) +2\right)}
                  {(n+1)(n+2)}
     \right]_{0}^{2n^{-1}\log(n)} \\
  &= -\frac{(1-2\log(n)/n)^{n}}{(n+1)(n+2)} O\left( [\log(n)]^2 +\log(n) +1\right) \\
  &\quad
    -\left( -\frac{2}{(n+1)(n+2)} \right) \\
  &= 2n^{-2}\left[1-O(1/n)\right]
    +O\left(n^{-4}[\log(n)]^2\right) \\
  &= 2n^{-2} + O(n^{-3}) .
  \qedhere
\end{align*}
}{%
to approximate the integral as $2n^{-2}+O(n^{-3})$. 
}
}{}
\Supplemental{\end{proof}}{}

\Supplemental{\subsection*{Proof of Theorem \ref{thm:cdferror}(\ref{thm:cdferror-unif})}}{}

\Supplemental{\begin{proof}}{}
\Supplemental{%
We again assume that the realized values of $\mathbf Y$ and $\boldsymbol\lambda$ satisfy Condition $\star (2\log(n))$.  By application of Lemma \ref{lem:den-i}(\ref{lem:den-i-star-prob},\ref{lem:den-i-star-prob-fixed}), this induces at most $O(n^{-2})$ error in our calculations, which is asymptotically negligible, as shown explicitly in the proof of Theorem \ref{thm:cdferror}(\ref{thm:cdferror-ptwise}), equation \eqref{eqn:star-error}. 
}{}

\Supplemental{%
The first result of this part is obtained by solving a first order condition and evaluating the expression obtained in part (\ref{thm:cdferror-ptwise}) at the solution $K=\sqrt{\mathcal{V_{\boldsymbol{\psi}}}}$.  To solve the FOC,
}{%
For part (\ref{thm:cdferror-unif}), the first result comes from the FOC
}
\begin{align*}
0 &= \D{}{K} \frac{K \exp\left\{ -K^2/(2\mathcal{V}_\psi)\right\}}
                  {\sqrt{2\pi\mathcal{V}_\psi^3}}
     \left[
       \sum_{j=1}^{J} \left( \frac{\psi_j^2\epsilon_j(1-\epsilon_j)}{\left[f\left(F^{-1}\left(u_j\right)\right)\right]^2} \right)
     \right]
     n^{-1} 
\Supplemental{%
\\
  &= \left[
      \frac{\exp\left\{ -K^2/(2\mathcal{V}_\psi)\right\}}
           {\sqrt{2\pi\mathcal{V}_\psi^3}}
      -(K^2/\mathcal{V}_\psi)
       \frac{\exp\left\{ -K^2/(2\mathcal{V}_\psi)\right\}}
            {\sqrt{2\pi\mathcal{V}_\psi^3}}
     \right]
     \left[
       \sum_{j=1}^{J} \left( \frac{\psi_j^2\epsilon_j(1-\epsilon_j)}{\left[f\left(F^{-1}\left(u_j\right)\right)\right]^2} \right)
     \right]
     n^{-1} , \\
1 &= K^2/\mathcal{V}_\psi, \\
K &= \sqrt{\mathcal{V}_\psi} . 
}{}
\end{align*}
\Supplemental{%
Plugging in,
\begin{align*}
\frac{\sqrt{\mathcal{V}_\psi} \exp\left\{ -\mathcal{V}_\psi/(2\mathcal{V}_\psi)\right\}}
                  {\sqrt{2\pi\mathcal{V}_\psi^3}}
     \left[
       \sum_{j=1}^{J} \left( \frac{\psi_j^2\epsilon_j(1-\epsilon_j)}{\left[f\left(F^{-1}\left(u_j\right)\right)\right]^2} \right)
     \right]
     n^{-1} 
&= \frac{e^{-1/2}}
        {\sqrt{2\pi\mathcal{V}_\psi^2}}
     \left[
       \sum_{j=1}^{J} \left( \frac{\psi_j^2\epsilon_j(1-\epsilon_j)}{\left[f\left(F^{-1}\left(u_j\right)\right)\right]^2} \right)
     \right]
     n^{-1} . 
\end{align*}
}{%
whose solution $K=\sqrt{\mathcal{V}_\psi}$ is plugged into the expression in Theorem \ref{thm:cdferror}(\ref{thm:cdferror-ptwise}). 
}

The additional result for $L^B$ in part (\ref{thm:cdferror-unif}) follows from the Dirichlet PDF approximation in Lemma \ref{lem:den}(\ref{lem:den-iii}).%
\Supplemental{%
\footnote{A similar approximation error seems to follow from using the remainder in \citet{Portnoy2012}, or from applying equation (5.1) in Theorem 6 of \citet{CsorgoRevesz1978}.} 
Specifically, continuing from \eqref{eqn:W-EL-CDF-K}, 
\begin{align*}
P\left( L^L < \mathbb{X}_{0} + n^{-1/2} K \right)
 &= P(\mathbb{W}_{\boldsymbol{\epsilon,\Lambda}} < K) 
   +O\left(n^{-3/2}[\log(n)]^3\right) \\
 &= \int_{-\infty}^{K} 
    \int\limits_{[0,2\log(n)/n]^J}
      \overbrace{f_{\mathbb{W}_{\boldsymbol{\epsilon,\Lambda}}|\boldsymbol{\Lambda}}\left( w \mid \boldsymbol{\lambda} \right)}^{\textrm{Lemma \ref{lem:den}(\ref{lem:den-iii})}} 
      \,dF_{\boldsymbol{\Lambda}}(\boldsymbol{\lambda})
      \,dw 
   +O\left(n^{-3/2}[\log(n)]^3\right) \\
 &= \int_{-\infty}^{K} \phi_{\mathcal{V}_\psi}(w) \left[ 1 + \overbrace{O\left(n^{-1/2}[\log(n)]^3\right)}^{\textrm{uniform}} \right] \,dw 
   +O\left(n^{-3/2}[\log(n)]^3\right) \\
 &= \left[ 1 + O\left(n^{-1/2}[\log(n)]^3\right) \right]
    \Phi_{\mathcal{V}_\psi}(K) 
   +O\left(n^{-3/2}[\log(n)]^3\right) . 
\end{align*}
Starting from the $L^B$ end,
\begin{align}\notag
\begin{split}
P\left( L^B < \mathbb{X}_{0} + n^{-1/2} K \right)
  &= P\left( n^{1/2}\left( L^B - \mathbb{X}_{0}\right) < K \right) \\
  &= P\left( 
       n^{1/2}\left( 
         \boldsymbol{\psi}' \tilde{Q}^B_X(\mathbf{u}) - \boldsymbol{\psi}' F^{-1}(\mathbf{u}) 
       \right) 
      < K 
     \right) \\
  &= P\left( 
       n^{1/2}\left( 
         \sum_{j=1}^{J} \psi_j \left[ F^{-1}\left(u_j+n^{-1/2}B(u_j)\right) - F^{-1}(u_j) \right] 
       \right) 
      < K 
     \right) \\
  &= P\left( 
       n^{1/2}\left( 
         \sum_{j=1}^{J} \psi_j \left[ n^{-1/2}B(u_j)Q'(u_j) 
                                     +n^{-1}B(u_j)^2 \frac{Q''(\tilde{u}_j)}{2} \right]
       \right) 
      < K 
     \right) 
\end{split}
\\&= P\left(
      \boldsymbol{\psi}' \underline{\mathcal{A}}^{-1} B(\mathbf{u})
      +(1/2)n^{-1/2} \sum_{j=1}^J \psi_j Q''(\tilde{u}_j) B(u_j)^2
      < K 
     \right) .
\label{eqn:LB-int1}
\end{align}
where $\tilde{u}_j$ is between $u_j$ and $u_j+n^{-1/2}B(u_j)$, using notation from the main text like $B(\cdot)$ for a standard Brownian bridge, and using notation $B(\mathbf{u})\equiv(B(u_1),\ldots,B(u_J))'$. 

If we consider a variant of \ref{cond:star} where all $|B(u_j)|<\sqrt{\log(n)}$, then since $B(u_j)\sim N(0,\sigma^2)$ with $\sigma^2=u_j(1-u_j)$, 
\begin{align*}
P\left( B(u_j) > \log(n)\right)
 &= \frac{1}{\sigma \sqrt{2\pi}} 
    \int_{\sqrt{\log(n)}}^{\infty} 
    e^{-x^2/(2\sigma^2)} \,dx \\
 &< \frac{\sigma}{\sqrt{2\pi}} 
    \int_{\sqrt{\log(n)}}^{\infty} 
    \frac{x}{\sigma^2\sqrt{\log(n)}} e^{-x^2/(2\sigma^2)} \,dx \\
 &= \frac{\sigma}{\sqrt{2\pi\log(n)}} 
    \left[ -e^{-x^2/(2\sigma^2)} \right]_{\sqrt{\log(n)}}^{\infty} \\
 &= \frac{\sigma}{\sqrt{2\pi\log(n)}} 
    \exp\{-[\sqrt{\log(n)}]^2/(2\sigma^2)\} \\
 &= \frac{\sigma}{\sqrt{2\pi\log(n)}} 
    \exp\{\log(1/n)/\overbrace{(2\sigma^2)}^{\le1/2}\} \\
\label{eqn:cond-star-B-prob}\refstepcounter{equation}\tag{\theequation}
 &= o\left(n^{-1}\right) . 
\end{align*}
Under this condition, $\tilde{u}_j\to u_j$, so \ref{a:hut-pf} bounds $Q''(\tilde{u}_j)$ in \eqref{eqn:LB-int1}.  

In \eqref{eqn:LB-int1}, note that the PDF of $\boldsymbol{\psi}' \underline{\mathcal{A}}^{-1} B(\mathbf{u})$ is $\phi_{\mathcal{V}_\psi}$. 
Also, given $|B(u_j)|<\sqrt{\log(n)}$, then $n^{-1/2}\psi_j Q''(\tilde{u}_j) B(u_j)^2=O(n^{-1/2}\log(n))$, and summing over $J$ such terms maintains the same order. 
Plugging this back into \eqref{eqn:LB-int1}, 
\begin{align}\notag
\Phi_{\mathcal{V}_\psi}\left( K - O\left(n^{-1/2}\log(n)\right) \right)
 &= \Phi_{\mathcal{V}_\psi}\left( K \right)
   -O\left(n^{-1/2}\log(n)\right)
    \phi_{\mathcal{V}_\psi}(\tilde{K}) \\
 &= \Phi_{\mathcal{V}_\psi}\left( K \right)
   -O\left(n^{-1/2}\log(n)\right) . 
\label{eqn:Phi-V-psi-K-approx}
\end{align}
Similar to the argument in \eqref{eqn:star-error}, this remainder is larger than the $o(n^{-1})$ in \eqref{eqn:cond-star-B-prob} from imposing the condition, so it dominates.  
Altogether,
\begin{align*}
& P\left( L^L < \mathbb{X}_{0} + n^{-1/2} K \right) \\
 &= \left[ 1 + O\left(n^{-1/2}[\log(n)]^3\right) \right]
    \Phi_{\mathcal{V}_\psi}(K) 
   +O\left(n^{-3/2}[\log(n)]^3\right) \\
 &= \left[ 1 + O\left(n^{-1/2}[\log(n)]^3\right) \right]
    \left[ P\left( L^B < \mathbb{X}_{0} + n^{-1/2} K \right) +\overbrace{o(n^{-1})}^{\textrm{\eqref{eqn:cond-star-B-prob}}} + \overbrace{O(n^{-1/2}\log(n))}^{\textrm{\eqref{eqn:Phi-V-psi-K-approx}}} \right] 
\\&\qquad+O\left(n^{-3/2}[\log(n)]^3\right) \\
 &= P\left( L^B < \mathbb{X}_{0} + n^{-1/2} K \right) 
   +O\left(n^{-1/2}[\log(n)]^3\right) . 
\qedhere
\end{align*}
}{%
}
\Supplemental{\end{proof}}{}

\subsection*{\Supplemental{Proof}{Sketch of proof} of Lemma \ref{lem:u1u2-approx}\label{sec:app-u1u2}}


\Supplemental{\begin{proof}}{}
\Supplemental{%
Let $u\equiv u^l(\alpha)$.  From \eqref{eqn:uhdef}, $u$ solves $P(B>p)=\alpha$ for $B\sim\beta\left(u(n+1),(1-u)(n+1)\right)$; equivalently, $F_B(p)=1-\alpha$.  We use a Cornish--Fisher-type expansion from (5.1) in \citet{Pratt1968}.  
\Supplemental{%
In his equation, $\gamma=[p-(1-p)]/\sigma=(2p-1)/\sqrt{np(1-p)}=O(n^{-1/2})$, so (in his notation) 
\begin{align*}
z &= z_0 -\frac{\gamma}{6}(z_0^2-1) +O\left(n^{-1}\right), 
\end{align*}
where $z_0$ and $\gamma$ depend on $n$, $p$, and $u$, and $\Phi(z)=F_B(p)=1-\alpha$ so that $z=z_{1-\alpha}$.  
}{}Let 
\[p-u = n^{-1/2}c_1 +n^{-1}c_2 +O(n^{-3/2}).\]
Using (2.2) from \citet{Pratt1968} and Table A from \citet{PeizerPratt1968}, 
\Supplemental{}{one can show}
\begin{align*}
\Supplemental{%
z_0 
    &= \frac{(1-u)(n+1)-(1/2)-n(1-p)}{\sqrt{np(1-p)}}
     = \frac{(1-u)(n+1)-(1/2)-n(1-p)}{\sqrt{np(1-p)}} \\
    &= \frac{\sqrt n(p-u)}{\sqrt{p(1-p)}} -(1/2)n^{-1/2}\frac{2u-1}{\sqrt{p(1-p)}} \\
    &= \frac{\sqrt n(p-u)}{\sqrt{p(1-p)}} -(1/2)n^{-1/2}\frac{2p-1}{\sqrt{p(1-p)}} +O(n^{-1}), \\
z_0^2 &= \frac{n(p-u)^2}{p(1-p)} + O(n^{-1/2}), \\ 
}{}
z_{1-\alpha}
\Supplemental{%
  &= z_0 -\frac{\gamma}{6}(z_0^2-1) +O\left(n^{-1}\right) \\
  &= \frac{\sqrt n(p-u)}{\sqrt{p(1-p)}} -(1/2)n^{-1/2}\frac{2p-1}{\sqrt{p(1-p)}} +O(n^{-1}) \\
  &\quad -\frac{2p-1}{6\sqrt{np(1-p)}}\left( \frac{n(p-u)^2}{p(1-p)} + O(n^{-1/2}) - 1 \right)
    +O(n^{-1}) \\
  &= \frac{\sqrt n(p-u)}{\sqrt{p(1-p)}} -(1/2)n^{-1/2}\frac{2p-1}{\sqrt{p(1-p)}} 
    -n^{-1/2} \frac{2p-1}{6\sqrt{p(1-p)}}\left( \frac{n(p-u)^2}{p(1-p)} - 1 \right)
    +O(n^{-1}) \\
  &= \frac{\sqrt n(p-u)}{\sqrt{p(1-p)}} 
    -n^{-1/2}\left[ \frac{3(2p-1)}{6\sqrt{p(1-p)}} 
                   +\frac{2p-1}{6\sqrt{p(1-p)}}\left( \frac{n(p-u)^2}{p(1-p)} - 1 \right) \right]
    +O(n^{-1}) \\
  &= \frac{\sqrt n(p-u)}{\sqrt{p(1-p)}} 
    -n^{-1/2} \frac{2p-1}{6\sqrt{p(1-p)}}\left( \frac{n(p-u)^2}{p(1-p)} + 2 \right)
    +O(n^{-1}) \\
  &= \frac{\sqrt n\left[ c_1/\sqrt n +c_2/n +O(n^{-3/2})\right]}{\sqrt{p(1-p)}}  \\
  &\quad
    -n^{-1/2} \frac{2p-1}{6\sqrt{p(1-p)}}\left( \frac{n\left[ c_1^2/n +O(n^{-3/2}) \right]}{p(1-p)} + 2 \right) 
    +O(n^{-1}) \\
  &= \frac{c_1 +c_2n^{-1/2} +O(n^{-1})}{\sqrt{p(1-p)}} 
    -n^{-1/2} \frac{2p-1}{6\sqrt{p(1-p)}}\left( \frac{c_1^2 +O(n^{-1/2}) }{p(1-p)} + 2 \right)
    +O(n^{-1}) \\
}{}
  &= \frac{c_1}{\sqrt{p(1-p)}} + n^{-1/2}\frac{c_2}{\sqrt{p(1-p)}}
    -n^{-1/2} \frac{2p-1}{6\sqrt{p(1-p)}}\left( \frac{c_1^2}{p(1-p)} + 2 \right)
    +O(n^{-1}) .
\end{align*}
\Supplemental{

To solve for $c_1$, we can set equal the first-order terms on the left- and right-hand sides:
\begin{align*}
z_{1-\alpha}
  &= \frac{c_1}{\sqrt{p(1-p)}} \implies
c_1  = z_{1-\alpha}\sqrt{p(1-p)} .
\end{align*}
Plugging in for $c_1$, we can solve for $c_2$ to zero the $n^{-1/2}$ terms:
\begin{align*}
0 &= n^{-1/2} \left[ \frac{c_2}{\sqrt{p(1-p)}}
                    -\frac{2p-1}{6\sqrt{p(1-p)}}\left( \frac{z_{1-\alpha}^2p(1-p)}{p(1-p)} + 2 \right)\right] 
\Supplemental{, \\}{\implies}
\Supplemental{%
\frac{c_2}{\sqrt{p(1-p)}}
  &= \frac{2p-1}{6\sqrt{p(1-p)}}\left( z_{1-\alpha}^2 + 2 \right), \\
}{}
c_2 \Supplemental{&}{}= \frac{2p-1}{6}(z_{1-\alpha}^2+2) .
\end{align*}
Altogether,
}{After solving for $c_1$ and $c_2$,}
\begin{align*}
u &= p - \frac{c_1}{\sqrt n} - \frac{c_2}{n} +O(n^{-3/2}) 
        = p - n^{-1/2}z_{1-\alpha}\sqrt{p(1-p)} - \frac{2p-1}{6n}(z_{1-\alpha}^2+2) +O(n^{-3/2}) .
\end{align*}

For $u=u^h(\alpha)$, we have $z_\alpha$ instead of $z_{1-\alpha}$, and using $-z_\alpha=z_{1-\alpha}$ yields the result. 
}{%
The results are based on the Cornish--Fisher-type expansion from \citet{Pratt1968} and \citet{PeizerPratt1968}, solving for the high-order constants. %
}%
\Supplemental{\end{proof}}{}

\subsection*{\Supplemental{Proof}{Sketch of proof} of Theorem \ref{thm:IDEAL-single}}

\Supplemental{\subsubsection*{Coverage probability}}{}

\Supplemental{\begin{proof}}{}
\Supplemental{Let }{For CP, let }%
$R_n=O\left(n^{-3/2}[\log(n)]^3\right)$ be the remainder from Theorem \ref{thm:cdferror}(\ref{thm:cdferror-ptwise}). 

For a lower one-sided CI,
\begin{align*}
u^h(\alpha) &= p+O(n^{-1/2}), \quad 
J =1 , \quad 
\epsilon_h
   = (n+1)u^h(\alpha) - \lfloor(n+1)u^h(\alpha)\rfloor, \\
\mathcal{V}_{\boldsymbol{\psi}} 
  &= \frac{u^h(\alpha)(1-u^h(\alpha))}{f(F^{-1}(u^h(\alpha)))^2},   \quad 
\mathbb{X}_{0}
   = F^{-1}(u^h(\alpha)),   \\
K &= n^{1/2}\left[F^{-1}(p)-F^{-1}\left(u^h(\alpha)\right)\right]  
   = -\frac{z_{1-\alpha} \sqrt{u^h(\alpha)(1-u^h(\alpha))}}{f(F^{-1}(u^h(\alpha)))} + O(n^{-1/2}) \\
  &= -z_{1-\alpha}\sqrt{\mathcal{V}_\psi} + O(n^{-1/2})  ,
\end{align*}
where the first and last lines use Lemma \ref{lem:u1u2-approx}, and the last line uses Assumption \ref{a:hut-pf}.  
Then, the rate of coverage probability error is 
\begin{align*}
& P\left( \hat{Q}^L_X\left(u^h(\alpha)\right) < Q(p) \right) 
\Supplemental{\\&}{}= P\left( \hat{Q}^L_X\left(u^h(\alpha)\right) < \mathbb{X}_{0} + n^{-1/2}K \right) \\
&= P\left( \hat{Q}^I_X\left(u^h(\alpha)\right) < \mathbb{X}_{0} + n^{-1/2}K \right) 
    +n^{-1}\frac{\epsilon_h(1-\epsilon_h)}{[f(F^{-1}(u^h(\alpha)))]^2}
     \frac{K \exp\{-K^2/(2\mathcal{V}_\psi)\}}{\sqrt{2\pi\mathcal{V}_\psi^3}}
    +R_n \\
\Supplemental{%
&= \overbrace{P\left( \hat{Q}^I_X\left(u^h(\alpha)\right) < Q(p) \right) }^{=\alpha\textrm{, by definition of }u^h(\alpha)} \\
&\quad
    +n^{-1}\frac{\epsilon_h(1-\epsilon_h)}{[f(F^{-1}(u^h(\alpha)))]^2}
     \left[ -\frac{z_{1-\alpha} \sqrt{u^h(\alpha)(1-u^h(\alpha))}}{f(F^{-1}(u^h(\alpha)))} + O(n^{-1/2}) \right] \\
&\qquad\times
     \frac{\exp\{-z_{1-\alpha}^2/2+O(n^{-1/2})\}}{\sqrt{2\pi}}
     \frac{\left[f(F^{-1}(u^h(\alpha)))\right]^3}{\left[u^h(\alpha)(1-u^h(\alpha))\right]^{3/2}}
    +R_n \\
&= \alpha
  -n^{-1}
   z_{1-\alpha}
   \frac{\epsilon_h(1-\epsilon_h)}{\sqrt{2\pi}u^h(\alpha)\left(1-u^h(\alpha)\right)}
   \exp\{-z_{1-\alpha}^2/2\}\left[1+O(n^{-1/2})\right] 
  +O(n^{-3/2})
  +R_n \\
}{}
\label{eqn:upper-non-CP}\refstepcounter{equation}\tag{\theequation}
&= \alpha
  -n^{-1}
   z_{1-\alpha}
   \frac{\epsilon_h(1-\epsilon_h)}{p(1-p)}
   \phi(z_{1-\alpha})
  +O(n^{-3/2})
  +R_n ,
\end{align*}
where $f(F^{-1}(u^h(\alpha)))$ is uniformly (for large enough $n$) bounded away from zero by \ref{a:hut-pf} since $u^h(\alpha)=p+O(n^{-1/2})\to p$. 
\Supplemental{

For the lower endpoint,
\begin{align*}
u^l(\alpha) &= p-O(n^{-1/2}), \quad 
J =1 , \quad 
\epsilon_\ell
   = (n+1)u^l(\alpha) - \lfloor(n+1)u^l(\alpha)\rfloor, \\
\mathcal{V}_{\boldsymbol{\psi}} 
  &= \frac{u^l(\alpha)(1-u^l(\alpha))}{f(F^{-1}(u^l(\alpha)))^2},   \quad 
\mathbb{X}_{0}
   = F^{-1}(u^l(\alpha)),   \\
K &= n^{1/2}\left[F^{-1}(p)-F^{-1}\left(u^l(\alpha)\right)\right]  
   = \frac{z_{1-\alpha} \sqrt{u^l(\alpha)(1-u^l(\alpha))}}{f(F^{-1}(u^l(\alpha)))} + O(n^{-1/2}) \\
  &= z_{1-\alpha}\sqrt{\mathcal{V}_\psi} + O(n^{-1/2})  ,
\end{align*}
where the first and last lines use Lemma \ref{lem:u1u2-approx}, and the last line uses Assumption \ref{a:hut-pf}. 
Then, 
coverage probability is 
\begin{align*}
& P\left( \hat{Q}^L_X\left(u^l(\alpha)\right) < Q(p) \right) 
\Supplemental{\\&}= P\left( \hat{Q}^L_X\left(u^l(\alpha)\right) < \mathbb{X}_{0} + n^{-1/2}K \right) \\
&= P\left( \hat{Q}^I_X\left(u^l(\alpha)\right) < \mathbb{X}_{0} + n^{-1/2}K \right) 
    +n^{-1}\frac{\epsilon_\ell(1-\epsilon_\ell)}{[f(F^{-1}(u^l(\alpha)))]^2}
     \frac{K \exp\{-K^2/(2\mathcal{V}_\psi)\}}{\sqrt{2\pi\mathcal{V}_\psi^3}}
    +R_n \\
\Supplemental{%
&= \overbrace{P\left( \hat{Q}^I_X\left(u^l(\alpha)\right) < Q(p) \right) }^{=1-\alpha\textrm{, by definition of }u^l(\alpha)} \\
&\quad
    +n^{-1}\frac{\epsilon_\ell(1-\epsilon_\ell)}{[f(F^{-1}(u^l(\alpha)))]^2}
     \left[ \frac{z_{1-\alpha} \sqrt{u^l(\alpha)(1-u^l(\alpha))}}{f(F^{-1}(u^l(\alpha)))} + O(n^{-1/2}) \right] \\
&\qquad\times
     \frac{\exp\{-z_{1-\alpha}^2/2+O(n^{-1/2})\}}{\sqrt{2\pi}}
     \frac{\left[f(F^{-1}(u^l(\alpha)))\right]^3}{\left[u^l(\alpha)(1-u^l(\alpha))\right]^{3/2}}
    +R_n \\
&= 1-\alpha
  +n^{-1}
   z_{1-\alpha}
   \frac{\epsilon_\ell(1-\epsilon_\ell)}{\sqrt{2\pi}u^l(\alpha)\left(1-u^l(\alpha)\right)}
   \exp\{-z_{1-\alpha}^2/2\}\left[1+O(n^{-1/2})\right] 
  +O(n^{-3/2})
  +R_n \\
}{}
\label{eqn:lower-CP}\refstepcounter{equation}\tag{\theequation}
&= 1-\alpha
  +n^{-1}
   z_{1-\alpha}
   \frac{\epsilon_\ell(1-\epsilon_\ell)}{p(1-p)}
   \phi(z_{1-\alpha})
  +O(n^{-3/2})
  +R_n , 
\end{align*}
where $R_n$ is as before and $f(F^{-1}(u^l(\alpha)))$ is uniformly bounded away from zero by \ref{a:hut-pf} since $u^l(\alpha)=p-O(n^{-1/2})\to p$. 
}{%
The argument for the lower endpoint is similar.%
}

\Supplemental{%
For two-sided coverage probability, 
\Supplemental{%
the events $\hat{Q}^L_X\left(u^h(\alpha)\right) < Q(p)$ and $\hat{Q}^L_X\left(u^l(\alpha)\right) > Q(p)$ are mutually exclusive since $u^h(\alpha)>u^l(\alpha)$ by construction, so the probability of their union is the sum of their probabilities.  Using 
}{%
using }%
\eqref{eqn:upper-non-CP} and \Supplemental{\eqref{eqn:lower-CP}}{its analog}, replacing $\alpha$ with $\alpha/2$,
\begin{align*}
& P\left( \hat{Q}^L_X\left(u^l(\alpha/2)\right) < Q(p) < \hat{Q}^L_X\left(u^h(\alpha/2)\right) \right) \\
&= \overbrace{P\left( \hat{Q}^L_X\left(u^l(\alpha/2)\right) < Q(p) \right)}^{\textrm{\Supplemental{\eqref{eqn:lower-CP}}{analog of \eqref{eqn:upper-non-CP}}}}
     - \overbrace{P\left( \hat{Q}^L_X\left(u^h(\alpha/2)\right) < Q(p) \right)}^{\textrm{\eqref{eqn:upper-non-CP}}} \\
\Supplemental{%
&= \overbrace{1 - \alpha/2
    +n^{-1}
     z_{1-\alpha/2}
     \frac{\epsilon_\ell(1-\epsilon_\ell)}{p(1-p)}
     \phi(z_{1-\alpha/2})
    +O(n^{-3/2}) +R_n}^{\textrm{\eqref{eqn:lower-CP}}}
\\&\quad
  -\overbrace{\left\{ \alpha/2
     -n^{-1}
      z_{1-\alpha/2}
      \frac{\epsilon_h(1-\epsilon_h)}{p(1-p)}
      \phi(z_{1-\alpha/2})
     +O(n^{-3/2})
     +R_n 
   \right\}}^{\textrm{\eqref{eqn:upper-non-CP}}} \\
}{}
&= 1 - \alpha
  +n^{-1}
   \left[
      z_{1-\alpha/2}
      \phi(z_{1-\alpha/2})
      \frac{\epsilon_h(1-\epsilon_h)+\epsilon_\ell(1-\epsilon_\ell)}{p(1-p)}
   \right]
  +O(n^{-3/2}) + R_n . 
\end{align*}
}{%
Two-sided CP comes directly from the two one-sided results, replacing $\alpha$ with $\alpha/2$. %
}%
\Supplemental{

For the yet-higher-order calibration, let 
\begin{gather*}
\tilde\alpha = \alpha + n^{-1}\frac{\epsilon_h(1-\epsilon_h)z_{1-\alpha}\phi(z_{1-\alpha})}{p(1-p)} , \\
z_{1-\tilde\alpha} = \Phi^{-1}(1-\tilde\alpha) = \Phi^{-1}(1-\alpha-O(1/n)) = z_{1-\alpha}-O(n^{-1}) , \\
\begin{split}
\tilde{u}^h &= u^h(\tilde\alpha) 
   = \overbrace{p + n^{-1/2}z_{1-\tilde\alpha}\sqrt{p(1-p)} - \frac{2p-1}{6n}(z_{1-\tilde\alpha}^2+2) + O(n^{-3/2})}^{\textrm{Lemma \ref{lem:u1u2-approx}}} \\
\Supplemental{%
  &= p + n^{-1/2}[z_{1-\alpha}-O(n^{-1})]\sqrt{p(1-p)} - \frac{2p-1}{6n}\left\{\left[z_{1-\alpha}-O(n^{-1})\right]^2+2\right\} + O(n^{-3/2}) \\
  &= p + n^{-1/2}z_{1-\alpha}\sqrt{p(1-p)} - \frac{2p-1}{6n}\left\{z_{1-\alpha}^2+2\right\} + O(n^{-3/2}) \\
}{}
  &= u^h(\alpha) - O(n^{-3/2})
   = p + O(n^{-1/2}) , \\
\end{split}\\
\tilde{\epsilon}_h = (n+1)\tilde{u}^h - \lfloor(n+1)\tilde{u}^h\rfloor 
  = (n+1)u^h(\alpha)-O(n^{-1/2}) - \lfloor(n+1)u^h(\alpha)-O(n^{-1/2})\rfloor . 
\end{gather*}
If $\lfloor(n+1)\tilde{u}^h\rfloor=\lfloor(n+1)u^h(\alpha)\rfloor$, then $\tilde{\epsilon}_h=\epsilon_h-O(n^{-1/2})$, 
so \eqref{eqn:upper-non-CP} becomes
\begin{align*}
& P\left( \hat{Q}^L_X\left(u^h(\tilde\alpha)\right) < Q(p) \right) \\
&= \tilde\alpha
  -n^{-1}
   z_{1-\tilde\alpha}
   \frac{\tilde\epsilon_h(1-\tilde\epsilon_h)}{p\left(1-p\right)}
   \phi(z_{1-\tilde\alpha})
  +O(n^{-3/2})
  +R_n \\
\label{eqn:upper-non-CP-calib-int}\refstepcounter{equation}\tag{\theequation}
\begin{split}
&= \alpha 
  +n^{-1}\frac{\epsilon_h(1-\epsilon_h)z_{1-\alpha}\phi(z_{1-\alpha})}{p(1-p)} \\
&\quad
  -n^{-1}
   \left[ z_{1-\alpha} + O(n^{-1}) \right]
   \frac{\left[\epsilon_h-O(n^{-1/2})\right]\left\{1-\left[\epsilon_h-O(n^{-1/2})\right]\right\}}{p\left(1-p\right)}
   \left[ \phi(z_{1-\alpha}) + O(n^{-1}) \right] \\
&\quad
  +O(n^{-3/2})
  +R_n 
\end{split}\\
\label{eqn:upper-non-CP-calib}\refstepcounter{equation}\tag{\theequation}
&= \alpha +O(n^{-3/2}) +R_n . 
\end{align*}
Alternatively, if $\lfloor(n+1)\tilde{u}^h\rfloor=\lfloor(n+1)u^h\rfloor-1$, then 
\begin{align*}
\epsilon_h 
  &= (n+1)u^h-\lfloor(n+1)u^h\rfloor < (n+1)(u^h-\tilde{u}^h) = O(n^{-1/2}) , \\
1-\tilde\epsilon_h 
  &= 1 - \left[ (n+1)\tilde{u}^h - \overbrace{\lfloor(n+1)\tilde{u}^h\rfloor}^{=\lfloor(n+1)u^h\rfloor-1} \right] 
\\&= \lfloor(n+1)u^h\rfloor-(n+1)\tilde{u}^h < (n+1)(u^h-\tilde{u}^h) = O(n^{-1/2}) , 
\end{align*}
so $\epsilon_h(1-\epsilon_h)=O(n^{-1/2})$ and $\tilde\epsilon_h(1-\tilde\epsilon_h)=O(n^{-1/2})$.  Consequently, the two terms in \eqref{eqn:upper-non-CP-calib-int} with leading $n^{-1}$ become $O(n^{-3/2})$, so the final result in \eqref{eqn:upper-non-CP-calib} remains the same. 
\Supplemental{%

For the lower endpoint,
\begin{gather*}
\tilde\alpha = \alpha + n^{-1}\frac{\epsilon_\ell(1-\epsilon_\ell)z_{1-\alpha}\phi(z_{1-\alpha})}{p(1-p)} , \\
z_{1-\tilde\alpha} = \Phi^{-1}(1-\tilde\alpha) = \Phi^{-1}(1-\alpha-O(1/n)) = z_{1-\alpha}-O(n^{-1}) , \\
\begin{split}
\tilde{u}^l &= u^l(\tilde\alpha) 
   = \overbrace{p - n^{-1/2}z_{1-\tilde\alpha}\sqrt{p(1-p)} - \frac{2p-1}{6n}(z_{1-\tilde\alpha}^2+2) + O(n^{-3/2})}^{\textrm{Lemma \ref{lem:u1u2-approx}}} \\
  &= p - n^{-1/2}[z_{1-\alpha}-O(n^{-1})]\sqrt{p(1-p)} - \frac{2p-1}{6n}\left\{\left[z_{1-\alpha}-O(n^{-1})\right]^2+2\right\} + O(n^{-3/2}) \\
  &= p - n^{-1/2}z_{1-\alpha}\sqrt{p(1-p)} - \frac{2p-1}{6n}\left\{z_{1-\alpha}^2+2\right\} + O(n^{-3/2}) \\
  &= u^l(\alpha) + O(n^{-3/2})
   = p - O(n^{-1/2}) , \\
\tilde{\epsilon}_\ell &= (n+1)\tilde{u}^l - \lfloor(n+1)\tilde{u}^l\rfloor 
  = (n+1)u^l(\alpha)-O(n^{-1/2}) - \lfloor(n+1)u^l(\alpha)-O(n^{-1/2})\rfloor . 
\end{split}
\end{gather*}
Similar to the case for the upper endpoint, if $\lfloor(n+1)\tilde{u}^l\rfloor=\lfloor(n+1)u^l(\alpha)\rfloor$, then $\tilde{\epsilon}_\ell=\epsilon_\ell+O(n^{-1/2})$, 
so \eqref{eqn:lower-CP} becomes
\begin{align*}
& P\left( \hat{Q}^L_X\left(u^l(\tilde\alpha)\right) < Q(p) \right) \\
&= 1-\tilde\alpha
  +n^{-1}
   z_{1-\tilde\alpha}
   \frac{\tilde\epsilon_h(1-\tilde\epsilon_h)}{p\left(1-p\right)}
   \phi(z_{1-\tilde\alpha})
  +O(n^{-3/2})
  +R_n \\
\label{eqn:lower-CP-calib-int}\refstepcounter{equation}\tag{\theequation}
\begin{split}
&= 1-\alpha 
  -n^{-1}\frac{\epsilon_\ell(1-\epsilon_\ell)z_{1-\alpha}\phi(z_{1-\alpha})}{p(1-p)} \\
&\quad
  +n^{-1}
   \left[ z_{1-\alpha} + O(n^{-1}) \right]
   \frac{\left[\epsilon_\ell-O(n^{-1/2})\right]\left\{1-\left[\epsilon_\ell-O(n^{-1/2})\right]\right\}}{p\left(1-p\right)}
   \left[ \phi(z_{1-\alpha}) + O(n^{-1}) \right] \\
&\quad
  +O(n^{-3/2})
  +R_n 
\end{split}\\
\label{eqn:lower-CP-calib}\refstepcounter{equation}\tag{\theequation}
&= 1-\alpha +O(n^{-3/2}) +R_n . 
\end{align*}
Alternatively, if $\lfloor(n+1)\tilde{u}^l\rfloor=\lfloor(n+1)u^l\rfloor+1$, then 
\begin{align*}
\tilde\epsilon_\ell 
  &= (n+1)\tilde{u}^l-\lfloor(n+1)\tilde{u}^l\rfloor < (n+1)(\tilde{u}^l-u^l) = O(n^{-1/2}) , \\
1-\epsilon_\ell
  &= 1 - \left[ (n+1)u^l - \overbrace{\lfloor(n+1)u^l\rfloor}^{=\lfloor(n+1)\tilde{u}^l\rfloor-1} \right] 
\\&= \lfloor(n+1)\tilde{u}^l\rfloor-(n+1)u^l < (n+1)(\tilde{u}^l-u^l) = O(n^{-1/2}) , 
\end{align*}
so $\epsilon_\ell(1-\epsilon_\ell)=O(n^{-1/2})$ and $\tilde\epsilon_\ell(1-\tilde\epsilon_\ell)=O(n^{-1/2})$.  Consequently, the two terms in \eqref{eqn:lower-CP-calib-int} with leading $n^{-1}$ become $O(n^{-3/2})$, so the final result in \eqref{eqn:lower-CP-calib} remains the same. 
}{%
The argument for the lower endpoint is symmetric. 
}

For two-sided calibration, let
\begin{align*}
\tilde\alpha_\ell/2 &= \frac{\alpha}{2} + n^{-1}z_{1-\alpha/2}\phi(z_{1-\alpha/2})\frac{\epsilon_\ell(1-\epsilon_\ell)}{p(1-p)}, \\
\tilde\alpha_h/2 &= \frac{\alpha}{2} + n^{-1}z_{1-\alpha/2}\phi(z_{1-\alpha/2})\frac{\epsilon_h(1-\epsilon_h)}{p(1-p)} . 
\end{align*}
Using \eqref{eqn:upper-non-CP-calib} and \Supplemental{\eqref{eqn:lower-CP-calib}}{an analogous lower endpoint expression}, the coverage probability is
\begin{align*}
& P\left( \hat{Q}^L_X\left(u^l(\tilde\alpha_\ell/2)\right) < Q(p) < \hat{Q}^L_X\left(u^h(\tilde\alpha_h/2)\right) \right) \\
&= \overbrace{P\left( \hat{Q}^L_X\left(u^l(\tilde\alpha_\ell/2)\right) < Q(p) \right)}^{\textrm{\Supplemental{\eqref{eqn:lower-CP-calib}}{analog of \eqref{eqn:upper-non-CP-calib}}}}
     - \overbrace{P\left( \hat{Q}^L_X\left(u^h(\tilde\alpha_h/2)\right) < Q(p) \right)}^{\textrm{\eqref{eqn:upper-non-CP-calib}}} \\
\Supplemental{
&= \overbrace{1 - \alpha/2 + O(n^{-3/2}) + R_n}^{\textrm{\Supplemental{\eqref{eqn:lower-CP-calib}}{analog of \eqref{eqn:upper-non-CP-calib}}}}
  -\overbrace{\left\{ \alpha/2
    +O(n^{-3/2})
    +R_n
   \right\}}^{\textrm{\eqref{eqn:upper-non-CP-calib}}} \\
}{}
&= 1 - \alpha
  +O(n^{-3/2}) + R_n .   \qedhere
\end{align*}
}{%
For the yet-higher-order calibration, the results follow from plugging in the proposed $\tilde\alpha$. 
}
\Supplemental{\end{proof}}{}

\Supplemental{\subsubsection*{Power}}{}


\Supplemental{\begin{proof}}{}
\Supplemental{%
We use the first-order equivalence of $\hat Q^L_X(u)$ and $\tilde Q^B_X(u)$ to show that power against deviations of magnitude $n^{-1/2}$ is the same as the power of the test based on asymptotic normality of the sample quantile. 

Consider one-sided power of the test of $H_0:Q(p)\ge D_n \equiv Q(p)+Kn^{-1/2}$ against $H_1:Q(p)<D_n$ with $K>0$.  That is, the true quantile function $Q(\cdot)$ is fixed while the hypothesis $D_n$ drifts toward $Q(p)$ (rather than a fixed $H_0$ and drifting true quantile). 
Type II error occurs if the lower one-sided CI contains $D_n$.  Using notation for $\tilde{Q}^B_X(\cdot)$ and related objects from the main text and applying Theorem \ref{thm:cdferror} and Lemma \ref{lem:u1u2-approx}, with $\tilde{u}$ between $u^h(\alpha)$ and $u^h(\alpha)+n^{-1/2}B(u^h(\alpha))$ and $\tilde{\tilde{u}}$ between $p$ and $u^h(\alpha)$, power is
\begin{align*}
\mathcal{P}^l_n(D_n) 
  &= P\left(D_n \not\in (-\infty, \hat Q^L_X(u^h(\alpha)))\right) 
\Supplemental{%
\\&= P\left(\hat Q^L_X(u^h(\alpha)) < D_n \right) \\
  &}{}= \overbrace{P\left(\hat Q^L_X(u^h(\alpha)) < Q(p)+Kn^{-1/2} \right)}^{\textrm{apply Theorem \ref{thm:cdferror}}} \\
  &= P\left(\tilde{Q}^B_X(u^h(\alpha)) < Q(p)+Kn^{-1/2} \right)
    +O\left(n^{-1/2}\log(n)\right) \\
  &= P\left(F^{-1}\left(\tilde{Q}^B_U(u^h(\alpha))\right) - F^{-1}(u^h(\alpha)) 
            < F^{-1}(p)-F^{-1}(u^h(\alpha)) +Kn^{-1/2} \right) \\
  &\quad
    +O\left(n^{-1/2}\log(n)\right) \\
  &= P\Biggl(\frac{\tilde{Q}^B_U(u^h(\alpha)) - u^h(\alpha)}{f(F^{-1}(u^h(\alpha)))} 
             +(1/2)Q''(\tilde{u})\left[\tilde{Q}^B_U(u^h(\alpha)) - u^h(\alpha)\right]^2 \\
  &\qquad\qquad
            < \frac{p-u^h(\alpha)}{f\left(F^{-1}\left(u^h(\alpha)\right)\right)} + (1/2)Q''\left(\tilde{\tilde{u}}\right)[p-u^h(\alpha)]^2 
             +Kn^{-1/2} \Biggr) \\
  &\quad
    +O\left(n^{-1/2}\log(n)\right) \\
\Supplemental{%
  &= P\Biggl(\frac{u^h(\alpha)+n^{-1/2}B(u^h(\alpha)) - u^h(\alpha)}{f(F^{-1}(u^h(\alpha)))} 
             +\frac{Q''(\tilde{u})}{2}\left[u^h(\alpha)+n^{-1/2}B(u^h(\alpha)) - u^h(\alpha)\right]^2 \\
  &\qquad\qquad
            < \frac{n^{-1/2}z_\alpha\sqrt{u^h(\alpha)(1-u^h(\alpha))}+O(n^{-1})}{f\left(F^{-1}\left(u^h(\alpha)\right)\right)} + \frac{Q''\left(\tilde{\tilde{u}}\right)}{2}O(n^{-1}) 
             +Kn^{-1/2} \Biggr) \\
  &\quad
    +O\left(n^{-1/2}\log(n)\right) \\
  &= P\Biggl(\frac{B(u^h(\alpha))}{f(F^{-1}(u^h(\alpha)))} 
             +(1/2)Q''(\tilde{u})n^{-1/2}\left[B(u^h(\alpha))\right]^2 \\
  &\qquad\qquad
            < \frac{z_\alpha\sqrt{u^h(\alpha)(1-u^h(\alpha))}}{f(F^{-1}(u^h(\alpha)))} 
             +O(n^{-1/2}) 
             +K \Biggr) \\
  &\quad
    +O\left(n^{-1/2}\log(n)\right) \\
}{}
  &= P\Biggl(\frac{B(u^h(\alpha))}{\sqrt{u^h(\alpha)(1-u^h(\alpha))}} 
             +n^{-1/2}(1/2)Q''(\tilde{u})\left[B(u^h(\alpha))\right]^2 \\
  &\qquad\qquad
            < z_\alpha
             +O(n^{-1/2}) 
             +K\frac{f(F^{-1}(u^h(\alpha)))}{\sqrt{u^h(\alpha)(1-u^h(\alpha))}} \Biggr) \\
  &\quad
    +O\left(n^{-1/2}\log(n)\right) , 
\end{align*}
where Assumption \ref{a:hut-pf} guarantees: $Q''(\tilde{\tilde{u}})$ is uniformly bounded since $u^h(\alpha)\to p$ implies $\tilde{\tilde{u}}\to p$, $O(n^{-1})/f(F^{-1}(u^h(\alpha)))=O(n^{-1})$ uniformly since $u^h(\alpha)\to p$, and $f(F^{-1}(u^h(\alpha)))<\infty$ uniformly, again since $u^h(\alpha)\to p$. 
By definition, $B(u^h(\alpha))\sim N\left(0,u^h(\alpha)[1-u^h(\alpha)]\right)$, so $Z\equiv B(u^h(\alpha))/\sqrt{u^h(\alpha)(1-u^h(\alpha))}\sim N(0,1)$.  Additionally, $\sqrt{u^h(\alpha)(1-u^h(\alpha))}=\sqrt{p(1-p)}+O(n^{-1/2})$ since $u^h(\alpha)=p+O(n^{-1/2})$ by Lemma \ref{lem:u1u2-approx}, and, additionally using \ref{a:hut-pf},
\begin{gather*}
Q(u^h(\alpha)) = Q(p)+(u^h(\alpha)-p)Q'(\tilde{u}) = Q(p) + O(n^{-1/2}), \\
\begin{split}
f(F^{-1}(u^h(\alpha)))
  &= f(F^{-1}(p)+O(n^{-1/2})) 
   = f(F^{-1}(p)) + O(n^{-1/2})f'(\overbrace{\tilde{q}}^{\to Q(p)}) \\
  &= f(F^{-1}(p)) +O(n^{-1/2}) . 
\end{split}
\end{gather*}

\Supplemental{
It remains to treat $Q''(\tilde{u})\left[B(u^h(\alpha))\right]^2$.  We consider a condition like \ref{cond:star} on $B(u^h(\alpha))$: 
\begin{equation*}
\left|\frac{B(u^h(\alpha))}{\sqrt{u^h(\alpha)[1-u^h(\alpha)]}}\right|<\sqrt{\log(n)} , 
\end{equation*}
which is equivalent to $|Z|<\sqrt{\log(n)}$ and was shown in \eqref{eqn:cond-star-B-prob} to have probability $1-o(n^{-1})$.  Under such a condition, $\left[B(u^h(\alpha))\right]^2=O(\log(n))$ and $\tilde{u}\to p$ so that $Q''(\tilde{u})$ is uniformly bounded by \ref{a:hut-pf}. 
Continuing from before, with $Z\sim N(0,1)$, 
}{%
Under the condition $|Z|<\sqrt{\log(n)}$, $\tilde{u}\to p$ so that \ref{a:hut-pf} applies; this condition can be shown to have probability $1-o(n^{-1})$.  Then,
}
\begin{align*}
\mathcal{P}^l_n(D_n) 
\Supplemental{
  &= P\Biggl(\frac{B(u^h(\alpha))}{\sqrt{u^h(\alpha)(1-u^h(\alpha))}} 
             +O(n^{-1/2})Q''(\tilde{u})\left[B(u^h(\alpha))\right]^2 \\
  &\qquad\qquad
            < z_\alpha
             +O(n^{-1/2}) 
             +K\frac{f(F^{-1}(u^h(\alpha)))}{\sqrt{u^h(\alpha)(1-u^h(\alpha))}} \Biggr) \\
  &\quad
    +O\left(n^{-1/2}\log(n)\right) \\
}{}
  &= P\Biggl(Z 
             +O(n^{-1/2})Q''(\tilde{u})\left[B(u^h(\alpha))\right]^2 \\
  &\qquad\qquad
            < z_\alpha
             +O(n^{-1/2}) 
             +K\frac{f(F^{-1}(p))}{\sqrt{p(1-p)}}+O(n^{-1/2}) \Biggr) \\
  &\quad
    +O\left(n^{-1/2}\log(n)\right) \\
\Supplemental{
  &= P\left(Z
            < z_\alpha
             +K\frac{f(F^{-1}(p))}{\sqrt{p(1-p)}}
             +O(n^{-1/2})Q''(\tilde{u})\left[B(u^h(\alpha))\right]^2 \right) \\
  &\quad
    +O(n^{-1/2})
    +O\left(n^{-1/2}\log(n)\right) \\
  &= P\left(Z < z_\alpha
               +K\frac{f(F^{-1}(p))}{\sqrt{p(1-p)}}
               +O(n^{-1/2})\overbrace{Q''(\tilde{u})}^{O(1)}\overbrace{\left[B(u^h(\alpha))\right]^2}^{O(\log(n))} 
            \mathrel{\bigg|} |Z|<\sqrt{\log(n)} \right) \\
  &\quad\times
     P\left( |Z|<\sqrt{\log(n)} \right) \\
  &\quad 
    +P\left(Z < z_\alpha
               +K\frac{f(F^{-1}(p))}{\sqrt{p(1-p)}}
               +O(n^{-1/2})Q''(\tilde{u})\left[B(u^h(\alpha))\right]^2 
            \mathrel{\bigg|} |Z|\ge\sqrt{\log(n)} \right) \\
  &\qquad\times
     P\left( |Z|\ge\sqrt{\log(n)} \right) \\
  &\quad
    +O(n^{-1/2})
    +O\left(n^{-1/2}\log(n)\right) \\
}{}
  &= 
         P\left(Z < z_\alpha +K\frac{f(F^{-1}(p))}{\sqrt{p(1-p)}} \textrm{ and }
               |Z|<\sqrt{\log(n)} \right)
      +O\left(n^{-1/2}\log(n)\right) \\
  &\quad 
    +O(1) o(n^{-1}) 
    +O(n^{-1/2})
    +O\left(n^{-1/2}\log(n)\right) \\
  &= \Phi\left(z_\alpha +K\frac{f(F^{-1}(p))}{\sqrt{p(1-p)}}\right)
    -\overbrace{\Phi\left(-\sqrt{\log(n)}\right)}^{o(n^{-1})}
    +O\left(n^{-1/2}\log(n)\right) \\
  &\to 
  \Phi\left(z_{\alpha} + \frac{K f(F^{-1}(p))}{\sqrt{p(1-p)}}\right) 
\Supplemental{,}{.} 
\end{align*}
\Supplemental{
where $\Phi(\cdot)$ is the standard normal CDF. 
Note that the result holds regardless of $K\ge0$ or $K\le0$. 
}{}

\Supplemental{
For the other one-sided direction, $H_0:Q(p)\le D_n\equiv Q(p)+Kn^{-1/2}$, the only changes are from $z_\alpha$ to $z_{1-\alpha}$, $u^h(\alpha)$ to $u^l(\alpha)$, and the direction of the inequality.  This leaves
\begin{align*}
\mathcal{P}^u_n(D_n) 
  &= P\left(Z > z_{1-\alpha} +K\frac{f(F^{-1}(p))}{\sqrt{p(1-p)}} \textrm{ and }
           |Z|<\sqrt{\log(n)} \right)
    +O\left(n^{-1/2}\log(n)\right) \\
  &= \Phi\left(z_\alpha -K\frac{f(F^{-1}(p))}{\sqrt{p(1-p)}}\right)
    -\overbrace{\Phi\left(-\sqrt{\log(n)}\right)}^{o(n^{-1})}
    +O\left(n^{-1/2}\log(n)\right) \\
  &\to 
  \Phi\left(z_{\alpha} - \frac{K f(F^{-1}(p))}{\sqrt{p(1-p)}}\right) . 
\end{align*}
As it should be, this is greater than $\alpha$ when $K<0$ and $D_n<Q(p)$, rather than when $K>0$ and $D_n>Q(p)$ like for $\mathcal{P}^l_n$. 

For a two-sided CI, using the two one-sided results,
\begin{align*}
\mathcal{P}^t_n(D_n) 
  &= P\left(D_n \not\in \left[ \hat{Q}^L_X(u^l(\alpha/2)), \hat Q^L_X(u^h(\alpha/2))\right] \right) \\
  &= P\left(\hat Q^L_X(u^h(\alpha/2)) < D_n \right) 
   + P\left(\hat Q^L_X(u^l(\alpha/2)) > D_n \right) \\
  &\to \Phi\left(z_{\alpha/2} + \frac{K f(F^{-1}(p))}{\sqrt{p(1-p)}}\right)
      +\Phi\left(z_{\alpha/2} - \frac{K f(F^{-1}(p))}{\sqrt{p(1-p)}}\right) . 
      \qedhere
\end{align*}
}{%
The other one-sided result follows symmetrically, and the two-sided result simply combines the two one-sided results. 
}
}{%
The results for power are derived using the normal approximation $\tilde{Q}^B_X(u)$ of $\hat{Q}^L_X(u)$, along with a first-order Taylor approximation and arguments that the remainder terms are negligible.  
}
\Supplemental{\end{proof}}{}

\subsection*{\Supplemental{Proof}{Sketch of proof} of Lemma \ref{lem:bias}\label{app:pf-bias}}

\Supplemental{
%
We start by taking the difference between identities for $Q_{Y|X}(p;C_h)$ and $Q_{Y|X}(p;0)$, where $x_0=0$.  Using the smoothness assumptions, an expansion of the difference may be taken.  The bias appears in the lowest-order term; some cancellation and rearrangement leads to the final expression.  Further manipulations to replace $Q_{Y|X}$ (and its derivatives) with $F_{Y|X}$ (and its derivatives) lead to the same expression as in \citet{BhattacharyaGangopadhyay1990} for $d=1$, $b=2$, as shown explicitly below. 

\subsubsection*{Case $b=2$, with $k_Q\ge2$ and $k_X\ge1$} 

Regardless of $d$ (other than the implicit dimension of $C_h$), by definition $Q_{Y|X}(p;C_h)$ satisfies
\begin{align}
\label{eqn:QYCh-def}
p &= \int_{C_h} \left\{\int_{-\infty}^{Q_{Y|X}(p;C_h)} f_{Y|X}(y;x)\,dy\right\}
           f_{X|C_h}(x)\,dx 
\end{align}
since $Y$ is (conditionally) continuous, where $f_{Y|X}(y;x)=f_{Y,X}(y,x)/f_X(x)$ is the conditional PDF of $Y$ given $X$ evaluated at $Y=y$ and $X=x$.  Similarly, $f_{X|C_h}(x)=f_X(x)/P(X\in C_h)$ 
is the conditional PDF of $X$ within $C_h$.  Also by definition, for any $x$, $Q_{Y|X}(p;x)$ satisfies
\begin{align}
\label{eqn:QYX-def}
p &= \int_{-\infty}^{Q_{Y|X}(p;x)} f_{Y|X}(y;x)\,dy . 
\end{align}

Decomposing the $\{\cdot\}$ term, \eqref{eqn:QYCh-def} becomes
\begin{align*}
p &= \int_{C_h} \Biggl\{\overbrace{\int_{-\infty}^{Q_{Y|X}(p;x)} f_{Y|X}(y;x)\,dy}^{=p\textrm{, by \eqref{eqn:QYX-def}}}
                      +\int_{Q_{Y|X}(p;x)}^{Q_{Y|X}(p;C_h)} f_{Y|X}(y;x)\,dy  \Biggr\}
           f_{X|C_h}(x)\,dx  \\
  &= \int_{C_h} \left\{p
                      +\int_{Q_{Y|X}(p;x)}^{Q_{Y|X}(p;C_h)} f_{Y|X}(y;x)\,dy  \right\}
           f_{X|C_h}(x) \,dx  \\
  &= p + \int_{C_h} \left\{ \int_{Q_{Y|X}(p;x)}^{Q_{Y|X}(p;C_h)} f_{Y|X}(y;x)\,dy  \right\}  f_{X|C_h}(x) \,dx , \\
0 &= \int_{C_h} \left\{ \int_{Q_{Y|X}(p;x)}^{Q_{Y|X}(p;C_h)} f_{Y|X}(y;x)\,dy  \right\}  f_{X|C_h}(x) \,dx .
\end{align*}
Our goal is to determine $Q_{Y|X}(p;C_h)-Q_{Y|X}(p;0)$, which is possible by expanding the right-hand side of the preceding equation. 

With $d=1$, $C_h=[-h,h]$ as in \eqref{eqn:Ch}.  Recall that the point of interest is $x_0=0$.  Using a change of variables $x=wh$ (so $dx=h\,dw$) and converting $f_{X|C_h}(\cdot)=f_X(\cdot)/P(X\in C_h)$,
\begin{align*}
0 &= \int_{C_h} \left\{ \int_{Q_{Y|X}(p;x)}^{Q_{Y|X}(p;C_h)} f_{Y|X}(y;x)\,dy  \right\}  f_{X|C_h}(x) \,dx \\
  &= \int_{-h}^{h} \left\{ \int_{Q_{Y|X}(p;x)}^{Q_{Y|X}(p;C_h)} f_{Y|X}(y;x)\,dy  \right\}  \frac{f_X(x)}{P(X\in C_h)} \,dx \\
  &= \frac{1}{P(X\in C_h)} 
     \int_{-1}^{1} \left\{ \int_{Q_{Y|X}(p;wh)}^{Q_{Y|X}(p;C_h)} f_{Y|X}(y;wh)\,dy  \right\}  f_X(wh) [h\,dw] \\
\label{eqn:bias-ID-d1}\refstepcounter{equation}\tag{\theequation}
  &= \frac{h}{P(X\in C_h)} \int_{-1}^1 \left\{ \int_{Q_{Y|X}(p;wh)}^{Q_{Y|X}(p;C_h)} f_{Y|X}(y;wh)\,dy  \right\}  f_X(wh) \,dw .
\end{align*}
With any $d\ge1$, so $x=(x_1,\ldots,x_d)'$ and $w=(w_1,\ldots,w_d)'$ are $d\times1$ column vectors (with $x'$ and $w'$ denoting the transposes, which are row vectors), we take $C_h=[-h,h]^d$ (a hypercube).  The change of variables yields
\begin{align*}
0
  &= \idotsint\limits_{[-h,h]^d} \left\{ \int_{Q_{Y|X}(p;x)}^{Q_{Y|X}(p;C_h)} f_{Y|X}(y;x)\,dy  \right\}  \frac{f_X(x)}{P(X\in C_h)} \,dx_1\cdots dx_d \\
  &= \frac{1}{P(X\in C_h)}
     \idotsint\limits_{[-1,1]^d} \left\{ \int_{Q_{Y|X}(p;wh)}^{Q_{Y|X}(p;C_h)} f_{Y|X}(y;wh)\,dy  \right\}  f_X(wh) [h\,dw_1]\cdots[h\,dw_d] \\
\label{eqn:bias-ID-d}\refstepcounter{equation}\tag{\theequation}
  &= \frac{h^d}{P(X\in C_h)}
     \idotsint\limits_{[-1,1]^d} \left\{ \int_{Q_{Y|X}(p;wh)}^{Q_{Y|X}(p;C_h)} f_{Y|X}(y;wh)\,dy  \right\}  f_X(wh) \,dw_1\cdots dw_d .
\end{align*}

With enough smoothness, a second-order expansion around $w=0$ can be taken of the integrand.  An explicit expression for the bias comes out of the zeroth-order term.  The first-order terms zero out, leaving the second-order terms of order $h^2$ to determine the bias. 

\Supplemental{%
Schematically, we are expanding 
\[ \left\{ \int_{a(x)}^b f(y;x)\,dy \right\} g(x) \] 
around $x=0$.  Let
\begin{equation*}
f^{(0,1)}(y;0) \equiv \Dz{f(y;x)}{x}, \quad
f^{(0,2)}(y;0) \equiv \DDz{f(y;x)}{x} , \quad
f^{(1,0)}(a;0) \equiv \Dz[a]{f(y;0)}{y} . 
\end{equation*}
With $d=1$, using the mean value theorem (MVT) result
\begin{equation}\label{eqn:integral-MVT}
\int_a^b f(x)\,dx = (b-a)f(\tilde{x})
\end{equation}
for some $\tilde{x}\in[a,b]$, 
as $x\to0$, 
\begin{align*}
& \int_{a(x)}^{b} f(y;x)\,dy\,g(x) \\
&= 
   \int_{a(0)}^{b} f(y;0)\,dy\,g(0)
  +\int_{a(0)}^{b} f(y;0)\,dy\,[g(x)-g(0)] \\
\label{eqn:bias-exp-alt-1}\refstepcounter{equation}\tag{\theequation}
&\quad
  +\int_{a(0)}^{b}\left[f(y;x)-f(y;0)\right]\,dy\,g(x)
  +\int_{a(x)}^{a(0)} f(y;x)\,dy\,g(x) \\
&= 
   \int_{a(0)}^{b} f(y;0)\,dy\,g(0)
  +\int_{a(0)}^{b} f(y;0)\,dy 
   \,\left\{xg'(0)+x\left[g'(\tilde{x})-g'(0)\right]\right\} \\
&\quad
  +\int_{a(0)}^{b}\left\{xf^{(0,1)}(y;0)+x\left[f^{(0,1)}(y;\tilde{x})-f^{(0,1)}(y;0)\right] \right\}\,dy \\
&\qquad\times
   \left\{ g(0) +xg'(0) + x\left[g'(\tilde{x})-g'(0)\right] \right\} \\
&\quad
  -\left\{
     \int\limits_{a(0)}^{a(x)} f(y;0)\,dy
    +\int_{a(0)}^{a(x)} \left[f(y;x)-f(y;0)\right]\,dy
   \right\} \\
&\qquad\times
   \left\{ g(0) +xg'(0) + x\left[g'(\tilde{x})-g'(0)\right] \right\} 
   \\
&= 
   \overbrace{\int_{a(0)}^{b} f(y;0)\,dy\,g(0)%
   }^{A} \\
&\quad
  +\overbrace{%
   \int_{a(0)}^{b} f(y;0)\,dy 
   \,\left\{xg'(0)\right\} %
   }^{B_1}
  +\overbrace{%
   \int_{a(0)}^{b} f(y;0)\,dy 
   \,\left\{x\left[g'(\tilde{x})-g'(0)\right]\right\}%
   }^{R_1}
   \\
&\quad
  +\overbrace{x \int_{a(0)}^{b} f^{(0,1)}(y;0)\,dy\,g(0)%
  }^{B_2} \\
&\quad
  +\overbrace{x \int_{a(0)}^{b} f^{(0,1)}(y;0)\,dy
   \,\left\{ xg'(0) \right\} %
   }^{R_2}
  +\overbrace{x \int_{a(0)}^{b} f^{(0,1)}(y;0)\,dy
   \,\left\{ x\left[g'(\tilde{x})-g'(0)\right] \right\}%
   }^{R_3} \\
&\quad
  +\overbrace{%
   \int_{a(0)}^{b}\left\{x\left[f^{(0,1)}(y;\tilde{x})-f^{(0,1)}(y;0)\right] \right\}\,dy%
  \,
  \left\{ g(0) +xg'(0) + x\left[g'(\tilde{x})-g'(0)\right] \right\} 
   }^{R_4} 
\\&\quad
   -\overbrace{\int_{a(0)}^{a(x)} f(y;0)\,dy
               \left\{ g(0) +xg'(0) + x\left[g'(\tilde{x})-g'(0)\right] \right\} %
    }^{D}
\\&\quad
  -\overbrace{%
   \int_{a(0)}^{a(x)} x f^{(0,1)}(y;\tilde{x}) \,dy
   \,g(0) %
   }^{C_1}
\\&\quad
  -\overbrace{\int_{a(0)}^{a(x)} \left[f(y;x)-f(y;0)\right]\,dy %
   \left\{ xg'(0) + x\left[g'(\tilde{x})-g'(0)\right] \right\} %
   }^{R_5} , 
\label{eqn:bias-exp-alt}\refstepcounter{equation}\tag{\theequation}
\end{align*}
where for term $D$,  
\begin{align*}
&\int_{a(0)}^{a(x)} f(y;0)\,dy \\
&= \int_{a(0)}^{a(x)} \left\{ f(a(0);0) +[y-a(0)]f^{(1,0)}(a(0);0) +(1/2)[y-a(0)]^2 f^{(2,0)}(\tilde{y};0) \right\}\,dy
\\&= \left[a(x)-a(0)\right] f(a(0);0)
    +f^{(1,0)}(a(0);0) \left[ (1/2)y^2-a(0)y\right]_{a(0)}^{a(x)} 
\\&\qquad
    +\left[ a(x)-a(0)\right] (1/2)[\tilde{y}-a(0)]^2 f^{(2,0)}(\tilde{y};0)
\label{eqn:TBD-pre1}\refstepcounter{equation}\tag{\theequation}
\\&= \left[ xa'(0) + (1/2)x^2a''(0)\right] f(a(0);0)
\\&\qquad
    +f^{(1,0)}(a(0);0) \left[ (1/2)a(x)^2-a(x)a(0)-(1/2)a(0)^2+a(0)^2\right]
    +O(x^3)
\\&= xa'(0)f(a(0);0) + (1/2)x^2a''(0)f(a(0);0)
\\&\qquad
    +f^{(1,0)}(a(0);0) (1/2)\overbrace{\left[ a(x)-a(0)\right]^2}^{=[xa'(0)+o(x)]^2}
    +O(x^3)
\\&= xa'(0)f(a(0);0) + (1/2)x^2a''(0)f(a(0);0)
    +f^{(1,0)}(a(0);0) (1/2)x^2[a'(0)]^2 
\\&\qquad
    +o(x^2)
    +O(x^3) 
,\\
D &= \left\{ 
     x f(a(0);0) a'(0)
     +(1/2)x^2 
      \left[ f^{(1,0)}(a(0);0)(a'(0))^2 + f(a(0);0) a''(0) \right]
     +o(x^2)
   \right\} 
\\&\quad\times   
   \left\{ g(0) +xg'(0) + x\left[g'(\tilde{x})-g'(0)\right] \right\} 
   . 
\end{align*}
}{}

For compactness, we define 
\begin{align*}
& \xi_p \equiv Q_{Y|X}(p;0), \\
& Q_{Y|X}^{(0,1)}(p;0) \equiv \Dz{Q_{Y|X}(p;x)}{x}, 
&& Q_{Y|X}^{(0,2)}(p;0) \equiv \DDz{Q_{Y|X}(p;x)}{x}, \\
& f_{Y|X}^{(0,1)}(y;0) \equiv \Dz{f_{Y|X}(y;x)}{x}, 
&& f_{Y|X}^{(1,0)}(\xi_p;0) \equiv \Dz[\xi_p]{f_{Y|X}(y;0)}{y}  .
\end{align*}

For the expansion of
\begin{align*}
\left\{ \int_{Q_{Y|X}(p;wh)}^{Q_{Y|X}(p;C_h)} f_{Y|X}(y;wh)\,dy  \right\}  f_X(wh)
\end{align*}
around $wh=0$, we use the earlier expansion and substitute 
\begin{gather*}
x = wh, \;
a(x) = Q_{Y|X}(p;wh), \;
g(x) = f_X(wh), \;
f(y;x) = f_{Y|X}(y;wh) , \;
b = Q_{Y|X}(p;C_h) . 
\end{gather*}
Equation \eqref{eqn:Qyx_02} below shows that the local smoothness (Definition \ref{def:smoothness}) of $f_{Y|X}(\cdot;\cdot)$ in $x$ is at least $s_Q-1$; directly, we see the existence of $Q_{Y|X}^{(0,2)}(p;x)$ implies the existence of $f_{Y|X}^{(0,1)}\left(Q_{Y|X}(p;x);x\right)$.  Thus, although the local smoothness is not assumed directly in the main text, we take it to be $s_Q-1$, with reference to \ref{a:Q}. 
Going back to \eqref{eqn:bias-exp-alt}, anticipating that the order of the bias is $b-a(0)=Q_{Y|X}(p;C_h)-Q_{Y|X}(p;0)=O(h^2)$, applying the MVT from \eqref{eqn:integral-MVT}, with $\tilde{y}$ and $\tilde{x}$ referring to values determined by the MVT (not necessarily the same across terms), 
\begin{align*}
R_1
&= 
\overbrace{%
         \overbrace{\left[Q_{Y|X}(p;C_h)-Q_{Y|X}(p;0)\right]}^{O(h^2)}
         \overbrace{f_{Y|X}(\tilde{y};0)}^{O(1)\textrm{ by \ref{a:GK2}}}%
       }^{\textrm{MVT}}
       (wh)
       \overbrace{\left[ f_X'(\tilde{w}h)-f_X'(0) \right]}^{O\left(h^{\gamma_X}\right)\textrm{ by \ref{a:f}}}
     = O\left( h^{3+\gamma_X} \right) 
       , \\
R_2
&=
(wh)^2
\overbrace{%
         \overbrace{\left[Q_{Y|X}(p;C_h)-Q_{Y|X}(p;0)\right]}^{O(h^2)}
         \overbrace{f_{Y|X}^{(0,1)}(\tilde{y};0)}^{O(1)\textrm{ by \ref{a:Q}, \eqref{eqn:Qyx_02}}}%
       }^{\textrm{MVT}}
       \overbrace{f_X'(0)}^{O\left(1\right)\textrm{ by \ref{a:f}, }s_X>1}
     = O\left( h^4 \right) 
       , \\
R_3
&=
O(h)
\overbrace{%
         \overbrace{\left[Q_{Y|X}(p;C_h)-Q_{Y|X}(p;0)\right]}^{O(h^2)}
         \overbrace{f_{Y|X}^{(0,1)}(\tilde{y};0)}^{O(1)\textrm{ by \ref{a:Q}, \eqref{eqn:Qyx_02}}}%
       }^{\textrm{MVT}}
       O(h)
       \overbrace{\left[ f_X'(\tilde{w}h)-f_X'(0) \right]}^{O\left(h^{\gamma_X}\right)\textrm{ by \ref{a:f}}}
     = O\left( h^{4+\gamma_X} \right) 
     , \\
R_4
&=
O(h)
\overbrace{%
         \overbrace{\left[Q_{Y|X}(p;C_h)-Q_{Y|X}(p;0)\right]}^{O(h^2)}
         \overbrace{f_{Y|X}^{(0,1)}(\tilde{y};h\tilde{w})-f_{Y|X}^{(0,1)}(\tilde{y};0)}^{O(h^{\gamma_Q})\textrm{ by \ref{a:Q}, \eqref{eqn:Qyx_02}}}%
       }^{\textrm{MVT}}
\overbrace{f_X(hw)}^{O(1)\textrm{ by \ref{a:f}}}
= O\left( h^{3+\gamma_Q} \right) 
, \\
C_1
&=
wh \overbrace{%
     \overbrace{\left[ Q_{Y|X}(p;wh)-Q_{Y|X}(p;0) \right]%
               }^{=whQ_{Y|X}^{(0,1)}(p;0)+o(h)\textrm{; \ref{a:Q}, }s_Q>1} 
     \overbrace{f_{Y|X}^{(0,1)}(\tilde{y};\tilde{x})%
               }^{=f_{Y|X}^{(0,1)}\left(Q_{Y|X}(p;0);0\right)+o(1)\textrm{; \eqref{eqn:Qyx_02}, \ref{a:Q}, }s_Q>1}%
   }^{\textrm{MVT}} 
f_X(0) \\
&= (wh)^2 Q_{Y|X}^{(0,1)}(p;0) f_{Y|X}^{(0,1)}\left(Q_{Y|X}(p;0);0\right) f_X(0)
  +o(h^2)
,\\
R_5
&=
\overbrace{%
  \overbrace{\left[ Q_{Y|X}(p;wh)-Q_{Y|X}(p;0) \right]%
            }^{=O(h)\textrm{; \ref{a:Q}, }s_Q\ge1} 
  \overbrace{wh f_{Y|X}^{(0,1)}(\tilde{y};\tilde{x})%
            }^{=O(h)\textrm{; \eqref{eqn:Qyx_02}, \ref{a:Q}, }s_Q>1}%
  }^{\textrm{MVT}} 
O(h)
 = O\left( h^3 \right) .
\end{align*}
For the $D$ term,
\begin{gather*}
wh \overbrace{\left[ f_X'(h\tilde{w}) - f_X'(0) \right]}^{o(1)\textrm{; \ref{a:f}, }s_X>1}
= o(h) , \\
\begin{split}
D
&= 
 \overbrace{wh f_{Y|X}\left( Q_{Y|X}(p;0); 0\right) Q_{Y|X}^{(0,1)}(p;0) f_X(0)}^{B_3}
\\&\quad
+\overbrace{(wh)^2 f_{Y|X}\left( Q_{Y|X}(p;0); 0\right) Q_{Y|X}^{(0,1)}(p;0) f_X'(0)}^{C_2}
\\&\quad
+\overbrace{\frac{(wh)^2}{2} 
   f_{Y|X}^{(1,0)}\left(Q_{Y|X}(p;0);0\right) \left[Q_{Y|X}^{(0,1)}(p;0)\right]^2 
 f_X(0)}^{C_3}
\\&\quad
+\overbrace{\frac{(wh)^2}{2} 
  f_{Y|X}\left(Q_{Y|X}(p;0);0\right) Q_{Y|X}^{(0,2)}(p;0) 
 f_X(0)}^{C_4}
\\&\quad
+\overbrace{o(h^2) +o(h^2) +O(h^3)}^{R_6} . 
\end{split}
\end{gather*}
Using these, \eqref{eqn:bias-exp-alt} reduces to
\begin{align*}
& A + B_1 + B_2 - C_1 - \overbrace{(B_3+C_2+C_3+C_4+R_6)}^{D}
+ \overbrace{R_1 + R_2 + R_3 + R_4 - R_5}^{o(h^2)} 
\\&= 
\overbrace{\left\{ \int_{\xi_p}^{Q_{Y|X}(p;C_h)} f_{Y|X}(y;0)\,dy  \right\}  f_X(0)}^{A}
\\&\quad
  +\overbrace{hw\int_{\xi_p}^{Q_{Y|X}(p;C_h)} f_{Y|X}(y;0)\,dy\, f_X'(0)}^{B_1}
  +\overbrace{hw\int_{\xi_p}^{Q_{Y|X}(p;C_h)} f_{Y|X}^{(0,1)}(y;0) dy\, f_X(0)}^{B_2}
\\&\quad
  -\overbrace{hw f_{Y|X}\left( Q_{Y|X}(p;0); 0\right) Q_{Y|X}^{(0,1)}(p;0) f_X(0)}^{B_3}
\\&\quad
  -\overbrace{h^2w^2 Q_{Y|X}^{(0,1)}(p;0) f_{Y|X}^{(0,1)}\left(Q_{Y|X}(p;0);0\right) f_X(0) +o(h^2)}^{C_1}
\\&\quad
  -\overbrace{h^2w^2 f_{Y|X}\left( Q_{Y|X}(p;0); 0\right) Q_{Y|X}^{(0,1)}(p;0) f_X'(0)}^{C_2}
\\&\quad
  -\overbrace{(1/2)h^2w^2 
 \left\{ 
   f_{Y|X}^{(1,0)}\left[Q_{Y|X}(p;0);0\right] \left(Q_{Y|X}^{(0,1)}(p;0)\right)^2 
  +f_{Y|X}\left(Q_{Y|X}(p;0);0\right) Q_{Y|X}^{(0,2)}(p;0) 
 \right\}
 f_X(0)}^{C_3+C_4}
\\&\quad
  +o(h^2) . 
\end{align*}

For $d=1$, noting $w\in[-1,1]$, this yields
\begin{gather*}
\label{eqn:bias-exp-d1}\refstepcounter{equation}\tag{\theequation}
\left\{ \int_{Q_{Y|X}(p;wh)}^{Q_{Y|X}(p;C_h)} f_{Y|X}(y;wh)\,dy  \right\}  f_X(wh)
= A + hwB - (1/2)h^2 w^2 C + o(h^2), \\
\begin{split}
A &= \left\{ \int_{\xi_p}^{Q_{Y|X}(p;C_h)} f_{Y|X}(y;0)\,dy  \right\}  f_X(0) , \\
B &= \Biggl\{ \int_{\xi_p}^{Q_{Y|X}(p;C_h)} f_{Y|X}(y;0)\,dy\, f_X'(0)
       + \int_{\xi_p}^{Q_{Y|X}(p;C_h)} f_{Y|X}^{(0,1)}(y;0) dy\, f_X(0)  \\
&\qquad- f_{Y|X}\left(\xi_p;0\right) Q_{Y|X}^{(0,1)}(p;0) f_X(0)  \Biggr\} , \\
\end{split}\\
\label{eqn:bias-C-approx}\refstepcounter{equation}\tag{\theequation}
\begin{split}
C &= \Biggl\{ 
       2 f_{Y|X}\left(\xi_p;0\right) Q_{Y|X}^{(0,1)}(p;0) f_X'(0) 
      +f_{Y|X}\left(\xi_p;0\right) Q_{Y|X}^{(0,2)}(p;0) f_X(0)  \\
&\qquad
      +2 f_{Y|X}^{(0,1)}(\xi_p;0) Q_{Y|X}^{(0,1)}(p;0) f_X(0)      
      +f_{Y|X}^{(1,0)}(\xi_p;0)\left[Q_{Y|X}^{(0,1)}(p;0)\right]^2 f_X(0) 
           \Biggr\}  . 
\end{split}
\end{gather*}

With $d\ge1$, let $\nabla h(t)$ denote the gradient of function $h(\cdot)$ evaluated at vector $t=(t_1,\ldots,t_d)'$, 
\begin{equation*}
\nabla h(t)
  \equiv \left. \left(\D{h(x)}{x_1},\ldots,\D{h(x)}{x_d} \right)' \right|_{x=t} , 
\end{equation*}
and we include a subscript $x$ if the function has multiple arguments: for function $f(y;x)$,
\begin{equation*}
\nabla_x f(y_0;x_0) 
  \equiv \left. \D{}{x}f(y;x) \right|_{y=y_0,x=x_0} .
\end{equation*}
Similarly, for Hessian notation, let
\begin{gather*}
\nabla^2 h(t)
  \equiv \left. \DD{}{x}h(x) \right|_{x=t} , \quad
\nabla^2_x f(y_0;x_0) 
  \equiv \left. \DD{}{x}f(y;x) \right|_{y=y_0,x=x_0} .
\end{gather*}
The expansion is analogous to \eqref{eqn:bias-exp-d1}:
\begin{gather*}
\label{eqn:bias-exp-d}\refstepcounter{equation}\tag{\theequation}
\left\{ \int_{Q_{Y|X}(p;wh)}^{Q_{Y|X}(p;C_h)} f_{Y|X}(y;wh)\,dy  \right\}  f_X(wh)
= A + hw' B - (1/2)h^2 w' C w +o\left(h^2\right), \\
\begin{split}
A &= \int_{Q_{Y|X}(p;0)}^{Q_{Y|X}(p;C_h)} f_{Y|X}(y;0) \,dy\, f_X(0) , \\
B &=\Biggl\{ 
      \int_{Q_{Y|X}(p;0)}^{Q_{Y|X}(p;C_h)} f_{Y|X}(y;0) \,dy\, \nabla f_X(0) 
      \\
 &\qquad
     +\left[
        \int_{Q_{Y|X}(p;0)}^{Q_{Y|X}(p;C_h)} \nabla_x f_{Y|X}(y;0) \,dy 
       -f_{Y|X}\left(Q_{Y|X}(p;0);0\right) \nabla_x Q_{Y|X}(p;0) 
      \right] f_X(0) 
    \Biggr\} , \\
-C&=\Biggl\{ 
     -2 f_{Y|X}\left(Q_{Y|X}(p;0);0\right) \nabla_x Q_{Y|X}(p;0) \nabla f_X(0)' 
\\&\qquad
     -2\nabla_x f_{Y|X}\left(Q_{Y|X}(p;0);0\right) \nabla_x Q_{Y|X}(p;0)' f_X(0)
\\&\qquad 
     -f_{Y|X}\left(Q_{Y|X}(p;0);0\right) \nabla^2_x Q_{Y|X}(p;0) f_X(0)
\\&\qquad 
             -\Dz[Q_{Y|X}(p;0)]{f_{Y|X}(y;0)}{y} 
              \nabla_x Q_{Y|X}(p;0) 
              \nabla_x Q_{Y|X}(p;0)' 
              f_X(0) 
    \Biggr\} 
    . 
\end{split}
\end{gather*}

Going back to 
\eqref{eqn:bias-ID-d}, we now must solve 
\begin{align*}
0 &= \frac{h^d}{P\left(X\in C_h\right)}
     \idotsint\limits_{[-1,1]^d} 
         \left[ A + hw' B - (1/2)h^2 w' C w +o\left(h^2\right) \right] 
     \,dw_1\cdots dw_d 
\end{align*}
for $A$, inside of which is the bias. 
To extract the bias from $A$, we need another expansion: 
\begin{align*}
\int_{\xi_p}^{Q_{Y|X}(p;C_h)} f_{Y|X}(y;0) \,dy  
  &= \int_{\xi_p}^{Q_{Y|X}(p;C_h)} 
     \left[f_{Y|X}(\xi_p;0) +f_{Y|X}^{(1,0)}(\tilde y;0) \left(y-\xi_p\right)\right] \,dy \\
  &= f_{Y|X}(\xi_p;0) \left[Q_{Y|X}(p;C_h) - \xi_p\right]
    +O\left(\left[Q_{Y|X}(p;C_h) - \xi_p\right]^2\right),
\end{align*}
where $\tilde y$ lies between the lower and upper limits of integration, so $\tilde{y}\to\xi_p$ and thus $f_{Y|X}^{(1,0)}(\tilde y;0)$ is uniformly bounded by Assumption \ref{a:GK2}.  Anticipating that the bias is $\left[Q_{Y|X}(p;C_h) - \xi_p\right]=O(h^{2})$, 
\begin{equation}\label{eqn:bias-A}
A  = f_{Y|X}(\xi_p;0) \left[Q_{Y|X}(p;C_h) - \xi_p\right] f_X(0)
    +O\left(h^{4}\right)  .
\end{equation}
Since there is no $w$ in $A$, $\int_{-1}^1A\,dw=A\int_{-1}^{1}dw=2A$; more generally, integrating over $[-1,1]^d$ yields $2^d A$. 

The $B$ term zeroes out since it is linear in $w$:
\begin{equation}\label{eqn:bias-B}
\idotsint\limits_{[-1,1]^d} hw'B\,dw
  = hB' \overbrace{\idotsint\limits_{[-1,1]^d} w\,dw}^{=0}
  = 0 .
\end{equation}

The $C$ term is quadratic in $w$; $C$ itself does not depend on $w$. 
With $d=1$, 
\[ \int_{-1}^1(1/2)h^2w^2C\,dw
  =(1/2)h^2 C \int_{-1}^{1}w^2\,dw
  =(1/2)h^2 C (2/3)
  =h^2C/3 . \]

With $d\ge1$, 
\begin{equation*}
\idotsint\limits_{[-1,1]^d} w' C w\, dw_1\cdots dw_d 
  = \idotsint\limits_{[-1,1]^d} \sum_{i=1}^{d}\sum_{k=1}^{d} C_{ik}w_i w_k \, dw_1\cdots dw_d . 
\end{equation*}
For the $i\ne k$ terms,
\begin{align*}
\idotsint\limits_{[-1,1]^d} C_{ik}w_i w_k \, dw_1\cdots dw_d
  &= C_{ik} 
     \idotsint\limits_{[-1,1]^{d-2}} 
     \int_{-1}^{1} w_i 
     \overbrace{\int_{-1}^{1} w_k\,dw_k}^{=0}\,dw_i\,dw_1\cdots dw_d
   = 0 . 
\end{align*}
For the $i=k$ terms,
\begin{align*}
\idotsint\limits_{[-1,1]^d} C_{ii}w_i^2 \, dw_1\cdots dw_d
  &= C_{ii} 
     \idotsint\limits_{[-1,1]^{d-1}} 
     \int_{-1}^{1} w_i^2\,dw_i
     \,dw_1\cdots dw_{i-1}\,dw_{i+1}\cdots dw_d \\
  &= C_{ii} 
     \idotsint\limits_{[-1,1]^{d-1}} 
     (2/3)
     \,dw_1\cdots dw_{i-1}\,dw_{i+1}\cdots dw_d \\
  &= C_{ii} (2/3) 2^{d-1}
   = 2^d C_{ii} / 3 . 
\end{align*}
Altogether,
\begin{equation}\label{eqn:bias-C}
\idotsint\limits_{[-1,1]^d} (h^2/2)w' C w\, dw_1\cdots dw_d 
   = (h^2/2) \sum_{i=1}^{d} 2^d C_{ii} / 3
   = (h^2/2)(2^d/3) \trp{C} . 
\end{equation}

Plugging results from \eqref{eqn:bias-A}, \eqref{eqn:bias-B}, and \eqref{eqn:bias-C} into \eqref{eqn:bias-ID-d}, for any $d$, 
\begin{align*}
0 
&= \frac{h^d}{P(X\in C_h)}
     \idotsint\limits_{[-1,1]^d} \left[ A + hw'B -(1/2)h^2 w'Cw +o(h^2) \right] \,dw_1\cdots dw_d , \\
0 &= 2^d A + 0 - (1/2)h^2 (2^d/3) \trp{C} +o(h^2) , \\
(h^2/6) \trp{C} + o(h^2)
&= A
 = f_{Y|X}(\xi_p;0) \left[Q_{Y|X}(p;C_h) - \xi_p\right] f_X(0)
    +O\left(h^{4}\right) , \\
Q_{Y|X}(p;C_h) - \xi_p
&= \frac{h^2}{6}
   \frac{\trp{C}}{f_{Y|X}(\xi_p;0) f_X(0)}
  +o(h^2) 
 = O(h^2) . 
\end{align*}

For $d=1$, plugging in for $C$ from \eqref{eqn:bias-C-approx}, 
\begin{align}\notag
& Q_{Y|X}(p;C_h) - \xi_p  \\ \notag
  &= h^2 \frac{C}{6 f_{Y|X}\left(\xi_p;0\right) f_X(0)} +o(h^2)  \\ \notag
  &= \frac{h^2}{6}   
     \biggl\{ 
        2 Q_{Y|X}^{(0,1)}(p;0) f_X'(0) / f_X(0)
       +Q_{Y|X}^{(0,2)}(p;0)
       +2 f_{Y|X}^{(0,1)}(\xi_p;0) Q_{Y|X}^{(0,1)}(p;0) / f_{Y|X}\left(\xi_p;0\right) \\
  &\qquad\quad 
       +f_{Y|X}^{(1,0)}(\xi_p;0)
        \left[Q_{Y|X}^{(0,1)}(p;0)\right]^2 
        / f_{Y|X}\left(\xi_p;0\right)
     \biggr\} +o(h^2)  .
\label{eqn:bias-orig-d1}
\end{align}

The bias expression in \eqref{eqn:bias-orig-d1} is equivalent to the bias in \citet[Thm.\ K1]{BhattacharyaGangopadhyay1990}, just in different notation.  Since $Y$ is (conditionally) continuous, by definition, for all $x$,
\begin{align*}
F_{Y|X}\left(Q_{Y|X}(p;x);x\right) &= p  .
\end{align*}

Differentiating once with respect to $x$, 
\begin{align*}
0 &= Q_{Y|X}^{(0,1)}(p;x)f_{Y|X}\left(Q_{Y|X}(p;x);x\right) 
    +F_{Y|X}^{(0,1)}\left(Q_{Y|X}(p;x);x\right), \\
\label{eqn:Qyx_01}\refstepcounter{equation}\tag{\theequation}
Q_{Y|X}^{(0,1)}(p;x) 
  &= -\frac{F_{Y|X}^{(0,1)}\left(Q_{Y|X}(p;x);x\right)}
           {f_{Y|X}\left(Q_{Y|X}(p;x);x\right)}  .
\end{align*}
Differentiating again with respect to $x$, 
\begin{align*}
0 &= Q_{Y|X}^{(0,2)}(p;x) f_{Y|X}\left(Q_{Y|X}(p;x);x\right)
    +F_{Y|X}^{(0,2)}\left(Q_{Y|X}(p;x);x\right)  \\
  &\quad 
    +Q_{Y|X}^{(0,1)}(p;x) f_{Y|X}^{(0,1)}\left(Q_{Y|X}(p;x);x\right)
    +\left[Q_{Y|X}^{(0,1)}(p;x)\right]^2 f_{Y|X}^{(1,0)}\left(Q_{Y|X}(p;x);x\right)  \\
  &\quad 
    +f_{Y|X}^{(0,1)}\left(Q_{Y|X}(p;x);x\right) Q_{Y|X}^{(0,1)}(p;x),   \\
Q_{Y|X}^{(0,2)}(p;x)
  &= -\frac{1}{f_{Y|X}\left(Q_{Y|X}(p;x);x\right)} \\
\label{eqn:Qyx_02}\refstepcounter{equation}\tag{\theequation}
  &\quad\times
     \biggl\{
       F_{Y|X}^{(0,2)}\left(Q_{Y|X}(p;x);x\right)
      +2 Q_{Y|X}^{(0,1)}(p;x)f_{Y|X}^{(0,1)}\left(Q_{Y|X}(p;x);x\right) \\
  &\qquad\quad
      +\left[Q_{Y|X}^{(0,1)}(p;x)\right]^2 f_{Y|X}^{(1,0)}\left(Q_{Y|X}(p;x);x\right)
     \biggr\}  .
\end{align*}

Plugging these substitutions into \eqref{eqn:bias-orig-d1}, 
\begin{align*}
\textrm{Bias}
  &= \frac{h^2}{6}
     \Biggl\{
       -\frac{2F_{Y|X}^{(0,1)}(\xi_p;0) f_X'(0)}{f_X(0) f_{Y|X}(\xi_p;0)} \\
  &\qquad\quad   
      -\frac{1}{f_{Y|X}(\xi_p;0)}
       \biggl\{
         F_{Y|X}^{(2,0)}(\xi_p;0) 
        -2F_{Y|X}^{(0,1)}(\xi_p;0) f_{Y|X}^{(0,1)}(\xi_p;0)/f_{Y|X}(\xi_p;0) \\
  &\qquad\qquad\qquad\qquad\qquad   
        +\left[F_{Y|X}^{(0,1)}(\xi_p;0)\right]^2 f_{Y|X}^{(1,0)}(\xi_p;0)/f_{Y|X}(\xi_p;0)
       \biggr\}     \\
  &\qquad\quad   
      -\frac{2 f_{Y|X}^{(0,1)}(\xi_p;0) F_{Y|X}^{(0,1)}(\xi_p;0)}
            {\left[f_{Y|X}(\xi_p;0)\right]^2}
      +\frac{f_{Y|X}^{(1,0)}(\xi_p;0) \left[F_{Y|X}^{(0,1)}(\xi_p;0)\right]^2}
            {\left[f_{Y|X}(\xi_p;0)\right]^3}
     \Biggr\}    \\
\Supplemental{%
  &= \frac{h^2}{6 [f_{Y|X}(\xi_p;0)]^3} \\
  &\quad\times  \Biggl\{
       -2 F_{Y|X}^{(0,1)}(\xi_p;0) [f_{Y|X}(\xi_p;0)]^2 f_X'(0)/f_X(0) 
       -F_{Y|X}^{(0,2)}(\xi_p;0) [f_{Y|X}(\xi_p;0)]^2  \\
  &\qquad\quad
       +2 F_{Y|X}^{(0,1)}(\xi_p;0) f_{Y|X}^{(0,1)}(\xi_p;0) f_{Y|X}(\xi_p;0)
       -[F_{Y|X}^{(0,1)}(\xi_p;0)]^2 f_{Y|X}^{(1,0)}(\xi_p;0)  \\
  &\qquad\quad
       -2 F_{Y|X}^{(0,1)}(\xi_p;0) f_{Y|X}^{(0,1)}(\xi_p;0) f_{Y|X}(\xi_p;0)
       +[F_{Y|X}^{(0,1)}(\xi_p;0)]^2 f_{Y|X}^{(1,0)}(\xi_p;0)
     \Biggr\}  \\}{}
  &= -h^2
      \frac{f_X(0) F_{Y|X}^{(0,2)}(\xi_p;0)
            +2 f_X'(0) F_{Y|X}^{(0,1)}(\xi_p;0)}
           {6 f_X(0) f_{Y|X}(\xi_p;0)}   .
\end{align*}
\Supplemental{This is equivalent to the bias in \citet[Thm.\ K1]{BhattacharyaGangopadhyay1990}.  To convert from their notation, note that their bandwidth is $h/2$ (and $(h/2)^2/6=h^2/24$), $x_0=0$, and 
\[ g(\xi)\equiv f_{Y|X}(\xi_p;0), \quad G_x(\xi\mid x_0)\equiv F_{Y|X}^{(0,1)}(\xi_p;0), \quad G_{xx}(\xi\mid x_0)\equiv F_{Y|X}^{(0,2)}(\xi_p;0). \]}{}

\subsubsection*{Case $b\le2$}

In the prior calculations, we used $k_X\ge1$ and $k_Q\ge2$, in addition to $k_Y\ge2$ from \ref{a:GK2}. 
By examining the original expansion, we can see how the order of the bias grows as $s_X$ and $s_Q$ decrease.  We also let $s_f=k_f+\gamma_f$ denote the local smoothness of $f_{Y|X}(y;x)$ with respect to $x$, noting that $s_f\ge s_Q-1$ from \eqref{eqn:Qyx_01} and \eqref{eqn:Qyx_02}. 

The expression for $A$ in \eqref{eqn:bias-A} is the same except for the remainder, which is the squared bias.  Anticipating that the bias is $B_h=O(h^{b})$, now
\begin{equation}\label{eqn:bias-A-ble2}
A  = f_{Y|X}(\xi_p;0) \left[Q_{Y|X}(p;C_h) - \xi_p\right] f_X(0)
    +O\left(h^{2b}\right) 
   = B_h O(1) + O(B_h^2) .
\end{equation}
Thus, the order of magnitude of $B_h$ equals the largest remaining term (outside $A$) integrated over $w\in[-1,1]^d$. 





Returning to the first step of the expansion in \eqref{eqn:bias-exp-alt-1},
\begin{align*}
&\int_{a(x)}^{b} f(y;x)\,dy\,g(x)
\\&= A
  +\int_{a(0)}^{b} f(y;0)\,dy\,[g(x)-g(0)]
\\&\quad
  +\int_{a(0)}^{b}\left[f(y;x)-f(y;0)\right]\,dy\,g(x)
  +\int_{a(x)}^{a(0)} f(y;x)\,dy\,g(x) 
\\&= A 
  +\overbrace{%
   \overbrace{\left[ b - a(0) \right]}^{B_h} 
   \overbrace{f(\tilde{y};0) }^{O(1)}
   \overbrace{\left[g(x)-g(0)\right]}^{xg'(0)+o(x)\textrm{ if }k_X=1\textrm{, }O(x^{s_X})\textrm{ if }k_X=0} 
   }^{o(B_h)}
\\&\quad
  +\overbrace{%
   \overbrace{\left[ b - a(0) \right]}^{B_h} 
   \overbrace{\left[ f(\tilde{y};x)-f(\tilde{y};0) \right]}^{O\left(x^{\min\{1,s_f\}}\right)} 
   O(1) 
   }^{o(B_h)}
\\&\quad
  -\int_{a(0)}^{a(x)} \left\{ f\left(a(0);0\right) + \left[y-a(0)\right]f^{(1,0)}\left(a(0);0\right) + (1/2)\left[y-a(0)\right]^2 f^{(2,0)}(\tilde{y};0) \right\} \,dy
\\&\qquad\times
  \left\{ g(0) +\left[g(x)-g(0)\right] \right\} 
\\&\quad
  -\int_{a(0)}^{a(x)} \left[ f\left(y;x\right)-f\left(y;0\right) \right] \,dy
   \,\left\{ g(0) +\overbrace{\left[g(x)-g(0)\right]}^{o(1)} \right\} 
\\&= A 
  -\int_{a(0)}^{a(x)} f\left(a(0);0\right)  \,dy
  \,g(0)
  -\int_{a(0)}^{a(x)} f\left(a(0);0\right)  \,dy
  \,\left[g(x)-g(0)\right] 
\\&\quad
  -\int_{a(0)}^{a(x)} \left[y-a(0)\right]f^{(1,0)}\left(a(0);0\right) \,dy
  \,\left\{ g(0) +\left[g(x)-g(0)\right] \right\} 
\\&\quad
  -\int_{a(0)}^{a(x)} (1/2)\left[y-a(0)\right]^2 f^{(2,0)}(\tilde{y};0) \,dy
  \,\left\{ g(0) +\left[g(x)-g(0)\right] \right\} 
\\&\quad
  -\left[a(x)-a(0)\right] 
   \left[ f\left(\tilde{y};x\right)-f\left(\tilde{y};0\right) \right] 
   \left\{ g(0) +o(1) \right\} 
\\&\quad
  +o(B_h) 
\\&= A 
  -\left[ a(x) - a(0)\right] \overbrace{f\left(a(0);0\right) g(0)}^{O(1)}
  -\left[ a(x) - a(0)\right] \overbrace{f\left(a(0);0\right)}^{O(1)}
   \left[g(x)-g(0)\right] 
\\&\quad
  -\int_{a(0)}^{a(x)} \left[y-a(0)\right]f^{(1,0)}\left(a(0);0\right) \,dy
  \,\left\{ g(0) +\left[g(x)-g(0)\right] \right\} 
\\&\quad
  -(1/2) \overbrace{\left[a(x)-a(0)\right] \left[\tilde{y}-a(0)\right]^2 }^{O\left([a(x)-a(0)]^3\right)}
  \overbrace{f^{(2,0)}(\tilde{\tilde{y}};0) }^{O(1)\textrm{ by \ref{a:GK2}}} 
  \overbrace{\left\{ g(0) +\left[g(x)-g(0)\right] \right\} }^{O(1)}
\\&\quad
  -\left[a(x)-a(0)\right] 
   \left[ f\left(\tilde{y};x\right)-f\left(\tilde{y};0\right) \right] 
   \left\{ g(0) +o(1) \right\} 
\\&\quad
  +o(B_h) 
\\&= A 
  -\left[ a(x) - a(0)\right] f\left(a(0);0\right) g(0)
  -\left[ a(x) - a(0)\right] f\left(a(0);0\right) \left[g(x)-g(0)\right] 
\\&\quad
  -f^{(1,0)}\left(a(0);0\right) \overbrace{\left[ (1/2)y^2-a(0)y \right]_{a(0)}^{a(x)}}^{(1/2)a(x)^2-a(x)a(0)-(1/2)a(0)^2+a(0)^2}
  \left\{ g(0) +\left[g(x)-g(0)\right] \right\} 
  -\overbrace{O\left([a(x)-a(0)]^3\right)}^{O(x^{3\min\{1,s_Q\}})}
\\&\quad
  -\left[a(x)-a(0)\right] 
   \left[ f\left(\tilde{y};x\right)-f\left(\tilde{y};0\right) \right] 
   \left\{ g(0) +o(1) \right\} 
\\&\quad
  +o(B_h) 
 . 
\label{eqn:bias-exp-ble2-1}\refstepcounter{equation}\tag{\theequation}
\end{align*}

From \eqref{eqn:bias-exp-ble2-1}, with $k_Q=2$ (implying $k_f\ge1$, from \eqref{eqn:Qyx_02}) and $k_X=1$, 
\begin{align*}
A
&=
   \left[ xa'(0)+(1/2)x^2a''(0)+o(x^2)\right] f\left(a(0);0\right) g(0)
  +\left[ xa'(0)+o(x)\right] f\left(a(0);0\right) \left[xg'(0)+o(x)\right] 
\\&\quad
  +f^{(1,0)}\left(a(0);0\right) 
   (1/2)[a(x)-a(0)]^2
  \left[ g(0) +o(1) \right]
 +O\left(x^3\right)
\\&\quad
  +\left[xa'(0)+o(x)\right] 
   \left[ xf^{(0,1)}(\overbrace{\tilde{y}}^{=a(\tilde{x})\to a(0)};0) +o(x) \right] 
   \left[ g(0) +o(1) \right] 
  +o(B_h)
\\&=
   xa'(0) f\left(a(0);0\right) g(0)
  +(1/2)x^2a''(0) f\left(a(0);0\right) g(0)
  +xa'(0) f\left(a(0);0\right) xg'(0) 
\\&\quad
  +f^{(1,0)}\left(a(0);0\right) 
   (1/2)[xa'(0)+o(x)]^2
   g(0) 
  +xa'(0)
   x\left[ f^{(0,1)}\left(a(0);0\right) +o(1)\right]
   g(0) 
\\&\quad
  +o\left(x^2\right)
  +o(B_h)
\\&=
   x a'(0) f\left(a(0);0\right) g(0)
\\&\quad
  +(1/2)x^2
   \Bigl\{ 
     a''(0) f\left(a(0);0\right) g(0)
    +2 a'(0) f\left(a(0);0\right) g'(0) 
\\&\quad\qquad\qquad\quad
    +f^{(1,0)}\left(a(0);0\right) [a'(0)]^2 g(0) 
    +2 a'(0) f^{(0,1)}\left(a(0);0\right) g(0) 
   \Bigr\}
\\&\quad
  +o\left(x^2\right)
  +o(B_h)
 ,
\end{align*}
matching \eqref{eqn:bias-exp-d1} (with some terms from $B$ now in the $o(B_h)$ remainder). 

From \eqref{eqn:bias-exp-ble2-1}, we can identify the orders of magnitude of the various terms, in terms of $s_Q$, $s_X$, and $s_f$.  Integrating $x=wh$ over $w\in[-1,1]$ complicates only the first term: if $k_Q=1$, then $a(x)-a(0)=xa'(0)+O(x^{s_Q})=O(x)$, but since $a'(0) \int_{-h}^{h} x\,dx=0$, the eventual result is $O(x^{s_Q})$, so we label it as such.  Thus,
\begin{align*}
A 
&= 
   \left[ a(x) - a(0)\right] \overbrace{f\left(a(0);0\right) g(0)}^{O(1)\textrm{ constant}}
  +\left[ a(x) - a(0)\right] \left[g(x)-g(0)\right] 
   \overbrace{ f\left(a(0);0\right)}^{O(1)}
\\&\quad
  +\overbrace{f^{(1,0)}\left(a(0);0\right) (1/2)}^{O(1)}
   [a(x)-a(0)]^2
   \overbrace{\left[ g(0) +o(1) \right] }^{O(1)}
\\&\quad
  +\left[a(x)-a(0)\right] 
   \left[ f\left(\tilde{y};x\right)-f\left(\tilde{y};0\right) \right] 
   \overbrace{\left[ g(0) +o(1) \right] }^{O(1)}
  +O\left(x^{3\min\{1,s_Q\}}\right) +o(B_h) 
\\&= 
   \overbrace{\left[ a(x) - a(0)\right] O(1)}^{O\left(x^{\min\{2,s_Q\}}\right)}
  +\overbrace{\left[ a(x) - a(0)\right]}^{O\left(x^{\min\{1,s_Q\}}\right)}
   \overbrace{\left[g(x)-g(0)\right] }^{O\left(x^{\min\{1,s_X\}}\right)}
   O(1)
\\&\quad
  +O(1) \overbrace{[a(x)-a(0)]^2}^{O\left(x^{2\min\{1,s_Q\}}\right)}
  +\overbrace{\left[a(x)-a(0)\right] }^{O\left(x^{\min\{1,s_Q\}}\right)}
   \overbrace{\left[ f\left(\tilde{y};x\right)-f\left(\tilde{y};0\right) \right] }^{O\left(x^{\min\{1,s_f\}}\right)}
   O(1)
  +O\left(x^{3\min\{1,s_Q\}}\right) +o(B_h) 
 . 
\end{align*}
Now, $3\min\{1,s_Q\}>\min\{2,s_Q\}$ and $2\min\{1,s_Q\}\ge\min\{2,s_Q\}$, so neither of those terms can dominate, nor can $o(B_h)$.  This leaves three terms to consider; the dominant term will be $O(x^b)$ with
\begin{align*}
b
&= \min\left\{ \min\{2,s_Q\}, 
              \min\{1,s_Q\}+\min\{1,s_X\}, 
              \min\{1,s_Q\}+\min\{1,s_f\} \right\} \\
&= \left\{ \begin{array}{cl}
              \min\left\{2,1+\min\{1,s_X\},1+\min\{1,s_f\}\right\}
              =1+\min\{1,s_X,s_f\} & \textrm{if }k_Q\ge2 \\
              \min\left\{s_Q,1+\min\{1,s_X\},1+\min\{1,s_f\}\right\}
              =\min\left\{s_Q,1+\min\{s_X,s_f\}\right\} & \textrm{if }k_Q=1 \\
              \min\left\{s_Q,s_Q+\min\{1,s_X\},s_Q+\min\{1,s_f\}\right\}
              =s_Q & \textrm{if }k_Q=0 \\
          \end{array} \right.
\\&= \min\left\{ 2, s_Q, 1+s_X, 1+s_f\right\}
. 
\end{align*}
Since $s_f\ge s_Q-1$ from \eqref{eqn:Qyx_01} and \eqref{eqn:Qyx_02}, then $1+s_f\ge s_Q$, so $b=\min\left\{ 2, s_Q, 1+s_X\right\}$. %
\qed
}{
Since the result is similar to other kernel bias results, and since the special case of $d=1$ and $b=2$ is already given in \citet{BhattacharyaGangopadhyay1990}, we leave the proof to the supplemental appendix and provide only a very brief sketch here.  The approach is to start from the definitions of $Q_{Y|X}(p;C_h)$ and $Q_{Y|X}(p;x)$,
\begin{align*}
p &= \int_{C_h} \left\{\int_{-\infty}^{Q_{Y|X}(p;C_h)} f_{Y|X}(y;x)\,dy\right\}
           f_{X|C_h}(x)\,dx, \quad
p = \int_{-\infty}^{Q_{Y|X}(p;x)} f_{Y|X}(y;x)\,dy, 
\quad \textrm{so} \\
0 &= \int_{C_h} \left\{ \int_{Q_{Y|X}(p;x)}^{Q_{Y|X}(p;C_h)} f_{Y|X}(y;x)\,dy  \right\}  f_{X|C_h}(x) \,dx .
\end{align*}
After a change of variables to $w=x/h$, an expansion around $w=0$ is taken, and the bias can be isolated.  If $b=2$, $k_Q\ge2$, and $k_X\ge1$, then a second-order expansion is justified; otherwise, the smoothness determines the order of both the expansion and the remainder.  
}

\subsection*{\Supplemental{Proof}{Sketch of proof} of Theorem \ref{thm:h-rate}\label{app:pf-h-rate}}

As in \citet{Chaudhuri1991}, we consider a deterministic bandwidth sequence, leaving treatment of a random (data-dependent) bandwidth to future work. 
\Supplemental{

Theorem \ref{thm:IDEAL-single} states the CPE for \citeposs{Hutson1999} CIs in terms of $n$.  In the conditional setting, we instead have the local sample size $N_n$, which is random (binomially distributed) and defined in \eqref{eqn:Ch}. 
\Supplemental{From \citet[proof of Thm.\ 3.1, p.\ 769]{Chaudhuri1991}, using \ref{a:sampling}, \ref{a:f}, and Bernstein's inequality, we can choose $c_1,c_2,c_3,c_4>0$ such that 
\[ P\left(A_n\right) \ge 1-c_3\exp(-c_4 n h^d) \] 
for all $n$, where $A_n\equiv \left\{c_1 n h^d\le N_n\le c_2 n h^d\right\}$, using a hypercube $C_h$ from Definition \ref{def:cube}.  
If the rate of $h$ leads to $\sum_{n=1}^{\infty} [1-P(A_n)]<\infty$, then the Borel--Cantelli Lemma gives $P(\liminf A_n)=1$.  That is, with probability 1, there exists $n_0$ such that for all $n>n_0$, $A_n$ is true and $N_n$ equals $nh^d$ times a constant in the interval $[c_1,c_2]$, so we may use $N_n\asymp nh^d$ in our asymptotic analysis.  
The optimal bandwidth rates derived here all satisfy \ref{a:h} and thus all satisfy the summability condition:
\begin{align*}
\sum_{n=1}^{\infty} 1-\left[ 1-c_3\exp(-c_4 n h^d) \right]
&=
c_3 \sum_{n=1}^{\infty} \exp(-c_4 \overbrace{n h^d}^{\gtrsim[\log(n)]^2\textrm{, \ref{a:h}(ii)}})
\\&\le c_3 \sum_{n=1}^{\infty} \exp\{-c_4 [\log(n)]^2\} ,
\end{align*}
which is finite by comparison with 
\[ \sum_{n=1}^{\infty} \exp\{-2\log(n)\}=\sum_{n=1}^{\infty} \exp\{\log(n^{-2})\}=\sum_{n=1}^{\infty} n^{-2} = \pi^2/6 \]
since $c_4[\log(n)]^2>2\log(n)$ for large enough $n$ for any $c_4>0$. 
This result might hold under dependent sampling (for $X$) by replacing Bernstein's inequality for independent data with that in \citet{FanEtAl2012}, but such is beyond our scope. 
}{Using the same argument as \citet[proof of Thm.\ 3.1, p.\ 769]{Chaudhuri1991}, invoking Bernstein's inequality and the Borel--Cantelli Lemma, $N_n$ is almost surely of order $nh^d$.} 

}{%
Whereas $n$ is a deterministic sequence, $N_n$ is random, but $N_n\stackrel{a.s.}{\asymp}nh^d$ as shown in \citet{Chaudhuri1991}. %
}%
Another difference with the unconditional case is that the local sample's distribution, $F_{Y|X}(\cdot;C_h)$, changes with $n$ (through $h$). 
The uniformity of the remainder term in Theorem \ref{thm:IDEAL-single} relies on the properties of the PDF in Assumption \ref{a:hut-pf}.  In the conditional case, we 
\Supplemental{must }{}show that these properties hold uniformly over the PDFs $f_{Y|X}(\cdot;C_h)$ as $h\to0$, for which we rely on \ref{a:f}, \ref{a:Q}, \ref{a:GK1}, and \ref{a:GK2}. 

\Supplemental{%
The analog of \ref{a:hut-pf}(i) is $f_{Y|X}\bigl(F_{Y|X}^{-1}(p;C_h);C_h\bigr)$ being uniformly bounded away from zero as $h\to0$%
\Supplemental{.  
From \ref{a:f} and \ref{a:GK1}, with $\tilde{x},\tilde{\tilde{x}}\in C_h$ determined by the MVT, 
\begin{align*}
&F_{Y|X}(t;C_h) 
= P\left( Y\le t\mid X\in C_h\right) 
 = \frac{\int_{C_h}\int_{-\infty}^{t}f_{Y|X}(y; x)\,dy\,dx}{P(X\in C_h)} 
,\\
&f_{Y|X}(t;C_h)
= \frac{d}{dt} F_{Y|X}(t;C_h)
= \frac{\int_{C_h}f_{Y|X}(t; x)\,dx}{P(X\in C_h)}
= \frac{\overbrace{(2h)^d f_{Y|X}(t; \tilde{x})}^{\textrm{MVT}}}{\int_{C_h}f_X(x)\,dx}
= \frac{(2h)^d f_{Y|X}(t; \tilde{x})}{\underbrace{(2h)^d f_X(\tilde{\tilde{x}})}_{\textrm{MVT}}} 
\\&\quad= \frac{f_{Y|X}(t; \tilde{x})}{f_X(\tilde{\tilde{x}})} 
\label{eqn:cond-unif-int2}\refstepcounter{equation}\tag{\theequation}
,\\
\label{eqn:cond-unif-int1}\refstepcounter{equation}\tag{\theequation}
&F_{Y|X}^{-1}(p;C_h) - F_{Y|X}^{-1}(p;\tilde{x})
= \overbrace{Q_{Y|X}(p;C_h)-Q_{Y|X}(p;0)}^{=B_h=O(h^b)\textrm{ (Lemma \ref{lem:bias})}} 
  +\overbrace{\left[Q_{Y|X}(p;0)-Q_{Y|X}(p;\tilde{x})\right]}^{=o(1)\textrm{ uniformly by \ref{a:Q}}}
,\\
&f_{Y|X}\left(F_{Y|X}^{-1}(p;C_h);C_h\right) 
\\&\quad= 
 \frac{f_{Y|X}\left(F_{Y|X}^{-1}(p;\tilde{x})+[F_{Y|X}^{-1}(p;C_h)-F_{Y|X}^{-1}(p;\tilde{x})] ; \tilde{x}\right)}{f_X(\tilde{\tilde{x}})}
\\&\quad= 
\frac{f_{Y|X}\left(F_{Y|X}^{-1}(p;\tilde{x}) ; \tilde{x}\right)%
         +\overbrace{\overbrace{[F_{Y|X}^{-1}(p;C_h)-F_{Y|X}^{-1}(p;\tilde{x})]}^{\textrm{uniformly $o(1)$ by \eqref{eqn:cond-unif-int1}}}
                     \overbrace{f^{(1,0)}_{Y|X}\left(\tilde{y};\tilde{x}\right)}^{\textrm{uniformly $O(1)$ by \ref{a:GK2}}}}^{\textrm{MVT}}}{f_X(\tilde{\tilde{x}})} 
\\&\quad= \frac{f_{Y|X}\left(F_{Y|X}^{-1}(p;\tilde{x}) ; \tilde{x}\right)%
         +o(1)}{f_X(\tilde{\tilde{x}})} 
, 
\end{align*}
where the first-order term in the numerator is uniformly bounded away from zero over any $\tilde{x}$ in a neighborhood of zero (and thus over $C_h$) by \ref{a:GK1}, the $o(1)$ term is uniform over $\tilde{x}$ in $C_h$ as shown (by \eqref{eqn:cond-unif-int1} and \ref{a:GK2}), and the denominator is uniformly bounded and bounded away from zero by \ref{a:f} over any $\tilde{\tilde{x}}$ in a neighborhood of zero (and thus over $C_h$). 

}{%
, which can be shown using Assumptions \ref{a:f}, \ref{a:Q}, \ref{a:GK1}, and \ref{a:GK2}, Lemma \ref{lem:bias}, $h\to0$, and the mean value theorem. 
}
The analog of \ref{a:hut-pf}(ii) is $f_{Y|X}''(\cdot;C_h)$ being uniformly continuous (in $y$) in a neighborhood of $F_{Y|X}^{-1}(p;C_h)$%
\Supplemental{.  
By Lemma \ref{lem:bias}, $F_{Y|X}^{-1}(p;C_h)-Q_{Y|X}(p;0)=O(h^b)\to0$, so it is sufficient to establish the property over a neighborhood of $Q_{Y|X}(p;0)$.  Continuing from \eqref{eqn:cond-unif-int2},
\begin{align*}
f_{Y|X}''(t;C_h) 
&=
\frac{f_{Y|X}^{(2,0)}(t; \tilde{x})}{f_X(\tilde{\tilde{x}})} 
.
\end{align*}
Again, the denominator is uniformly bounded and bounded away from zero by \ref{a:f} over any $\tilde{\tilde{x}}$ in a neighborhood of zero.  The numerator is uniformly bounded and continuous by \ref{a:GK2} for all $t$ in a neighborhood of $Q_{Y|X}(p;0)$ (as desired) and $\tilde{x}$ in a neighborhood of the origin (which it is, since $\tilde{x}\in C_h\to0$). 

In sum: when applying Theorem \ref{thm:IDEAL-single} to inference on $Q_{Y|X}(p;C_h)$ in the conditional setting, the remainder term is uniform due to $h\to0$ and sufficient smoothness of the conditional distribution near the point of interest. 
}{%
, which can be shown using \ref{a:f}, \ref{a:GK2}, Lemma \ref{lem:bias}, $h\to0$, and the mean value theorem. 
}
}{}

\Supplemental{\subsubsection*{One-sided case}}{}

In the lower one-sided case, let $\hat Q_{Y|C_h}^L(u_h)$ be the \citet{Hutson1999} upper endpoint, with notation analogous to Section \ref{sec:cdf-err}, with $u_h=u^h(\alpha)$. 
\Supplemental{This is the linearly interpolated $u_h(N_n+1)$th order statistic from the $N_n$ values of $Y_i$ with $X_i\in C_h$. 

}{}%
The CP of the lower one-sided CI is
\begin{align}
P\left(Q_{Y|X}(p;0) < \hat Q_{Y|C_h}^L(u_h) \right)  
\Supplemental{
\notag
  &= P\left(Q_{Y|X}(p;C_h) < \hat Q_{Y|C_h}^L(u_h) \right) \\
  &\quad\notag
    +\left[ P\left(Q_{Y|X}(p;0) < \hat Q_{Y|C_h}^L(u_h) \right)
           -P\left(Q_{Y|X}(p;C_h) < \hat Q_{Y|C_h}^L(u_h) \right)\right] \\*
}{}
\label{eqn:cpe-lower}
  &= 1-\alpha +\textrm{CPE}_\textrm{U} +\textrm{CPE}_\textrm{Bias},  
\end{align}
where $\textrm{CPE}_\textrm{U}$ is CPE due to the unconditional method and $\textrm{CPE}_\textrm{Bias}$ comes from the bias:
\begin{align*}
\textrm{CPE}_\textrm{U}
 &\equiv P\left(Q_{Y|X}(p;C_h) < \hat Q_{Y|C_h}^L(u_h) \right) 
         -(1-\alpha)
  = O\left(N_n^{-1}\right)  , \\
\textrm{CPE}_\textrm{Bias}
  &\equiv P\left(Q_{Y|X}(p;0) < \hat Q_{Y|C_h}^L(u_h) \right)
           -P\left(Q_{Y|X}(p;C_h) < \hat Q_{Y|C_h}^L(u_h) \right) . 
\end{align*}
\Supplemental{%

As in the main text, define $B_h\equiv Q_{Y|X}(p;C_h)-Q_{Y|X}(p;0)$.  From Lemma \ref{lem:bias}, $B_h=O(h^b)$ with $b\equiv\min\{s_Q,s_X+1,2\}$. 

For $\textrm{CPE}_\textrm{Bias}$, we apply Lemma \ref{lem:den}\Supplemental{ as in the steps leading to \eqref{eqn:W-EL-CDF-K}}{}.  
\Supplemental{%
We assume that Condition $\star(2\log(n))$ holds.  By Lemma \ref{lem:den-i}(\ref{lem:den-i-star-prob},\ref{lem:den-i-star-prob-fixed}), this introduces only $O(N_n^{-2})$ error, as shown in \eqref{eqn:star-error}, which is negligible in all following results. The points of evaluation also satisfy $\star(2\log(n))$ since they are all within $O(N_n^{-1/2})$ of $\mathbb{X}_0$. 
Here, $u=u_h$, $L^L=\hat Q_{Y|C_h}^L(u_h)$, $\mathbb{X}_0=Q_{Y|X}(u_h;C_h)$, and $\mathcal{V}_\psi=u_h(1-u_h)/[f_{Y|X}(Q_{Y|X}(u_h;C_h);C_h)]^2$, with $\xi_p \equiv Q_{Y|X}(p;0)$ as in Lemma \ref{lem:bias}, so
}{%
Eventually, 
}
\begin{align*}
\textrm{CPE}_\textrm{Bias}
\Supplemental{
&=
P\left(Q_{Y|X}(p;0) < \hat Q_{Y|C_h}^L(u_h) \right)
           -P\left(Q_{Y|X}(p;C_h) < \hat Q_{Y|C_h}^L(u_h) \right)
\\&=
 1-P\left( L^L < \xi_p\right) 
-\left[ 1- P\left( L^L<Q_{Y|X}(p;C_h) \right) \right]
\\&=
 P\left( N_n^{1/2}\left[L^L-\mathbb{X}_0\right] 
        <N_n^{1/2}\left[Q_{Y|X}(p;C_h)-\mathbb{X}_0\right] \right)
\\&\quad
-P\left( N_n^{1/2}\left[L^L-\mathbb{X}_0\right] 
        <N_n^{1/2}\left[\xi_p-\mathbb{X}_0\right] \right)
\\&=
 P\Bigl( \mathbb{W}_{\epsilon,\Lambda} 
        +\overbrace{O\left(N_n^{-3/2}[\log(N_n)]^3\right)}^{\textrm{\ref{lem:den}(\ref{lem:den-bd2})}}
        <N_n^{1/2}\left[Q_{Y|X}(p;C_h)-\mathbb{X}_0\right] \Bigr)
\\&\quad
-P\Bigl( \mathbb{W}_{\epsilon,\Lambda} 
        +\overbrace{O\left(N_n^{-3/2}[\log(N_n)]^3\right)}^{\textrm{\ref{lem:den}(\ref{lem:den-bd2})}}
        <N_n^{1/2}\left[\xi_p-\mathbb{X}_0\right] \Bigr)
\\&=
 P\left( \mathbb{W}_{\epsilon,\Lambda} 
        <N_n^{1/2}\left[Q_{Y|X}(p;C_h)-\mathbb{X}_0\right] \right)
-P\left( \mathbb{W}_{\epsilon,\Lambda} 
        <N_n^{1/2}\left[\xi_p-\mathbb{X}_0\right] \right)
\\&\quad
+\overbrace{[f_{\mathbb{W}_{\epsilon,\Lambda}}(\tilde{K_1})
  +f_{\mathbb{W}_{\epsilon,\Lambda}}(\tilde{K_2})]}^{O(1)\textrm{ by \ref{lem:den}(\ref{lem:den-iii})}}
  O\left( N_n^{-3/2}[\log(N_n)]^3 \right) 
\\&=
\int_{N_n^{1/2}\left[\xi_p-\mathbb{X}_0\right]}^{N_n^{1/2}\left[Q_{Y|X}(p;C_h)-\mathbb{X}_0\right]}
\int_0^1 
    f_{\mathbb{W}_{\epsilon,\Lambda}}(w\mid\lambda)
    f_{\Lambda}(\lambda)
    \,d\lambda \,dw
 +O\left( N_n^{-3/2}[\log(N_n)]^3 \right) 
\\&=
\int_{N_n^{1/2}\left[\xi_p-\mathbb{X}_0\right]}^{N_n^{1/2}\left[Q_{Y|X}(p;C_h)-\mathbb{X}_0\right]}
    \overbrace{\phi_{\mathcal{V}_\psi}(w)
    \left[1+O(N_n^{-1/2}[\log(N_n)]^3)\right]}^{\textrm{Lemma \ref{lem:den}(\ref{lem:den-iii})}}
    \,dw
 +O\left( N_n^{-3/2}[\log(N_n)]^3 \right) 
\\&=
\left\{N_n^{1/2}\left[Q_{Y|X}(p;C_h)-\mathbb{X}_0\right]
- N_n^{1/2}\left[\xi_p-\mathbb{X}_0\right] \right\}
    \phi_{\mathcal{V}_\psi}(\tilde{w})
    \left[1+O(N_n^{-1/2}[\log(N_n)]^3)\right]
\\&\quad
 +O\left( N_n^{-3/2}[\log(N_n)]^3 \right) 
\\&=
N_n^{1/2}\overbrace{\left[Q_{Y|X}(p;C_h)-\xi_p\right] }^{B_h}
\\&\quad\times
\overbrace{%
\Bigl[ \overbrace{\phi_{\mathcal{V}_\psi}\left( N_n^{1/2}\left[Q_{Y|X}(p;C_h)-\mathbb{X}_0\right] \right)}^{O(1)}
      +\overbrace{\phi_{\mathcal{V}_\psi}'(\tilde{\tilde{w}})}^{O(1)}
       \overbrace{\left\{\tilde{w} - N_n^{1/2}\left[Q_{Y|X}(p;C_h)-\mathbb{X}_0\right] \right\}}^{\le N_n^{1/2}\left[Q_{Y|X}(p;C_h)-\xi_p\right]=O\left(N_n^{1/2}B_h\right)}
\Bigr]
}^{\textrm{MVT expansion of }\phi_{\mathcal{V}_\psi}(\tilde{w})\textrm{ around }w=N_n^{1/2}\left[Q_{Y|X}(p;C_h)-\mathbb{X}_0\right]}
\\&\quad\times
    \left[1+O(N_n^{-1/2}[\log(N_n)]^3)\right]
\\&\quad
 +O\left( N_n^{-3/2}[\log(N_n)]^3 \right) 
\\&=
N_n^{1/2} B_h
\phi_{\mathcal{V}_\psi}%
\Bigl( \overbrace{N_n^{1/2}(p-u_h) }^{\textrm{use Lemma \ref{lem:u1u2-approx}}}
       \overbrace{Q^{(1,0)}_{Y|X}(u_h;C_h)}^{O(1)\textrm{ by \ref{a:GK1}}}
      +\overbrace{N_n^{1/2}(1/2)(p-u_h)^2}^{O(N_n^{-1/2})} \overbrace{Q^{(2,0)}_{Y|X}(\tilde{p};C_h)}^{O(1)\textrm{ by \ref{a:GK1}, \ref{a:GK2}}} \Bigr)
\\&\quad
 +O\Bigl( N_n^{-3/2}[\log(N_n)]^3 + N_n B_h^2 
         +B_h [\log(N_n)]^3 
   \Bigr) 
\\&=
N_n^{1/2} B_h
\overbrace{%
\phi_{\mathcal{V}_\psi}%
\Bigl( -z_{1-\alpha}
        \overbrace{\sqrt{u_h(1-u_h)} / f_{Y|X}\left(Q_{Y|X}(u_h;C_h);C_h\right)}^{=\sqrt{\mathcal{V}_\psi}}
       +O(N_n^{-1/2}) \Bigr)%
}^{\textrm{use }\phi_V(x)=V^{-1/2}\phi(x/\sqrt{V})}
\\&\quad
 +O\left( N_n^{-3/2}[\log(N_n)]^3 + N_n B_h^2 + B_h [\log(N_n)]^3 \right) 
\\&=
N_n^{1/2} B_h
\left[ \mathcal{V}_\psi^{-1/2} \phi( z_{1-\alpha}) + O(N_n^{-1/2}) \right] 
\\&\quad
 +O\left( N_n^{-3/2}[\log(N_n)]^3 + N_n B_h^2 + B_h [\log(N_n)]^3 \right) 
\\
}{}
&=
\label{eqn:cpe-bias-lower}\refstepcounter{equation}\tag{\theequation}
N_n^{1/2} B_h
\frac{f_{Y|X}\left(Q_{Y|X}(u_h;C_h);C_h\right)}{\sqrt{u_h(1-u_h)}}
\phi( z_{1-\alpha} )
\\&\quad
 +O\left( N_n^{-3/2}[\log(N_n)]^3 + N_n B_h^2 + B_h [\log(N_n)]^3 \right) 
\\&= O\left( (nh^d)^{1/2} h^b \right)
   = \overbrace{O\left( n^{1/2} h^{b+d/2} \right)}^{o(1)\textrm{ by \ref{a:h}(i')}}
,
\end{align*}
using $N_n\stackrel{a.s.}{\asymp}nh^d$ in the last line. 
\Supplemental{
To show the remainder terms are smaller-order, 
first, $B_h[\log(N_n)]^3=o(B_h N_n^{1/2})$ is immediate.  
Second, since $N_n^{1/2}B_h=o(1)$, then 
\begin{equation*}
N_n B_h^2 = (N_n^{1/2}B_h)(N_n^{1/2}B_h) = N_n^{1/2}B_h o(1) = o(N_n^{1/2}B_h) . 
\end{equation*}
Third, the CPE-optimal bandwidth below will set (up to log terms) $N_n^{-1}\asymp N_n^{1/2}B_h$, so $N_n^{-3/2}[\log(N_n)]^3\asymp B_h[\log(N_n)]^3 = o(B_h N_n^{1/2})$. 
The expression in \eqref{eqn:cpe-bias-lower} holds for $B_h>0$ (leading to over-coverage) or $B_h<0$ (under-coverage). 
}{}

The result in \eqref{eqn:cpe-bias-lower} can also be derived using Theorem \ref{thm:cdferror} instead of Lemma \ref{lem:den}; see the Supplemental Appendix for details.  

The dominant terms of $\textrm{CPE}_\textrm{U}$ and $\textrm{CPE}_\textrm{Bias}$ are thus respectively $O(N_n^{-1})$ and $O(N_n^{1/2}h^b)$.  
\Supplemental{
These are both sharp except in the special case of $u_h(N_n+1)$ being an integer or of $f_X(0)F_{Y|X}^{(0,2)}(\xi_p;0)+2f_X'(0)F_{Y|X}^{(0,1)}(\xi_p;0)=0$ when $k_Q\ge2$ and $k_X\ge1$ (or similar conditions otherwise); the following assumes we are not in such a special case. 
}{}
The term $\textrm{CPE}_\textrm{U}$ in Theorem \ref{thm:IDEAL-single} is always positive.  
The optimal $h$ sets equal the orders of magnitude $N_n^{-1}\asymp N_n^{1/2}h^b$, leading to
\begin{align*}
(nh^d)^{-1} &\asymp (nh^d)^{1/2} h^b \implies 
(nh^d)^{-3/2} \asymp h^b \implies
n^{-3/2} \asymp h^{b+3d/2} \implies
n^{-3/(2b+3d)} \asymp h .
\end{align*}
The resulting CPE is
\begin{align*}
O\left( N_n^{-1} \right)
&= O\left( (nh^d)^{-1} \right)
 = O\left( n^{-1} n^{3d/(2b+3d)} \right)
 = O\left( n^{(3d-2b-3d)/(2b+3d)} \right)
 = O\left( n^{-2b/(2b+3d)} \right) .
\end{align*}

If the calibrated method from Theorem \ref{thm:IDEAL-single} is used, then $\textrm{CPE}_\textrm{U}=O\left(N_n^{-3/2}[\log(N_n)]^3\right)$, while $\textrm{CPE}_\textrm{Bias}$ is unchanged.  Ignoring the $\log$ for simplicity, the optimal bandwidth solves
\begin{align*}
(nh^d)^{-3/2} &\asymp (nh^d)^{1/2} h^b \implies 
(nh^d)^{-2} \asymp h^b \implies
n^{-2} \asymp h^{b+2d} \implies
n^{-2/(b+2d)} \asymp h .
\end{align*}
Since $\log(N_n)=\log(nh^d)=\log(n^{1-2d/(b+2d)})\propto\log(n)$, this bandwidth leads to CPE
\begin{align*}
  O\left(N_n^{-3/2}[\log(N_n)]^3\right)
&=O\left(n^{-3/2} h^{-3d/2} [\log(n)]^3\right) 
\Supplemental{
=O\left(n^{-3/2} n^{3d/(b+2d)} [\log(n)]^3\right) 
\\&
 =O\left(n^{(-3b-6d+6d)/(2b+4d)}[\log(n)]^3\right) 
}{}
 =O\left(n^{-3b/(2b+4d)}[\log(n)]^3\right) 
. 
\end{align*}
\Supplemental{
We must also double-check that none of the remainder terms in \eqref{eqn:cpe-bias-lower} is strictly larger than this calibrated CPE (in order of magnitude).  Two of the remainder terms were shown to be smaller than $\textrm{CPE}_\textrm{Bias}=O(N_n^{1/2}B_h)$, so the same arguments hold here. 
For the third remainder term, $O\left(N_n^{-3/2}[\log(N_n)]^3\right)$ now equals $\textrm{CPE}_\textrm{U}$, but it will not affect the optimal rates since it is not strictly larger than $\textrm{CPE}_\textrm{U}$. 
}{}

\Supplemental{
In the upper one-sided case, let $u_\ell\equiv u^\ell(\alpha)$ instead of $u_h$, with $\hat Q_{Y|C_h}^L(u_\ell)$ the lower endpoint.  Similar to before, $L^L=\hat Q_{Y|C_h}^L(u_\ell)$, $\mathbb{X}_0=Q_{Y|X}(u_\ell;C_h)$, and $\mathcal{V}_\psi=u_\ell(1-u_\ell)/[f_{Y|X}(Q_{Y|X}(u_\ell;C_h);C_h)]^2$. 
Coverage probability is 
\begin{align*}
&P\left( \hat{Q}^L_{Y|C_h}(u_\ell) < Q_{Y|X}(p;0) \right) 
\\&=
\overbrace{P\left( \hat{Q}^L_{Y|C_h}(u_\ell) < Q_{Y|X}(p;C_h) \right)}^{=1-\alpha+\textrm{CPE}_\textrm{U}}
\\&\quad
+ \overbrace{%
  P\left( \hat{Q}^L_{Y|C_h}(u_\ell) < Q_{Y|X}(p;0) \right)
 -P\left( \hat{Q}^L_{Y|C_h}(u_\ell) < Q_{Y|X}(p;C_h) \right)%
  }^{\textrm{CPE}_\textrm{Bias}} , 
\end{align*}
so
\begin{align*}
\textrm{CPE}_\textrm{Bias}
&= 
  P\left( \hat{Q}^L_{Y|C_h}(u_\ell) < Q_{Y|X}(p;0) \right)
 -P\left( \hat{Q}^L_{Y|C_h}(u_\ell) < Q_{Y|X}(p;C_h) \right) 
\\&=
  P\left( N_n^{1/2} [L^L-\mathbb{X}_0] < N_n^{1/2}[Q_{Y|X}(p;0)-\mathbb{X}_0] \right)
\\&\quad
 -P\left( N_n^{1/2} [L^L-\mathbb{X}_0] < N_n^{1/2}[Q_{Y|X}(p;C_h)-\mathbb{X}_0] \right) 
\\&=
 P\left( \mathbb{W}_{\epsilon,\Lambda} < N_n^{1/2}[Q_{Y|X}(p;0)-\mathbb{X}_0] \right)
-P\left( \mathbb{W}_{\epsilon,\Lambda} < N_n^{1/2}[Q_{Y|X}(p;C_h)-\mathbb{X}_0] \right)
\\&\quad
+\overbrace{O\left( N_n^{-3/2}[\log(N_n)]^3 \right)}^{\textrm{Lemma \ref{lem:den}}}
\\&=
\int_{N_n^{1/2}\left[Q_{Y|X}(p;C_h)-\mathbb{X}_0\right]}^{N_n^{1/2}\left[\xi_p-\mathbb{X}_0\right]}
    \overbrace{\phi_{\mathcal{V}_\psi}(w)
    \left[1+O(N_n^{-1/2}[\log(N_n)]^3)\right]}^{\textrm{Lemma \ref{lem:den}(\ref{lem:den-iii})}}
    \,dw
 +O\left( N_n^{-3/2}[\log(N_n)]^3 \right) 
\\&=
N_n^{1/2}\overbrace{\left[\xi_p-Q_{Y|X}(p;C_h)\right]}^{=-B_h}
\left[ \phi_{\mathcal{V}_\psi}\left( N_n^{1/2}\left[Q_{Y|X}(p;C_h)-Q_{Y|X}(u_\ell;C_h)\right] \right)
      +O\left(N_n^{1/2}B_h\right) \right]
\\&\quad\times
    \left[1+O(N_n^{-1/2}[\log(N_n)]^3)\right]
\\&\quad
 +O\left( N_n^{-3/2}[\log(N_n)]^3 \right) 
\\&=
- N_n^{1/2} B_h \mathcal{V}_\psi^{-1/2} \phi(z_{1-\alpha})
+O\left( N_n^{-3/2}[\log(N_n)]^3 +N_n B_h^2 +B_h[\log(N_n)]^3 \right)
\\
\label{eqn:cpe-bias-upper}\refstepcounter{equation}\tag{\theequation}
\begin{split}
&=
-N_n^{1/2} B_h
\frac{f_{Y|X}\left(Q_{Y|X}(u_\ell;C_h);C_h\right)}{\sqrt{u_\ell(1-u_\ell)}}
\phi( z_{1-\alpha} )
\\&\quad
 +O\left( N_n^{-3/2}[\log(N_n)]^3 + N_n B_h^2 + B_h[\log(N_n)]^3 \right) 
\end{split}
\\&=
O\left( N_n^{1/2} B_h \right) = O\left( n^{1/2} h^{b+d/2} \right) 
. 
\end{align*}
Alternatively, again, the result in \eqref{eqn:cpe-bias-upper} can be derived using Theorem \ref{thm:cdferror} instead of Lemma \ref{lem:den}; see \eqref{eqn:cpe-bias-ell} below.  
Opposite before, $B_h>0$ now contributes under-coverage and $B_h<0$ over-coverage, but the order of $\textrm{CPE}_\textrm{Bias}$ is the same.  Since $\textrm{CPE}_\textrm{U}$ is of the same order and sign as for the upper endpoint, the resulting rates of CPE-optimal bandwidth and overall CPE are also the same. 
}{
The upper one-sided case follows similarly, except that $B_h>0$ contributes under-coverage and $B_h<0$ over-coverage, but the rates are all the same. 
}
}{%
Using Lemmas \ref{lem:bias} and \ref{lem:den}, or alternatively Theorem \ref{thm:cdferror}, one can show $\textrm{CPE}_\textrm{Bias}=O(N_n^{1/2}h^b)$.  Then, one can solve for the $h$ that equates the orders of $\textrm{CPE}_\textrm{U}$ and $\textrm{CPE}_\textrm{Bias}$, i.e.,\ so that $N_n^{-1}\asymp N_n^{1/2}h^b$, using $N_n\asymp nh^d$. 
}

\Supplemental{\subsubsection*{Two-sided case}}{}

With two-sided inference, the lower and upper endpoints have opposite bias effects.  For the median, the dominant terms of these effects cancel completely.  For other quantiles, there is a partial, order-reducing cancellation.  
\Supplemental{

Below, we apply Theorem \ref{thm:cdferror}. 
Let $F_{\beta,\ell}(\cdot)$ denote the $\beta\left(u_\ell(N_n+1),(1-u_\ell)(N_n+1)\right)$ distribution's CDF, which is the distribution of $F_{Y|X}\left(\tilde Q_{Y|C_h}^I(u_\ell);C_h\right)$, and let $f_{\beta,\ell}(\cdot)$ and $f_{\beta,\ell}'(\cdot)$ denote the corresponding PDF and PDF derivative. 
Let $\mathbb{X}_{0\ell}\equiv Q_{Y|X}(u_\ell;C_h)$, 
$L^L_\ell\equiv \hat Q_{Y|C_h}^L(u_\ell)$, 
$L^I_\ell\equiv \tilde{Q}_{Y|C_h}^I(u_\ell)$, $\mathcal{V}_\ell=u_\ell(1-u_\ell)/[f_{Y|X}(Q_{Y|X}(u_\ell;C_h);C_h)]^2\asymp1$ (since \ref{a:GK1} uniformly bounds the denominator away from zero and $u_\ell=p-O(N_n^{-1/2})$), 
and 
$\epsilon_\ell\equiv u_\ell(N_n+1)-\lfloor u_\ell(N_n+1)\rfloor$. 
Letting 
\begin{align*}
K_{\ell1} &\equiv N_n^{1/2}\left[ Q_{Y|X}(p;C_h) - Q_{Y|X}(u_\ell;C_h)\right] \asymp 1, \\
K_{\ell2} &\equiv N_n^{1/2}\left[ \xi_p - Q_{Y|X}(u_\ell;C_h)\right]
           = K_{\ell1} - N_n^{1/2}B_h , 
\end{align*}
then
\begin{align*}
&
P\left( \hat Q_{Y|C_h}^L(u_\ell) < Q_{Y|X}(p;0) \right)
\\&=
\overbrace{P\left( L^L_\ell < \mathbb{X}_{0\ell} + N_n^{-1/2}K_{\ell2} \right)}^{\textrm{apply Theorem \ref{thm:cdferror}}}
\\&=
 P\left( \tilde{Q}_{Y|C_h}^I(u_\ell) < Q_{Y|X}(p;0) \right)
+N_n^{-1}
 \frac{K_{\ell2}\exp\{-K_{\ell2}^2/(2\mathcal{V}_\ell)\}}{\sqrt{2\pi\mathcal{V}_\ell^3}}
 \frac{\epsilon_\ell(1-\epsilon_\ell)}{\underbrace{\left[f_{Y|X}\left(Q_{Y|X}(u_\ell;C_h);C_h\right)\right]^2}_{\textrm{\ref{a:GK1}: uniformly bdd $>0$}}}
\\&\quad
+O\left(N_n^{-3/2}[\log(N_n)]^3\right)
\\&=
 P\left( \tilde{Q}_{Y|C_h}^I(u_\ell) < Q_{Y|X}(p;0) \right)
+N_n^{-1}
 \frac{K_{\ell1}\exp\{-K_{\ell1}^2/(2\mathcal{V}_\ell)\}}{\sqrt{2\pi\mathcal{V}_\ell^3}}
 \frac{\epsilon_\ell(1-\epsilon_\ell)}{\left[f_{Y|X}\left(Q_{Y|X}(u_\ell;C_h);C_h\right)\right]^2}
\\&\quad
+\overbrace{O\left( N_n^{-1} N_n^{1/2}B_h \right)}^{o(N_n^{-1})\textrm{ by \ref{a:h}(i')}}
+\overbrace{O\left(N_n^{-3/2}[\log(N_n)]^3\right)}^{o(N_n^{-1})}
. 
\intertext{For comparison, the endpoint CDF evaluated at the biased $Q_{Y|X}(p;C_h)$ is}
&P\left( \hat Q_{Y|C_h}^L(u_\ell) < Q_{Y|X}(p;C_h) \right)
\\&=
\overbrace{P\left( L^L_\ell < \mathbb{X}_{0\ell} + N_n^{-1/2}K_{\ell1} \right)}^{\textrm{apply Theorem \ref{thm:cdferror}}}
\\&=
 P\left( \tilde{Q}_{Y|C_h}^I(u_\ell) < Q_{Y|X}(p;C_h) \right)
+N_n^{-1}
 \frac{K_{\ell1}\exp\{-K_{\ell1}^2/(2\mathcal{V}_\ell)\}}{\sqrt{2\pi\mathcal{V}_\ell^3}}
 \frac{\epsilon_\ell(1-\epsilon_\ell)}{\left[f_{Y|X}\left(Q_{Y|X}(u_\ell;C_h);C_h\right)\right]^2}
\\&\quad
+O\left(N_n^{-3/2}[\log(N_n)]^3\right)
. 
\intertext{For notational compactness, let}
%
%
&D \equiv
-B_h f_{Y|X}\left(Q_{Y|X}(p;C_h);C_h\right)
             +\frac{1}{2}(-B_h)^2 f^{(1,0)}_{Y|X}\left(Q_{Y|X}(p;C_h);C_h\right)
             +\frac{1}{6}(-B_h)^3 f^{(2,0)}_{Y|X}(\tilde{Q};C_h) 
\\&\quad= 
O(B_h)
, 
\intertext{so that}&
 P\left( \hat Q_{Y|C_h}^L(u_\ell) < Q_{Y|X}(p;0) \right)
-P\left( \hat Q_{Y|C_h}^L(u_\ell) < Q_{Y|X}(p;C_h) \right)
\\&=
 P\left( \tilde{Q}_{Y|C_h}^I(u_\ell) < Q_{Y|X}(p;0) \right)
-P\left( \tilde{Q}_{Y|C_h}^I(u_\ell) < Q_{Y|X}(p;C_h) \right)
\\&\quad
+O\left( N_n^{-1/2} B_h + N_n^{-3/2}[\log(N_n)]^3 \right)
\\&=
 P\left( F_{Y|X}\left(\tilde{Q}_{Y|C_h}^I(u_\ell);C_h\right) 
        <F_{Y|X}\left(Q_{Y|X}(p;0);C_h\right) \right)
\\&\quad
-P\left( F_{Y|X}\left(\tilde{Q}_{Y|C_h}^I(u_\ell);C_h\right) 
        <F_{Y|X}\left(Q_{Y|X}(p;C_h);C_h\right) \right)
+O\left( N_n^{-1/2} B_h + N_n^{-3/2}[\log(N_n)]^3 \right)
\\&=
 F_{\beta,\ell}\left( F_{Y|X}\left(Q_{Y|X}(p;C_h)-B_h;C_h\right) \right)
-F_{\beta,\ell}(p)
+O\left( N_n^{-1/2} B_h + N_n^{-3/2}[\log(N_n)]^3 \right)
\\&=
 F_{\beta,\ell}\left( 
             p +D
               \right)
-F_{\beta,\ell}(p)
+O\left( N_n^{-1/2} B_h + N_n^{-3/2}[\log(N_n)]^3 \right)
\\&=
 F_{\beta,\ell}(p) 
+D f_{\beta,\ell}(p) 
+(1/2) D^2 f'_{\beta,\ell}(p)
+(1/6) D^3 \overbrace{f''_{\beta,\ell}(\tilde{p})}^{O(N_n^{3/2})}
-F_{\beta,\ell}(p)
\\&\quad
+O\left( N_n^{-1/2} B_h + N_n^{-3/2}[\log(N_n)]^3 \right)
\\&=
\left[ -B_h f_{Y|X}\left(Q_{Y|X}(p;C_h);C_h\right)
       +(1/2) B_h^2 f^{(1,0)}_{Y|X}\left(Q_{Y|X}(p;C_h);C_h\right)
 \right]
 \overbrace{f_{\beta,\ell}(p) }^{O(N_n^{1/2})}
 +O\left( B_h^3 N_n^{1/2} \right)
\\&\quad
+(1/2) B_h^2 \left[f_{Y|X}\left(Q_{Y|X}(p;C_h);C_h\right)\right]^2 
 \overbrace{f'_{\beta,\ell}(p)}^{O(N_n)}
 +O\left( N_n B_h^3 \right)
+\overbrace{O\left( B_h^3 N_n^{3/2} \right) }^{=o(B_h^2 N_n)\textrm{ by \ref{a:h}(i')}}
\\&\quad
+O\left( N_n^{-1/2} B_h + N_n^{-3/2}[\log(N_n)]^3 \right)
\\
\label{eqn:cpe-bias-ell}\refstepcounter{equation}\tag{\theequation}
\begin{split}
&=
-B_h f_{Y|X}\left(Q_{Y|X}(p;C_h);C_h\right) f_{\beta,\ell}(p)
+(1/2) B_h^2 f^{(1,0)}_{Y|X}\left(Q_{Y|X}(p;C_h);C_h\right) f_{\beta,\ell}(p)
\\&\quad
+\frac{1}{2} B_h^2 \left[f_{Y|X}\left(Q_{Y|X}(p;C_h);C_h\right)\right]^2 
 f'_{\beta,\ell}(p)
\\&\quad
+O\left( N_n^{-1/2} B_h + N_n^{-3/2}[\log(N_n)]^3 +B_h^3 N_n^{3/2} \right)
,
\end{split}
\end{align*}
where $\tilde{Q}$ is between $Q_{Y|X}(p;C_h)\to\xi_p$ and $Q_{Y|X}(p;C_h)-B_h\to\xi_p$, and $\tilde{p}$ is between $p$ and $p+D\to p$, both determined by the mean value theorem (i.e.,\ in the Lagrange form of the remainder from Taylor's theorem). 

For the upper endpoint, we essentially replace $\ell$ with $h$ in the subscripts and repeat the above calculations. 
Let $F_{\beta,h}(\cdot)$ denote the CDF of the $\beta\left(u_h(N_n+1),(1-u_h)(N_n+1)\right)$ distribution, which is the distribution of $F_{Y|X}\left(\tilde Q_{Y|C_h}^I(u_h);C_h\right)$, and let $f_{\beta,h}(\cdot)$ and $f_{\beta,h}'(\cdot)$ denote the corresponding PDF and PDF derivative. 
Let $\mathbb{X}_{0h}\equiv Q_{Y|X}(u_h;C_h)$, 
$L^L_h\equiv \hat Q_{Y|C_h}^L(u_h)$, 
$L^I_h\equiv \tilde{Q}_{Y|C_h}^I(u_h)$, $\mathcal{V}_h=u_h(1-u_h)/[f_{Y|X}(Q_{Y|X}(u_h;C_h);C_h)]^2\asymp1$ (since \ref{a:GK1} uniformly bounds the denominator away from zero and $u_h=p+O(N_n^{-1/2})$), 
and 
$\epsilon_h\equiv u_h(N_n+1)-\lfloor u_h(N_n+1)\rfloor$. 
Letting  
\begin{align*}
K_{h1} &\equiv N_n^{1/2}\left[ Q_{Y|X}(p;C_h) - Q_{Y|X}(u_h;C_h)\right] \asymp 1, \\
K_{h2} &\equiv N_n^{1/2}\left[ \xi_p - Q_{Y|X}(u_h;C_h)\right]
           = K_{h1} - N_n^{1/2}B_h , 
\end{align*}
then 
\begin{align*}
&
P\left( \hat Q_{Y|C_h}^L(u_h) < Q_{Y|X}(p;0) \right)
\\&=
\overbrace{P\left( L^L_h < \mathbb{X}_{0h} + N_n^{-1/2}K_{h2} \right)}^{\textrm{apply Theorem \ref{thm:cdferror}}}
\\&=
 P\left( \tilde{Q}_{Y|C_h}^I(u_h) < Q_{Y|X}(p;0) \right)
+N_n^{-1}
 \frac{K_{h2}\exp\{-K_{h2}^2/(2\mathcal{V}_h)\}}{\sqrt{2\pi\mathcal{V}_h^3}}
 \frac{\epsilon_h(1-\epsilon_h)}{\underbrace{\left[f_{Y|X}\left(Q_{Y|X}(u_h;C_h);C_h\right)\right]^2}_{\textrm{\ref{a:GK1}: uniformly bdd $>0$}}}
\\&\quad
+O\left(N_n^{-3/2}[\log(N_n)]^3\right)
\\&=
 P\left( \tilde{Q}_{Y|C_h}^I(u_h) < Q_{Y|X}(p;0) \right)
+N_n^{-1}
 \frac{K_{h1}\exp\{-K_{h1}^2/(2\mathcal{V}_h)\}}{\sqrt{2\pi\mathcal{V}_h^3}}
 \frac{\epsilon_h(1-\epsilon_h)}{\left[f_{Y|X}\left(Q_{Y|X}(u_h;C_h);C_h\right)\right]^2}
\\&\quad
+\overbrace{O\left( N_n^{-1} N_n^{1/2}B_h \right)}^{o(N_n^{-1})\textrm{ by \ref{a:h}(i')}}
+\overbrace{O\left(N_n^{-3/2}[\log(N_n)]^3\right)}^{o(N_n^{-1})}
,\\&
P\left( \hat Q_{Y|C_h}^L(u_h) < Q_{Y|X}(p;C_h) \right)
\\&=
\overbrace{P\left( L^L_h < \mathbb{X}_{0h} + N_n^{-1/2}K_{h1} \right)}^{\textrm{apply Theorem \ref{thm:cdferror}}}
\\&=
 P\left( \tilde{Q}_{Y|C_h}^I(u_h) < Q_{Y|X}(p;C_h) \right)
+N_n^{-1}
 \frac{K_{h1}\exp\{-K_{h1}^2/(2\mathcal{V}_h)\}}{\sqrt{2\pi\mathcal{V}_h^3}}
 \frac{\epsilon_h(1-\epsilon_h)}{\left[f_{Y|X}\left(Q_{Y|X}(u_h;C_h);C_h\right)\right]^2}
\\&\quad
+O\left(N_n^{-3/2}[\log(N_n)]^3\right)
,\\&
D \equiv -B_h f_{Y|X}\left(Q_{Y|X}(p;C_h);C_h\right)
             +\frac{1}{2}(-B_h)^2 f^{(1,0)}_{Y|X}\left(Q_{Y|X}(p;C_h);C_h\right)
             +\frac{1}{6}(-B_h)^3 f^{(2,0)}_{Y|X}(\tilde{Q};C_h)
\\&\quad
 = O(B_h)
,\\&
 P\left( \hat Q_{Y|C_h}^L(u_h) < Q_{Y|X}(p;0) \right)
-P\left( \hat Q_{Y|C_h}^L(u_h) < Q_{Y|X}(p;C_h) \right)
\\&=
 P\left( \tilde{Q}_{Y|C_h}^I(u_h) < Q_{Y|X}(p;0) \right)
-P\left( \tilde{Q}_{Y|C_h}^I(u_h) < Q_{Y|X}(p;C_h) \right)
\\&\quad
+O\left( N_n^{-1/2} B_h + N_n^{-3/2}[\log(N_n)]^3 \right)
\\&=
 P\left( F_{Y|X}\left(\tilde{Q}_{Y|C_h}^I(u_h);C_h\right) 
        <F_{Y|X}\left(Q_{Y|X}(p;0);C_h\right) \right)
\\&\quad
-P\left( F_{Y|X}\left(\tilde{Q}_{Y|C_h}^I(u_h);C_h\right) 
        <F_{Y|X}\left(Q_{Y|X}(p;C_h);C_h\right) \right)
+O\left( N_n^{-1/2} B_h + N_n^{-3/2}[\log(N_n)]^3 \right)
\\&=
 F_{\beta,h}\left( F_{Y|X}\left(Q_{Y|X}(p;C_h)-B_h;C_h\right) \right)
-F_{\beta,h}(p)
+O\left( N_n^{-1/2} B_h + N_n^{-3/2}[\log(N_n)]^3 \right)
\\&=
 F_{\beta,h}\left( 
             p +D
               \right)
\\&\quad
-F_{\beta,h}(p)
+O\left( N_n^{-1/2} B_h + N_n^{-3/2}[\log(N_n)]^3 \right)
\\&=
 F_{\beta,h}(p) 
+D f_{\beta,h}(p) 
+(1/2) D^2 f'_{\beta,h}(p)
+(1/6) D^3 \overbrace{f''_{\beta,h}(\tilde{p})}^{O(N_n^{3/2})}
-F_{\beta,h}(p)
\\&\quad
+O\left( N_n^{-1/2} B_h + N_n^{-3/2}[\log(N_n)]^3 \right)
\\&=
\left[ -B_h f_{Y|X}\left(Q_{Y|X}(p;C_h);C_h\right)
       +(1/2) B_h^2 f^{(1,0)}_{Y|X}\left(Q_{Y|X}(p;C_h);C_h\right)
 \right]
 \overbrace{f_{\beta,h}(p) }^{O(N_n^{1/2})}
 +O\left( B_h^3 N_n^{1/2} \right)
\\&\quad
+(1/2) B_h^2 \left[f_{Y|X}\left(Q_{Y|X}(p;C_h);C_h\right)\right]^2 
 \overbrace{f'_{\beta,h}(p)}^{O(N_n)}
 +O\left( N_n B_h^3 \right)
+\overbrace{O\left( B_h^3 N_n^{3/2} \right) }^{=o(B_h^2 N_n)\textrm{ by \ref{a:h}(i')}}
\\&\quad
+O\left( N_n^{-1/2} B_h + N_n^{-3/2}[\log(N_n)]^3 \right)
\\&=
-B_h f_{Y|X}\left(Q_{Y|X}(p;C_h);C_h\right) f_{\beta,h}(p)
+(1/2) B_h^2 f^{(1,0)}_{Y|X}\left(Q_{Y|X}(p;C_h);C_h\right) f_{\beta,h}(p)
\\&\quad
+\frac{1}{2} B_h^2 \left[f_{Y|X}\left(Q_{Y|X}(p;C_h);C_h\right)\right]^2 
 f'_{\beta,h}(p)
\\&\quad
+O\left( N_n^{-1/2} B_h + N_n^{-3/2}[\log(N_n)]^3 +B_h^3 N_n^{3/2} \right)
\label{eqn:cpe-bias-h}\refstepcounter{equation}\tag{\theequation}
,
\end{align*}
where $\tilde{Q}$ is between $Q_{Y|X}(p;C_h)\to\xi_p$ and $Q_{Y|X}(p;C_h)-B_h\to\xi_p$, and $\tilde{p}$ is between $p$ and $p+D\to p$, both determined by the mean value theorem (i.e.,\ in the Lagrange form of the remainder from Taylor's theorem). 

The overall two-sided CP is
\begin{align*}
P & \left( \hat Q_{Y|C_h}^L(u_\ell) < Q_{Y|X}(p;0) < \hat Q_{Y|C_h}^L(u_h) \right)  
\\&= 
1 - P\left( \hat Q_{Y|C_h}^L(u_\ell) > Q_{Y|X}(p;0) \right)
  - P\left( \hat Q_{Y|C_h}^L(u_h)    < Q_{Y|X}(p;0) \right)  
\\&= 
1 - P\left( \hat Q_{Y|C_h}^L(u_\ell) > Q_{Y|X}(p;C_h) \right) 
\\&\quad 
 +\left[ P\left( \hat Q_{Y|C_h}^L(u_\ell) > Q_{Y|X}(p;C_h) \right)
        -P\left( \hat Q_{Y|C_h}^L(u_\ell) > Q_{Y|X}(p;0) \right) \right] 
\\&\quad
 -P\left( \hat Q_{Y|C_h}^L(u_h) < Q_{Y|X}(p;C_h) \right) 
\\&\quad  
 +\left[ P\left( \hat Q_{Y|C_h}^L(u_h) < Q_{Y|X}(p;C_h) \right)
        -P\left( \hat Q_{Y|C_h}^L(u_h) < Q_{Y|X}(p;0) \right) \right] 
\\&= 
1-\alpha+\textrm{CPE}_\textrm{U}
\\&\quad
  +\overbrace{P\left( \hat{Q}^L_{Y|C_h}(u_\ell) < \xi_p \right)
  -P\left( \hat{Q}^L_{Y|C_h}(u_\ell) < Q_{Y|X}(p;C_h) \right) }^{\textrm{\eqref{eqn:cpe-bias-ell}}}
\\&\quad
  -\overbrace{\left[P\left( \hat{Q}^L_{Y|C_h}(u_h) < \xi_p \right)
  -P\left( \hat{Q}^L_{Y|C_h}(u_h) < Q_{Y|X}(p;C_h) \right)\right] }^{\textrm{\eqref{eqn:cpe-bias-h}}}
\\&= 
1-\alpha+\textrm{CPE}_\textrm{U}
\\&\quad
  +\Bigl\{
  -B_h f_{Y|X}\left(Q_{Y|X}(p;C_h);C_h\right) f_{\beta,\ell}(p)
+(1/2) B_h^2 f^{(1,0)}_{Y|X}\left(Q_{Y|X}(p;C_h);C_h\right) f_{\beta,\ell}(p)
\\&\qquad\quad
+\frac{1}{2} B_h^2 \left[f_{Y|X}\left(Q_{Y|X}(p;C_h);C_h\right)\right]^2 
 f'_{\beta,\ell}(p)
   \Bigr\}
\\&\quad
  -\Bigl\{
  -B_h f_{Y|X}\left(Q_{Y|X}(p;C_h);C_h\right) f_{\beta,h}(p)
+(1/2) B_h^2 f^{(1,0)}_{Y|X}\left(Q_{Y|X}(p;C_h);C_h\right) f_{\beta,h}(p)
\\&\qquad\quad
+\frac{1}{2} B_h^2 \left[f_{Y|X}\left(Q_{Y|X}(p;C_h);C_h\right)\right]^2 
 f'_{\beta,h}(p)
   \Bigr\}
\\&\quad
+O\left( N_n^{-1/2} B_h + N_n^{-3/2}[\log(N_n)]^3 +B_h^3 N_n^{3/2} \right)
\\&= 
1-\alpha+\textrm{CPE}_\textrm{U}
\\&\quad
  +\overbrace{B_h f_{Y|X}\left(Q_{Y|X}(p;C_h);C_h\right) 
   \overbrace{\left[ f_{\beta,h}(p) - f_{\beta,\ell}(p) \right]}^{=0\textrm{ if }p=1/2\textrm{; else }=O(1)\textrm{ by \eqref{eqn:beta-PDF-diff-approx}}}}^{=0\textrm{ if $p=1/2$; else }=O(B_h^2)=o(N_n^{-1/2}B_h)\textrm{ by \ref{a:h}(i')}}
\\&\quad
  +\overbrace{(1/2) B_h^2 f^{(1,0)}_{Y|X}\left(Q_{Y|X}(p;C_h);C_h\right)
   \overbrace{\left[ f_{\beta,\ell}(p) - f_{\beta,h}(p) \right]}^{O(1)\textrm{ by \eqref{eqn:beta-PDF-diff-approx}}}}^{=O(B_h^2)=o(N_n^{-1/2}B_h)\textrm{ by \ref{a:h}(i')}}
\\&\quad
+\frac{1}{2} B_h^2 \left[f_{Y|X}\left(Q_{Y|X}(p;C_h);C_h\right)\right]^2 
 \bigl[ \overbrace{f'_{\beta,\ell}(p)}^{<0} - \overbrace{f'_{\beta,h}(p)}^{>0} \bigr]
\\&\quad
+O\left( N_n^{-1/2} B_h + N_n^{-3/2}[\log(N_n)]^3 +B_h^3 N_n^{3/2} \right)
\\&= 
1-\alpha+\overbrace{\textrm{CPE}_\textrm{U}}^{O(N_n^{-1})}
\\&\quad
  +\overbrace{B_h f_{Y|X}\left(Q_{Y|X}(p;C_h);C_h\right) 
   \left[ f_{\beta,h}(p) - f_{\beta,\ell}(p) \right]}^{O(B_h)\textrm{, or zero if }p=1/2}
\\&\quad
-\overbrace{\frac{1}{2} B_h^2 \left[f_{Y|X}\left(Q_{Y|X}(p;C_h);C_h\right)\right]^2 
 \overbrace{\bigl[ f'_{\beta,h}(p) - f'_{\beta,\ell}(p) \bigr]}^{O(N_n)}}^{O(B_h^2 N_n)}
\\&\quad
+O\left( N_n^{-1/2} B_h + N_n^{-3/2}[\log(N_n)]^3 +B_h^3 N_n^{3/2} \right)
.
\label{eqn:cpe-bias-two}\refstepcounter{equation}\tag{\theequation}
\end{align*}

}{%
The calculations, which use Theorem \ref{thm:cdferror}, are extensive and thus left to the supplemental appendix.  
Ultimately, it can be shown that two-sided CP is $1-\alpha$ plus terms of $O(N_n^{-1})$, $O(B_h)$, and $O(B_h^2 N_n)$, in addition to smaller-order remainders. 
With the new CPE terms, one can again solve for the $h$ that sets the orders equal. 
}%
\Supplemental{%
If $p=1/2$, then the $B_h$ term is identically zero, and the CPE-optimal bandwidth equates the remaining two dominant terms\Supplemental{ in \eqref{eqn:cpe-bias-two}}{}:
\begin{align*}
         B_h^2 N_n   \asymp N_n^{-1}
\implies h^{2b} nh^d \asymp (nh^d)^{-1}
\Supplemental{
\implies h^{2(b+d)}  \asymp n^{-2}
}{}
\implies h \asymp n^{-1/(b+d)} 
.
\end{align*}
The resulting CPE is
\begin{align*}
O\left( N_n^{-1} \right)
= O\left( [n n^{-d/(b+d)}]^{-1}\right)
= O(n^{-b/(b+d)} )
. 
\end{align*}
Even if $p\ne1/2$, the fact that $B_h^2 N_n\asymp N_n^{-1}$ implies $B_h\asymp N_n^{-1}$, so the $B_h$ term does not affect the CPE-optimal bandwidth rate or resulting CPE. 

With the calibrated method, $\textrm{CPE}_\textrm{U}=O\bigl( N_n^{-3/2}[\log(N_n)]^3\bigr)$. 
With $p=1/2$, the CPE-optimal bandwidth equates (ignoring the $\log$)
\begin{align*}
         B_h^2 N_n   \asymp N_n^{-3/2}
&\implies h^{2b} nh^d \asymp (nh^d)^{-3/2}
\Supplemental{
\implies h^{2b+5d/2}  \asymp n^{-5/2}
\implies h^{4b+5d}  \asymp n^{-5}
\\&
}{}
\implies h \asymp n^{-5/(4b+5d)} 
.
\end{align*}
The resulting CPE is
\begin{align*}
O\left( N_n^{-3/2}[\log(N_n)]^3 \right)
&= O\left( n^{-3/2} n^{(15d/2)/(4b+5d)} [\log(n)]^3\right)
\Supplemental{
= O\left( n^{(-12b-15d+15d)/(8b+10d)} [\log(n)]^3 \right)
\\&
}{}
= O\left( n^{-6b/(4b+5d)} [\log(n)]^3 \right)
. 
\end{align*}
\Supplemental{
We also double-check that the remainder terms in \eqref{eqn:cpe-bias-two} are not larger than $\textrm{CPE}_\textrm{U}$ or $\textrm{CPE}_\textrm{Bias}$ with the CPE-optimal bandwidth. 
By \ref{a:h}(i'), $B_h N_n^{1/2}=o(1)$, so 
\[ B_h^3 N_n^{3/2} = B_h^2 N_n (B_h N_n^{1/2}) = o(B_h^2 N_n) . \]
With $h \asymp n^{-5/(4b+5d)}$,
\begin{align*}
B_h N_n^{-1/2}
&= O\left( h^{b} n^{-1/2} h^{-d/2} \right)
= O\left( n^{(-5b-2b-5d/2)/(4b+5d)} n^{(5d/2)/(4b+5d)} \right)
= O\left( n^{-7b/(4b+5d)} \right)
\\&= o\left( n^{-6b/(4b+5d)} [\log(n)]^3 \right) 
. 
\end{align*}
}{}

With $p\ne1/2$ and the calibrated method, note that $B_h^2 N_n\asymp N_n^{-3/2}$ implies $B_h\asymp N_n^{-5/4}\gtrsim B_h^2 N_n$, so the $B_h$ term must be accounted for.  Consequently,
\begin{align*}
          B_h        \asymp N_n^{-3/2}
&\implies h^{b}      \asymp (nh^d)^{-3/2}
\Supplemental{
 \implies h^{b+3d/2} \asymp n^{-3/2}
 \implies h^{b+3d}  \asymp n^{-3}
\\&
}{}
\implies h \asymp n^{-3/(b+3d)} 
.
\end{align*}
\Supplemental{
Since $B_h \asymp N_n^{-3/2}$, the $B_h^2 N_n\asymp N_n^{-3}N_n=N_n^{-2}$ term is smaller-order. 
The other remainder terms are also smaller-order: $B_h^3N_n^{3/2}\asymp N_n^{-3}$ and $B_h N_n^{-1/2}\asymp N_n^{-2}$. 
}{}%
The resulting CPE with $h \asymp n^{-3/(b+3d)}$ is
\begin{align*}
O\left( N_n^{-3/2}[\log(N_n)]^3 \right)
&= O\left( n^{-3/2} n^{(9d/2)/(b+3d)} [\log(n)]^3\right)
\Supplemental{
= O\left( n^{(-3b-9d+9d)/(2b+6d)} [\log(n)]^3 \right)
\\&
}{}
= O\left( n^{-3b/(2b+6d)} [\log(n)]^3 \right)
. 
\end{align*}
}{}

\Supplemental{


\section{Details on methods and theoretical properties}

\subsection{Steps to implement unconditional method}

To summarize unconditional quantile CI construction, as in \citet{Hutson1999} and implemented in our code:
\begin{enumerate}[1.]
 \item Parameters: determine the sample size $n$, quantile of interest $p\in(0,1)$, and coverage level $1-\alpha$.  Optionally, use the calibrated $\alpha$ in \eqref{eqn:hutson-CI-lower-calibrated}, \eqref{eqn:hutson-CI-upper}, or \eqref{eqn:hutson-CI-2s-calibrated}.  
 \item Endpoint index computation: solve for $u^h$ (for two-sided or lower one-sided CIs) and/or $u^l$ (for two-sided or upper one-sided CIs) from \eqref{eqn:uhdef}, using \eqref{eqn:QUIdist}.  For a two-sided CI, replace $\alpha$ in \eqref{eqn:uhdef} with $\alpha/2$.  
 \item CI construction: compute the lower and/or upper endpoints, respectively $\hat Q^L_X(u^l)$ and $\hat Q^L_X(u^h)$, using \eqref{eqn:QXLdef}. 
\end{enumerate}

\subsection{CPE comparison with other methods: details}\label{sec:app-CPE-comp}

As will be seen, Theorem \ref{thm:h-rate} implies that for the most common values of dimension $d$ and most plausible values of smoothness $s_Q$, our method is more accurate than inference based on asymptotic normality (or, equivalently, unsmoothed bootstrap) with a local polynomial estimator.  Only our uncalibrated method is compared here; our calibrated method is even more accurate. 

The only opportunity for normality to yield smaller CPE is to greatly reduce bias by using a very large local polynomial if $s_Q$ is large.  However, our method has smaller CPE when $d=1$ or $d=2$ even if $s_Q=\infty$. 

In finite samples, the normality approach may have additional error from estimation of the asymptotic variance, which includes the probability density of the error term at zero as discussed in \citet[pp.\ 764--766]{Chaudhuri1991}.  Additionally, large local polynomials cannot perform well unless there is a large local sample size.  

For the local polynomial approach, the CPE-optimal bandwidth and optimal CPE can be derived from the results in \citet{Chaudhuri1991}.  We balance the CPE from the bias with additional CPE from the Bahadur remainder from his Theorem 3.3(ii).  The CPE from bias is of order $N_n^{1/2}h^{s_Q}$.  \citet[Thm.\ 3.3]{Chaudhuri1991} gives a Bahadur-type expansion of the local polynomial quantile regression estimator that has remainder $R_n\sqrt{N_n}=O(N_n^{-1/4})$ (up to $\log$ terms) as in \citet{Bahadur1966}, but recently \citet{Portnoy2012} has shown that the CPE is nearly $O(N_n^{-1/2})$ in such cases.  
Solving $N_n^{1/2}h^{s_Q}\asymp N_n^{-1/2}\asymp (nh^d)^{-1/2}$ yields $h^*\asymp n^{-1/(s_Q+d)}$ and optimal CPE (nearly) $O(\sqrt{N_n}B_n)=O(N_n^{-1/2})=O(n^{-s_Q/(2s_Q+2d)})$.  

  
As illustrated in the left panel of Figure \ref{fig:CPE-comp}, if $s_Q=2$ (one Lipschitz-continuous derivative), then the optimal CPE from asymptotic normality is nearly $O(n^{-2/(4+2d)})$, which is always larger than our method's CPE.  
With $d=1$, this is $n^{-1/3}$, significantly larger than our $n^{-2/3}$ (one-sided: $n^{-4/7}$).  With $d=2$, $n^{-1/4}$ is larger than our $n^{-1/2}$ (one-sided: $n^{-2/5}$).  It remains larger for all $d$, even for one-sided inference, since the bias is the same for both methods and the unconditional $L$-statistic inference is more accurate than normality. 

From another perspective, as in the right panel of Figure \ref{fig:CPE-comp}: what amount of smoothness (and local polynomial degree) is needed for asymptotic normality to match our method's CPE?  
For the most common cases of $d\in\{1,2\}$ (one-sided: $d=1$), normality is worse even with infinite smoothness. 
With $d=3$ (one-sided: $d=2$), normality needs $s_Q\ge12$ to match our CPE.  If $n$ is large, then maybe such a high degree ($k_Q\ge11$) local polynomial will be appropriate, but often it is not.  Interaction terms are required, so an $11$th-degree polynomial has $\sum_{T=d-1}^{k_Q+d-1} \binom{T}{d-1} = 364$ terms. 
As $d\to\infty$, the required smoothness approaches $s_Q=4$ (one-sided: $8/3$) from above, though again the number of terms in the local polynomial grows with $d$ as well as $k_Q$ and may thus still be prohibitive in finite samples.  

Basic bootstraps claim no refinement over asymptotic normality, so the foregoing discussion also applies to such methods.  Without using smoothed or $m$-out-of-$n$ bootstrap, which require additional smoothing parameter choice and computation time, Studentization offers no CPE improvement.  \citet{GangopadhyaySen1990} examine the bootstrap percentile method for a uniform kernel (or nearest neighbor) estimator with $d=1$ and $s_Q=2$, but they only show first-order consistency.  Based on their (2.14), the optimal error seems to be $O(n^{-1/3})$ when $h\asymp n^{-1/3}$ (matching our derivation from \citet{Chaudhuri1991}) to balance the bias and remainder terms, improved by \citet{Portnoy2012} from $O(n^{-2/11})$ CPE when $h\asymp n^{-3/11}$.

\subsection{Joint confidence intervals and a possible uniform band}

Beyond the pointwise CIs in Definition \ref{def:method}, the usual Bonferroni adjustment $\alpha/m$ gives joint CIs over $m$ different values of $x_0$.  If the $C_h$ are mutually exclusive, then the local samples are independent due to Assumption \ref{a:sampling}, and the adjustment can be refined to $1-(1-\alpha)^{1/m}$, which maintains exact joint coverage (not conservative like Bonferroni).  

The number of such mutually exclusive windows grows as $n\to\infty$, which may be helpful for constructing a uniform confidence band.  \citet[\S4.3]{HorowitzLee2012} suggest forming uniform confidence bands for nonparametric functions by considering $L$ grid points with $L\to\infty$ as $n\to\infty$, interpolating the corresponding $L$ CIs (possessing joint $1-\alpha$ coverage).  They also consider restrictions on smoothness or monotonicity that can replace $L\to\infty$.  If our Assumption \ref{a:Q} is strengthened to a uniform bound on the second derivative (with respect to $X$) of the conditional quantile function, then the local bandwidths have a common asymptotic rate, $h$, and there can be order $1/h$ mutually exclusive local samples.  
The error from a linear approximation of the conditional quantile function over a span of $h$ is $O(h^2)$.  This is smaller than the length of the CIs, which is order $1/\sqrt{nh^d}=h$ when using the CPE-optimal bandwidth rate $h\asymp n^{-1/(2+d)}$ for two-sided CIs from Theorem \ref{thm:h-rate} with $b=2$. 

The growing number of CIs is not problematic given independent sampling, but the fact that $\alpha\to0$ (at rate $h$) has not been rigorously treated.  We leave such proof to future work, suggesting due caution in applying the uniform band until then. 
As discussed in Section \ref{sec:sim}, a \citet{Hotelling1939} tube-based calibration of $\alpha$ also appears to yield reasonable uniform confidence bands, but it has even less formal justification.

\section{Plug-in bandwidth}\label{sec:h-plugin}

\subsection{Plug-in bandwidth: overview and setup}

The dominant terms of $\textrm{CPE}_\textrm{U}$ and $\textrm{CPE}_\textrm{Bias}$ have known expressions, not simply rates.  Theorem \ref{thm:IDEAL-single} gives an exact expression for the $O(N_n^{-1})$ $\textrm{CPE}_\textrm{U}$ term. 
\Supplemental{This could be used to reduce the CPE to $O\left(N_n^{-3/2}[\log(N_n)]^3\right)$ via analytic calibration, but we use it to determine the precise value of the optimal bandwidth.  In theory, it is better to use the analytic calibration and an ad hoc bandwidth of the proper rate.  In practice, it is helpful to know that the $O(N_n^{-1})$ CPE term only leads to over-coverage, which implies that smaller $h$ is always more conservative, and the detrimental effect of an ad hoc bandwidth can be significant in smaller samples.}{Instead of using this for calibration, which achieves smaller theoretical CPE, we use it to pin down the optimal bandwidth value since the bandwidth is a crucial determinant of finite-sample CPE.}  CPE from bias is found in \eqref{eqn:cpe-bias-lower}, \eqref{eqn:cpe-bias-upper}, and \eqref{eqn:cpe-bias-two}.  Our bandwidth sets the dominant terms of $\textrm{CPE}_\textrm{U}$ and $\textrm{CPE}_\textrm{Bias}$ to sum to zero if possible, or else minimizes their sum. 
 
\Supplemental{While we know $\textrm{CPE}_\textrm{U}>0$, we do not always know the sign of the bias.  Specifically, when the rate-limiting term of $B_h$ is determined by H\"older continuity (of $Q_{Y|X}(\cdot;\cdot)$ or $f_X(\cdot)$), we lose track of the sign at that step.  However, when $k_Q\ge2$ and $k_X\ge1$, the H\"older continuity terms all end up in the remainder while the rate-limiting terms are signed derivatives.}{}


The plug-in bandwidth is derived for $d=1$ and $b=2$.  Let $\phi(\cdot)$ be the standard normal PDF, $z_{1-\alpha}$ be the $(1-\alpha)$-quantile of the standard normal distribution, $u_h>p$ (or $u_\ell<p$) be the quantile determining the high (low) endpoint of a lower (upper) one-sided CI, $\epsilon_h\equiv (N_n+1) u_h-\lfloor (N_n+1) u_h\rfloor$, and $\epsilon_\ell\equiv (N_n+1) u_\ell-\lfloor (N_n+1) u_\ell\rfloor$.   Let $I_H$ denote a $100(1-\alpha)\%$ CI constructed using \citeposs{Hutson1999} method on a univariate data sample of size $N_n$.  The coverage probability of a lower one-sided $I_H$ (with the upper one-sided result substituting $\epsilon_\ell$ for $\epsilon_h$) is given in Theorem \ref{thm:IDEAL-single}:
\begin{align}
\label{eqn:CPE-hut-1s}
P\{F^{-1}(p)\in I_H\}
  &= 1-\alpha 
    +N_n^{-1}z_{1-\alpha} \frac{\epsilon_h(1-\epsilon_h)}{p(1-p)} \phi(z_{1-\alpha})
    +O\left(N_n^{-3/2}[\log(N_n)]^3\right)  ,
\intertext{or for a two-sided CI,}
\label{eqn:CPE-hut-2s}
\begin{split}
P\{F^{-1}(p)\in I_H\}
  &= 1-\alpha
    +N_n^{-1} z_{1-\alpha/2} 
     \frac{\epsilon_h(1-\epsilon_h)+\epsilon_\ell(1-\epsilon_\ell)}{p(1-p)} 
     \phi(z_{1-\alpha/2})
\\&\quad
   +O\left(N_n^{-3/2}[\log(N_n)]^3\right)  .
\end{split}
\end{align}
In either case there is $O(N_n^{-1})$ over-coverage.

Since $f''_X(\cdot)$ is uniformly bounded in a neighborhood of zero, then 
\begin{align*}
P_C &= \int_{C_h}f_X(x)\,\mathrm{d}x
     = \int_{C_h}[f_X(0)+xf'_X(0)+(1/2)x^2f''_X(\tilde x)] \,\mathrm{d}x 
     = 2hf_X(0)+O(h^3) . 
\end{align*}
Since $N_n\sim\textrm{Binomial}(n,P_C)$, $N_n=nP_C+o_p(1)$, so we approximate
\begin{align}
\label{eqn:N-approx}
N_n
  &\doteq nP_C
   = 2nhf_X(0) +O(h^3)  .
\end{align}

\subsection{Plug-in bandwidth: one-sided}

For a one-sided CI, $\textrm{CPE}_\textrm{U}$ is in \eqref{eqn:CPE-hut-1s}, and $\textrm{CPE}_\textrm{Bias}$ is in \eqref{eqn:cpe-bias-lower} or \eqref{eqn:cpe-bias-upper}\Supplemental{ for lower and upper one-sided CIs, respectively}{}.  Since $\textrm{CPE}_\textrm{U}>0$ always, when $\textrm{CPE}_\textrm{Bias}<0$, with $\epsilon\in\{\epsilon_h,\epsilon_\ell\}$ and correspondingly $u\in\{u_h,u_\ell\}$, the optimal $h$ equates 
\begin{align}
\label{eqn:plugin-1s-equality}
\begin{split}
N_n^{-1} 
    z_{1-\alpha} \frac{\epsilon(1-\epsilon)}{p(1-p)} \phi(z_{1-\alpha}) 
  &= h^2
      \left|
      \frac{f_X(0) F_{Y|X}^{(0,2)}(\xi_p;0)
            +2 f_X'(0) F_{Y|X}^{(0,1)}(\xi_p;0)}
           {6 f_X(0) f_{Y|X}(\xi_p;0)}  
      \right|
\\&\quad\times
        N_n^{1/2} \frac{f_{Y|X}\left(Q_{Y|X}(u;C_h);C_h\right)}{\sqrt{u(1-u)}} \phi(z_{1-\alpha})
.
\end{split}
\end{align}
%
Plugging in $N_n\doteq 2nhf_X(0)$ and $u=p+O(N_n^{-1/2})$, the optimal $h$ is now an explicit function of known values and estimable objects.  
To avoid iteration, we plug in $\epsilon_h=\epsilon_\ell=0.2$ as a rule of thumb.  (More precisely, $\epsilon=0.5-(1/2)\sqrt{1-4\sqrt2/9}\approx0.20$ makes the constants cancel.) 
Let $\hat h_{++}$ be the plug-in bandwidth if both $\textrm{CPE}_\textrm{U}>0$ and $\textrm{CPE}_\textrm{Bias}>0$, and $\hat h_{+-}$ if $\textrm{CPE}_\textrm{Bias}<0$. 
Up to smaller-order terms,
\begin{align*}
[2n\hat h_{+-}f_X(0)]^{-1}z_{1-\alpha}\frac{\epsilon(1-\epsilon)}{p(1-p)}\phi(z_{1-\alpha})
  &= \hat h_{+-}^2
      \frac{f_X(0) F_{Y|X}^{(0,2)}(\xi_p;0)
            +2 f_X'(0) F_{Y|X}^{(0,1)}(\xi_p;0)}
           {6 f_X(0) f_{Y|X}(\xi_p;0)} \\
  &\quad\times
     [2n\hat h_{+-}f_X(0)]^{1/2}  
     [p(1-p)]^{-1/2} \phi(z_{1-\alpha})
     f_{Y|X}(\xi_p;0),  
\end{align*}
\Supplemental{%
\begin{align*}
\hat h_{+-}^{7/2}
  &= \frac{n^{-3/2} 2^{-3/2} z_{1-\alpha} \epsilon(1-\epsilon)/\sqrt{p(1-p)}}
     { \sqrt{f_X(0)} \left\{f_X(0) F_{Y|X}^{(0,2)}(\xi_p;0) +2 f_X'(0) F_{Y|X}^{(0,1)}(\xi_p;0) \right\} / 6}  , \\ 
\hat h_{+-}
  &= n^{-3/7} 
     \left(
       \frac{z_{1-\alpha}}
            {3 \sqrt{p(1-p)f_X(0)} \left\{f_X(0) F_{Y|X}^{(0,2)}(\xi_p;0) +2 f_X'(0) F_{Y|X}^{(0,1)}(\xi_p;0) \right\}}
     \right)^{2/7}  , \\
\hat h_{++} 
  &= \hat h_{+-} [-2d/(2b+d)]^{2/(2b+3d)}
   \approx -0.77\hat h_{+-}  .
\end{align*}
}{and solving for $\hat h$ (and plugging in estimators of the unknown objects) yields the expression in the text.} 
The extra coefficient on $\hat h_{++}$ is due to taking the first-order condition of the CPE minimization problem rather than simply setting the two types of CPE to sum to zero.  These $\hat h$ hold for both lower and upper one-sided inference.

\subsection{Plug-in bandwidth: two-sided, \texorpdfstring{$p=1/2$}{p equals 1/2}}

For two-sided inference with $p=1/2$, the $B_h$ term in $\textrm{CPE}_\textrm{Bias}$ zeroes out because $f_{\beta,h}(p)=f_{\beta,\ell}(p)$ exactly.  This happens because $u_h-p=p-u_\ell$ exactly, so $u_\ell=1-u_h$, and the corresponding beta PDFs are reflections of each other about $p=1/2$: $f_{\beta,h}(x)=f_{\beta,\ell}(1-x)$ since $\beta(u_\ell(N_n+1),(1-u_\ell)(N_n+1))=\beta((1-u_h)(N_n+1),u_h(N_n+1))$. 

Consequently, $\textrm{CPE}_\textrm{Bias}<0$ always, while $\textrm{CPE}_\textrm{U}>0$, 
so the optimal $h$ sets $\textrm{CPE}_\textrm{Bias}+\textrm{CPE}_\textrm{U}\doteq0$. 
Applying \eqref{eqn:dfX-dphi} with $k=\lfloor u_h(N_n+1)\rfloor=u_h(N_n+1)+O(N_n^{-1})$ or $\lfloor u_\ell(N_n+1)\rfloor=u_\ell(N_n+1)+O(N_n^{-1})$, using $u_h-p=N_n^{-1/2}z_{1-\alpha/2}\sqrt{u_h(1-u_h)}+O(N_n^{-1})$ and $p-u_\ell=N_n^{-1/2}z_{1-\alpha/2}\sqrt{u_\ell(1-u_\ell)}+O(N_n^{-1})$ from Lemma \ref{lem:u1u2-approx},
\begin{align*}
f'_{\beta,h}(p) 
&= 
\frac{N_n}{u_h(1-u_h)} 
\phi'\overbrace{\left( \frac{p-u_h}{\sqrt{u_h(1-u_h)/N_n}} \right)}^{\textrm{Lemma \ref{lem:u1u2-approx}}}
+O\left( N_n^{1/2}[\log(N_n)]^4 \right)
\\&= 
\frac{N_n}{p(1-p)} 
\phi'\left( -z_{1-\alpha/2} \right)
+O\left( N_n^{1/2}[\log(N_n)]^4 \right)
\\&= 
\frac{N_n}{p(1-p)} 
z_{1-\alpha/2}
\phi\left( z_{1-\alpha/2} \right)
+O\left( N_n^{1/2}[\log(N_n)]^4 \right)
,\\
f'_{\beta,\ell}(p) 
&= 
\frac{N_n}{u_\ell(1-u_\ell)} 
\phi'\overbrace{\left( \frac{p-u_\ell}{\sqrt{u_\ell(1-u_\ell)/N_n}} \right)}^{\textrm{Lemma \ref{lem:u1u2-approx}}}
+O\left( N_n^{1/2}[\log(N_n)]^4 \right)
\\&= 
\frac{N_n}{p(1-p)} 
\phi'\left( z_{1-\alpha/2} \right)
+O\left( N_n^{1/2}[\log(N_n)]^4 \right)
\\&= 
-\frac{N_n}{p(1-p)} 
 z_{1-\alpha/2}
 \phi\left( z_{1-\alpha/2} \right)
+O\left( N_n^{1/2}[\log(N_n)]^4 \right)
,\\
f'_{\beta,h}(p) - f'_{\beta,\ell}(p) 
&= 
2\frac{N_n}{p(1-p)} 
 z_{1-\alpha/2}
 \phi\left( z_{1-\alpha/2} \right)
+O\left( N_n^{1/2}[\log(N_n)]^4 \right)
.
\label{eqn:fpbetah-fpbetaell}\refstepcounter{equation}\tag{\theequation}
\end{align*}

Using \eqref{eqn:CPE-hut-2s}, \eqref{eqn:cpe-bias-two}, and \eqref{eqn:fpbetah-fpbetaell}, 
\begin{align*}
& N_n^{-1} 
    z_{1-\alpha/2}
    \frac{\epsilon_h(1-\epsilon_h)+\epsilon_\ell(1-\epsilon_\ell)}{p(1-p)}
    \phi(z_{1-\alpha/2})   \\
  &\quad= (1/2) \overbrace{\hat h^4 \left(
     \frac{f_X(0) F_{Y|X}^{(0,2)}(\xi_p;0) + 2 f_X'(0) F_{Y|X}^{(0,1)}(\xi_p;0)}
          {6 f_X(0) f_{Y|X}(\xi_p;0)}
     \right)^2}^{B_h^2\textrm{, Lemma \ref{lem:bias}}}  \\ 
  &\qquad\times
    \Bigl\{ \overbrace{z_{1-\alpha/2}N_n\phi(z_{1-\alpha/2})
     2[p(1-p)]^{-1}}^{\textrm{\eqref{eqn:fpbetah-fpbetaell}}} \left[f_{Y|X}(\xi_p;0)\right]^2 \Bigr\} , \\
& [2n\hat hf_X(0)]^{-2} 
    2 \epsilon(1-\epsilon) \\
  &\quad= \frac{\hat h^4}{36 f_X(0)^2 f_{Y|X}(\xi_p;0)}
     \left(f_X(0) F_{Y|X}^{(0,2)}(\xi_p;0) + 2 f_X'(0) F_{Y|X}^{(0,1)}(\xi_p;0)\right)^2
       ,  \\
\label{eqn:h-plugin-median-app}\refstepcounter{equation}\tag{\theequation}
\hat h
  &= n^{-1/3} 
     \left| f_X(0) F_{Y|X}^{(0,2)}(\xi_p;0) + 2 f_X'(0) F_{Y|X}^{(0,1)}(\xi_p;0) \right|^{-1/3} ,
\end{align*}
again using rule-of-thumb $\epsilon=0.2$ and rounding down a leading $1.19$ to one.

\subsection{Plug-in bandwidth: two-sided, \texorpdfstring{$p\ne1/2$}{p not 1/2}}

For two-sided inference with $p\ne1/2$, the dominant bias terms do not cancel completely, but the difference is of smaller order than the $O(N_n^{1/2})$ in the one-sided case.  From \eqref{eqn:beta-PDF-diff-approx},
\begin{align*}
f_{\beta,h}(p)-f_{\beta,\ell}(p)
  &= -(1/3) z_{1-\alpha/2} \phi(z_{1-\alpha/2}) \frac{2p-1}{p(1-p)}
     +O(N_n^{-1/2})  
   = O(1) . 
\end{align*}
The $B_h$ term is driven by skewness: it is zero if $p=1/2$, 
positive for $p<1/2$, and negative for $p>1/2$. 
The CPE terms are now of orders $N_n^{-1}\asymp n^{-1}h^{-1}$, $B_h\asymp h^2$, and $N_n B_h^2\asymp nh^5$.  This implies that $h^*\asymp n^{-1/3}$, so all CPE terms are $O(n^{-2/3})$.  

Using \eqref{eqn:cpe-bias-two}, \eqref{eqn:fpbetah-fpbetaell}, and \eqref{eqn:beta-PDF-diff-approx}, up to smaller-order terms (hence $\doteq$), $\textrm{CPE}_\textrm{Bias}$ is
\begin{align*}
\textrm{CPE}_\textrm{Bias}
&\doteq 
  B_h f_{Y|X}\left(Q_{Y|X}(p;C_h);C_h\right) 
   \overbrace{\left[ f_{\beta,h}(p) - f_{\beta,\ell}(p) \right]}^{\textrm{\eqref{eqn:beta-PDF-diff-approx}}}
\\&\quad
-\frac{1}{2} B_h^2 \left[f_{Y|X}\left(Q_{Y|X}(p;C_h);C_h\right)\right]^2 
 \overbrace{\bigl[ f'_{\beta,h}(p) - f'_{\beta,\ell}(p) \bigr]}^{\textrm{\eqref{eqn:fpbetah-fpbetaell}}}
\\&\doteq 
  B_h f_{Y|X}\left(Q_{Y|X}(p;C_h);C_h\right) 
   \left[ -(1/3) z_{1-\alpha/2} \phi(z_{1-\alpha/2}) \frac{2p-1}{p(1-p)} \right]
\\&\quad
-\frac{1}{2} B_h^2 \left[f_{Y|X}\left(Q_{Y|X}(p;C_h);C_h\right)\right]^2 
 2\frac{N_n}{p(1-p)} 
 z_{1-\alpha/2}
 \phi\left( z_{1-\alpha/2} \right)
\\&\doteq 
\frac{f_{Y|X}(\xi_p;0)z_{1-\alpha/2}\phi(z_{1-\alpha/2})}{p(1-p)}
\left[ B_h(1-2p)/3 - B_h^2 N_n f_{Y|X}(\xi_p;0) \right] 
.
\label{eqn:CPE-bias-2s-general}\refstepcounter{equation}\tag{\theequation}  
\end{align*}
With large enough $h$, \eqref{eqn:CPE-bias-2s-general} is always negative (unless $B_h=0$ exactly) since if $h\gtrsim h^*\asymp n^{-1/3}$ then $N_n B_h^2\asymp nh^5 \gtrsim h^2\asymp B_h$, so there always exists an optimal $h$ that cancels the dominant terms of $\textrm{CPE}_\textrm{U}$ and $\textrm{CPE}_\textrm{Bias}$.  
 This $\hat h$ solves 
\begin{align*}
-N_n^{-1}
  & z_{1-\alpha/2} 
     \frac{\epsilon_h(1-\epsilon_h)+\epsilon_\ell(1-\epsilon_\ell)}{p(1-p)} 
     \phi(z_{1-\alpha/2})  \\
  &= \frac{f_{Y|X}(\xi_p;0)z_{1-\alpha/2}\phi(z_{1-\alpha/2})}{p(1-p)}
    \left[ B_h(1-2p)/3 - B_h^2 N_n f_{Y|X}(\xi_p;0) \right] . 
\end{align*}
Continuing to solve for $\hat h$ using \eqref{eqn:N-approx} and then Lemma \ref{lem:bias}, 
\begin{align*}
- & [2n\hat hf_X(0)]^{-1} 
    2\epsilon(1-\epsilon) \\
  &= f_{Y|X}(\xi_p;0)
     \left\{ \hat h^2(B_h/h^2)(1-2p)/3
            -\hat h^4(B_h^2/h^4) [2n\hat hf_X(0)] f_{Y|X}(\xi_p;0) \right\}, \\
\big\{n^{-1} & [f_X(0)]^{-1} \epsilon(1-\epsilon) \big\}  \\
  &= \hat h^6 \left\{2 n \left[f_{Y|X}(\xi_p;0)\right]^2 f_X(0) B_h^2/h^4 \right\}
    -\hat h^3 \left\{f_{Y|X}(\xi_p;0) [(1-2p)/3] B_h/h^2 \right\} , \\
0 &= (\hat h^3)^2 \{na\} - (\hat h^3) \{b\} -\{c/n\}, \\
\hat h^3 &= \frac{b\pm\sqrt{b^2+4ac}}{2an}   , \\
\hat h
\Supplemental{%
  &= n^{-1/3} \left( \frac{b+\sqrt{b^2+4ac}}{2a} \right)^{1/3}  \\}{}
\Supplemental{%
  &= n^{-1/3} \left[ \frac
       {\frac{B_h}{|B_h|} (1-2p)/3
        +\sqrt{(1-2p)^2/9 
               +8 \epsilon(1-\epsilon)} }
       {4 f_{Y|X}(\xi_p;0) f_X(0) |B_h|/h^2 } 
     \right]^{1/3}  \\}{}
  &= n^{-1/3} \left( \frac
       { (B_h/|B_h|) (1-2p)
        +\sqrt{(1-2p)^2  +4} }
       {2 \left|   f_X(0)  F_{Y|X}^{(0,2)}(\xi_p;0)
                +2 f_X'(0) F_{Y|X}^{(0,1)}(\xi_p;0) \right| } 
     \right)^{1/3}  ,
\end{align*}
setting $\epsilon$ to match \eqref{eqn:h-plugin-median-app} when $p=1/2$.

\subsection{Plug-in bandwidth: beta PDF difference approximation}

The following approximates the beta PDF difference $f_{\beta,h}(p)-f_{\beta,\ell}(p)$, which are the respective PDFs of the $\beta\left(u_h(N_n+1),(1-u_h)(N_n+1)\right)$ and $\beta\left(u_\ell(N_n+1),(1-u_\ell)(N_n+1)\right)$ distributions, evaluated at $p$, referring to notation in the proof of Theorem \ref{thm:h-rate} (two-sided case).  For either $u=u_h$ or $u=u_\ell$, from the proof of Lemma \ref{lem:den-i}, where here $\Delta x_1=p$, $\Delta x_2=1-p$, $\Delta k_1=(n+1)u$, $\Delta k_2=(n+1)(1-u)$, $J=1$, and $u-p=O(n^{-1/2})$, 
\begin{align*}
\log(f_{\Delta X}(\Delta x))
  &= K + h(\Delta x)+O(n^{-1}), \\
K &= (J/2)\log(n/2\pi)+\frac{1}{2}\sum_{j=1}^{J+1}\log(n/[\Delta k_j-1]) , \\
e^K \Supplemental{&= \sqrt{\frac{n}{2\pi}}\sqrt{\frac{n}{(n+1)u-1}}\sqrt{\frac{n}{(n+1)(1-u)-1}} \\}{}
  &= \left[1+O(n^{-1})\right] \sqrt{\frac{n}{2\pi}}\sqrt{\frac{1}{u(1-u)}} , \\
h(\Delta x)
  &= -(1/2)(n-J)^2 \sum_{j=1}^{J+1}\frac{\left(\Delta x_j-[\Delta k_j-1]/[n-J]\right)^2}{\Delta k_j-1} \\
  &\quad
     +(1/3)(n-J)^3 \sum_{j=1}^{J+1}\frac{\left(\Delta x_j-[\Delta k_j-1]/[n-J]\right)^3}{(\Delta k_j-1)^2}
     +O(n^{-1}) \\
\Supplemental{%
  &= -(1/2)n^2 \Bigl[ \frac{(p-u)^2+2(p-u)(1-2u)/n+O(n^{-2})}{nu}  \\
  &\qquad\qquad\qquad+\frac{(1-p-(1-u))^2+2(u-p)(2u-1)/n+O(n^{-2})}{n(1-u)} \Bigr]    \\
  &\quad
     +(1/3)n^3 \left[ (p-u)^3/(nu)^2 + (1-p-(1-u))^3/(n(1-u))^2 \right] 
     +O(n^{-1}) \\
}{}
\label{eqn:h-delta-x}\refstepcounter{equation}\tag{\theequation}
  &= -(1/2) n \left[ (u-p)^2/[u(1-u)]  +2(u-p)(2u-1)/[nu(1-u)] \right] \\
  &\quad
     +(1/3) n (p-u)^3 [u^{-2}-(1-u)^{-2}]
     +O(n^{-1}) .
\end{align*}

First, note that $e^K$ will put the $\sqrt n$ coefficient in front, so the $O(n^{-1})$ remainder in $h(\Delta x)$ is required in order to yield $O(n^{-1/2})$ overall.  Even within $e^K$, we must keep track of $u=p\pm n^{-1/2}z\sqrt{p(1-p)}+O(n^{-1})$, not just $u\approx p$, where for brevity $z\equiv z_{1-\alpha/2}$.  So
\begin{align*}
[u_h(1-u_h)]^{-1/2}
\Supplemental{%
  &= \left[ p+n^{-1/2}z\sqrt{p(1-p)}+O(n^{-1}) \right]^{-1/2} \\
  &\quad\times \left[ 1-p-n^{-1/2}z\sqrt{p(1-p)}+O(n^{-1}) \right]^{-1/2} \\
  &= \left[
          p(1-p) +n^{-1/2}z\sqrt{p(1-p)}(1-2p)
         \right]^{-1/2} \\
  &= [p(1-p)]^{-1/2} -(1/2)[p(1-p)]^{-3/2}n^{-1/2}z\sqrt{p(1-p)}(1-2p) +O(n^{-1})  \\
}{}
  &= [p(1-p)]^{-1/2} -(1/2)n^{-1/2}z\frac{1-2p}{p(1-p)} +O(n^{-1})  , \\
e^K(u_h)
  &= \frac{n^{1/2}}{\sqrt{2\pi p(1-p)}} 
     \left[ 1 +(1/2)n^{-1/2}z\frac{2p-1}{\sqrt{p(1-p)}} \right]
    +O(n^{-1/2}) , \\
e^K(u_\ell)
  &= \frac{n^{1/2}}{\sqrt{2\pi p(1-p)}} 
     \left[ 1 -(1/2)n^{-1/2}z\frac{2p-1}{\sqrt{p(1-p)}} \right]
    +O(n^{-1/2}) .
\end{align*}

Second, for $u_h$, let $u_h=p+n^{-1/2}c_1+n^{-1}c_2+O(n^{-3/2})$ from Lemma \ref{lem:u1u2-approx}.  Then, $(u-p)^2 = n^{-1}c_1^2 +2n^{-3/2}c_1c_2 +O(n^{-2})$ and $(u-p)^3=n^{-3/2}c_1^3+O(n^{-2})$.  Plugging into \eqref{eqn:h-delta-x} and using the small $x$ approximations $e^x=1+x+O(x^2)$ and $1/(1+x)=1-x+O(x^2)$,
\begin{align*}
h(\Delta x)
\Supplemental{%
  &= -\frac{1}{2}\frac{n\left[n^{-1}c_1^2+2n^{-3/2}c_1c_2+O(n^{-2}) 
                        +2\left(n^{-1/2}c_1+O(n^{-1})\right)\left(2p-1+O(n^{-1/2})\right)/n\right]}{u(1-u)} \\
  &\quad
     +\frac{1}{3}\frac{(2u-1)n[n^{-3/2}c_1^3+O(n^{-2})]}{u^2(1-u)^2}
     +O(n^{-1}) \\
  &= -\frac{1}{2}\frac{c_1^2}{u(1-u)}
     -\frac{1}{2}\frac{2n^{-1/2}c_1c_2}{u(1-u)}
     -\frac{1}{2}\frac{2n^{-1/2}c_1(2p-1)}{u(1-u)}
     +\frac{1}{3}\frac{(2p-1)n^{-1/2}c_1^3}{p^2(1-p)^2}
     +O(n^{-1}) \\
  &= -\frac{1}{2}\frac{c_1^2}{p(1-p) -n^{-1/2}z\sqrt{p(1-p)}(2p-1) +O(n^{-1})}
     -n^{-1/2}\frac{c_1(c_2+2p-1)}{p(1-p)} \\
  &\quad
     +n^{-1/2}\frac{1}{3}\frac{(2p-1)c_1^3}{p^2(1-p)^2}
     +O(n^{-1}) \\
  &= -\frac{1}{2}\frac{c_1^2}{p(1-p)}\left[ 1 +n^{-1/2}z(2p-1)/\sqrt{p(1-p)} \right] \\
  &\quad
     -n^{-1/2}\left( \frac{c_1(c_2+2p-1)}{p(1-p)}
                    -\frac{(2p-1)c_1^3}{3p^2(1-p)^2} \right)
     +O(n^{-1}) \\
}{}
  &= -\frac{1}{2}\frac{c_1^2}{p(1-p)}
     -n^{-1/2}\left( \frac{c_1(c_2+2p-1)}{p(1-p)}
                    -\frac{(2p-1)c_1^3}{3p^2(1-p)^2}
                    +\frac{c_1^2 z (2p-1)}{2[p(1-p)]^{3/2}}\right)
     +O(n^{-1}) \\
\Supplemental{%
  &= -\frac{1}{2}\frac{z^2 p(1-p)}{p(1-p)} \\
  &\quad
     -n^{-1/2}\biggl( \frac{z\sqrt{p(1-p)} \left(6(2p-1)-(2p-1)(z^2+2)\right)}{6p(1-p)} \\
  &\qquad\qquad\qquad
                     -\frac{(2p-1)z^3[p(1-p)]^{3/2}}{3p^2(1-p)^2}
                     +\frac{z^2 p(1-p) z (2p-1)}{2[p(1-p)]^{3/2}}\biggr)
     +O(n^{-1}) \\
  &= -z^2/2
     -n^{-1/2}\left(-\frac{(2p-1)(z^3-4z)}{6\sqrt{p(1-p)}} 
                    -\frac{2 (2p-1)z^3}{6\sqrt{p(1-p)}}
                    +\frac{3 z^3 (2p-1)}{6\sqrt{p(1-p)}} \right)
     +O(n^{-1}) \\
  &= -z^2/2
     -n^{-1/2}\frac{2p-1}{6\sqrt{p(1-p)}} \left( -(z^3-4z) -2z^3 +3 z^3 \right)
     +O(n^{-1}) \\
}{}
  &= -z^2/2 - 2 n^{-1/2} z \frac{2p-1}{3\sqrt{p(1-p)}} +O(n^{-1}) , \\
\exp&\{h(\Delta x)\}
  = e^{-z^2/2}
     \left[1 - (2/3)n^{-1/2} z \frac{2p-1}{\sqrt{p(1-p)}} \right] 
    +O(n^{-1}) ,
\intertext{while for $u_\ell$,}
\exp&\{h(\Delta x)\}
  = e^{-z^2/2}
     \left[1 + (2/3)n^{-1/2} z \frac{2p-1}{\sqrt{p(1-p)}} \right] 
    +O(n^{-1}) .
\end{align*}

Altogether,
\begin{align*}
f_{\beta,h}(p)
\Supplemental{%
  &= \left( \frac{n^{1/2}}{\sqrt{2\pi p(1-p)}} 
     \left[ 1 +(1/2)n^{-1/2}z\frac{2p-1}{\sqrt{p(1-p)}} \right]
    +O(n^{-1/2}) \right) \\
  &\quad\times
    \left( e^{-z^2/2}
     \left[ 1 - (2/3) n^{-1/2} z \frac{2p-1}{\sqrt{p(1-p)}} \right]
    +O(n^{-1}) \right) \\
}{}
  &= \frac{n^{1/2}}{\sqrt{p(1-p)}} \phi(z)
     \left[ 1 - (1/6) n^{-1/2} z \frac{2p-1}{\sqrt{p(1-p)}} \right]
    +O(n^{-1/2}) , \\
f_{\beta,\ell}(p)
  &= \frac{n^{1/2}}{\sqrt{p(1-p)}} \phi(z)
     \left[ 1 + (1/6) n^{-1/2} z \frac{2p-1}{\sqrt{p(1-p)}} \right]
    +O(n^{-1/2}) , \\
\Supplemental{%
f_{\beta,h}(p) &- f_{\beta,\ell}(p) \\
  &= \frac{n^{1/2}}{\sqrt{p(1-p)}} \phi(z)
     \left[ -(1/3)  n^{-1/2} z \frac{2p-1}{\sqrt{p(1-p)}} \right]
    +O(n^{-1/2}) \\
}{%
f_{\beta,h}(p) - f_{\beta,\ell}(p) 
}
\label{eqn:beta-PDF-diff-approx}\refstepcounter{equation}\tag{\theequation}
  &= -(1/3) z \phi(z) \frac{2p-1}{p(1-p)}
     +O(n^{-1/2})  .
\end{align*}
}{}


\Supplemental{%

\section{Additional empirical results}

The $L$-statistic joint CIs (``L-stat'') in Figure \ref{fig:emp2} are the same as those for $p=0.5$ and $p=0.9$ in Figure \ref{fig:emp}, but now a nonparametric (instead of quadratic) conditional quantile estimate is shown, along with joint CIs from \citet{FanLiu2016}.  The L-stat CIs are generally shorter, but there can be exceptions, as seen especially in the bottom right graph.  Of course, shorter is not better if coverage probability is sacrificed; we explore both properties in the simulations of Section \ref{sec:sim}.

\begin{figure}[htbp]
 \centering
\hspace*{\fill}%
 \includegraphics[clip=true,trim=10 40 30 70,width=0.45\textwidth]{2015_08_16_Figures/quantile_inf_np_ex_jt_food_comparison_sub8528.pdf}%
 \hfill%
 \includegraphics[clip=true,trim=40 40 0 70,width=0.45\textwidth]{2015_08_16_Figures/quantile_inf_np_ex_jt_fuel_comparison_sub8528.pdf}%
\hspace*{\fill}%
 \break%
\hspace*{\fill}%
 \includegraphics[clip=true,trim=10 5 30 80,width=0.45\textwidth]{2015_08_16_Figures/quantile_inf_np_ex_jt_clothing_comparison_sub8528.pdf}%
 \hfill%
 \includegraphics[clip=true,trim=40 5 0 80,width=0.45\textwidth]{2015_08_16_Figures/quantile_inf_np_ex_jt_alc_comparison_sub8528.pdf}%
\hspace*{\fill}%
 \caption{\label{fig:emp2}Joint (over the nine expenditure levels) $90\%$ confidence intervals for quantile Engel curves, $p=0.5$ and $p=0.9$: food (top left), fuel (top right), clothing (bottom left), and alcohol (bottom right).}
\end{figure}

\section{Additional simulation results and details}

\subsection{Unconditional simulations}

\begin{table}[htbp]
\caption{\label{tab:sim-un3}Coverage probability and median CI length, as in Table \ref{tab:sim-un1}.}
\centering
\begin{tabular}[c]{cccccccc}
 $n$  &   $p$   &       $F$       & Method & CP & Too low & Too high & Length \\
\hline
$ 25$ & $0.2  $ &          Normal & L-stat & 0.956 & 0.023 & 0.021 & 1.18 \\
$ 25$ & $0.2  $ &          Normal & BH     & 0.967 & 0.028 & 0.006 & 1.33 \\
$ 25$ & $0.2  $ &          Normal & Norm   & 0.944 & 0.011 & 0.045 & 1.13 \\
$ 25$ & $0.2  $ &          Normal & K15    & 0.952 & 0.006 & 0.042 & 1.42 \\
$ 25$ & $0.2  $ &          Normal & BStsym & 0.944 & 0.022 & 0.034 & 1.34 \\[2pt]
$ 25$ & $0.2  $ &          Cauchy & L-stat & 0.960 & 0.024 & 0.016 & 5.05 \\
$ 25$ & $0.2  $ &          Cauchy & BH     & 0.965 & 0.028 & 0.006 & 7.87 \\
$ 25$ & $0.2  $ &          Cauchy & Norm   & 0.909 & 0.014 & 0.077 & 2.43 \\
$ 25$ & $0.2  $ &          Cauchy & K15    & 0.946 & 0.001 & 0.053 & 4.00 \\
$ 25$ & $0.2  $ &          Cauchy & BStsym & 0.959 & 0.005 & 0.036 & 5.18 \\
$ 25$ & $0.2  $ &         Uniform & L-stat & 0.952 & 0.022 & 0.026 & 0.29 \\[2pt]
$ 25$ & $0.2  $ &         Uniform & BH     & 0.963 & 0.026 & 0.010 & 0.30 \\
$ 25$ & $0.2  $ &         Uniform & Norm   & 0.960 & 0.006 & 0.034 & 0.33 \\
$ 25$ & $0.2  $ &         Uniform & K15    & 0.953 & 0.015 & 0.032 & 0.39 \\
$ 25$ & $0.2  $ &         Uniform & BStsym & 0.930 & 0.042 & 0.028 & 0.35 \\
\hline
\hline
\end{tabular}
\end{table}

Table \ref{tab:sim-un3} shows a case with $p=0.2$ (away from the median) where BH can be computed.  The patterns among the normal, K15, and BStsym CIs are similar to before.  L-stat and BH both attain $95\%$ CP in each case, but L-stat is much closer to equal-tailed, and L-stat is shorter or the same length.  

\begin{table}[htbp]
\caption{\label{tab:sim-calib}Coverage probability and median CI length, as in Table \ref{tab:sim-un1}; $10^5$ replications.  ``L-stat'' uses equation \eqref{eqn:hutson-CI-2s}; ``Calib'' uses \eqref{eqn:hutson-CI-2s-calibrated}.}
\centering
\begin{tabular}[c]{cccccccc}
 $n$  &   $p$   &       $F$       & Method & CP & Too low & Too high & Length \\
\hline
$ 10$ & $0.35 $ &          Normal & L-stat & 0.959 & 0.023 & 0.018 & 1.73 \\
$ 10$ & $0.35 $ &          Normal & Calib  & 0.952 & 0.025 & 0.023 & 1.64 \\
$ 10$ & $0.40 $ &          Normal & L-stat & 0.965 & 0.020 & 0.015 & 1.69 \\
$ 10$ & $0.40 $ &          Normal & Calib  & 0.943 & 0.030 & 0.027 & 1.50 \\
$ 10$ & $0.45 $ &          Normal & L-stat & 0.953 & 0.023 & 0.024 & 1.62 \\
$ 10$ & $0.45 $ &          Normal & Calib  & 0.948 & 0.027 & 0.025 & 1.59 \\
$ 10$ & $0.50 $ &          Normal & L-stat & 0.962 & 0.019 & 0.019 & 1.64 \\
$ 10$ & $0.50 $ &          Normal & Calib  & 0.942 & 0.029 & 0.029 & 1.47 \\
\hline
\hline
\end{tabular}
\end{table}

Table \ref{tab:sim-calib} compares the original two-sided CI in equation \eqref{eqn:hutson-CI-2s} with the calibrated (``Calib'') version in \eqref{eqn:hutson-CI-2s-calibrated}.  Even with small $n$, the true differences are small, so we use a larger number of replications ($10^5$) to reduce simulation error.  By construction, the calibrated $\alpha$ is always larger than the original $\alpha$, which is why the Calib CI is always shorter and has smaller CP.  Similarly, while the dominant $n^{-1}$ term in the L-stat CPE is always positive (i.e.,\ over-coverage), the dominant term in the Calib CPE may be positive or negative, although the worst under-coverage in Table \ref{tab:sim-calib} is still $94.3\%$.  Since the calibration is separate for the upper and lower endpoints, it also makes the Calib CI (slightly) more equal-tailed.

Tables \ref{tab:sim-un-full1}--\ref{tab:sim-un-full5} show additional results for unconditional quantile CI simulations.

\begin{table}[htbp]
\caption{\label{tab:sim-un-full1}Coverage probability and median CI length, $1-\alpha=0.95$; various $n$, $p$, and distributions of $X_i$ ($F$) shown in table.  ``Too high'' is the proportion of simulations in which the lower endpoint was above the true value, $F^{-1}(p)$, and ``too low'' is the proportion when the upper endpoint was below $F^{-1}(p)$.}
\centering
\begin{tabular}[c]{cccccccc}
 $n$  &   $p$   &       $F$       & Method & CP & Too low & Too high & Length \\
\hline
$ 25$ & $0.2  $ &          Normal & L-stat & 0.956 & 0.023 & 0.021 & 1.18 \\
$ 25$ & $0.2  $ &          Normal & BH     & 0.967 & 0.028 & 0.006 & 1.33 \\
$ 25$ & $0.2  $ &          Normal & Norm   & 0.944 & 0.011 & 0.045 & 1.13 \\
$ 25$ & $0.2  $ &          Normal & K15    & 0.952 & 0.006 & 0.042 & 1.42 \\
$ 25$ & $0.2  $ &          Normal & BStsym & 0.944 & 0.022 & 0.034 & 1.34 \\
\hline
$ 25$ & $0.2  $ &          Cauchy & L-stat & 0.960 & 0.024 & 0.016 & 5.05 \\
$ 25$ & $0.2  $ &          Cauchy & BH     & 0.965 & 0.028 & 0.006 & 7.87 \\
$ 25$ & $0.2  $ &          Cauchy & Norm   & 0.909 & 0.014 & 0.077 & 2.43 \\
$ 25$ & $0.2  $ &          Cauchy & K15    & 0.946 & 0.001 & 0.053 & 4.00 \\
$ 25$ & $0.2  $ &          Cauchy & BStsym & 0.959 & 0.005 & 0.036 & 5.18 \\
\hline
$ 25$ & $0.2  $ &         Uniform & L-stat & 0.952 & 0.022 & 0.026 & 0.29 \\
$ 25$ & $0.2  $ &         Uniform & BH     & 0.963 & 0.026 & 0.010 & 0.30 \\
$ 25$ & $0.2  $ &         Uniform & Norm   & 0.960 & 0.006 & 0.034 & 0.33 \\
$ 25$ & $0.2  $ &         Uniform & K15    & 0.953 & 0.015 & 0.032 & 0.39 \\
$ 25$ & $0.2  $ &         Uniform & BStsym & 0.930 & 0.042 & 0.028 & 0.35 \\
\hline
$ 25$ & $0.2  $ &     Exponential & L-stat & 0.952 & 0.023 & 0.025 & 0.38 \\
$ 25$ & $0.2  $ &     Exponential & BH     & 0.964 & 0.027 & 0.009 & 0.38 \\
$ 25$ & $0.2  $ &     Exponential & Norm   & 0.989 & 0.001 & 0.010 & 0.61 \\
$ 25$ & $0.2  $ &     Exponential & K15    & 0.960 & 0.018 & 0.023 & 0.50 \\
$ 25$ & $0.2  $ &     Exponential & BStsym & 0.933 & 0.044 & 0.023 & 0.45 \\
\hline
$ 25$ & $0.2  $ &      Log-normal & L-stat & 0.954 & 0.023 & 0.023 & 0.49 \\
$ 25$ & $0.2  $ &      Log-normal & BH     & 0.966 & 0.027 & 0.006 & 0.50 \\
$ 25$ & $0.2  $ &      Log-normal & Norm   & 0.993 & 0.001 & 0.007 & 0.79 \\
$ 25$ & $0.2  $ &      Log-normal & K15    & 0.962 & 0.014 & 0.024 & 0.63 \\
$ 25$ & $0.2  $ &      Log-normal & BStsym & 0.938 & 0.039 & 0.023 & 0.55 \\
\hline
\hline
\end{tabular}
\end{table}

\begin{table}[htbp]
\caption{\label{tab:sim-un-full2}Coverage probability and median CI length, $1-\alpha=0.95$; various $n$, $p$, and distributions of $X_i$ ($F$) shown in table.  ``Too high'' is the proportion of simulations in which the lower endpoint was above the true value, $F^{-1}(p)$, and ``too low'' is the proportion when the upper endpoint was below $F^{-1}(p)$.}
\centering
\begin{tabular}[c]{cccccccc}
 $n$  &   $p$   &       $F$       & Method & CP & Too low & Too high & Length \\
\hline
$ 99$ & $0.05 $ &          Normal & L-stat & 0.960 & 0.024 & 0.016 & 0.92 \\
$ 99$ & $0.05 $ &          Normal & BH     &   NA &   NA &   NA &   NA \\
$ 99$ & $0.05 $ &          Normal & Norm   & 0.922 & 0.017 & 0.062 & 0.76 \\
$ 99$ & $0.05 $ &          Normal & K15    & 0.977 & 0.011 & 0.012 & 1.26 \\
$ 99$ & $0.05 $ &          Normal & BStsym & 0.947 & 0.023 & 0.031 & 1.05 \\
\hline
$ 99$ & $0.95 $ &          Normal & L-stat & 0.960 & 0.017 & 0.023 & 0.92 \\
$ 99$ & $0.95 $ &          Normal & BH     &   NA &   NA &   NA &   NA \\
$ 99$ & $0.95 $ &          Normal & Norm   & 0.921 & 0.061 & 0.017 & 0.76 \\
$ 99$ & $0.95 $ &          Normal & K15    & 0.975 & 0.013 & 0.012 & 1.28 \\
$ 99$ & $0.95 $ &          Normal & BStsym & 0.947 & 0.032 & 0.021 & 1.06 \\
\hline
\hline
\end{tabular}
\end{table}

\begin{table}[htbp]
\caption{\label{tab:sim-un-full3}Coverage probability and median CI length, $1-\alpha=0.95$; various $n$, $p$, and distributions of $X_i$ ($F$) shown in table.  ``Too high'' is the proportion of simulations in which the lower endpoint was above the true value, $F^{-1}(p)$, and ``too low'' is the proportion when the upper endpoint was below $F^{-1}(p)$.}
\centering
\begin{tabular}[c]{cccccccc}
 $n$  &   $p$   &       $F$       & Method & CP & Too low & Too high & Length \\
\hline
$ 99$ & $0.037$ &          Normal & L-stat & 0.951 & 0.023 & 0.026 & 1.02 \\
$ 99$ & $0.037$ &          Normal & BH     &   NA &   NA &   NA &   NA \\
$ 99$ & $0.037$ &          Normal & Norm   & 0.925 & 0.016 & 0.059 & 0.83 \\
$ 99$ & $0.037$ &          Normal & K15    & 0.970 & 0.009 & 0.021 & 1.55 \\
$ 99$ & $0.037$ &          Normal & BStsym & 0.950 & 0.020 & 0.030 & 1.20 \\
\hline
$ 99$ & $0.037$ &          Cauchy & L-stat & 0.950 & 0.022 & 0.028 & 39.37 \\
$ 99$ & $0.037$ &          Cauchy & BH     &   NA &   NA &   NA &   NA \\
$ 99$ & $0.037$ &          Cauchy & Norm   & 0.784 & 0.082 & 0.134 & 18.90 \\
$ 99$ & $0.037$ &          Cauchy & K15    & 0.957 & 0.002 & 0.041 & 36.55 \\
$ 99$ & $0.037$ &          Cauchy & BStsym & 0.961 & 0.002 & 0.037 & 48.77 \\
\hline
$ 99$ & $0.037$ &         Uniform & L-stat & 0.951 & 0.024 & 0.026 & 0.07 \\
$ 99$ & $0.037$ &         Uniform & BH     &   NA &   NA &   NA &   NA \\
$ 99$ & $0.037$ &         Uniform & Norm   & 0.990 & 0.000 & 0.010 & 0.12 \\
$ 99$ & $0.037$ &         Uniform & K15    & 0.963 & 0.028 & 0.009 & 0.11 \\
$ 99$ & $0.037$ &         Uniform & BStsym & 0.924 & 0.053 & 0.022 & 0.08 \\
\hline
$ 99$ & $0.037$ &     Exponential & L-stat & 0.956 & 0.020 & 0.024 & 0.07 \\
$ 99$ & $0.037$ &     Exponential & BH     &   NA &   NA &   NA &   NA \\
$ 99$ & $0.037$ &     Exponential & Norm   & 0.999 & 0.000 & 0.002 & 0.17 \\
$ 99$ & $0.037$ &     Exponential & K15    & 0.965 & 0.026 & 0.009 & 0.12 \\
$ 99$ & $0.037$ &     Exponential & BStsym & 0.932 & 0.046 & 0.022 & 0.09 \\
\hline
$ 99$ & $0.037$ &      Log-normal & L-stat & 0.951 & 0.022 & 0.026 & 0.15 \\
$ 99$ & $0.037$ &      Log-normal & BH     &   NA &   NA &   NA &   NA \\
$ 99$ & $0.037$ &      Log-normal & Norm   & 0.994 & 0.000 & 0.006 & 0.23 \\
$ 99$ & $0.037$ &      Log-normal & K15    & 0.972 & 0.015 & 0.013 & 0.24 \\
$ 99$ & $0.037$ &      Log-normal & BStsym & 0.938 & 0.036 & 0.026 & 0.18 \\
\hline
\hline
\end{tabular}
\end{table}

\begin{table}[htbp]
\caption{\label{tab:sim-un-full4}Coverage probability and median CI length, $1-\alpha=0.95$; various $n$, $p$, and distributions of $X_i$ ($F$) shown in table.  ``Too high'' is the proportion of simulations in which the lower endpoint was above the true value, $F^{-1}(p)$, and ``too low'' is the proportion when the upper endpoint was below $F^{-1}(p)$.}
\centering
\begin{tabular}[c]{cccccccc}
 $n$  &   $p$   &       $F$       & Method & CP & Too low & Too high & Length \\
\hline
$  9$ & $0.5  $ &          Normal & L-stat & 0.956 & 0.021 & 0.024 & 1.72 \\
$  9$ & $0.5  $ &          Normal & BH     & 0.959 & 0.019 & 0.022 & 1.76 \\
$  9$ & $0.5  $ &          Normal & Norm   & 0.908 & 0.047 & 0.045 & 1.58 \\
$  9$ & $0.5  $ &          Normal & K15    & 0.983 & 0.008 & 0.008 & 2.35 \\
$  9$ & $0.5  $ &          Normal & BStsym & 0.947 & 0.026 & 0.026 & 2.00 \\
\hline
$  9$ & $0.5  $ &          Cauchy & L-stat & 0.954 & 0.021 & 0.024 & 3.49 \\
$  9$ & $0.5  $ &          Cauchy & BH     & 0.957 & 0.020 & 0.023 & 3.59 \\
$  9$ & $0.5  $ &          Cauchy & Norm   & 0.958 & 0.022 & 0.020 & 2.66 \\
$  9$ & $0.5  $ &          Cauchy & K15    & 0.995 & 0.003 & 0.002 & 4.82 \\
$  9$ & $0.5  $ &          Cauchy & BStsym & 0.972 & 0.014 & 0.015 & 3.40 \\
\hline
$  9$ & $0.5  $ &         Uniform & L-stat & 0.957 & 0.021 & 0.022 & 0.57 \\
$  9$ & $0.5  $ &         Uniform & BH     & 0.960 & 0.019 & 0.021 & 0.59 \\
$  9$ & $0.5  $ &         Uniform & Norm   & 0.861 & 0.068 & 0.071 & 0.53 \\
$  9$ & $0.5  $ &         Uniform & K15    & 0.971 & 0.014 & 0.015 & 0.78 \\
$  9$ & $0.5  $ &         Uniform & BStsym & 0.936 & 0.031 & 0.033 & 0.74 \\
\hline
$  9$ & $0.5  $ &     Exponential & L-stat & 0.953 & 0.022 & 0.025 & 1.41 \\
$  9$ & $0.5  $ &     Exponential & BH     & 0.956 & 0.021 & 0.024 & 1.44 \\
$  9$ & $0.5  $ &     Exponential & Norm   & 0.889 & 0.076 & 0.035 & 1.17 \\
$  9$ & $0.5  $ &     Exponential & K15    & 0.976 & 0.018 & 0.006 & 1.92 \\
$  9$ & $0.5  $ &     Exponential & BStsym & 0.946 & 0.037 & 0.017 & 1.71 \\
\hline
$  9$ & $0.5  $ &      Log-normal & L-stat & 0.955 & 0.020 & 0.024 & 1.91 \\
$  9$ & $0.5  $ &      Log-normal & BH     & 0.958 & 0.019 & 0.023 & 1.95 \\
$  9$ & $0.5  $ &      Log-normal & Norm   & 0.898 & 0.078 & 0.024 & 1.53 \\
$  9$ & $0.5  $ &      Log-normal & K15    & 0.981 & 0.016 & 0.003 & 2.62 \\
$  9$ & $0.5  $ &      Log-normal & BStsym & 0.950 & 0.037 & 0.013 & 2.29 \\
\hline
\hline
\end{tabular}
\end{table}

\begin{table}[htbp]
\caption{\label{tab:sim-un-full5}Coverage probability and median CI length, $1-\alpha=0.95$; various $n$, $p$, and distributions of $X_i$ ($F$) shown in table.  ``Too high'' is the proportion of simulations in which the lower endpoint was above the true value, $F^{-1}(p)$, and ``too low'' is the proportion when the upper endpoint was below $F^{-1}(p)$.}
\centering
\begin{tabular}[c]{cccccccc}
 $n$  &   $p$   &       $F$       & Method & CP & Too low & Too high & Length \\
\hline
$ 25$ & $0.5  $ &          Normal & L-stat & 0.953 & 0.022 & 0.025 & 0.99 \\
$ 25$ & $0.5  $ &          Normal & BH     & 0.955 & 0.021 & 0.024 & 1.00 \\
$ 25$ & $0.5  $ &          Normal & Norm   & 0.942 & 0.028 & 0.030 & 1.02 \\
$ 25$ & $0.5  $ &          Normal & K15    & 0.971 & 0.014 & 0.015 & 1.19 \\
$ 25$ & $0.5  $ &          Normal & BStsym & 0.942 & 0.028 & 0.030 & 1.13 \\
\hline
$ 25$ & $0.5  $ &          Cauchy & L-stat & 0.951 & 0.024 & 0.024 & 1.44 \\
$ 25$ & $0.5  $ &          Cauchy & BH     & 0.953 & 0.023 & 0.024 & 1.46 \\
$ 25$ & $0.5  $ &          Cauchy & Norm   & 0.979 & 0.010 & 0.011 & 1.59 \\
$ 25$ & $0.5  $ &          Cauchy & K15    & 0.988 & 0.006 & 0.006 & 1.83 \\
$ 25$ & $0.5  $ &          Cauchy & BStsym & 0.951 & 0.024 & 0.025 & 1.47 \\
\hline
$ 25$ & $0.5  $ &         Uniform & L-stat & 0.953 & 0.022 & 0.025 & 0.37 \\
$ 25$ & $0.5  $ &         Uniform & BH     & 0.954 & 0.021 & 0.025 & 0.37 \\
$ 25$ & $0.5  $ &         Uniform & Norm   & 0.908 & 0.046 & 0.046 & 0.35 \\
$ 25$ & $0.5  $ &         Uniform & K15    & 0.963 & 0.018 & 0.020 & 0.44 \\
$ 25$ & $0.5  $ &         Uniform & BStsym & 0.937 & 0.031 & 0.032 & 0.45 \\
\hline
$ 25$ & $0.5  $ &     Exponential & L-stat & 0.953 & 0.024 & 0.023 & 0.79 \\
$ 25$ & $0.5  $ &     Exponential & BH     & 0.954 & 0.024 & 0.022 & 0.80 \\
$ 25$ & $0.5  $ &     Exponential & Norm   & 0.924 & 0.056 & 0.020 & 0.75 \\
$ 25$ & $0.5  $ &     Exponential & K15    & 0.968 & 0.022 & 0.010 & 0.96 \\
$ 25$ & $0.5  $ &     Exponential & BStsym & 0.941 & 0.039 & 0.020 & 0.93 \\
\hline
\hline
\end{tabular}
\end{table}

\subsection{Conditional simulations: implementation details}

For L-stat, for estimating $f_X(x_0)$ in order to calculate the plug-in bandwidth, any kernel density estimator will suffice. We use \texttt{kde} from package \texttt{ks} \citep{R.ks}, with the Gaussian-based bandwidth.  From the same package, \texttt{kdde} estimates $f_X'(x_0)$.  Both functions work for up to $d=6$-dimensional $X$. 
The objects $F_{Y|X}^{(0,1)}(\xi_p;x_0)$ and $F_{Y|X}^{(0,2)}(\xi_p;x_0)$ are estimated using local cubic (mean) regression of $1\{Y_i\le\xi_p\}$ on $X$. 

The first other method (``rqss'') is directly available in the popular \texttt{quantreg} package in R \citep{R.quantreg}.  The regression quantile smoothing spline function \texttt{rqss} is used as on page 10 of its vignette, with the \citet{Schwarz1978} information criterion (SIC) for model selection; pointwise CIs come from \texttt{predict.rqss}, and uniform bands from a \citet{Hotelling1939} tube approach coded in \texttt{plot.rqss}.  
 
The second method is a local polynomial following \citet{Chaudhuri1991} but with bootstrapped standard errors (``boot''). 
For this method, we adjust our method's plug-in bandwidth by a factor of $n^{1/12}$ to yield the local cubic CPE-optimal bandwidth rate. 
The function \texttt{rq} performs the local cubic point estimation, and \texttt{summary.rq} generates the bootstrap standard errors; we use $299$ replications.  
  
The third method uses the asymptotic normality of a local linear estimator, using results and ideas from \citet{QuYoon2015}, although they are concerned more with uniform (over quantiles) than pointwise inference.  They suggest using the MSE-optimal bandwidth (Corollary 1) and a conservative type of bias correction for the CIs (Remark 7) that increases the length of the CI.  We tried a uniform kernel, but show only results for the version with a Gaussian kernel (``QYg'') since it was usually better. 

The fourth method is from Section 3.1 in \citet{FanLiu2016}, based on a symmetrized $k$-NN estimator using a bisquare kernel (``FLb'').  We use the same code from their simulations (graciously provided to us).\footnote{This differs somewhat from the description in the text, most notably by an additional factor of $0.4$ in the bandwidth.}  Interestingly, although in principle they are just slightly undersmoothing the MSE-optimal bandwidth, the result is very close to the CPE-optimal bandwidth for the sample sizes considered. 

Figure \ref{fig:true-fn} shows the conditional median function from \eqref{eqn:DGP-rqss}. 

\begin{figure}[thbp]
 \centering
 \includegraphics[clip=true,trim=0 30 30 50,width=0.7\textwidth]
	{2014_07_23_Figures/F010_true_fn.pdf}
 \caption{\label{fig:true-fn}Conditional median function from \eqref{eqn:DGP-rqss}.}
\end{figure}

\subsection{Conditional simulations: additional results}

Figure \ref{fig:unifPWR} shows uniform power curves, as discussed in the main text.  These were evaluated at $231$ different $x_0$ values. 
\begin{figure}[thbp]
  \centering
  \includegraphics[clip=true,trim=15 45 10 58,width=0.445\textwidth]
    {2015_08_27_Figures/qinfnp_2015_08_27_unifPWR_F010_a05_n400_p50_numx231_trimx04_nrep1000.pdf}
  \includegraphics[clip=true,trim=51 45 10 58,width=0.406\textwidth]
    {2015_08_27_Figures/qinfnp_2015_08_27_unifPWR_F012_a05_n400_p50_numx231_trimx04_nrep1000.pdf} \\
  \includegraphics[clip=true,trim=15 45 10 58,width=0.445\textwidth]
    {2015_08_27_Figures/qinfnp_2015_08_27_unifPWR_F014_a05_n400_p50_numx231_trimx04_nrep1000.pdf}
  \includegraphics[clip=true,trim=51 45 10 58,width=0.406\textwidth]
    {2015_08_27_Figures/qinfnp_2015_08_27_unifPWR_F016_a05_n400_p50_numx231_trimx04_nrep1000.pdf} 
  \caption{\label{fig:unifPWR}Uniform power curves, analogous to column 3 of Figure \ref{fig:ptCP-1}, with normal (top left), $t_3$ (top right), Cauchy (bottom left), and centered $\chi^2_3$ (bottom right) $U_i$.}
\end{figure}

Figure \ref{fig:ptPWR2-2} is the same as Figure \ref{fig:ptPWR2-1} but with heteroskedasticity. 
\begin{figure}[thbp]
  \centering 
  \includegraphics[clip=true,trim=15 15 10 58,width=0.445\textwidth]
    {2015_08_16_Figures/qinfnp_2015_08_16_ptPWR2_F011_a05_n400_p50_numx47_trimx04_nrep1000.pdf}
  \includegraphics[clip=true,trim=51 15 10 58,width=0.406\textwidth]
    {2015_08_16_Figures/qinfnp_2015_08_16_ptPWR2_F013_a05_n400_p50_numx47_trimx04_nrep1000.pdf}
\\
  \includegraphics[clip=true,trim=15 15 10 58,width=0.445\textwidth]
    {2015_08_16_Figures/qinfnp_2015_08_16_ptPWR2_F015_a05_n400_p50_numx47_trimx04_nrep1000.pdf}
  \includegraphics[clip=true,trim=51 15 10 58,width=0.406\textwidth]
    {2015_08_16_Figures/qinfnp_2015_08_16_ptPWR2_F017_a05_n400_p50_numx47_trimx04_nrep1000.pdf}
  \caption{\label{fig:ptPWR2-2}Pointwise power by $X$; same as Figure \ref{fig:ptPWR2-1} but with $\sigma(x)=0.2(1+x)$.}
\end{figure}

\begin{figure}[htbp]
  \centering
  \includegraphics[clip=true,trim=15 15 30 25,width=0.32\textwidth]
    {2015_08_26_Figures/qinfnp_2015_08_26_ptCP_F010_a05_n400_p25_numx47_trimx04_nrep1000.pdf}
  \includegraphics[clip=true,trim=15 15 30 25,width=0.32\textwidth]
    {2015_08_26_Figures/qinfnp_2015_08_26_ptPWR2_F010_a05_n400_p25_numx47_trimx04_nrep1000.pdf}
  \includegraphics[clip=true,trim=15 15 30 25,width=0.32\textwidth]
    {2015_08_26_Figures/qinfnp_2015_08_26_jtPWR_F010_a05_n400_p25_numx47_trimx04_nrep1000.pdf}
  \caption{\label{fig:p25}Same as Figures \ref{fig:ptCP-1} and \ref{fig:ptPWR2-1} but with $p=0.25$ and homoskedastic normal $U_i$.} 
\end{figure}

Figure \ref{fig:p25} shows that L-stat continues to perform well even with $p=0.25$ instead of $p=0.5$ (still with $n=400$, $\alpha=0.05$).  We show this in Figure \ref{fig:p25} with normal $U_i$ centered to have $P(U_i<0)=p$. 
Pointwise CP (left panel) is similar to $p=0.5$, except FLb is somewhat worse.  
At $x_0$ where all methods have good coverage, pointwise power (middle panel) is generally highest for L-stat, although QYg is highest at points (but near zero at others).  
The L-stat joint CIs are the only ones to achieve correct coverage, as seen in the right panel; next-closest is boot with 86.2\% CP.  The joint L-stat test is also least biased.  Due to its equal-tailed nature, the L-stat power curve is less steep for negative deviations (where the data are sparser) but still roughly parallel to QYg and boot, and it is more steep for positive deviations.  

\begin{figure}[htbp]
  \centering 
  \includegraphics[clip=true,trim=15 15 30 25,width=0.32\textwidth]
    {2015_08_26_Figures/qinfnp_2015_08_26_ptCP_F014_a05_n400_p25_numx47_trimx04_nrep1000.pdf}
  \includegraphics[clip=true,trim=15 15 30 25,width=0.32\textwidth]
    {2015_08_26_Figures/qinfnp_2015_08_26_ptPWR2_F014_a05_n400_p25_numx47_trimx04_nrep1000.pdf}
  \includegraphics[clip=true,trim=15 15 30 25,width=0.32\textwidth]
    {2015_08_26_Figures/qinfnp_2015_08_26_jtPWR_F014_a05_n400_p25_numx47_trimx04_nrep1000.pdf}
  \caption{\label{fig:p25-Cauchy}Same as Figures \ref{fig:ptCP-1} and \ref{fig:ptPWR2-1} but with $p=0.25$ and homoskedastic Cauchy $U_i$.}
\end{figure}
Figure \ref{fig:p25-Cauchy} shows $p=0.25$ and (recentered) Cauchy errors.  Only boot maintains correct pointwise and joint CP.  L-stat has some under-coverage, but its CP only falls below $93\%$ at 10 of 47 points, and the lowest CP is still above $88\%$.  In contrast, rqss, FLb, and QYg have severe under-coverage at many $x_0$. 
The cost of boot's conservative coverage is pointwise power (middle panel) near 5\% everywhere, compared with around 20\% for L-stat.  
For joint testing (right panel), L-stat and boot share some slight (below $10\%$ rejection rates) size distortion and bias, while other methods are quite size distorted and/or biased.  In terms of power, L-stat is better than boot for positive deviations and vice-versa for negative.  However, the difference is more dramatic for positive deviations.  
Although L-stat has only 38\% power against the $-0.2$ deviation where boot has 80\% power, boot has almost trivial 8.5\% power against the $+0.15$ deviation where L-stat has 95.6\% power.


\begin{table}[htbp]
\caption{\label{tab:comp-time}Computation time in seconds (rounded), including bandwidth selection.  $X_i$ and $Y_i$ both iid $\textrm{Unif}(0,1)$, $X_i\independent Y_i$, $p=0.5$; $n$ and number of $x_0$ (spaced between lower and upper quartiles of $X$) shown in table.  Run on 3.2GHz Intel i5 processor with 8GB RAM.}
\centering
\begin{tabular}[c]{lrrrrrrr}
 && \multicolumn{6}{c}{sample size, $n$} \\
\cline{3-8}
Method & \multicolumn{1}{c}{\#$x_0$} & \multicolumn{1}{c}{$400$} & \multicolumn{1}{c}{$1000$} & \multicolumn{1}{c}{$4000$} 
                 & \multicolumn{1}{c}{$10^4$} & \multicolumn{1}{c}{$10^5$} & \multicolumn{1}{c}{$10^6$}  \\
\hline
L-stat      &   10 &    0 &    0 &    0 &    0 &    2 &    23 \\
Local cubic &   10 &    0 &    0 &    1 &    1 &   11 &   276 \\
rqss        &   10 &    1 &    5 &  191 &  n/a &  n/a &   n/a \\
L-stat      &  100 &    1 &    1 &    1 &    2 &   10 &   110 \\
Local cubic &  100 &    1 &    2 &    4 &    9 &  141 &  4501 \\
rqss        &  100 &    1 &    5 &  189 &  n/a &  n/a &   n/a \\
L-stat      & 1000 &    6 &    6 &    9 &   14 &   86 &   969 \\
Local cubic & 1000 &   11 &   15 &   37 &   80 & 1169 & 31994 \\
rqss        & 1000 &    1 &    5 &  191 &  n/a &  n/a &   n/a \\
\hline
\end{tabular}
\end{table}

Our method also has a computational advantage.  While rqss is the fastest with $n=400$ and scales best with the number of $x_0$, it slows to taking a few minutes with $n=4000$.  On the machine used, there was not enough memory to compute \texttt{rqss} for $n\ge10^4$ (tried to allocate $9.9$Gb, but only 8 available).  
The local cubic bootstrap scales worse with the number of $x_0$ than \texttt{rqss}, but scales better with $n$.  By only considering local samples of size $N_n$, it also scales much better than bootstrapping a global estimator using all $n$ observations.  Still, the reliance on resampling makes it slower than L-stat.  For example, to examine $100$ values of $x_0$ when $n=10^5$, the local cubic bootstrap takes over two minutes, whereas L-stat takes under 10 seconds.  See Table \ref{tab:comp-time} for other examples. 

}{}

\nocite{Kumaraswamy1980,Robbins1955,Jones2002,BhattacharyaGangopadhyay1990,Hutson1999,FanEtAl2012,Chaudhuri1991,Pratt1968,PeizerPratt1968,Muir1960,ShorackWellner1986,DasGupta2000,Bickel1967,Wilks1962,KaplanSun2016}
\Supplemental{\pagebreak\singlespacing\bibliographystyle{chicago}\bibliography{GoldmanKaplan}}{}

\begin{thebibliography}{}

\bibitem[\protect\citeauthoryear{Abrevaya}{Abrevaya}{2001}]{Abrevaya2001}
Abrevaya, J. (2001).
\newblock The effects of demographics and maternal behavior on the distribution
  of birth outcomes.
\newblock {\em Empirical Economics\/}~{\em 26\/}(1), 247--257.

\bibitem[\protect\citeauthoryear{Alan, Crossley, Grootendorst, and Veall}{Alan
  et~al.}{2005}]{AlanEtAl2005}
Alan, S., T.~F. Crossley, P.~Grootendorst, and M.~R. Veall (2005).
\newblock Distributional effects of `general population' prescription drug
  programs in {Canada}.
\newblock {\em Canadian Journal of Economics\/}~{\em 38\/}(1), 128--148.

\bibitem[\protect\citeauthoryear{Andrews and Guggenberger}{Andrews and
  Guggenberger}{2014}]{AndrewsGuggenberger2014}
Andrews, D. W.~K. and P.~Guggenberger (2014).
\newblock A conditional-heteroskedasticity-robust confidence interval for the
  autoregressive parameter.
\newblock {\em Review of Economics and Statistics\/}~{\em 96\/}(2), 376--381.

\bibitem[\protect\citeauthoryear{Banks, Blundell, and Lewbel}{Banks
  et~al.}{1997}]{BanksEtAl1997}
Banks, J., R.~Blundell, and A.~Lewbel (1997).
\newblock Quadratic {Engel} curves and consumer demand.
\newblock {\em Review of Economics and Statistics\/}~{\em 79\/}(4), 527--539.

\bibitem[\protect\citeauthoryear{Beran and Hall}{Beran and
  Hall}{1993}]{BeranHall1993}
Beran, R. and P.~Hall (1993).
\newblock Interpolated nonparametric prediction intervals and confidence
  intervals.
\newblock {\em Journal of the Royal Statistical Society: Series B (Statistical
  Methodology)\/}~{\em 55\/}(3), 643--652.

\bibitem[\protect\citeauthoryear{Bhattacharya and Gangopadhyay}{Bhattacharya
  and Gangopadhyay}{1990}]{BhattacharyaGangopadhyay1990}
Bhattacharya, P.~K. and A.~K. Gangopadhyay (1990).
\newblock Kernel and nearest-neighbor estimation of a conditional quantile.
\newblock {\em Annals of Statistics\/}~{\em 18\/}(3), 1400--1415.

\bibitem[\protect\citeauthoryear{Bickel}{Bickel}{1967}]{Bickel1967}
Bickel, P.~J. (1967).
\newblock Some contributions to the theory of order statistics.
\newblock In {\em Proceedings of the Fifth Berkeley Symposium on Mathematical
  Statistics and Probability, Volume 1: Statistics}. The Regents of the
  University of California.

\bibitem[\protect\citeauthoryear{Buchinsky}{Buchinsky}{1994}]{Buchinsky1994}
Buchinsky, M. (1994).
\newblock Changes in the {U.S.} wage structure 1963--1987: Application of
  quantile regression.
\newblock {\em Econometrica\/}~{\em 62\/}(2), 405--458.

\bibitem[\protect\citeauthoryear{Chamberlain}{Chamberlain}{1994}]{Chamberlain1994}
Chamberlain, G. (1994).
\newblock Quantile regression, censoring, and the structure of wages.
\newblock In {\em Advances in Econometrics: Sixth World Congress}, Volume~2,
  pp.\  171--209.

\bibitem[\protect\citeauthoryear{Chaudhuri}{Chaudhuri}{1991}]{Chaudhuri1991}
Chaudhuri, P. (1991).
\newblock Nonparametric estimates of regression quantiles and their local
  {Bahadur} representation.
\newblock {\em Annals of Statistics\/}~{\em 19\/}(2), 760--777.

\bibitem[\protect\citeauthoryear{Chen and Hall}{Chen and
  Hall}{1993}]{ChenHall1993}
Chen, S.~X. and P.~Hall (1993).
\newblock Smoothed empirical likelihood confidence intervals for quantiles.
\newblock {\em Annals of Statistics\/}~{\em 21\/}(3), 1166--1181.

\bibitem[\protect\citeauthoryear{DasGupta}{DasGupta}{2000}]{DasGupta2000}
DasGupta, A. (2000).
\newblock Best constants in {Chebyshev} inequalities with various applications.
\newblock {\em Metrika\/}~{\em 51\/}(3), 185--200.

\bibitem[\protect\citeauthoryear{David and Nagaraja}{David and
  Nagaraja}{2003}]{DavidNagaraja2003}
David, H.~A. and H.~N. Nagaraja (2003).
\newblock {\em Order Statistics\/} (3rd ed.).
\newblock New York: Wiley.

\bibitem[\protect\citeauthoryear{Deaton}{Deaton}{1997}]{Deaton1997}
Deaton, A. (1997).
\newblock {\em The analysis of household surveys: a microeconometric approach
  to development policy}.
\newblock Baltimore: The Johns Hopkins University Press.

\bibitem[\protect\citeauthoryear{Donald, Hsu, and Barrett}{Donald
  et~al.}{2012}]{DonaldEtAl2012}
Donald, S.~G., Y.-C. Hsu, and G.~F. Barrett (2012).
\newblock Incorporating covariates in the measurement of welfare and
  inequality: methods and applications.
\newblock {\em The Econometrics Journal\/}~{\em 15\/}(1), C1--C30.

\bibitem[\protect\citeauthoryear{Engel}{Engel}{1857}]{Engel1857}
Engel, E. (1857).
\newblock Die productions- und consumtionsverh\"altnisse des k\"onigreichs
  sachsen.
\newblock {\em Zeitschrift des Statistischen Bureaus des K\"oniglich
  S\"achsischen, Ministerium des Inneren\/}~{\em 8--9}, 1--54.

\bibitem[\protect\citeauthoryear{Fan, Grama, and Liu}{Fan
  et~al.}{2012}]{FanEtAl2012}
Fan, X., I.~Grama, and Q.~Liu (2012).
\newblock Hoeffding's inequality for supermartingales.
\newblock {\em Stochastic Processes and their Applications\/}~{\em 122\/}(10),
  3545--3559.

\bibitem[\protect\citeauthoryear{Fan and Liu}{Fan and Liu}{2016}]{FanLiu2016}
Fan, Y. and R.~Liu (2016).
\newblock A direct approach to inference in nonparametric and semiparametric
  quantile models.
\newblock {\em Journal of Econometrics\/}~{\em 191\/}(1), 196--216.

\bibitem[\protect\citeauthoryear{Ferguson}{Ferguson}{1973}]{Ferguson1973}
Ferguson, T.~S. (1973).
\newblock A {Bayesian} analysis of some nonparametric problems.
\newblock {\em Annals of Statistics\/}~{\em 1\/}(2), 209--230.

\bibitem[\protect\citeauthoryear{Fisher}{Fisher}{1932}]{Fisher1932}
Fisher, R.~A. (1932).
\newblock {\em Statistical Methods for Research Workers\/} (4th ed.).
\newblock Edinburg: Oliver and Boyd.

\bibitem[\protect\citeauthoryear{Goldman and Kaplan}{Goldman and
  Kaplan}{2016a}]{GoldmanKaplan2015c}
Goldman, M. and D.~M. Kaplan (2016a).
\newblock Evenly sensitive {KS-type} inference on distributions.
\newblock Working paper, available at
  \url{http://faculty.missouri.edu/~kaplandm}.

\bibitem[\protect\citeauthoryear{Goldman and Kaplan}{Goldman and
  Kaplan}{2016b}]{GoldmanKaplan2014b}
Goldman, M. and D.~M. Kaplan (2016b).
\newblock Nonparametric inference on conditional quantile differences, linear
  combinations, and vectors, using {$L$-statistics}.
\newblock Working paper, available at
  \url{http://faculty.missouri.edu/~kaplandm}.

\bibitem[\protect\citeauthoryear{Hall and Sheather}{Hall and
  Sheather}{1988}]{HallSheather1988}
Hall, P. and S.~J. Sheather (1988).
\newblock On the distribution of a {Studentized} quantile.
\newblock {\em Journal of the Royal Statistical Society: Series B (Statistical
  Methodology)\/}~{\em 50\/}(3), 381--391.

\bibitem[\protect\citeauthoryear{Ho and Lee}{Ho and Lee}{2005a}]{HoLee2005a}
Ho, Y. H.~S. and S.~M.~S. Lee (2005a).
\newblock Calibrated interpolated confidence intervals for population
  quantiles.
\newblock {\em Biometrika\/}~{\em 92\/}(1), 234--241.

\bibitem[\protect\citeauthoryear{Ho and Lee}{Ho and Lee}{2005b}]{HoLee2005b}
Ho, Y. H.~S. and S.~M.~S. Lee (2005b).
\newblock Iterated smoothed bootstrap confidence intervals for population
  quantiles.
\newblock {\em Annals of Statistics\/}~{\em 33\/}(1), 437--462.

\bibitem[\protect\citeauthoryear{Hogg}{Hogg}{1975}]{Hogg1975}
Hogg, R. (1975).
\newblock Estimates of percentile regression lines using salary data.
\newblock {\em Journal of the American Statistical Association\/}~{\em
  70\/}(349), 56--59.

\bibitem[\protect\citeauthoryear{Horowitz and Lee}{Horowitz and
  Lee}{2012}]{HorowitzLee2012}
Horowitz, J.~L. and S.~Lee (2012).
\newblock Uniform confidence bands for functions estimated nonparametrically
  with instrumental variables.
\newblock {\em Journal of Econometrics\/}~{\em 168\/}(2), 175--188.

\bibitem[\protect\citeauthoryear{Hotelling}{Hotelling}{1939}]{Hotelling1939}
Hotelling, H. (1939).
\newblock Tubes and spheres in $n$-space and a class of statistical problems.
\newblock {\em American Journal of Mathematics\/}~{\em 61}, 440--460.

\bibitem[\protect\citeauthoryear{Hutson}{Hutson}{1999}]{Hutson1999}
Hutson, A.~D. (1999).
\newblock Calculating nonparametric confidence intervals for quantiles using
  fractional order statistics.
\newblock {\em Journal of Applied Statistics\/}~{\em 26\/}(3), 343--353.

\bibitem[\protect\citeauthoryear{Jones}{Jones}{2002}]{Jones2002}
Jones, M.~C. (2002).
\newblock On fractional uniform order statistics.
\newblock {\em Statistics \& Probability Letters\/}~{\em 58\/}(1), 93--96.

\bibitem[\protect\citeauthoryear{Kaplan}{Kaplan}{2014}]{Kaplan2014}
Kaplan, D.~M. (2014).
\newblock Nonparametric inference on quantile marginal effects.
\newblock Working paper, available at
  \url{http://faculty.missouri.edu/~kaplandm}.

\bibitem[\protect\citeauthoryear{Kaplan}{Kaplan}{2015}]{Kaplan2015}
Kaplan, D.~M. (2015).
\newblock Improved quantile inference via fixed-smoothing asymptotics and
  {Edgeworth} expansion.
\newblock {\em Journal of Econometrics\/}~{\em 185\/}(1), 20--32.

\bibitem[\protect\citeauthoryear{Kaplan and Sun}{Kaplan and
  Sun}{2016}]{KaplanSun2016}
Kaplan, D.~M. and Y.~Sun (2016).
\newblock Smoothed estimating equations for instrumental variables quantile
  regression.
\newblock {\em Econometric Theory\/}~{\em XX\/}(XX), XX--XX.
\newblock Forthcoming.

\bibitem[\protect\citeauthoryear{Koenker}{Koenker}{2012}]{R.quantreg}
Koenker, R. (2012).
\newblock {\em quantreg: Quantile Regression}.
\newblock R package version 4.81.

\bibitem[\protect\citeauthoryear{Kumaraswamy}{Kumaraswamy}{1980}]{Kumaraswamy1980}
Kumaraswamy, P. (1980).
\newblock A generalized probability density function for double-bounded random
  processes.
\newblock {\em Journal of Hydrology\/}~{\em 46\/}(1--2), 79--88.

\bibitem[\protect\citeauthoryear{Lehmann}{Lehmann}{1951}]{Lehmann1951}
Lehmann, E.~L. (1951).
\newblock A general concept of unbiasedness.
\newblock {\em Annals of Mathematical Statistics\/}~{\em 22\/}(4), 587--592.

\bibitem[\protect\citeauthoryear{Li and Racine}{Li and
  Racine}{2007}]{LiRacine2007}
Li, Q. and J.~S. Racine (2007).
\newblock {\em Nonparametric econometrics: Theory and practice}.
\newblock Princeton University Press.

\bibitem[\protect\citeauthoryear{Manning, Blumberg, and Moulton}{Manning
  et~al.}{1995}]{ManningEtAl1995}
Manning, W., L.~Blumberg, and L.~Moulton (1995).
\newblock The demand for alcohol: the differential response to price.
\newblock {\em Journal of Health Economics\/}~{\em 14\/}(2), 123--148.

\bibitem[\protect\citeauthoryear{Muir}{Muir}{1960}]{Muir1960}
Muir, T. (1960).
\newblock {\em A Treatise on the Theory of Determinants}.
\newblock Dover Publications.

\bibitem[\protect\citeauthoryear{Neyman}{Neyman}{1937}]{Neyman1937}
Neyman, J. (1937).
\newblock {\guillemotright}{Smooth} test{\guillemotright} for goodness of fit.
\newblock {\em Skandinavisk Aktuarietidskrift\/}~{\em 20\/}(3--4), 149--199.

\bibitem[\protect\citeauthoryear{{Office for National Statistics and Department
  for Environment, Food and Rural Affairs}}{{Office for National Statistics and
  Department for Environment, Food and Rural Affairs}}{2012}]{UKLCFS2012}
{Office for National Statistics and Department for Environment, Food and Rural
  Affairs} (2012).
\newblock {Living Costs and Food Survey}.
\newblock {2nd Edition. Colchester, Essex: UK Data Archive.
  \url{http://dx.doi.org/10.5255/UKDA-SN-7472-2}}.

\bibitem[\protect\citeauthoryear{Pearson}{Pearson}{1933}]{Pearson1933}
Pearson, K. (1933).
\newblock On a method of determining whether a sample of size n supposed to
  have been drawn from a parent population having a known probability integral
  has probably been drawn at random.
\newblock {\em Biometrika\/}~{\em 25}, 379--410.

\bibitem[\protect\citeauthoryear{Peizer and Pratt}{Peizer and
  Pratt}{1968}]{PeizerPratt1968}
Peizer, D.~B. and J.~W. Pratt (1968).
\newblock A normal approximation for binomial, {$F$}, beta, and other common,
  related tail probabilities, {I}.
\newblock {\em Journal of the American Statistical Association\/}~{\em
  63\/}(324), 1416--1456.

\bibitem[\protect\citeauthoryear{Polansky and Schucany}{Polansky and
  Schucany}{1997}]{PolanskySchucany1997}
Polansky, A.~M. and W.~R. Schucany (1997).
\newblock Kernel smoothing to improve bootstrap confidence intervals.
\newblock {\em Journal of the Royal Statistical Society: Series B (Statistical
  Methodology)\/}~{\em 59\/}(4), 821--838.

\bibitem[\protect\citeauthoryear{Polonik and Yao}{Polonik and
  Yao}{2002}]{PolonikYao2002}
Polonik, W. and Q.~Yao (2002).
\newblock Set-indexed conditional empirical and quantile processes based on
  dependent data.
\newblock {\em Journal of Multivariate Analysis\/}~{\em 80\/}(2), 234--255.

\bibitem[\protect\citeauthoryear{Pratt}{Pratt}{1968}]{Pratt1968}
Pratt, J.~W. (1968).
\newblock A normal approximation for binomial, {$F$}, beta, and other common,
  related tail probabilities, {II}.
\newblock {\em Journal of the American Statistical Association\/}~{\em
  63\/}(324), 1457--1483.

\bibitem[\protect\citeauthoryear{Qu and Yoon}{Qu and Yoon}{2015}]{QuYoon2015}
Qu, Z. and J.~Yoon (2015).
\newblock Nonparametric estimation and inference on conditional quantile
  processes.
\newblock {\em Journal of Econometrics\/}~{\em 185\/}(1), 1--19.

\bibitem[\protect\citeauthoryear{R{\'e}nyi}{R{\'e}nyi}{1953}]{Renyi1953}
R{\'e}nyi, A. (1953).
\newblock On the theory of order statistics.
\newblock {\em Acta Mathematica Hungarica\/}~{\em 4\/}(3), 191--231.

\bibitem[\protect\citeauthoryear{Robbins}{Robbins}{1955}]{Robbins1955}
Robbins, H. (1955).
\newblock A remark on {Stirling's} formula.
\newblock {\em The American Mathematical Monthly\/}~{\em 62\/}(1), 26--29.

\bibitem[\protect\citeauthoryear{Ruppert, Wand, and Carroll}{Ruppert
  et~al.}{2003}]{RuppertEtAl2003}
Ruppert, D., M.~P. Wand, and R.~J. Carroll (2003).
\newblock {\em Semiparametric Regression}.
\newblock Cambridge Series in Statistical and Probabilistic Mathematics.
  Cambridge University Press.

\bibitem[\protect\citeauthoryear{Shorack}{Shorack}{1972}]{Shorack1972}
Shorack, G.~R. (1972).
\newblock Convergence of quantile and spacings processes with applications.
\newblock {\em Annals of Mathematical Statistics\/}~{\em 43\/}(5), 1400--1411.

\bibitem[\protect\citeauthoryear{Shorack and Wellner}{Shorack and
  Wellner}{1986}]{ShorackWellner1986}
Shorack, G.~R. and J.~A. Wellner (1986).
\newblock {\em Empirical Processes with Applications to Statistics}.
\newblock New York: John Wiley \& Sons.

\bibitem[\protect\citeauthoryear{Stigler}{Stigler}{1977}]{Stigler1977}
Stigler, S.~M. (1977).
\newblock Fractional order statistics, with applications.
\newblock {\em Journal of the American Statistical Association\/}~{\em
  72\/}(359), 544--550.

\bibitem[\protect\citeauthoryear{Thompson}{Thompson}{1936}]{Thompson1936}
Thompson, W.~R. (1936).
\newblock On confidence ranges for the median and other expectation
  distributions for populations of unknown distribution form.
\newblock {\em Annals of Mathematical Statistics\/}~{\em 7\/}(3), 122--128.

\bibitem[\protect\citeauthoryear{Wilks}{Wilks}{1962}]{Wilks1962}
Wilks, S.~S. (1962).
\newblock {\em Mathematical Statistics}.
\newblock New York: Wiley.

\end{thebibliography}

\end{document}